\pdfoutput=1  %% COMMENT OUT IF "lualatex" IS USED

\documentclass[12pt,a4paper]{article}

\usepackage{amsfonts,amssymb,amsthm,amsmath,latexsym}
\usepackage{eqsection,indent}
\usepackage{subeqnarray, eepicemu, url,cite,bm}
\usepackage{multirow}  % See http://en.wikibooks.org/wiki/LaTeX/Tables
\usepackage{tensor}  % To help do 2F1, 2\phi1, etc
\usepackage{graphicx,graphbox}
\usepackage{tikz,xcolor}
\usepackage{float}  % http://tex.stackexchange.com/questions/2275/keeping-tables-figures-close-to-where-they-are-mentioned
\usepackage[scr=dutchcal,scrscaled=1.05]{mathalfa}

\usepackage{color}
\def\blue{\textcolor{blue}}
\def\red{\textcolor{red}}
\def\green{\textcolor{green}}

\date{December 14, 2022}   %% REMEMBER TO PUT FINAL DATE HERE!!!!!
\begin{document}

\title{Classical continued fractions \\
       for some multivariate polynomials \\
       \hspace*{-7mm} \hbox{generalizing the Genocchi and median Genocchi numbers}
      }

\author{ \\
      \hspace*{-1cm}
      {\large Bishal Deb${}^1$ and Alan D.~Sokal${}^{1,2}$}
   \\[5mm]
     \hspace*{-1.3cm}
      \normalsize
           ${}^1$Department of Mathematics, University College London,
                    London WC1E 6BT, UK   \\[1mm]
     \hspace*{-2.9cm}
      \normalsize
           ${}^2$Department of Physics, New York University,
                    New York, NY 10003, USA     \\[1mm]
       \\[5mm]
     \hspace*{-0.5cm}
     {\tt bishal.deb.19@ucl.ac.uk}, {\tt sokal@nyu.edu}  \\[1cm]
}

\maketitle
\thispagestyle{empty}   % Suppress page number on front page.

\begin{abstract}
A D-permutation is a permutation of $[2n]$ satisfying
$2k-1 \le \sigma(2k-1)$ and $2k \ge \sigma(2k)$ for all $k$;
they provide a combinatorial model
for the Genocchi and median Genocchi numbers.
We find Stieltjes-type and Thron-type continued fractions
for some multivariate polynomials that enumerate D-permutations
with respect to a very large (sometimes infinite) number
of simultaneous statistics
that measure cycle status, record status, crossings and nestings.
\end{abstract}

\bigskip
\noindent
{\bf Key Words:}
Genocchi numbers, median Genocchi numbers,
permutation, D-permutation, D-semiderangement, D-derangement, D-cycle,
continued fraction, J-fraction, S-fraction, T-fraction,
Dyck path, almost-Dyck path, Schr\"oder path.
%% Hankel-total positivity (???).

\bigskip
\bigskip
\noindent
{\bf Mathematics Subject Classification (MSC 2020) codes:}
05A19 (Primary);
05A05, 05A15, 05A30, 11B68, 30B70 (Secondary).

%%\clearpage
\vspace*{1cm}

\newtheorem{theorem}{Theorem}[section]
\newtheorem{proposition}[theorem]{Proposition}
\newtheorem{lemma}[theorem]{Lemma}
\newtheorem{corollary}[theorem]{Corollary}
\newtheorem{definition}[theorem]{Definition}
\newtheorem{conjecture}[theorem]{Conjecture}
\newtheorem{question}[theorem]{Question}
\newtheorem{problem}[theorem]{Problem}
\newtheorem{openproblem}[theorem]{Open Problem}
\newtheorem{example}[theorem]{Example}

\renewcommand{\theenumi}{\alph{enumi}}
\renewcommand{\labelenumi}{(\theenumi)}
\def\eop{\hbox{\kern1pt\vrule height6pt width4pt
depth1pt\kern1pt}\medskip}
\def\prf{\par\noindent{\bf Proof.\enspace}\rm}
\def\rmk{\par\medskip\noindent{\bf Remark\enspace}\rm}

\newcommand{\textbfit}[1]{\textbf{\textit{#1}}}

\newcommand{\bigdash}{%
\smallskip\begin{center} \rule{5cm}{0.1mm} \end{center}\smallskip}

% DO NOT USE \par WITHIN A \footnote;  USE \safepar INSTEAD.
\newcommand{\safepar}{ {\protect\hfill\protect\break\hspace*{5mm}} }

\newcommand{\be}{\begin{equation}}
\newcommand{\ee}{\end{equation}}
\newcommand{\<}{\langle}
\renewcommand{\>}{\rangle}
\newcommand{\widebar}{\overline}
\def\reff#1{(\protect\ref{#1})}
\def\spose#1{\hbox to 0pt{#1\hss}}
\def\ltapprox{\mathrel{\spose{\lower 3pt\hbox{$\mathchar"218$}}
    \raise 2.0pt\hbox{$\mathchar"13C$}}}
\def\gtapprox{\mathrel{\spose{\lower 3pt\hbox{$\mathchar"218$}}
    \raise 2.0pt\hbox{$\mathchar"13E$}}}
\def\textprime{${}^\prime$}
\def\proof{\par\medskip\noindent{\sc Proof.\ }}
\def\firstproof{\par\medskip\noindent{\sc First Proof.\ }}
\def\secondproof{\par\medskip\noindent{\sc Second Proof.\ }}
\def\alternateproof{\par\medskip\noindent{\sc Alternate Proof.\ }}
\def\algebraicproof{\par\medskip\noindent{\sc Algebraic Proof.\ }}
\def\combinatorialproof{\par\medskip\noindent{\sc Combinatorial Proof.\ }}
\def\proofof#1{\bigskip\noindent{\sc Proof of #1.\ }}
\def\firstproofof#1{\bigskip\noindent{\sc First Proof of #1.\ }}
\def\secondproofof#1{\bigskip\noindent{\sc Second Proof of #1.\ }}
\def\thirdproofof#1{\bigskip\noindent{\sc Third Proof of #1.\ }}
\def\algebraicproofof#1{\bigskip\noindent{\sc Algebraic Proof of #1.\ }}
\def\combinatorialproofof#1{\bigskip\noindent{\sc Combinatorial Proof of #1.\ }}
\def\sketchofproof{\par\medskip\noindent{\sc Sketch of proof.\ }}
\renewcommand{\qed}{ $\square$ \bigskip}
\newcommand{\myendremark}{ $\blacksquare$ \bigskip}
\def\half{ {1 \over 2} }
\def\third{ {1 \over 3} }
\def\twothird{ {2 \over 3} }
\def\smfrac#1#2{{\textstyle{#1\over #2}}}
\def\smhalf{ {\smfrac{1}{2}} }
\newcommand{\real}{\mathop{\rm Re}\nolimits}
\renewcommand{\Re}{\mathop{\rm Re}\nolimits}
\newcommand{\imag}{\mathop{\rm Im}\nolimits}
\renewcommand{\Im}{\mathop{\rm Im}\nolimits}
\newcommand{\sgn}{\mathop{\rm sgn}\nolimits}
\newcommand{\tr}{\mathop{\rm tr}\nolimits}
\newcommand{\supp}{\mathop{\rm supp}\nolimits}
\newcommand{\disc}{\mathop{\rm disc}\nolimits}
\newcommand{\diag}{\mathop{\rm diag}\nolimits}
\newcommand{\tridiag}{\mathop{\rm tridiag}\nolimits}
\newcommand{\AZ}{\mathop{\rm AZ}\nolimits}
\newcommand{\NC}{\mathop{\rm NC}\nolimits}
\newcommand{\PF}{{\rm PF}}
\newcommand{\rk}{\mathop{\rm rk}\nolimits}
\newcommand{\perm}{\mathop{\rm perm}\nolimits}
\def\hboxscript#1{ {\hbox{\scriptsize\em #1}} }
\renewcommand{\emptyset}{\varnothing}
\newcommand{\eqdef}{\stackrel{\rm def}{=}}

\newcommand{\restrict}{\upharpoonright}

\newcommand{\compinv}{{\langle -1 \rangle}}   % Compositional inverse

\newcommand{\scra}{{\mathcal{A}}}
\newcommand{\scrb}{{\mathcal{B}}}
\newcommand{\scrc}{{\mathcal{C}}}
\newcommand{\bfscra}{{\bm{\mathcal{A}}}}
\newcommand{\bfscrb}{{\bm{\mathcal{B}}}}
\newcommand{\bfscrc}{{\bm{\mathcal{C}}}}
\newcommand{\bfscrap}{{\bm{\mathcal{A}'}}}
\newcommand{\bfscrbp}{{\bm{\mathcal{B}'}}}
\newcommand{\bfscrcp}{{\bm{\mathcal{C}'}}}
\newcommand{\bfscrapp}{{\bm{\mathcal{A}''}}}
\newcommand{\bfscrbpp}{{\bm{\mathcal{B}''}}}
\newcommand{\bfscrcpp}{{\bm{\mathcal{C}''}}}
\newcommand{\scrd}{{\mathcal{D}}}
\newcommand{\scre}{{\mathcal{E}}}
\newcommand{\scrf}{{\mathcal{F}}}
\newcommand{\scrg}{{\mathcal{G}}}
\newcommand{\scrgg}{{\mathscr{g}}}  % Lower-case calligraphic g (using Dutchcal)
\newcommand{\scrh}{{\mathcal{H}}}
\newcommand{\scri}{{\mathcal{I}}}
\newcommand{\scrj}{{\mathcal{J}}}
\newcommand{\scrk}{{\mathcal{K}}}
\newcommand{\scrl}{{\mathcal{L}}}
\newcommand{\scrm}{{\mathcal{M}}}
\newcommand{\scrn}{{\mathcal{N}}}
\newcommand{\scro}{{\mathcal{O}}}
\newcommand\scroo{
  \mathchoice
    {{\scriptstyle\mathcal{O}}}% \displaystyle
    {{\scriptstyle\mathcal{O}}}% \textstyle
    {{\scriptscriptstyle\mathcal{O}}}% \scriptstyle
    {\scalebox{0.6}{$\scriptscriptstyle\mathcal{O}$}}%\scriptscriptstyle
  }
%% Taken from https://tex.stackexchange.com/questions/191059/how-to-get-a-small-letter-version-of-mathcalo
\newcommand{\scrp}{{\mathcal{P}}}
\newcommand{\scrq}{{\mathcal{Q}}}
\newcommand{\scrr}{{\mathcal{R}}}
\newcommand{\scrs}{{\mathcal{S}}}
\newcommand{\scrss}{{\mathscr{s}}}  % Lower-case calligraphic s (using Dutchcal)
\newcommand{\scrt}{{\mathcal{T}}}
\newcommand{\scrv}{{\mathcal{V}}}
\newcommand{\scrw}{{\mathcal{W}}}
\newcommand{\scrz}{{\mathcal{Z}}}
\newcommand{\SP}{{\mathcal{SP}}}
\newcommand{\ST}{{\mathcal{ST}}}

\newcommand{\bfa}{{\mathbf{a}}}
\newcommand{\bfb}{{\mathbf{b}}}
\newcommand{\bfc}{{\mathbf{c}}}
\newcommand{\bfd}{{\mathbf{d}}}
\newcommand{\bfe}{{\mathbf{e}}}
\newcommand{\bfh}{{\mathbf{h}}}
\newcommand{\bfj}{{\mathbf{j}}}
\newcommand{\bfi}{{\mathbf{i}}}
\newcommand{\bfk}{{\mathbf{k}}}
\newcommand{\bfl}{{\mathbf{l}}}
\newcommand{\bfL}{{\mathbf{L}}}
\newcommand{\bfm}{{\mathbf{m}}}
\newcommand{\bfn}{{\mathbf{n}}}
\newcommand{\bfp}{{\mathbf{p}}}
\newcommand{\bfr}{{\mathbf{r}}}
\newcommand{\bfu}{{\mathbf{u}}}
\newcommand{\bfv}{{\mathbf{v}}}
\newcommand{\bfw}{{\mathbf{w}}}
\newcommand{\bfx}{{\mathbf{x}}}
\newcommand{\bfy}{{\mathbf{y}}}
\newcommand{\bfz}{{\mathbf{z}}}
\renewcommand{\k}{{\mathbf{k}}}
\newcommand{\n}{{\mathbf{n}}}
\newcommand{\vv}{{\mathbf{v}}}
\newcommand{\bv}{{\mathbf{v}}}
\newcommand{\w}{{\mathbf{w}}}
\newcommand{\x}{{\mathbf{x}}}
\newcommand{\y}{{\mathbf{y}}}
\newcommand{\cc}{{\mathbf{c}}}
\newcommand{\zero}{{\mathbf{0}}}
\newcommand{\one}{{\mathbf{1}}}
\newcommand{\bmm}{{\mathbf{m}}}

\newcommand{\ahat}{{\widehat{a}}}
\newcommand{\vhat}{{\widehat{v}}}
\newcommand{\yhat}{{\widehat{y}}}
\newcommand{\phat}{{\widehat{p}}}
\newcommand{\qhat}{{\widehat{q}}}
\newcommand{\Zhat}{{\widehat{Z}}}
\newcommand{\myhat}{{\widehat{\;}}}
\newcommand{\vtilde}{{\widetilde{v}}}
\newcommand{\ytilde}{{\widetilde{y}}}

\newcommand{\C}{{\mathbb C}}
\newcommand{\D}{{\mathbb D}}
\newcommand{\Z}{{\mathbb Z}}
\newcommand{\N}{{\mathbb N}}
\newcommand{\Q}{{\mathbb Q}}
\newcommand{\PP}{{\mathbb P}}
\newcommand{\R}{{\mathbb R}}
\newcommand{\RR}{{\mathbb R}}
\newcommand{\E}{{\mathbb E}}

\newcommand{\Sym}{{\mathfrak{S}}}
\newcommand{\SymB}{{\mathfrak{B}}}
\newcommand{\Alt}{{\mathrm{Alt}}}
\newcommand{\dperm}{{\mathfrak{D}}}
\newcommand{\dcycle}{{\mathfrak{DC}}}

\newcommand{\germanA}{{\mathfrak{A}}}
\newcommand{\germanB}{{\mathfrak{B}}}
\newcommand{\germanQ}{{\mathfrak{Q}}}
\newcommand{\germanh}{{\mathfrak{h}}}

\newcommand{\myle}{\preceq}
\newcommand{\myge}{\succeq}
\newcommand{\mygt}{\succ}

\newcommand{\B}{{\sf B}}
\newcommand{\OB}{B^{\rm ord}}
\newcommand{\OS}{{\sf OS}}
\newcommand{\OO}{{\sf O}}
\newcommand{\OSP}{{\sf OSP}}
\newcommand{\Eu}{{\sf Eu}}
\newcommand{\ERR}{{\sf ERR}}
\newcommand{\sfB}{{\sf B}}
\newcommand{\sfD}{{\sf D}}
\newcommand{\sfE}{{\sf E}}
\newcommand{\sfG}{{\sf G}}
\newcommand{\sfJ}{{\sf J}}
\newcommand{\sfP}{{\sf P}}
\newcommand{\sfQ}{{\sf Q}}
\newcommand{\sfS}{{\sf S}}
\newcommand{\sfT}{{\sf T}}
\newcommand{\sfW}{{\sf W}}
\newcommand{\sfMV}{{\sf MV}}
\newcommand{\AMV}{{\sf AMV}}
\newcommand{\BM}{{\sf BM}}
\newcommand{\emIB}{B^{\rm irr}}
\newcommand{\emIP}{P^{\rm irr}}
\newcommand{\emOB}{B^{\rm ord}}
\newcommand{\emCB}{B^{\rm cyc}}
\newcommand{\emSC}{P^{\rm cyc}}

\newcommand{\lev}{{\rm lev}}
\newcommand{\stat}{{\rm stat}}
\newcommand{\cyc}{{\rm cyc}}
\newcommand{\mysteryone}{{\rm mys1}}
\newcommand{\mysterytwo}{{\rm mys2}}
\newcommand{\Asc}{{\rm Asc}}
\newcommand{\asc}{{\rm asc}}
\newcommand{\Des}{{\rm Des}}
\newcommand{\des}{{\rm des}}
\newcommand{\Exc}{{\rm Exc}}

\newcommand{\EArec}{{\rm EArec}}
\newcommand{\earec}{{\rm earec}}
\newcommand{\recarec}{{\rm recarec}}
\newcommand{\erec}{{\rm erec}}
\newcommand{\nonrec}{{\rm nonrec}}
\newcommand{\nrar}{{\rm nrar}}
\newcommand{\ereccval}{{\rm ereccval}}
\newcommand{\ereccdrise}{{\rm ereccdrise}}
\newcommand{\eareccpeak}{{\rm eareccpeak}}
\newcommand{\eareccdfall}{{\rm eareccdfall}}
\newcommand{\eareccval}{{\rm eareccval}}
\newcommand{\ereccpeak}{{\rm ereccpeak}}
\newcommand{\rar}{{\rm rar}}
\newcommand{\evenrar}{{\rm evenrar}}
\newcommand{\oddrar}{{\rm oddrar}}
\newcommand{\nrcpeak}{{\rm nrcpeak}}
\newcommand{\nrcval}{{\rm nrcval}}
\newcommand{\nrcdrise}{{\rm nrcdrise}}
\newcommand{\nrcdfall}{{\rm nrcdfall}}
\newcommand{\nrfix}{{\rm nrfix}}
\newcommand{\Evenfix}{{\rm Evenfix}}
\newcommand{\evenfix}{{\rm evenfix}}
\newcommand{\Oddfix}{{\rm Oddfix}}
\newcommand{\oddfix}{{\rm oddfix}}
\newcommand{\evennrfix}{{\rm evennrfix}}
\newcommand{\oddnrfix}{{\rm oddnrfix}}
\newcommand{\Cpeak}{{\rm Cpeak}}
\newcommand{\cpeak}{{\rm cpeak}}
\newcommand{\Cval}{{\rm Cval}}
\newcommand{\cval}{{\rm cval}}
\newcommand{\Cdasc}{{\rm Cdasc}}
\newcommand{\cdasc}{{\rm cdasc}}
\newcommand{\Cddes}{{\rm Cddes}}
\newcommand{\cddes}{{\rm cddes}}
\newcommand{\Cdrise}{{\rm Cdrise}}
\newcommand{\cdrise}{{\rm cdrise}}
\newcommand{\Cdfall}{{\rm Cdfall}}
\newcommand{\cdfall}{{\rm cdfall}}

\newcommand{\maxpeak}{{\rm maxpeak}}
\newcommand{\nmaxpeak}{{\rm nmaxpeak}}
\newcommand{\minval}{{\rm minval}}
\newcommand{\nminval}{{\rm nminval}}

\newcommand{\exc}{{\rm exc}}
\newcommand{\excee}{{\rm excee}}
\newcommand{\exceo}{{\rm exceo}}
\newcommand{\excoe}{{\rm excoe}}
\newcommand{\excoo}{{\rm excoo}}
\newcommand{\erecexcoe}{{\rm erecexcoe}}
\newcommand{\erecexcoo}{{\rm erecexcoo}}
\newcommand{\nrexcoe}{{\rm nrexcoe}}
\newcommand{\nrexcoo}{{\rm nrexcoo}}
\newcommand{\aexc}{{\rm aexc}}
\newcommand{\aexcee}{{\rm aexcee}}
\newcommand{\aexceo}{{\rm aexceo}}
\newcommand{\aexcoe}{{\rm aexcoe}}
\newcommand{\aexcoo}{{\rm aexcoo}}
\newcommand{\earecaexcee}{{\rm earecaexcee}}
\newcommand{\earecaexceo}{{\rm earecaexceo}}
\newcommand{\nraexcee}{{\rm nraexcee}}
\newcommand{\nraexceo}{{\rm nraexceo}}
\newcommand{\Fix}{{\rm Fix}}
\newcommand{\fix}{{\rm fix}}
\newcommand{\fixe}{{\rm fixe}}
\newcommand{\fixo}{{\rm fixo}}
\newcommand{\rare}{{\rm rare}}
\newcommand{\raro}{{\rm raro}}
\newcommand{\nrfixe}{{\rm nrfixe}}
\newcommand{\nrfixo}{{\rm nrfixo}}
\newcommand{\xee}{x_{\rm ee}}
\newcommand{\xeo}{x_{\rm eo}}
\newcommand{\uee}{u_{\rm ee}}
\newcommand{\ueo}{u_{\rm eo}}
\newcommand{\yoo}{y_{\rm oo}}
\newcommand{\yoe}{y_{\rm oe}}
\newcommand{\voo}{v_{\rm oo}}
\newcommand{\voe}{v_{\rm oe}}
\newcommand{\zo}{z_{\rm o}}
\newcommand{\ze}{z_{\rm e}}
\newcommand{\wo}{w_{\rm o}}
\newcommand{\we}{w_{\rm e}}
\newcommand{\so}{s_{\rm o}}
\newcommand{\se}{s_{\rm e}}

\newcommand{\Wex}{{\rm Wex}}
\newcommand{\wex}{{\rm wex}}
\newcommand{\lrmax}{{\rm lrmax}}
\newcommand{\rlmax}{{\rm rlmax}}
\newcommand{\Rec}{{\rm Rec}}
\newcommand{\rec}{{\rm rec}}
\newcommand{\Arec}{{\rm Arec}}
\newcommand{\arec}{{\rm arec}}
\newcommand{\ERec}{{\rm ERec}}
\newcommand{\Val}{{\rm Val}}
\newcommand{\val}{{\rm val}}
\newcommand{\dasc}{{\rm dasc}}
\newcommand{\ddes}{{\rm ddes}}
\newcommand{\inv}{{\rm inv}}
\newcommand{\maj}{{\rm maj}}
\newcommand{\rs}{{\rm rs}}
\newcommand{\cross}{{\rm cr}}
\newcommand{\crosshat}{{\widehat{\rm cr}}}
\newcommand{\nest}{{\rm ne}}
\newcommand{\ucross}{{\rm ucross}}
\newcommand{\ucrosscval}{{\rm ucrosscval}}
\newcommand{\ucrosscpeak}{{\rm ucrosscpeak}}
\newcommand{\ucrosscdrise}{{\rm ucrosscdrise}}
\newcommand{\lcross}{{\rm lcross}}
\newcommand{\lcrosscpeak}{{\rm lcrosscpeak}}
\newcommand{\lcrosscval}{{\rm lcrosscval}}
\newcommand{\lcrosscdfall}{{\rm lcrosscdfall}}
\newcommand{\unest}{{\rm unest}}
\newcommand{\unestcval}{{\rm unestcval}}
\newcommand{\unestcpeak}{{\rm unestcpeak}}
\newcommand{\unestcdrise}{{\rm unestcdrise}}
\newcommand{\lnest}{{\rm lnest}}
\newcommand{\lnestcpeak}{{\rm lnestcpeak}}
\newcommand{\lnestcval}{{\rm lnestcval}}
\newcommand{\lnestcdfall}{{\rm lnestcdfall}}
\newcommand{\ujoin}{{\rm ujoin}}
\newcommand{\ljoin}{{\rm ljoin}}
\newcommand{\psnest}{{\rm psnest}}
\newcommand{\upsnest}{{\rm upsnest}}
\newcommand{\lpsnest}{{\rm lpsnest}}
\newcommand{\epsnest}{{\rm epsnest}}
\newcommand{\opsnest}{{\rm opsnest}}
\newcommand{\rodd}{{\rm rodd}}
\newcommand{\reven}{{\rm reven}}
\newcommand{\lodd}{{\rm lodd}}
\newcommand{\leven}{{\rm leven}}
\newcommand{\sg}{{\rm sg}}
\newcommand{\bl}{{\rm bl}}
\newcommand{\tran}{{\rm tr}}
\newcommand{\area}{{\rm area}}
\newcommand{\ret}{{\rm ret}}
\newcommand{\peaks}{{\rm peaks}}
\newcommand{\hl}{{\rm hl}}
\newcommand{\sll}{{\rm sl}}
\newcommand{\negg}{{\rm neg}}
\newcommand{\imp}{{\rm imp}}
\newcommand{\osg}{{\rm osg}}
\newcommand{\ons}{{\rm ons}}
\newcommand{\isg}{{\rm isg}}
\newcommand{\ins}{{\rm ins}}
\newcommand{\LL}{{\rm LL}}
\newcommand{\height}{{\rm ht}}
\newcommand{\as}{{\rm as}}

\newcommand{\ba}{{\bm{a}}}
\newcommand{\bahat}{{\widehat{\bm{a}}}}
\newcommand{\bb}{{\bm{b}}}
\newcommand{\bc}{{\bm{c}}}
\newcommand{\bchat}{{\widehat{\bm{c}}}}
\newcommand{\bd}{{\bm{d}}}
\newcommand{\bee}{{\bm{e}}}
\newcommand{\beh}{{\bm{eh}}}
\newcommand{\bff}{{\bm{f}}}
\newcommand{\bg}{{\bm{g}}}
\newcommand{\bh}{{\bm{h}}}
\newcommand{\bll}{{\bm{\ell}}}
\newcommand{\bp}{{\bm{p}}}
\newcommand{\br}{{\bm{r}}}
\newcommand{\bs}{{\bm{s}}}
\newcommand{\bu}{{\bm{u}}}
\newcommand{\bw}{{\bm{w}}}
\newcommand{\bx}{{\bm{x}}}
\newcommand{\by}{{\bm{y}}}
\newcommand{\bz}{{\bm{z}}}
\newcommand{\bA}{{\bm{A}}}
\newcommand{\bB}{{\bm{B}}}
\newcommand{\bC}{{\bm{C}}}
\newcommand{\bE}{{\bm{E}}}
\newcommand{\bF}{{\bm{F}}}
\newcommand{\bG}{{\bm{G}}}
\newcommand{\bH}{{\bm{H}}}
\newcommand{\bI}{{\bm{I}}}
\newcommand{\bJ}{{\bm{J}}}
\newcommand{\bM}{{\bm{M}}}
\newcommand{\bN}{{\bm{N}}}
\newcommand{\bP}{{\bm{P}}}
\newcommand{\bQ}{{\bm{Q}}}
\newcommand{\bR}{{\bm{R}}}
\newcommand{\bS}{{\bm{S}}}
\newcommand{\bT}{{\bm{T}}}
\newcommand{\bW}{{\bm{W}}}
\newcommand{\bX}{{\bm{X}}}
\newcommand{\bY}{{\bm{Y}}}
\newcommand{\bIB}{{\bm{B}^{\rm irr}}}
\newcommand{\bOB}{{\bm{B}^{\rm ord}}}
\newcommand{\bOS}{{\bm{OS}}}
\newcommand{\bERR}{{\bm{ERR}}}
\newcommand{\bSP}{{\bm{SP}}}
\newcommand{\bMV}{{\bm{MV}}}
\newcommand{\bBM}{{\bm{BM}}}
\newcommand{\balpha}{{\bm{\alpha}}}
\newcommand{\balphapre}{{\bm{\alpha}^{\rm pre}}}
\newcommand{\bbeta}{{\bm{\beta}}}
\newcommand{\bgamma}{{\bm{\gamma}}}
\newcommand{\bdelta}{{\bm{\delta}}}
\newcommand{\bkappa}{{\bm{\kappa}}}
\newcommand{\bmu}{{\bm{\mu}}}
\newcommand{\bomega}{{\bm{\omega}}}
\newcommand{\bsigma}{{\bm{\sigma}}}
\newcommand{\btau}{{\bm{\tau}}}
\newcommand{\bphi}{{\bm{\phi}}}
\newcommand{\bpsi}{{\bm{\psi}}}
\newcommand{\bzeta}{{\bm{\zeta}}}
\newcommand{\bone}{{\bm{1}}}
\newcommand{\bzero}{{\bm{0}}}

\newcommand{\sfa}{{{\sf a}}}
\newcommand{\sfb}{{{\sf b}}}
\newcommand{\sfc}{{{\sf c}}}
\newcommand{\sfd}{{{\sf d}}}
\newcommand{\sfe}{{{\sf e}}}
\newcommand{\sff}{{{\sf f}}}
\newcommand{\sfg}{{{\sf g}}}
\newcommand{\sfh}{{{\sf h}}}
\newcommand{\sfi}{{{\sf i}}}
\newcommand{\bsfa}{{\mbox{\textsf{\textbf{a}}}}}
\newcommand{\bsfb}{{\mbox{\textsf{\textbf{b}}}}}
\newcommand{\bsfc}{{\mbox{\textsf{\textbf{c}}}}}
\newcommand{\bsfd}{{\mbox{\textsf{\textbf{d}}}}}
\newcommand{\bsfe}{{\mbox{\textsf{\textbf{e}}}}}
\newcommand{\bsff}{{\mbox{\textsf{\textbf{f}}}}}
\newcommand{\bsfg}{{\mbox{\textsf{\textbf{g}}}}}
\newcommand{\bsfh}{{\mbox{\textsf{\textbf{h}}}}}
\newcommand{\bsfi}{{\mbox{\textsf{\textbf{i}}}}}

\newcommand{\Cbar}{{\overline{C}}}
\newcommand{\Dbar}{{\overline{D}}}
\newcommand{\dbar}{{\overline{d}}}
\def\Ctilde{{\widetilde{C}}}
\def\Ftilde{{\widetilde{F}}}
\def\Gtilde{{\widetilde{G}}}
\def\Htilde{{\widetilde{H}}}
\def\Ptilde{{\widetilde{P}}}
\def\Chat{{\widehat{C}}}
\def\ctilde{{\widetilde{c}}}
\def\zbar{{\overline{Z}}}
\def\pitilde{{\widetilde{\pi}}}
\def\omegahat{{\widehat{\omega}}}

\newcommand{\sech}{{\rm sech}}

\newcommand{\sinv}{\sigma^{-1}}

%
% Jacobian and Dixonian elliptic functions
%
\newcommand{\sn}{{\rm sn}}
\newcommand{\cn}{{\rm cn}}
\newcommand{\dn}{{\rm dn}}
\newcommand{\sm}{{\rm sm}}
\newcommand{\cm}{{\rm cm}}

%
% Commands for hypergeometric series
%
\newcommand{\zfz}{ {{}_0 \! F_0} }
\newcommand{\zfo}{ {{}_0  F_1} }
\newcommand{\ofz}{ {{}_1 \! F_0} }
\newcommand{\ofo}{ {{}_1 \! F_1} }
\newcommand{\oft}{ {{}_1 \! F_2} }
%\newcommand{\tfo}{ {{}_2 \! F_1} }

%
% Hypergeometric functions, using "tensor" package
%
\newcommand{\FHyper}[2]{ {\tensor[_{#1 \!}]{F}{_{#2}}\!} }
\newcommand{\FHYPER}[5]{ {\FHyper{#1}{#2} \!\biggl(
   \!\!\begin{array}{c} #3 \\[1mm] #4 \end{array}\! \bigg|\, #5 \! \biggr)} }
\newcommand{\tfo}{ {\FHyper{2}{1}} }
\newcommand{\tfz}{ {\FHyper{2}{0}} }
\newcommand{\threefz}{ {\FHyper{3}{0}} }
\newcommand{\FHYPERbottomzero}[3]{ {\FHyper{#1}{0} \hspace*{-0mm}\biggl(
   \!\!\begin{array}{c} #2 \\[1mm] \hbox{---} \end{array}\! \bigg|\, #3 \! \biggr)} }
\newcommand{\FHYPERtopzero}[3]{ {\FHyper{0}{#1} \hspace*{-0mm}\biggl(
   \!\!\begin{array}{c} \hbox{---} \\[1mm] #2 \end{array}\! \bigg|\, #3 \! \biggr)} }

\newcommand{\phiHyper}[2]{ {\tensor[_{#1}]{\phi}{_{#2}}} }
\newcommand{\psiHyper}[2]{ {\tensor[_{#1}]{\psi}{_{#2}}} }
\newcommand{\PhiHyper}[2]{ {\tensor[_{#1}]{\Phi}{_{#2}}} }
\newcommand{\PsiHyper}[2]{ {\tensor[_{#1}]{\Psi}{_{#2}}} }
\newcommand{\phiHYPER}[6]{ {\phiHyper{#1}{#2} \!\left(
   \!\!\begin{array}{c} #3 \\ #4 \end{array}\! ;\, #5, \, #6 \! \right)\!} }
\newcommand{\psiHYPER}[6]{ {\psiHyper{#1}{#2} \!\left(
   \!\!\begin{array}{c} #3 \\ #4 \end{array}\! ;\, #5, \, #6 \! \right)} }
\newcommand{\PhiHYPER}[5]{ {\PhiHyper{#1}{#2} \!\left(
   \!\!\begin{array}{c} #3 \\ #4 \end{array}\! ;\, #5 \! \right)\!} }
\newcommand{\PsiHYPER}[5]{ {\PsiHyper{#1}{#2} \!\left(
   \!\!\begin{array}{c} #3 \\ #4 \end{array}\! ;\, #5 \! \right)\!} }
\newcommand{\zerophizero}{ {\phiHyper{0}{0}} }
\newcommand{\ophizero}{ {\phiHyper{1}{0}} }
\newcommand{\zphio}{ {\phiHyper{0}{1}} }
\newcommand{\ophio}{ {\phiHyper{1}{1}} }
\newcommand{\tphio}{ {\phiHyper{2}{1}} }
\newcommand{\tphiz}{ {\phiHyper{2}{0}} }
\newcommand{\tPhio}{ {\PhiHyper{2}{1}} }
\newcommand{\opsio}{ {\psiHyper{1}{1}} }

%
% Variants of \binom  (defined using the AMS "genfrac" command)
%
\newcommand{\stirlingsubset}[2]{\genfrac{\{}{\}}{0pt}{}{#1}{#2}}
\newcommand{\stirlingcycle}[2]{\genfrac{[}{]}{0pt}{}{#1}{#2}}
\newcommand{\assocstirlingsubset}[3]{{\genfrac{\{}{\}}{0pt}{}{#1}{#2}}_{\! \ge #3}}
\newcommand{\genstirlingsubset}[4]{{\genfrac{\{}{\}}{0pt}{}{#1}{#2}}_{\! #3,#4}}
\newcommand{\irredstirlingsubset}[2]{{\genfrac{\{}{\}}{0pt}{}{#1}{#2}}^{\!\rm irr}}
\newcommand{\euler}[2]{\genfrac{\langle}{\rangle}{0pt}{}{#1}{#2}}
\newcommand{\eulergen}[3]{{\genfrac{\langle}{\rangle}{0pt}{}{#1}{#2}}_{\! #3}}
\newcommand{\eulersecond}[2]{\left\langle\!\! \euler{#1}{#2} \!\!\right\rangle}
\newcommand{\eulersecondgen}[3]{{\left\langle\!\! \euler{#1}{#2} \!\!\right\rangle}_{\! #3}}
\newcommand{\binomvert}[2]{\genfrac{\vert}{\vert}{0pt}{}{#1}{#2}}
\newcommand{\binomsquare}[2]{\genfrac{[}{]}{0pt}{}{#1}{#2}}
\newcommand{\doublebinom}[2]{\left(\!\! \binom{#1}{#2} \!\!\right)}

% Array for subscripts

\newenvironment{sarray}{
             \textfont0=\scriptfont0
             \scriptfont0=\scriptscriptfont0
             \textfont1=\scriptfont1
             \scriptfont1=\scriptscriptfont1
             \textfont2=\scriptfont2
             \scriptfont2=\scriptscriptfont2
             \textfont3=\scriptfont3
             \scriptfont3=\scriptscriptfont3
           \renewcommand{\arraystretch}{0.7}
           \begin{array}{l}}{\end{array}}

\newenvironment{scarray}{
             \textfont0=\scriptfont0
             \scriptfont0=\scriptscriptfont0
             \textfont1=\scriptfont1
             \scriptfont1=\scriptscriptfont1
             \textfont2=\scriptfont2
             \scriptfont2=\scriptscriptfont2
             \textfont3=\scriptfont3
             \scriptfont3=\scriptscriptfont3
           \renewcommand{\arraystretch}{0.7}
           \begin{array}{c}}{\end{array}}

% Circled math symbols:
% From http://latex-community.org/forum/viewtopic.php?f=44&t=22367

%\usepackage{tikz}
\newcommand*\circled[1]{\tikz[baseline=(char.base)]{
  \node[shape=circle,draw,inner sep=1pt] (char) {#1};}}
\newcommand{\ostar}{{\circledast}}
\newcommand{\ostarN}{{\,\circledast_{\vphantom{\dot{N}}N}\,}}
\newcommand{\ostarPsi}{{\,\circledast_{\vphantom{\dot{\Psi}}\Psi}\,}}
\newcommand{\starN}{{\,\ast_{\vphantom{\dot{N}}N}\,}}
\newcommand{\starpsi}{{\,\ast_{\vphantom{\dot{\bpsi}}\!\bpsi}\,}}
\newcommand{\starone}{{\,\ast_{\vphantom{\dot{1}}1}\,}}
\newcommand{\startwo}{{\,\ast_{\vphantom{\dot{2}}2}\,}}
\newcommand{\starinfty}{{\,\ast_{\vphantom{\dot{\infty}}\infty}\,}}
\newcommand{\starT}{{\,\ast_{\vphantom{\dot{T}}T}\,}}

%% For scaling equations (uses "graphicx" package):  see
%% http://tex.stackexchange.com/questions/60453/reducing-font-size-in-equation
\newcommand*{\Scale}[2][4]{\scalebox{#1}{$#2$}}

\newcommand*{\Scaletext}[2][4]{\scalebox{#1}{#2}} %% THIS DOESN'T SEEM TO WORK

\clearpage

\tableofcontents

\clearpage

\section{Introduction}
 
Our purpose here is to present continued fractions
for some multivariate polynomials that generalize
either the Genocchi numbers \cite[A110501]{OEIS}\footnote{
   The Genocchi numbers appear already in Euler's book
   {\em Foundations of Differential Calculus,
    with Applications to Finite Analysis and Series}\/,
   first published in 1755 \cite[paragraphs~181 and 182]{Euler_1755};
   this book is E212 in Enestr\"om's \cite{Enestrom_13} catalogue.
   These numbers were revisited by Genocchi \cite{Genocchi_1852} in 1852.
   The beautiful survey article of Viennot \cite{Viennot_81}
   contains a wealth of useful information.

   Our notation for the Genocchi and median Genocchi numbers
   is nonstandard but, we believe, sensible and logical.
   Later we will present a translation dictionary with respect to
   (the plethora of) previous notations:
   see footnotes~\ref{footnote.Gn} and \ref{footnote.Hn}.
}
\be
   (g_n)_{n \ge 0}
   \;=\;
   1, 1, 3, 17, 155, 2073, 38227, 929569, 28820619, 1109652905,
   %% 51943281731,
   \ldots
 \label{eq.genocchi}
\ee
or the median Genocchi numbers \cite[A005439]{OEIS}
\be
   (h_n)_{n \ge 0}
   \;=\;
   1, 1, 2, 8, 56, 608, 9440, 198272, 5410688, 186043904,
   %% 7867739648,
   \ldots
   \;.
 \label{eq.median}
\ee
The present paper can be viewed as a sequel
to a recent article by Zeng and one of us \cite{Sokal-Zeng_masterpoly}
in which we presented Stieltjes-type and Jacobi-type continued fractions
for some ``master polynomials''
that enumerate permutations, set partitions or perfect matchings
with respect to a large (sometimes infinite) number of independent statistics.
These polynomials systematize what we think of as the ``linear family'':
namely, sequences in which
the Stieltjes continued-fraction coefficients $(\alpha_n)_{n \ge 1}$
grow linearly in $n$.
More precisely, in the simplest case
\cite[Theorem~2.1]{Sokal-Zeng_masterpoly}
the even and odd coefficients grow affinely in $n$:
\begin{subeqnarray}
   \alpha_{2k-1}  & = &  x + (k-1) u \\
   \alpha_{2k}    & = &  y + (k-1) v
 \label{eq.linear_family}
\end{subeqnarray}
When $x=y=u=v=1$, these coefficients $\alpha_{2k-1} = \alpha_{2k} = k$
correspond to Euler's \cite[section~21]{Euler_1760}
continued fraction for the sequence $a_n = n!$;
so it is natural to expect that the resulting polynomials $P_n(x,y,u,v)$
can be interpreted as enumerating permutations of $[n]$
with respect to some suitable statistics.
The purpose of \cite{Sokal-Zeng_masterpoly}
was to exhibit explicitly those statistics,
and then to present some generalizations involving more refined statistics.

In the present paper we take one step up in complexity, to consider the
``quadratic family'', in~which the $(\alpha_n)_{n \ge 1}$
grow quadratically in $n$.
For instance, we could consider
\begin{subeqnarray}
   \alpha_{2k-1}  & = &  [x_1 + (k-1) u_1] \, [x_2 + (k-1) u_2] \\
   \alpha_{2k}    & = &  [y_1 + (k-1) v_1] \, [y_2 + (k-1) v_2]
 \label{eq.quadratic_family}
\end{subeqnarray}
With all parameters set to 1,
these coefficients $\alpha_{2k-1} = \alpha_{2k} = k^2$
correspond to the continued fraction
\cite[eq.~(9.7)]{Viennot_81} \cite[p.~V-15]{Viennot_83}
for the median Genocchi numbers;
so it is natural to seek a combinatorial model
that is enumerated by the median Genocchi numbers.
Many such models are known.
%% we briefly review a few of them in Appendix~\ref{app.models}.
In this paper we shall focus on a class of permutations of $[2n]$
called {\em D-permutations}\/ \cite{Lazar_22,Lazar_20},
which are defined by imposing some constraints
concerning the parity (even/odd) of excedances and anti-excedances.
Recall that an index $i$ in a permutation $\sigma$
is called an {\em excedance}\/ in case $i < \sigma(i)$,
an {\em anti-excedance}\/ in case $i > \sigma(i)$,
and a {\em fixed point}\/ in case $i = \sigma(i)$.
A permutation is called a \textbfit{D-permutation}
in case it contains no even excedances and no odd anti-excedances.
Equivalently, a permutation is a D-permutation in case
$2k-1 \le \sigma(2k-1)$ and $2k \ge \sigma(2k)$ for all $k$.
Let us say also that a permutation is an
{\em e-semiderangement}\/ (resp.\ {\em o-semiderangement}\/)
in case in contains no even (resp.~odd) fixed points;
it is a {\em derangement}\/ in case it contains no fixed points at all.
A D-permutation that is also an
e-semiderangement (resp.\ o-semiderangement, derangement)
will be called a \textbfit{D-e-semiderangement}
(resp.\ \textbfit{D-o-semiderangement}, \textbfit{D-derangement}).\footnote{
   In the past, D-o-semiderangements have been called
   {\em Genocchi permutations}\/
   \cite{Randrianarivony_96b,Randrianarivony_97,Han_99a}
   or {\em excedance-alternating permutations}\/
   \cite{Ehrenborg_00b};
   D-e-semiderangements \cite[p.~316, Corollaire~1]{Dumont_74}
   have been called {\em Dumont permutations}\/ \cite{Lazar_22,Lazar_20,Pan_21}
   or {\em Dumont permutations of the second kind}\/
   \cite{Burstein_06,Burstein_21};
   D-derangements have been called
   {\em Dumont derangements}\/ \cite{Lazar_22,Lazar_20,Pan_21}.
   D-permutations that are not semiderangements were apparently
   first considered in the recent work of
   Lazar and Wachs \cite{Lazar_22,Lazar_20}.

   Note also that the involution of $\Sym_{2n}$ defined by
   $\sigma \mapsto R \circ \sigma \circ R$,
   where $R(i) = 2n+1-i$ is the reversal map,
   takes D-permutations into D-permutations
   and interchanges e-semiderangements with o-semiderangements;
   it therefore yields a bijection between
   D-e-semiderangements and D-o-semiderangements.
   Therefore, any result about one of the two types of D-semiderangements
   can be expressed equivalently in terms of the other.
}
A D-permutation that contains exactly one cycle is called a \textbfit{D-cycle}. Notice that a D-cycle is also a D-derangement.
Let $\dperm_{2n}$
(resp.~$\dperm^{\rm e}_{2n}, \dperm^{\rm o}_{2n}, \dperm^{\rm eo}_{2n}, \dcycle_{2n}$)
denote the set of all D-permutations
(resp.\ D-e-semiderangements, D-o-semiderangements, D-derangements, D-cycles) of $[2n]$.
For instance,
\begin{subeqnarray}
   \dperm_2  & = & \{ 12 ,\, 21^{\rm eo} \}  \\[1mm]
   \dperm_4  & = & \{ 1234 ,\, 1243 ,\, 2134 ,\, 2143^{\rm eo} ,\,
                     3142^{\rm eo} ,\, 3241^{\rm o} ,\, 4132^{\rm e} ,\, 4231 \}\\[1mm]
\dcycle_2 & = & \{21\}\\[1mm]
\dcycle_4 & = & \{3142\}
%% For a D-permutation of [4], \sigma(2) must be 1 or 2,
%%    and \sigma(3) must be 3 or 4; for each such combination,
%%    there are then two choices for the other two entries.
\end{subeqnarray}
where $^{\rm e}$ denotes e-semiderangements that are not derangements,
$^{\rm o}$ denotes o-semi\-de\-range\-ments that are not derangements,
and $^{\rm eo}$ denotes derangements.

\smallskip

\begin{quote}
\small
{\bf Remark.}
Three natural variants of this setup lead to nothing new:
\begin{itemize}
   \item[1)] In a D-permutation of $[2n+1]$, $2n+1$ must be a fixed point,
       and the rest is a D-permutation of $[2n]$.
   \item[2)] Suppose we define an {\em anti-D-permutation}\/ to be one in which
       $2k-1 \ge \sigma(2k-1)$ and $2k \le \sigma(2k)$ for all $k$.
       Then, in an anti-D-permutation of $[2n]$,
       $1$ and $2n$ must be fixed points,
       and the rest is, after renumbering, a D-permutation of $[2n-2]$.
   \item[3)] In an anti-D-permutation of $[2n+1]$,
       $1$ must be a fixed point,
       and the rest is, after renumbering, a D-permutation of $[2n]$.
\end{itemize}
So there is no loss of generality in studying only D-permutations of $[2n]$.
\myendremark
\end{quote}

\smallskip

It is known \cite{Dumont_74,Dumont_94,Lazar_22,Lazar_20}
--- and we will recover as part of our work --- that
\begin{subeqnarray}
   |\dperm_{2n}|    & = &  h_{n+1}  \\[0.5mm]
   |\dperm^{\rm e}_{2n}|  \;=\; |\dperm^{\rm o}_{2n}|  & = &  g_n \\[1mm]
   |\dperm^{\rm eo}_{2n}|  & = &  h_n\\[0.5mm]
    |\dcycle_{2n}| & = & g_{n-1}
\end{subeqnarray}
%%%%{\bf Define also D-cycles $\dcycle_{2n}$, which satisfy
%%%%   $|\dcycle_{2n}| = g_{n-1}$???
%%%%   This formula should come as a corollary of our second T-fraction
%%%%   that includes $\lambda^{\cyc(\sigma)}$.}
%% {\bf Ehrenborg and Steigr\'{\i}msson \cite[p.~594]{Ehrenborg_00b}
%%    say that $h_n$ is also the number of D-derangements of $[2n+1]$.
%%    So maybe we should define these things also for odd length.}
%% NO!!!  A D-permutation of [2n+1] must have 2n+1 as a fixed point
%%    (so it *cannot* be a D-derangement!!),
%%    and the rest is a D-permutation of [2n].
%%    So there is *no* point in considering D-permutations of [2n+1]!
This suggests that we can obtain continued fractions
for multivariate polynomials enumerating
D-permutations, D-semiderangements, D-derangements or D-cycles
by generalizing the known continued fractions
for the Genocchi and median Genocchi numbers,
analogously to what was done in \cite{Sokal-Zeng_masterpoly}
by generalizing Euler's continued fraction for the factorials.
That is, indeed, what we shall do in this paper.
These continued fractions will be of Stieltjes and Thron types;
they can also be transformed by contraction into Jacobi type.\footnote{
   We call these {\em classical}\/ continued fractions,
   in order to contrast them with the recently-introduced
   {\em branched}\/ continued fractions \cite{latpath_SRTR},
   which are believed to apply to certain higher-order generalizations
   of the Genocchi numbers \cite[Conjecture~16.1]{latpath_SRTR}.
}
Our principal results will be:
\begin{itemize}
   \item[1)] A Thron-type continued fraction in 12~variables
      (Theorem~\ref{thm.Tfrac.first}) that enumerates D-permutations
      with respect to the parity-refined record-and-cycle classification
      (defined in Section~\ref{subsec.statistics.1}).
      Specializations of this continued fraction
      give Stieltjes-type continued fractions that enumerate
      D-semiderangements and D-derangements.
   \item[2)] A Thron-type continued fraction in 22~variables
      (Theorem~\ref{thm.Tfrac.first.pq}) that generalizes the previous
      T-fraction by including four pairs of $(p,q)$-variables
      that count crossings and nestings
      and a pair of variables that count pseudo-nestings of fixed points
      (defined in Section~\ref{subsec.statistics.2}).
   \item[3)] A Thron-type continued fraction in six infinite families
      of indeterminates (Theorem~\ref{thm.Tfrac.first.master})
      that generalizes the preceding two by further refining
      the counting of crossings, nestings and pseudo-nestings.
\end{itemize}
We call these the three versions of our ``first T-fraction''.
Already the first version (Theorem~\ref{thm.Tfrac.first})
contains several known continued fractions
for the Genocchi and median Genocchi numbers as special cases.
We also have three variant versions of the first T-fraction
(Theorems~\ref{thm.Tfrac.first.master.variant} and
 \ref{thm.Tfrac.first.pq.variant})
using slightly different statistics.

We will also have three versions of a ``second T-fraction''
(Theorems~\ref{thm.Tfrac.second}, \ref{thm.Tfrac.second.pq} and
\ref{thm.Tfrac.second.master})
that includes the counting of cycles,
at the expense of partially renouncing the counting of records;
as a corollary we obtain a continued fraction for D-cycles
(Corollary~\ref{cor.Tfrac.second.dcycle}).
The second T-fractions look less symmetrical than the first ones;
this defect seems to be inherent in including the counting of cycles
(just as in \cite{Sokal-Zeng_masterpoly}).
The proofs of both the first and second T-fractions
will be based on bijections from D-permutations to labeled Schr\"oder paths.

% \bigskip
% 
% {\bf WHAT ELSE TO SAY?????}

\bigskip

The plan of this paper is as follows:
In Section~\ref{sec.prelim} we recall some needed definitions and facts
concerning continued fractions, the Genocchi and median Genocchi numbers,
and permutation statistics.
In Section~\ref{sec.first} we state the three versions of the first T-fraction
and note some of their corollaries,
and in Section~\ref{sec.second} we do the same for the second T-fraction.
In Section~\ref{sec.prelimproofs} we recall how continued fractions
can be proven by bijection to labeled Dyck, Motzkin or Schr\"oder paths.
In Section~\ref{sec.proofs.1} we prove the first T-fraction
by a bijection that combines ideas of
Randrianarivony \cite{Randrianarivony_97}
and Foata--Zeilberger \cite{Foata_90}
together with some new ingredients.
In Section~\ref{sec.proofs.2} we prove the second T-fraction
by a variant bijection that analogously combines ideas of
Randrianarivony \cite{Randrianarivony_97}
and Biane \cite{Biane_93}.\footnote{
   Randrianarivony's proof \cite[Section~6]{Randrianarivony_97}
   is already motivated by ideas of Foata and Zeilberger \cite{Foata_90},
   as Randrianarivony himself points out
   \cite[pp.~78, 88]{Randrianarivony_97}.
   By contrast, our work here seems to be the first to apply a
   Biane-like \cite{Biane_93} bijection
   to D-permutations, D-semiderangements or D-derangements.
}
We conclude (Section~\ref{sec.final})
with some brief remarks on the relation of our work to
\cite{Sokal-Zeng_masterpoly}.

%% {\bf Also Appendix?????}

\section{Preliminaries}   \label{sec.prelim}

In this section we begin by explaining briefly
the types of continued fractions that will be employed
(Section~\ref{subsec.contfrac})
and recalling some contraction and transformation formulae
(Sections~\ref{subsec.contraction} and \ref{subsec.transformation}).
We then define the Euler, Genocchi and median Genocchi numbers
and recall their continued fractions
(Sections~\ref{subsec.euler}--\ref{subsec.median}).
Finally, we define some permutation statistics
that will play a central role in this paper
(Sections~\ref{subsec.statistics.1} and \ref{subsec.statistics.2}).
%% {\bf Also briefly discuss total positivity (Section~\ref{subsec.TP})?????}

\subsection{Classical continued fractions: S-fractions, J-fractions
   and T-fractions}
   \label{subsec.contfrac}

If $(a_n)_{n \ge 0}$ is a sequence of combinatorial numbers or polynomials
with $a_0 = 1$, it is often fruitful to seek to express its
ordinary generating function as a continued fraction.
The most commonly studied types of continued fractions
are Stieltjes-type (S-fractions),
\be
   \sum_{n=0}^\infty a_n t^n
   \;=\;
   \cfrac{1}{1 - \cfrac{\alpha_1 t}{1 - \cfrac{\alpha_2 t}{1 - \cdots}}}
   \label{def.Stype}
   \;\;,
\ee
and Jacobi-type (J-fractions),
\be
   \sum_{n=0}^\infty a_n t^n
   \;=\;
   \cfrac{1}{1 - \gamma_0 t - \cfrac{\beta_1 t^2}{1 - \gamma_1 t - \cfrac{\beta_2 t^2}{1 - \cdots}}}
   \label{def.Jtype}
   \;\;.
\ee
A less commonly studied type of continued fraction is 
the Thron-type (T-fraction):
\be
   \sum_{n=0}^\infty a_n t^n
   \;=\;
   \cfrac{1}{1 - \delta_1 t - \cfrac{\alpha_1 t}{1 - \delta_2 t - \cfrac{\alpha_2 t}{1 - \cdots}}}
   \label{def.Ttype}
   \;\;.
\ee
(Both sides of all these expressions are to be interpreted as
formal power series in the indeterminate $t$.)
This line of investigation goes back at least to
Euler \cite{Euler_1760,Euler_1788},
but it gained impetus following Flajolet's \cite{Flajolet_80}
seminal discovery that any S-fraction (resp.\ J-fraction)
can be interpreted combinatorially as a generating function
for Dyck (resp.\ Motzkin) paths with suitable weights for each rise and fall
(resp.\ each rise, fall and level step).
More recently, several authors
\cite{Fusy_15,Oste_15,Josuat-Verges_18,Sokal_totalpos,Elvey-Price-Sokal_wardpoly}
have found a similar combinatorial interpretation
of the general T-fraction:
namely, as a generating function for Schr\"oder paths with suitable
weights for each rise, fall and long level step.
These interpretations will be reviewed in
Section~\ref{subsec.prelimproofs.1} below.
% There are now literally dozens of sequences $(a_n)_{n \ge 0}$
% of combinatorial numbers or polynomials for which
% a continued-fraction expansion of the type
% \reff{def.Stype}, \reff{def.Jtype} or \reff{def.Ttype}
% is explicitly known.
%% {\bf Give some references??  Is there any review article???}

%% {\bf Mention branched continued fractions in a separate subsection?????}

\subsection{Contraction formulae}
   \label{subsec.contraction}

The formulae for even and odd contraction of an S-fraction
to an equivalent J-fraction are well known:
see e.g.~\cite[Lemmas~1 and 2]{Dumont_94b} \cite[Lemma~1]{Dumont_95}
for very simple algebraic proofs,
and see \cite[pp.~V-31--V-32]{Viennot_83}
for enlightening combinatorial proofs
%% based on Flajolet's \cite{Flajolet_80} combinatorial interpretation
%% of S-fractions (resp.\ J-fractions)
%% as generating functions for Dyck (resp.\ Motzkin) paths
%% with height-dependent weights.
based on grouping pairs of steps in a Dyck path.
Here we will need an extension of these formulae
to suitable subclasses of T-fractions \cite{Sokal_totalpos}:

\begin{proposition}[Even contraction for T-fractions with
   $\delta_2 = \delta_4 = \delta_6 = \ldots = 0$]
   \label{prop.contraction_even.Ttype}
We have
\begin{eqnarray}
    & &  \hspace*{-1.5cm}
   \cfrac{1}{1 - \delta_1 t - \cfrac{\alpha_1 t}{1 - \cfrac{\alpha_2 t}{1 - \delta_3 t -  \cfrac{\alpha_3 t}{1- \cdots}}}}
  \;\;=\;  \nonumber \\[2mm]
  & &  \hspace*{-7mm}
   \cfrac{1}{1 - (\alpha_1 + \delta_1) t - \cfrac{\alpha_1 \alpha_2 t^2}{1 - (\alpha_2 + \alpha_3 + \delta_3) t - \cfrac{\alpha_3 \alpha_4 t^2}{1 - (\alpha_4 + \alpha_5 + \delta_5) t - \cfrac{\alpha_5 \alpha_6 t^2}{1- \cdots}}}}
 \;\;.
   \qquad
 \label{eq.contraction_even.Ttype}
\end{eqnarray}
%% or in other words $T(t;\balpha,\bdelta_{\rm odd}) = J(t;\bbeta,\bgamma)$
%% where
That is, the T-fraction on the left-hand side of
\reff{eq.contraction_even.Ttype}
equals the J-fraction with coefficients
\begin{subeqnarray}
   \gamma_0  & = &  \alpha_1 + \delta_1   \\
   \gamma_n  & = &  \alpha_{2n} + \alpha_{2n+1} + \delta_{2n+1}
                            \quad\hbox{for $n \ge 1$} \\
   \beta_n  & = &  \alpha_{2n-1} \alpha_{2n}
 \label{eq.contraction_even.Ttype.coeffs}
\end{subeqnarray}
Here \reff{eq.contraction_even.Ttype}/\reff{eq.contraction_even.Ttype.coeffs}
holds as an identity in $\Z[\balpha,\bdelta_{\rm odd}][[t]]$,
where $\bdelta_{\rm odd} = (\delta_1,0,\delta_3,0,\ldots)$.
%
% {\bf I should also check whether this was known in the classical literature,
%   e.g.\ Perron or Jones--Thron.}
\end{proposition}

\begin{proposition}[Odd contraction for T-fractions with $\delta_1=\delta_3=\delta_5 = \ldots = 0$]
   \label{prop.contraction_odd.Ttype}
We have
\begin{eqnarray}
   & &  \hspace*{-15mm}
   \cfrac{1}{1 - \cfrac{\alpha_1 t}{1 - \delta_2 t - \cfrac{\alpha_2 t}{1 -  \cfrac{\alpha_3 t}{1- \delta_4 t - \cfrac{\alpha_4 t}{1 - \cdots}}}}}
   \;\;=\;  \nonumber \\[3mm]
   & &
   1 \:+\:
   \cfrac{\alpha_1 t}{1 - (\alpha_1 + \alpha_2 + \delta_2) t - \cfrac{\alpha_2 \alpha_3 t^2}{1 - (\alpha_3 + \alpha_4 + \delta_4) t - \cfrac{\alpha_4 \alpha_5 t^2}{1 - \cdots}}}
 \;\;.
   \qquad
 \label{eq.contraction_odd.Ttype}
\end{eqnarray}
%% or in other words
%% $T(t;\balpha,\bdelta_{\rm even}) = 1 \,+\, \alpha_1 t \, J(t;\bbeta,\bgamma)$
%% where
That is, the T-fraction on the left-hand side of
\reff{eq.contraction_odd.Ttype}
equals $1$ plus $\alpha_1 t$ times the J-fraction with coefficients
\begin{subeqnarray}
   \gamma_n  & = &  \alpha_{2n+1} + \alpha_{2n+2} + \delta_{2n+2} \\
   \beta_n  & = &  \alpha_{2n} \alpha_{2n+1}
 \label{eq.contraction_odd.coeffs.Ttype}
\end{subeqnarray}
Here \reff{eq.contraction_odd.Ttype}/\reff{eq.contraction_odd.coeffs.Ttype}
holds as an identity in $\Z[\balpha,\bdelta_{\rm even}][[t]]$,
where $\bdelta_{\rm even} = (0,\delta_2,0,\delta_4,\ldots)$.
\end{proposition}

Both the algebraic and the combinatorial proofs of the contraction formulae
for S-fractions can be easily generalized \cite{Sokal_totalpos}
to prove Propositions~\ref{prop.contraction_even.Ttype}
and \ref{prop.contraction_odd.Ttype}.

One consequence of Proposition~\ref{prop.contraction_even.Ttype}
is that if two T-fractions with only odd deltas ---
say, one with coefficients $(\balpha,\bdelta_{\rm odd})$
and the other with coefficients $(\balpha',\bdelta'_{\rm odd})$ ---
give rise by contraction to the same J-fraction $(\bbeta,\bgamma)$,
then they must be equal.
In some cases this principle can be used to transform
a T-fraction into an S-fraction (that is, $\bdelta'_{\rm odd} = \bzero$):
see, for instance, Corollary~\ref{cor.Sfrac.cyc} below.

By combining Propositions~\ref{prop.contraction_even.Ttype}
and \ref{prop.contraction_odd.Ttype}, we obtain:

\begin{corollary}[Combining odd and even contraction]
   \label{cor.contraction.Ttype}
If
\be
   \sum_{n=0}^\infty a_n t^n
   \;=\;
   \cfrac{1}{1 - \cfrac{\alpha_1 t}{1 - \delta_2 t - \cfrac{\alpha_2 t}{1 -  \cfrac{\alpha_3 t}{1- \delta_4 t - \cfrac{\alpha_4 t}{1 - \cdots}}}}}
   \;,
\ee
then
\be
   \sum_{n=0}^\infty a_{n+1} t^n
   \;=\;
   \cfrac{a_1}{1 - \delta'_1 t - \cfrac{\alpha'_1 t}{1 - \cfrac{\alpha'_2 t}{1 - \delta'_3 t -  \cfrac{\alpha'_3 t}{1- \cdots}}}}
\ee
whenever
\begin{subeqnarray}
   \alpha_1 + \alpha_2 + \delta_2  & = &  \alpha'_1 + \delta'_1 \\
   \alpha_{2n+1} + \alpha_{2n+2} + \delta_{2n+2}
         & = & \alpha'_{2n} + \alpha'_{2n+1} + \delta'_{2n+1}
                      \quad\hbox{for $n \ge 1$}  \\
   \alpha_{2n} \alpha_{2n+1}  & = &  \alpha'_{2n-1} \alpha'_{2n}
 \label{eq.cor.contraction.Ttype}
\end{subeqnarray}
\end{corollary}

\subsection{Transformation formula}
   \label{subsec.transformation}

We now prove a useful transformation formula for T-fractions.
First, a lemma:

\begin{lemma}
   \label{lemma.augment.Tfrac}
Let $R$ be a commutative ring, and let $f(t), g(t) \in R[[t]]$.
Then
\be
   1  \:+\: {t f(t) \over 1 - t f(t) - t g(t)}
   \;=\;
   \cfrac{1}{1 - \cfrac{t f(t)}{1 - t g(t)}}
 \label{eq.lemma.augment.Tfrac}
\ee
as an identity in $R[[t]]$.
\end{lemma}

\proof
Trivial: both sides equal $\dfrac{1 - tg(t)}{1 - tf(t) - tg(t)}$.
\qed

\begin{proposition}[Augmentation/restriction of T-fraction]
\hfill\break
   \label{prop.augment.Tfrac}
\be
   1 \:+\: \cfrac{\delta_1 t}{1 - \delta_1 t - \cfrac{\alpha_1 t}{1 - \delta_2 t - \cfrac{\alpha_2 t}{1 - \ldots}}}
   \;\;=\;\;
   \cfrac{1}{1 - \cfrac{\delta_1 t}{1 - \cfrac{\alpha_1 t}{1 - \delta_2 t - \cfrac{\alpha_2 t}{1 - \ldots}}}}
   \;\:.
   \quad
   \label{eq.prop.augment.Tfrac}
\ee
\end{proposition}

\proof
Use the lemma with $f(t) = \delta_1$ and
$g(t) = \displaystyle \cfrac{\alpha_1}{1 - \delta_2 t - \cfrac{\alpha_2 t}{1 - \ldots}}\;\,$.
\qed

Reading the identity \reff{eq.prop.augment.Tfrac} from left to right,
it says that if the ogf of a sequence $\ba = (a_0,a_1,a_2,\ldots)$
with $a_0 = 1$
is given by a T-fraction with coefficients $\balpha$ and $\bdelta$,
then the ogf of the ``augmented'' sequence
$\ba' = (1, \delta_1 a_0, \delta_1 a_1, \delta_1 a_2, \ldots)$
is given by a T-fraction with coefficients $\balpha'$ and $\bdelta'$, where
\begin{subeqnarray}
   \alpha'_1 & = & \delta_1  \\
   \alpha'_n & = & \alpha_{n-1} \quad\hbox{for $n \ge 2$} \\
   \delta'_1 = \delta'_2 & = & 0  \\
   \delta'_n & = & \delta_{n-1} \quad\hbox{for $n \ge 3$}
\end{subeqnarray}
In particular, if $\delta_2 = \delta_3 = \ldots = 0$,
then the T-fraction on the right-hand side is an S-fraction.
Of course, this transformation gives something interesting
only when $\delta_1 \neq 0$.

Alternatively, reading the identity \reff{eq.prop.augment.Tfrac}
from right to left,
it says that if the ogf of a sequence $\ba' = (a'_0,a'_1,a'_2,\ldots)$
with $a'_0 = 1$
is given by a T-fraction with coefficients $\balpha'$ and $\bdelta'$,
where $\delta'_1 = \delta'_2 = 0$ (and of course $\alpha'_1 \neq 0$),
then the ogf of the ``restricted'' sequence
$\ba = (a_1/\alpha'_1, a_2/\alpha'_1,\ldots)$
is given by a T-fraction with coefficients $\balpha$ and $\bdelta$, where
\begin{subeqnarray}
   \alpha_n & = & \alpha'_{n+1} \\
   \delta_1 & = & \alpha'_1 \\
   \delta_n & = & \delta'_{n+1} \quad\hbox{for $n \ge 2$}
\end{subeqnarray}

% \begin{proposition}[Transformation of special T-fraction with $\delta_2=\delta_4=\delta_6=\ldots=0$
% to S-fraction]
% \hfill\break
%    \label{prop.transform.special.Tfrac}
% \be
%    \cfrac{1}{1 - \delta_1 t - \cfrac{\alpha_1 t}{1 - \cfrac{\alpha_2 t}{1 - \cfrac{\alpha_3 t}{1-\ldots}}}}
%    \;\;=\;\;
%    \cfrac{1}{1 - \cfrac{\alpha_1' t}{1 - \cfrac{\alpha_2' t}{1 - \cfrac{\alpha_3' t}{1 - \ldots}}}}
%    \;\:
%    \quad
%    \label{eq.prop.transform.special.Tfrac}
% \ee
% where 
% \begin{subeqnarray}
%    \alpha_1' & = & \alpha_1 +\delta_1 \\
%    \alpha_{2n-1}' & = & \alpha_{2n-1} + \alpha_{2n-2} - \alpha_{2n-2}' \quad\hbox{for $n \ge 2$} \\
%    \alpha_{2n}' & = & \dfrac{\alpha_{2n-1}\alpha_{2n}}{\alpha_{2n-1}'}
% \end{subeqnarray} 
% \end{proposition}
% 
% {\bf What is the underlying ring?????????????}
% 
% \begin{proof} {\bf Use even contraction formula for T-fractions with even deltas $0$ and even contraction formula for S-fractions!!!}
% \end{proof}

\subsection{Euler numbers and augmented Euler numbers}
   \label{subsec.euler}

The \textbfit{Euler numbers} \cite[A000111]{OEIS}
\be
   (E_n)_{n \ge 0}
   \;=\;
   1, 1, 1, 2, 5, 16, 61, 272, 1385, 7936, 50521, 353792, 2702765, \ldots
\ee
are positive integers defined by the exponential generating function
\be
   \sec t + \tan t  \;=\;  \sum_{n=0}^\infty E_n \, {t^n \over n!}
   \;.
 \label{def.euler}
\ee
The $(E_{2n})_{n \ge 0}$ are also called \textbfit{secant numbers}
\cite[A000364]{OEIS},
and the $(E_{2n+1})_{n \ge 0}$ are called \textbfit{tangent numbers}
\cite[A000182]{OEIS}.
Their ordinary generating functions have the S-fraction representations
\be
   \sum_{n=0}^\infty E_{2n} \, t^n
   \;=\;
   \cfrac{1}{1 - \cfrac{1^2 t}{1 - \cfrac{2^2 t}{1 -  \cfrac{3^2 t}{1- \cdots}}}}
 \label{eq.cfrac.secant}
\ee
with coefficients $\alpha_n = n^2$,
and
\be
   \sum_{n=0}^\infty E_{2n+1} \, t^n
   \;=\;
   \cfrac{1}{1 - \cfrac{1\cdot2 t}{1 - \cfrac{2\cdot3 t}{1 -  \cfrac{3\cdot4 t}{1- \cdots}}}}
 \label{eq.cfrac.tangent}
\ee
with coefficients $\alpha_n = n(n+1)$.\footnote{
   These formulae were found by Stieltjes \cite[p.~H9]{Stieltjes_1889} in 1889
   %% (the first one explicitly, the second implicitly)
   and by Rogers \cite[p.~77]{Rogers_07} in 1907.
   They were given beautiful combinatorial proofs by
   Flajolet \cite{Flajolet_80} in 1980.
}
See \cite{Foata_71,Viennot_81,Kuznetsov_94,Stanley_10,Sokal_eulernumbers}
for further information on the Euler numbers,
including their combinatorial interpretations and moment representations.

We define the \textbfit{augmented Euler numbers} \cite[A065619]{OEIS} by
\be
   E^\sharp_n  \;\eqdef\; (n+1) \, E_n
   \;.
 \label{def.augeuler}
\ee
They have the exponential generating function
\be
   t \, (\sec t + \tan t)
   \;=\;
   \sum_{n=0}^\infty E^\sharp_n \, {t^{n+1} \over (n+1)!}
   \;.
 \label{eq.egf_augeuler}
\ee
The $(E^\sharp_{2n})_{n \ge 0}$ are the
\textbfit{augmented secant numbers} \cite[A009843]{OEIS},
and the $(E^\sharp_{2n+1})_{n \ge 0}$ are the
\textbfit{augmented tangent numbers} \cite[A009752]{OEIS}.

%% {\bf Give continued fractions?????}

\subsection{Genocchi numbers}
   \label{subsec.genocchi}

The \textbfit{Genocchi numbers} \cite[A110501]{OEIS}\footnote{
   Our $g_n$ is usually written by combinatorialists as $G_{2n+2}$.
   However, many texts ---
   particularly older ones,
   or those in the analysis and special-functions literature ---
   define the ({\em signed}\/) Genocchi numbers $(G_n)_{n \ge 0}$
   by \cite[24.15.1]{NIST}
   $$
      {2t \over e^t + 1}
      \;=\;
      \sum_{n=1}^\infty G_n \, {t^n \over n!}
      \;,
   $$
   which leads to $G_0 = 0$, $G_1 = 1$, $G_n = 0$ for odd $n \ge 3$,
   and $G_{2n+2} = (-1)^{n+1} g_n$.

   {\bf Warning:} Our $g_n$ is denoted $g_{n+1}$ by
   Lazar and Wachs \cite{Lazar_22,Lazar_20}
   and Eu, Fu, Lai and Lo \cite{Eu_21}.
 \label{footnote.Gn}
}
% \cite[A110501]{OEIS}\footnote{
%    The Genocchi numbers appear already in Euler's book
%    {\em Foundations of Differential Calculus,
%     with Applications to Finite Analysis and Series}\/,
%    first published in 1755 \cite[paragraphs~181 and 182]{Euler_1755};
%    this book is E212 in Enestr\"om's \cite{Enestrom_13} catalogue.
%    These numbers were revisited by Genocchi \cite{Genocchi_1852} in 1852.
%    The beautiful survey article of Viennot \cite{Viennot_81}
%    contains a wealth of useful information.
% }
\be
   (g_n)_{n \ge 0}
   \;=\;
   1, 1, 3, 17, 155, 2073, 38227, 929569, 28820619, 1109652905,
   %% 51943281731,
   \ldots
\ee
%% $(g_n)_{n \ge 0}$ [cf.\ \reff{eq.genocchi}]
are odd positive integers \cite{Lucas_1877,Barsky_81,Han_18}
\cite[pp.~217--218]{Foata_08}
defined by the exponential generating function
\be
   t \, \tan(t/2)
   \;=\;
   \sum_{n=0}^\infty g_n \, {t^{2n+2} \over (2n+2)!}
   \;.
\ee
They are therefore rescaled versions of the augmented tangent numbers:
\be
   g_n  \;=\; 4^{-n} \, (n+1) \, E_{2n+1}
        \;=\; 2^{-(2n+1)} \, E^\sharp_{2n+1}
   \;.
 \label{eq.gn.augmented_tangent}
\ee
% where
% \be
%    (E_{2n+1})_{n \ge 0}
%    \;=\;
%    1, 2, 16, 272, 7936, 353792, 22368256, 1903757312, 209865342976, \ldots
% \ee
% are the {\em tangent numbers}\/ \cite[A000182]{OEIS} defined by
% \be
%    \tan t
%    \;=\;
%    \sum_{n=0}^\infty E_{2n+1} \, {t^{2n+1} \over (2n+1)!}
%    \;.
% \ee

The ordinary generating function of the Genocchi numbers
has a classical S-fraction expansion
\cite[eq.~(7.5)]{Viennot_81}
\cite[p.~V-9]{Viennot_83}
\cite[eqns.~(1.4) and (3.9)]{Dumont_94b}
\be
   \sum_{n=0}^\infty g_n \, t^n
   \;=\;
   \cfrac{1}{1 - \cfrac{1 \cdot 1 t}{1 - \cfrac{1 \cdot 2t}{1 - \cfrac{2 \cdot 2t}{1- \cfrac{2 \cdot 3 t}{1-\cdots}}}}}
  \label{eq.genocchi.Sfrac}
\ee
with coefficients
\be
   \alpha_{2k-1} \:=\: k^2 \,,\quad
   \alpha_{2k} \:=\: k(k+1)  \;.
  \label{eq.genocchi.Sfrac.weights}
\ee
%% Otherwise put, the alphas are products of successive pairs
%% of the ``pre-alphas'' $\balpha^{\rm pre} = 1,1,2,2,3,3,4,4,\ldots\:$.
% It then follows from the contraction formulae
% (see Section~\ref{subsec.contraction})
% that the once-shifted Genocchi numbers $(g_{n+1})_{n \ge 0}$
% have a J-fraction
% \be
%    \sum_{n=0}^\infty g_{n+1} \, t^n
%    \;=\;
%    \cfrac{1}{1 - 3 t - \cfrac{8 t^2}{1 - 10 t - \cfrac{54 t^2}{1 - 21 t -\cfrac{192 t^2}{1 - \cdots}}}}
%    \label{eq.genocchi.Jfrac}
% \ee
% with coefficients
% \be
%    \gamma_n \:=\: (n+1)(2n+3) \,,\quad
%    \beta_n \:=\: n (n+1)^3
% \ee
% and also a T-fraction
It then follows from Proposition~\ref{prop.augment.Tfrac}
that the once-shifted Genocchi numbers $(g_{n+1})_{n \ge 0}$
have a T-fraction
\be
   \sum_{n=0}^\infty g_{n+1} \, t^n
   \;=\;
   \cfrac{1}{1 - t - \cfrac{1 \cdot 2 t}{1 - \cfrac{2 \cdot 2 t}{1 - \cfrac{2 \cdot 3 t}{1- \cfrac{3 \cdot 3 t}{1- \cdots}}}}}
   \label{eq.genocchi.Tfrac}
\ee
with coefficients
\be
   \alpha_{2k-1} \:=\: k(k+1) \,,\quad
   \alpha_{2k} \:=\: (k+1)^2 \,,\quad
   \delta_1 \:=\: 1 \,,\quad
   \delta_n \:=\: 0 \;\hbox{for $n \ge 2$}
   \;.
   \label{eq.genocchi.Tfrac.weights}
\ee
%{\bf Do they have a computable S-fraction also (with rational coefficients)???}
% In my notebook contfrac.nb (on the Mac),
%   Section "Genocchi numbers and Dumont permutations",
%   subsection "Genocchi numbers (OEIS A110501) and shifted Genocchi numbers",
%   I compute the S-fraction.  It seems to be related to the harmonic numbers,
%   but I don't have any explicit guess.

\medskip

{\bf Remark.}
Some generalizations of
\reff{eq.genocchi.Sfrac}/\reff{eq.genocchi.Sfrac.weights},
incorporating additional parameters,
are known in a variety of algebraic or combinatorial models:
see \cite[Section~6]{Dumont_86}
\cite{Dumont_95b}
\cite[eq.~(3.3) and Corollaire~8]{Zeng_96}
\cite[Th\'eor\`eme~3 and Proposition~13]{Randrianarivony_96}
\cite[Proposition~10]{Randrianarivony_96b}
\cite[Th\'eor\`eme~1.2]{Randrianarivony_97}
\cite[eq.~(5) and Th\'eor\`eme~2]{Han_99a}.
See also \cite[Corollaire~3]{Han_99a}
for a generalization of the T-fraction
\reff{eq.genocchi.Tfrac}/\reff{eq.genocchi.Tfrac.weights}.\footnote{
   There is a typographical error in \cite[eq.~(12)]{Han_99a}:
   on the left-hand side, $t^n$ should be $t^{n-1}$.
}
\myendremark

\subsection{Median Genocchi numbers}
   \label{subsec.median}

An alternative way of defining the Genocchi numbers is
by the \textbfit{Seidel recurrence}
%% https://tex.stackexchange.com/questions/308697/braces-around-both-sides-of-cases-with-text
%% https://tex.stackexchange.com/questions/240868/how-can-i-write-cases-with-grouping-braces-on-both-left-and-right-side-and-spli
\be
   s_{n,k}
   \;=\;
   \left\{\!\!
      \begin{array}{cl}
         \sum\limits_{j=0}^k s_{n-1,j}  & \text{if $n$ is even} \\
                                        &                       \\[-2mm]
         \sum\limits_{j=k}^{\lfloor n/2 \rfloor} s_{n-1,j}
                                        & \text{if $n$ is odd}
      \end{array}
   \!\!\right\}
   \quad
   \hbox{for $n \ge 1$ and $0 \le k \le \lfloor n/2 \rfloor$}
\ee
with $s_{0,0} = 1$
and $s_{n,k} = 0$ for $k > \lfloor n/2 \rfloor$.
The \textbfit{Seidel triangle} \cite[A014781]{OEIS} begins
\medskip
%
% Made in continued_fraction.nb using my program maketable2TRI
%
\begin{table}[H]
\centering
\small
\begin{tabular}{c|rrrrrr|r}
$n \setminus k$  & 0 & 1 & 2 & 3 & 4 & 5 & Row sums \\
\hline
0 & 1 &  &  &  &  &  & 1  \\
1 & 1 &  &  &  &  &  & 1  \\
2 & 1 & 1 &  &  &  &  & 2  \\
3 & 2 & 1 &  &  &  &  & 3  \\
4 & 2 & 3 & 3 &  &  &  & 8  \\
5 & 8 & 6 & 3 &  &  &  & 17  \\
6 & 8 & 14 & 17 & 17 &  &  & 56  \\
7 & 56 & 48 & 34 & 17 &  &  & 155  \\
8 & 56 & 104 & 138 & 155 & 155 &  & 608  \\
9 & 608 & 552 & 448 & 310 & 155 &  & 2073  \\
10 & 608 & 1160 & 1608 & 1918 & 2073 & 2073 & 9440  \\
11 & 9440 & 8832 & 7672 & 6064 & 4146 & 2073 & 38227  \\
\end{tabular}
\end{table}
\medskip
\noindent
It can be shown \cite{Dumont_80b}
that the Genocchi numbers occur on the ``diagonals'' and ``subdiagonals''
of the Seidel triangle,
\be
   g_n  \;=\; s_{2n,n}  \;=\;  s_{2n,n-1}   \;=\;  s_{2n+1,n}
\ee
(where we set $s_{0,-1} = 1$);
they are also the odd row sums:
\be
   \sum_{k=0}^n s_{2n+1,k}  \;=\;  g_{n+1}
   \;.
\ee
The \textbfit{median Genocchi numbers}
(or \textbfit{Genocchi medians} for short) \cite[A005439]{OEIS}\footnote{
   Our $h_n$ is usually written by combinatorialists as $H_{2n+1}$.

   {\bf Warning:}  Lazar and Wachs' \cite{Lazar_22,Lazar_20} $h_n$
   equals our $h_{n+1}$.
   Pan and Zeng's \cite{Pan_21} $h_n$ is our $h_{n+1}$ divided by $2^n$.
 \label{footnote.Hn}
}
are defined by the zeroth column of the Seidel triangle,
\be
   h_n  \;\eqdef\;  s_{2n,0}  \;=\;  s_{2n-1,0}
 \label{def.genocchimedians}
\ee
(where we set $s_{-1,0} = 1$);
they are also the even row sums:
\be
   \sum_{k=0}^n s_{2n,k}  \;=\;  h_{n+1}
   \;.
\ee
The Genocchi medians can alternatively be defined by
\cite[p.~63]{Han_99b}
%% \cite{?????}
\be
   h_n  \;=\;  \sum_{i=0}^{n-1} (-1)^i \, \binom{n}{2i+1} \, g_{n-1-i}
   \;.
 \label{eq.hn.binomgn}
\ee
They are also rescaled versions of the
1-binomial transform of the augmented secant numbers
\cite[Corollary~3]{Dumont_94b}:
\be
   h_n  \;=\; 4^{-n} \sum_{k=0}^n \binom{n}{k} \, (2k+1) E_{2k}
        \;=\; 4^{-n} \sum_{k=0}^n \binom{n}{k} \, E^\sharp_{2k}
 \label{eq.hn.augmented_secant}
\ee
[compare to \reff{eq.gn.augmented_tangent}].

The Genocchi medians
\be
   (h_n)_{n \ge 0}
   \;=\;
   1, 1, 2, 8, 56, 608, 9440, 198272, 5410688, 186043904, 7867739648,
   \ldots
\ee
%% $(h_n)_{n \ge 0}$ [cf.\ \reff{eq.median}]
do not have any known exponential generating function.
However, their ordinary generating function has a nice
classical S-fraction expansion
\cite[eq.~(9.7)]{Viennot_81}
\cite[p.~V-15]{Viennot_83}
\cite[eqns.~(1.5) and (3.8)]{Dumont_94b}:
\be
   \sum_{n=0}^\infty h_n \, t^n
   \;=\;
   \cfrac{1}{1 - \cfrac{1 t}{1 - \cfrac{1 t}{1 - \cfrac{4t}{1- \cfrac{4 t}{1-\cdots}}}}}
  \label{eq.mediangenocchi.Sfrac}
\ee
with coefficients
\be
   \alpha_{2k-1} \:=\: \alpha_{2k} \:=\: k^2
   \;.
  \label{eq.mediangenocchi.Sfrac.weights}
\ee
It then follows from Corollary~\ref{cor.contraction.Ttype}
that the once-shifted median Genocchi numbers $(h_{n+1})_{n \ge 0}$
have an S-fraction
\be
   \sum_{n=0}^\infty h_{n+1} \, t^n
   \;=\;
   \cfrac{1}{1 - \cfrac{2 t}{1 - \cfrac{2 t}{1 - \cfrac{6t}{1- \cfrac{6 t}{1-\cdots}}}}}
  \label{eq.mediangenocchi.Sfrac2}
\ee
with coefficients
\be
   \alpha_{2k-1} \:=\: \alpha_{2k} \:=\: k(k+1)
   \;.
  \label{eq.mediangenocchi.Sfrac2.weights}
\ee
Moreover, from Proposition~\ref{prop.augment.Tfrac} they also a T-fraction
\be
   \sum_{n=0}^\infty h_{n+1} \, t^n
   \;=\;
   \cfrac{1}{1 - t - \cfrac{1 t}{1 - \cfrac{4 t}{1 - \cfrac{4 t}{1- \cfrac{9 t}{1- \cdots}}}}}
   \label{eq.mediangenocchi.Tfrac}
\ee
with coefficients
\be
   \alpha_{2k-1} \:=\: k^2 \,,\quad
   \alpha_{2k} \:=\: (k+1)^2 \,,\quad
   \delta_1 \:=\: 1 \,,\quad
   \delta_n \:=\: 0 \;\hbox{for $n \ge 2$}
   \;.
   \label{eq.mediangenocchi.Tfrac.weights}
\ee

Finally, let us define, for future reference,
the sequence $(h^\flat_{n+1})_{n \ge 0}$
corresponding to the S-fraction underlying \reff{eq.mediangenocchi.Tfrac}:
\be
   \sum_{n=0}^\infty h^\flat_{n+1} \, t^n
   \;=\;
   \cfrac{1}{1 - \cfrac{1 t}{1 - \cfrac{4 t}{1 - \cfrac{4 t}{1- \cfrac{9 t}{1- \cdots}}}}}
   \label{eq.mediangenocchiflat.Sfrac}
\ee
with coefficients
\be
   \alpha_{2k-1} \:=\: k^2 \,,\quad
   \alpha_{2k} \:=\: (k+1)^2
   \;.
   \label{eq.mediangenocchiflat.Sfrac.weights}
\ee
This sequence begins
\be
   (h^\flat_{n+1})_{n \ge 0}
   \;=\;
   1, 1, 5, 41, 493, 8161, 178469, 4998905, 174914077, 7487810257,
     %% 385307632469,
     \ldots
\ee
and cannot be found, at present, in \cite{OEIS}.
In Section~\ref{subsec.first.statements}
we will give its combinatorial interpretation.
Using Proposition~\ref{prop.augment.Tfrac}, we also get a T-fraction
for the once-shifted sequence $(h^\flat_{n+2})_{n \ge 0}$:
\be
   \sum_{n=0}^\infty h^\flat_{n+2} \, t^n
   \;=\;
   \cfrac{1}{1 - t - \cfrac{4 t}{1 - \cfrac{4 t}{1- \cfrac{9 t}{1- \cdots}}}}
   \label{eq.mediangenocchiflat.Tfrac}
\ee
with coefficients
\be
   \alpha_{2k-1} \:=\: (k+1)^2 \,,\quad
   \alpha_{2k} \:=\: (k+1)^2 \,,\quad
   \delta_1 \:=\: 1 \,,\quad
   \delta_n \:=\: 0 \;\hbox{for $n \ge 2$}
   \;.
   \label{eq.mediangenocchiflat.Tfrac.weights}
\ee

%% {\bf What else to say???????????}

\medskip

{\bf Remark.}
Some generalizations of
\reff{eq.mediangenocchi.Sfrac}/\reff{eq.mediangenocchi.Sfrac.weights}
or \reff{eq.mediangenocchi.Sfrac2}/\reff{eq.mediangenocchi.Sfrac2.weights},
incorporating additional parameters,
are known in a variety of algebraic or combinatorial models:
see \cite[Section~6]{Dumont_86}
\cite{Dumont_95b}
\cite[eq.~(3.3) and Corollaire~8]{Zeng_96}
\cite[Th\'eor\`eme~3 and Proposition~13]{Randrianarivony_96}
\cite[Proposition~10 and Corollary~13]{Randrianarivony_96b}
\cite[Theorem~0.1]{Feigin_12}
\cite[Corollary~6]{Pan_21}.
\myendremark

\subsection{Permutation statistics: The record-and-cycle classification}
     \label{subsec.statistics.1}

We now define some permutation statistics that will play a central role
in what follows.

Given a permutation $\sigma \in \Sym_N$, an index $i \in [N]$ is called an
\begin{itemize}
   \item {\em excedance}\/ (exc) if $i < \sigma(i)$;
   \item {\em anti-excedance}\/ (aexc) if $i > \sigma(i)$;
   \item {\em fixed point}\/ (fix) if $i = \sigma(i)$.
\end{itemize}
Clearly every index $i$ belongs to exactly one of these three types;
we call this the \textbfit{excedance classification}.
We also say that $i$ is a {\em weak excedance}\/ if $i \le \sigma(i)$,
and a {\em weak anti-excedance}\/ if $i \ge \sigma(i)$.

A more refined classification is as follows:
an index $i \in [N]$ is called a
\begin{itemize}
   \item {\em cycle peak}\/ (cpeak) if $\sigma^{-1}(i) < i > \sigma(i)$;
   \item {\em cycle valley}\/ (cval) if $\sigma^{-1}(i) > i < \sigma(i)$;
   \item {\em cycle double rise}\/ (cdrise) if $\sigma^{-1}(i) < i < \sigma(i)$;
   \item {\em cycle double fall}\/ (cdfall) if $\sigma^{-1}(i) > i > \sigma(i)$;
   \item {\em fixed point}\/ (fix) if $\sigma^{-1}(i) = i = \sigma(i)$.
\end{itemize}
Clearly every index $i$ belongs to exactly one of these five types;
we refer to this classification as the \textbfit{cycle classification}.
Obviously, excedance = cycle valley or cycle double rise,
and anti-excedance = cycle peak or cycle double fall.

On the other hand, an index $i \in [N]$ is called a
\begin{itemize}
   \item {\em record}\/ (rec) (or {\em left-to-right maximum}\/)
         if $\sigma(j) < \sigma(i)$ for all $j < i$
      [note in particular that the indices 1 and $\sigma^{-1}(N)$
       are always records];
   \item {\em antirecord}\/ (arec) (or {\em right-to-left minimum}\/)
         if $\sigma(j) > \sigma(i)$ for all $j > i$
      [note in particular that the indices $N$ and $\sigma^{-1}(1)$
       are always antirecords];
   \item {\em exclusive record}\/ (erec) if it is a record and not also
         an antirecord;
   \item {\em exclusive antirecord}\/ (earec) if it is an antirecord
         and not also a record;
   \item {\em record-antirecord}\/ (rar) (or {\em pivot}\/)
      if it is both a record and an antirecord;
   \item {\em neither-record-antirecord}\/ (nrar) if it is neither a record
      nor an antirecord.
\end{itemize}
Every index $i$ thus belongs to exactly one of the latter four types;
we refer to this classification as the \textbfit{record classification}.

The record and cycle classifications of indices are related as follows:
\begin{quote}
\begin{itemize}
   \item[(a)]  Every record is a weak excedance,
      and every exclusive record is an excedance.
   \item[(b)]  Every antirecord is a weak anti-excedance,
      and every exclusive antirecord is an anti-excedance.
   \item[(c)]  Every record-antirecord is a fixed point.
\end{itemize}
\end{quote}
Therefore, by applying the record and cycle classifications simultaneously,
we obtain 10~disjoint categories \cite{Sokal-Zeng_masterpoly}:
\begin{itemize}
   \item ereccval:  exclusive records that are also cycle valleys;
   \item ereccdrise:  exclusive records that are also cycle double rises;
   \item eareccpeak:  exclusive antirecords that are also cycle peaks;
   \item eareccdfall:  exclusive antirecords that are also cycle double falls;
   \item rar:  record-antirecords (these are always fixed points);
   \item nrcpeak:  neither-record-antirecords that are also cycle peaks;
   \item nrcval:  neither-record-antirecords that are also cycle valleys;
   \item nrcdrise:  neither-record-antirecords that are also cycle double rises;
   \item nrcdfall:  neither-record-antirecords that are also cycle double falls;
   \item nrfix:  neither-record-antirecords that are also fixed points.
\end{itemize}
Clearly every index $i$ belongs to exactly one of these 10~types;
we call this the \textbfit{record-and-cycle classification}.

\bigskip

Now suppose that $\sigma$ is a D-permutation.
Then the cycle classification of a non-fixed-point index $i$
is equivalent to recording the parities of $\sinv(i)$ and $i$:
% 
% \begin{itemize}
%    \item {\em cycle peak}\/: $\sigma^{-1}(i) < i > \sigma(i)$
%             $\iff$ $\sinv(i)$ odd, $i$ even
%    \item {\em cycle valley}\/: $\sigma^{-1}(i) > i < \sigma(i)$
%             $\iff$ $\sinv(i)$ even, $i$ odd
%    \item {\em cycle double rise}\/: $\sigma^{-1}(i) < i < \sigma(i)$
%             $\iff$ $\sinv(i)$ odd, $i$ odd
%    \item {\em cycle double fall}\/:  $\sigma^{-1}(i) > i > \sigma(i)$
%             $\iff$ $\sinv(i)$ even, $i$ even
% \end{itemize}
% 
\be
   \begin{aligned}
   & \bullet\; \hbox{{\em cycle peak}\/: $\sigma^{-1}(i) < i > \sigma(i)$
            $\implies$ $\sinv(i)$ odd, $i$ even} \\[2mm]
   & \bullet\; \hbox{{\em cycle valley}\/: $\sigma^{-1}(i) > i < \sigma(i)$
            $\implies$ $\sinv(i)$ even, $i$ odd} \\[2mm]
   & \bullet\; \hbox{{\em cycle double rise}\/: $\sigma^{-1}(i) < i < \sigma(i)$
            $\implies$ $\sinv(i)$ odd, $i$ odd} \\[2mm]
   & \bullet\; \hbox{{\em cycle double fall}\/:  $\sigma^{-1}(i) > i > \sigma(i)$
            $\implies$ $\sinv(i)$ even, $i$ even} \hspace*{1.7cm} \\
   \end{aligned}
 \label{eq.parities.1}
\ee
For fixed points, by contrast,
we will need to explicitly record the parity of $i$,
by distinguishing even and odd fixed points:
% 
% \begin{itemize}
%   \item {\em even fixed point}\/ (evenfix):  $\sigma^{-1}(i) = i = \sigma(i)$
%          and $i$ is even
%   \item {\em odd fixed point}\/ (oddfix):  $\sigma^{-1}(i) = i = \sigma(i)$
%          and $i$ is odd
% \end{itemize}
% 
\be
   \begin{aligned}
   & \bullet\; \hbox{{\em even fixed point}\/ (evenfix):  $\sigma^{-1}(i) = i = \sigma(i)$ is even}  \\[2mm]
   & \bullet\; \hbox{{\em odd fixed point}\/ (oddfix):  $\sigma^{-1}(i) = i = \sigma(i)$ is odd}   \hspace*{3.7cm} \\
   \end{aligned}
 \label{eq.parities.2}
\ee
We therefore refine the record-and-cycle classification
by distinguishing even and odd fixed points:
\begin{itemize}
   \item evenrar:  even record-antirecords (these are always fixed points);
   \item oddrar:  odd record-antirecords (these are always fixed points);
   \item evennrfix:  even neither-record-antirecords that are also fixed points;
   \item oddnrfix:  odd neither-record-antirecords that are also fixed points.
\end{itemize}
This leads to the \textbfit{parity-refined record-and-cycle classification},
in which each index $i$ belongs to exactly one of 12 types.
More precisely, each even index $i$ belongs to exactly one of the 6 types
\begin{quote}
   eareccpeak, nrcpeak, eareccdfall, nrcdfall, evenrar, evennrfix,
\end{quote}
while each odd index $i$ belongs to exactly one of the 6 types
\begin{quote}
   ereccval, nrcval, ereccdrise, nrcdrise, oddrar, oddnrfix.
\end{quote}

\subsection{Permutation statistics: Crossings and nestings}
     \label{subsec.statistics.2}

We now define (following \cite{Sokal-Zeng_masterpoly})
some permutation statistics that count
\textbfit{crossings} and \textbfit{nestings}.

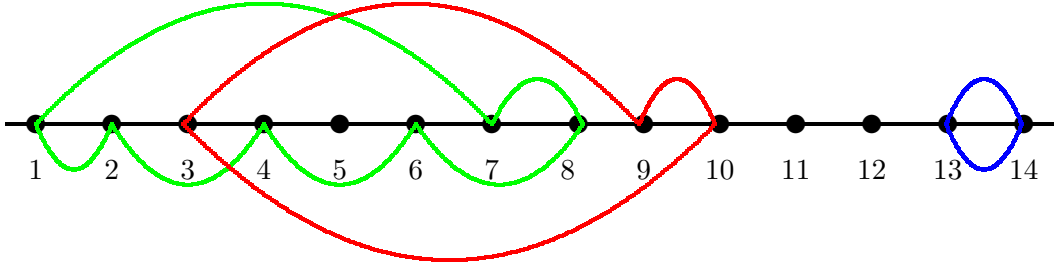
\begin{figure}[t]
\centering
\vspace*{4cm}
\begin{picture}(100,0)(140, -45)
\setlength{\unitlength}{2mm}
\linethickness{.5mm}
\put(-2,0){\line(1,0){69}}
\put(0,0){\circle*{1,3}}\put(0,0){\makebox(0,-6)[c]{\small 1}}
\put(5,0){\circle*{1,3}}\put(5,0){\makebox(0,-6)[c]{\small 2}}
\put(10,0){\circle*{1,3}}\put(10,0){\makebox(0,-6)[c]{\small 3}}
\put(15,0){\circle*{1,3}}\put(15,0){\makebox(0,-6)[c]{\small 4}}
\put(20,0){\circle*{1,3}}\put(20,0){\makebox(0,-6)[c]{\small 5}}
\put(25,0){\circle*{1,3}}\put(25,0){\makebox(0,-6)[c]{\small 6}}
\put(30,0){\circle*{1,3}}\put(30,0){\makebox(0,-6)[c]{\small 7}}
\put(35,0){ \circle*{1,3}}\put(35,0){\makebox(0,-6)[c]{\small 8}}
\put(40,0){\circle*{1,3}}\put(40,0){\makebox(0,-6)[c]{\small 9}}
\put(45,0){\circle*{1,3}}\put(45,0){\makebox(0,-6)[c]{\small 10}}
\put(50,0){\circle*{1,3}}\put(50,0){\makebox(0,-6)[c]{\small 11}}
\put(55,0){\circle*{1,3}}\put(55,0){\makebox(0,-6)[c]{\small 12}}
\put(60,0){\circle*{1,3}}\put(60,0){\makebox(0,-6)[c]{\small 13}}
\put(65,0){\circle*{1,3}}\put(65,0){\makebox(0,-6)[c]{\small 14}}
\green{\qbezier(0,0)(15,16)(30,0)
\qbezier(30,0)(33,6)(36,0)
\qbezier(36,0)(30.5,-8)(25,0)
\qbezier(25,0)(20,-8)(15,0)
\qbezier(15,0)(10,-8)(5,0)
\qbezier(5,0)(2.5,-6)(0,0)
}
\red{\qbezier(9,0)(24,16)(39,0)
\qbezier(39,0)(41.5,6)(44,0)
\qbezier(44,0)(26.5,-18)(9,0)
}
\blue{\qbezier(58.5,0)(61,6)(63.5,0)
\qbezier(63.5,0)(61,-6)(58.5,0)
}
\end{picture}
\caption{
   An example of a D-permutation
   $\sigma = 7\, 1\, 9\, 2\, 5\, 4\, 8\, 6\, 10\, 3\, 11\, 12\, 14\, 13\,
           = (1,7,8,6,4,2)\,(3,9,10)\,(5)\,(11)\,(12)\,(13,14) \in \dperm_{14}$.
 \label{fig.pictorial}
 \vspace*{7mm}
}
\end{figure}

First we associate to each permutation $\sigma \in \Sym_N$
a pictorial representation (Figure~\ref{fig.pictorial})
by placing vertices $1,2,\ldots,N$ along a horizontal axis
and then drawing an arc from $i$ to $\sigma(i)$
above (resp.\ below) the horizontal axis
in case $\sigma(i) > i$ [resp.\ $\sigma(i) < i$];
if $\sigma(i) = i$ we do not draw any arc.
Each vertex thus has either
out-degree = in-degree = 1 (if it is not a fixed point) or
out-degree = in-degree = 0 (if it is a fixed point).
Of course, the arrows on the arcs are redundant,
because the arrow on an arc above (resp.\ below) the axis
always points to the right (resp.\ left);
we therefore omit the arrows for simplicity.

We then say that a quadruplet $i < j < k < l$ forms an
\begin{itemize}
   \item {\em upper crossing}\/ (ucross) if $k = \sigma(i)$ and $l = \sigma(j)$;
   \item {\em lower crossing}\/ (lcross) if $i = \sigma(k)$ and $j = \sigma(l)$;
   \item {\em upper nesting}\/  (unest)  if $l = \sigma(i)$ and $k = \sigma(j)$;
   \item {\em lower nesting}\/  (lnest)  if $i = \sigma(l)$ and $j = \sigma(k)$.
\end{itemize}
We also consider some ``degenerate'' cases with $j=k$,
by saying that a triplet $i < j < l$ forms an
\begin{itemize}
   \item {\em upper joining}\/ (ujoin) if $j = \sigma(i)$ and $l = \sigma(j)$
      [i.e.\ the index $j$ is a cycle double rise];
   \item {\em lower joining}\/ (ljoin) if $i = \sigma(j)$ and $j = \sigma(l)$
      [i.e.\ the index $j$ is a cycle double fall];
   \item {\em upper pseudo-nesting}\/ (upsnest)
      if $l = \sigma(i)$ and $j = \sigma(j)$;
   \item {\em lower pseudo-nesting}\/ (lpsnest)
      if $i = \sigma(l)$ and $j = \sigma(j)$.
\end{itemize}
These are clearly degenerate cases of crossings and nestings, respectively.
See Figure~\ref{fig.crossnest}.
Note that $\upsnest(\sigma) = \lpsnest(\sigma)$ for all $\sigma$,
since for each fixed point~$j$,
the number of pairs $(i,l)$ with $i < j < l$ such that $l = \sigma(i)$
has to equal the number of such pairs with $i = \sigma(l)$;
we therefore write these two statistics simply as
\be
   \psnest(\sigma) \;\eqdef\; \upsnest(\sigma) \;=\;  \lpsnest(\sigma)
   \;.
\ee
And of course $\ujoin = \cdrise$ and $\ljoin = \cdfall$.

\begin{figure}[p]
\centering
\begin{picture}(30,15)(145, 10)
\setlength{\unitlength}{1.5mm}
\linethickness{.5mm}
\put(2,0){\line(1,0){28}}
\put(5,0){\circle*{1,3}}\put(5,0){\makebox(0,-6)[c]{\small $i$}}
\put(12,0){\circle*{1,3}}\put(12,0){\makebox(0,-6)[c]{\small $j$}}
\put(19,0){\circle*{1,3}}\put(19,0){\makebox(0,-6)[c]{\small $k$}}
\put(26,0){\circle*{1,3}}\put(26,0){\makebox(0,-6)[c]{\small $l$}}
\red{\qbezier(5,0)(12,10)(19,0)}
\blue{\qbezier(11,0)(18,10)(25,0)}
\put(15,-6){\makebox(0,-6)[c]{\small Upper crossing}}
%%%%%%%%%%%%%%%%%%%%
\put(43,0){\line(1,0){28}}
\put(47,0){\circle*{1,3}}\put(47,0){\makebox(0,-6)[c]{\small $i$}}
\put(54,0){\circle*{1,3}}\put(54,0){\makebox(0,-6)[c]{\small $j$}}
\put(61,0){\circle*{1,3}}\put(61,0){\makebox(0,-6)[c]{\small $k$}}
\put(68,0){\circle*{1,3}}\put(68,0){\makebox(0,-6)[c]{\small $l$}}
\red{\qbezier(47,0)(54,-10)(61,0)}
\blue{\qbezier(53,0)(60,-10)(67,0)}
\put(57,-6){\makebox(0,-6)[c]{\small Lower crossing}}
\end{picture}
\\[3.5cm]
\begin{picture}(30,15)(145, 10)
\setlength{\unitlength}{1.5mm}
\linethickness{.5mm}
\put(2,0){\line(1,0){28}}
\put(5,0){\circle*{1,3}}\put(5,0){\makebox(0,-6)[c]{\small $i$}}
\put(12,0){\circle*{1,3}}\put(12,0){\makebox(0,-6)[c]{\small $j$}}
\put(19,0){\circle*{1,3}}\put(19,0){\makebox(0,-6)[c]{\small $k$}}
\put(26,0){\circle*{1,3}}\put(26,0){\makebox(0,-6)[c]{\small $l$}}
\red{\qbezier(5,0)(15.5,10)(26,0)}
\blue{\qbezier(11,0)(14.5,5)(18,0)}
\put(15,-6){\makebox(0,-6)[c]{\small Upper nesting}}
%%%%%%%%%%%%%%%%%%%%
\put(43,0){\line(1,0){28}}
\put(47,0){\circle*{1,3}}\put(47,0){\makebox(0,-6)[c]{\small $i$}}
\put(54,0){\circle*{1,3}}\put(54,0){\makebox(0,-6)[c]{\small $j$}}
\put(61,0){\circle*{1,3}}\put(61,0){\makebox(0,-6)[c]{\small $k$}}
\put(68,0){\circle*{1,3}}\put(68,0){\makebox(0,-6)[c]{\small $l$}}
\red{\qbezier(47,0)(57,-10)(68,0)}
\blue{\qbezier(53,0)(56.5,-5)(60.5,0)}
\put(57,-6){\makebox(0,-6)[c]{\small Lower nesting}}
\end{picture}
\\[3.5cm]
\begin{picture}(30,15)(145, 10)
\setlength{\unitlength}{1.5mm}
\linethickness{.5mm}
\put(2,0){\line(1,0){28}}
\put(5,0){\circle*{1,3}}\put(5,0){\makebox(0,-6)[c]{\small $i$}}
\put(15.5,0){\circle*{1,3}}\put(15.5,0){\makebox(0,-6)[c]{\small $j$}}
%\put(19,0){\circle*{1,3}}\put(19,0){\makebox(0,-6)[c]%{\small $k$}}
\put(26,0){\circle*{1,3}}\put(26,0){\makebox(0,-6)[c]{\small $l$}}
\red{\qbezier(5,0)(10.5,10)(15.5,0)}
\blue{\qbezier(14.75,0)(19.5,10)(25.25,0)}
\put(15,-6){\makebox(0,-6)[c]{\small Upper joining}}
%%%%%%%%%%%%%%%%%%%%
\put(43,0){\line(1,0){28}}
\put(47,0){\circle*{1,3}}\put(47,0){\makebox(0,-6)[c]{\small $i$}}
\put(57.5,0){\circle*{1,3}}\put(57.5,0){\makebox(0,-6)[c]{\small $j$}}
%\put(51,0){\circle*{1,3}}\put(51,0){\makebox(0,-6)[c]{\small $k$}}
\put(68,0){\circle*{1,3}}\put(68,0){\makebox(0,-6)[c]{\small $l$}}
\red{\qbezier(47,0)(52,-10)(57.25,0)}
\blue{\qbezier(56.75,0)(62,-10)(67,0)}
%\blue{\qbezier(43,0)(46.5,-5)(50.5,0)}
\put(57,-6){\makebox(0,-6)[c]{\small Lower joining}}
\end{picture}
\\[3.5cm]
\begin{picture}(30,15)(145, 10)
\setlength{\unitlength}{1.5mm}
\linethickness{.5mm}
\put(2,0){\line(1,0){28}}
\put(5,0){\circle*{1,3}}\put(5,0){\makebox(0,-6)[c]{\small $i$}}
\put(15.5,0){\circle*{1,3}}\put(15.5,0){\makebox(0,-6)[c]{\small $j$}}
%\put(19,0){\circle*{1,3}}\put(19,0){\makebox(0,-6)[c]%{\small $k$}}
\put(26,0){\circle*{1,3}}\put(26,0){\makebox(0,-6)[c]{\small $l$}}
\red{\qbezier(5,0)(15.5,10)(26,0)}
%%\blue{\qbezier(14.25,0)(15,3)(14.75,0)}
\put(15,-6){\makebox(0,-6)[c]{\small Upper pseudo-nesting}}
%%%%%%%%%%%%%%%%%%%%
\put(43,0){\line(1,0){28}}
\put(47,0){\circle*{1,3}}\put(47,0){\makebox(0,-6)[c]{\small $i$}}
\put(57.5,0){\circle*{1,3}}\put(57.5,0){\makebox(0,-6)[c]{\small $j$}}
%\put(51,0){\circle*{1,3}}\put(51,0){\makebox(0,-6)[c]{\small $k$}}
\put(68,0){\circle*{1,3}}\put(68,0){\makebox(0,-6)[c]{\small $l$}}
\red{\qbezier(47,0)(57.5,-10)(68,0)}
%%\blue{\qbezier(56.75,0)(56.8,-3)(56.9,0)}
%\blue{\qbezier(53,0)(56.5,-5)(60.5,0)}
\put(57,-6){\makebox(0,-6)[c]{\small Lower pseudo-nesting}}
\end{picture}
\vspace*{3cm}
\caption{
   Crossing, nesting, joining and pseudo-nesting.
 \label{fig.crossnest}
}
\end{figure}

If $\sigma$ is a D-permutation, then its diagram has a special property:
all arrows emanating from odd (resp.~even) vertices
are upper (resp.~lower) arrows.
Otherwise put,
the leftmost (resp.~rightmost) vertex of an upper (resp.~lower) arc
is always odd (resp.~even).
It follows that in an upper crossing or nesting $i < j < k < l$,
the indices $i$ and $j$ must be odd;
and in a lower crossing or nesting $i < j < k < l$,
the indices $k$ and $l$ must be even.
Similar comments apply to upper and lower joinings and pseudo-nestings.

We can further refine the four crossing/nesting categories
by examining more closely the status of the inner index ($j$ or $k$)
whose {\em outgoing}\/ arc belonged to the crossing or nesting:
we say that a quadruplet $i < j < k < l$ forms an
\begin{itemize}
   \item {\em upper crossing of type cval}\/ (ucrosscval)
       if $k = \sigma(i)$ and $l = \sigma(j)$ and ${\sigma^{-1}(j) > j}$;
   \item {\em upper crossing of type cdrise}\/ (ucrosscdrise)
       \hbox{if $k = \sigma(i)$ and $l = \sigma(j)$ and ${\sigma^{-1}(j) < j}$;}
   \item {\em lower crossing of type cpeak}\/ (lcrosscpeak)
       if $i = \sigma(k)$ and $j = \sigma(l)$ and ${\sigma^{-1}(k) < k}$;
   \item {\em lower crossing of type cdfall}\/ (lcrosscdfall)
       if $i = \sigma(k)$ and $j = \sigma(l)$ and ${\sigma^{-1}(k) > k}$;
   \item {\em upper nesting of type cval}\/  (unestcval)
       if $l = \sigma(i)$ and $k = \sigma(j)$ and ${\sigma^{-1}(j) > j}$;
   \item {\em upper nesting of type cdrise}\/  (unestcdrise)
       if $l = \sigma(i)$ and $k = \sigma(j)$ and ${\sigma^{-1}(j) < j}$;
   \item {\em lower nesting of type cpeak}\/  (lnestcpeak)
       if $i = \sigma(l)$ and $j = \sigma(k)$ and ${\sigma^{-1}(k) < k}$;
   \item {\em lower nesting of type cdfall}\/  (lnestcdfall)
       if $i = \sigma(l)$ and $j = \sigma(k)$ and ${\sigma^{-1}(k) > k}$.
\end{itemize}
See Figure~\ref{fig.refined_crossnest}.
Please note that for the ``upper'' quantities
the distinguished index
(i.e.\ the one for which we examine both $\sigma$ and $\sigma^{-1}$)
is in second position ($j$),
while for the ``lower'' quantities
the distinguished index is in third position ($k$).

\begin{figure}[p]
\centering
\begin{picture}(30,15)(145, 10)
\setlength{\unitlength}{1.5mm}
\linethickness{.5mm}
\put(2,0){\line(1,0){28}}
\put(5,0){\circle*{1,3}}\put(5,0){\makebox(0,-6)[c]{\small $i$}}
\put(12,0){\circle*{1,3}}\put(12,0){\makebox(0,-6)[c]{\small $j$}}
\put(19,0){\circle*{1,3}}\put(19,0){\makebox(0,-6)[c]{\small $k$}}
\put(26,0){\circle*{1,3}}\put(26,0){\makebox(0,-6)[c]{\small $l$}}
\red{\qbezier(5,0)(12,10)(19,0)}
\blue{\qbezier(11,0)(18,10)(25,0)}
\blue{\qbezier(10.5,0)(11.75,-1)(13,-2)}
\put(15,-6){\makebox(0,-6)[c]{\small Upper crossing of type cval}}
%%%%%%%%%%%%%%%%%%%%
\put(43,0){\line(1,0){28}}
\put(47,0){\circle*{1,3}}\put(47,0){\makebox(0,-6)[c]{\small $i$}}
\put(54,0){\circle*{1,3}}\put(54,0){\makebox(0,-6)[c]{\small $j$}}
\put(61,0){\circle*{1,3}}\put(61,0){\makebox(0,-6)[c]{\small $k$}}
\put(68,0){\circle*{1,3}}\put(68,0){\makebox(0,-6)[c]{\small $l$}}
\red{\qbezier(47,0)(54,10)(61,0)}
\blue{\qbezier(53,0)(60,10)(67,0)}
\blue{\qbezier(50,2)(51.25,1)(52.5,0)}
\put(57,-6){\makebox(0,-6)[c]{\small Upper crossing of type cdrise}}
\end{picture}
\\[3.5cm]
\begin{picture}(30,15)(145, 10)
\setlength{\unitlength}{1.5mm}
\linethickness{.5mm}
\put(2,0){\line(1,0){28}}
\put(5,0){\circle*{1,3}}\put(5,0){\makebox(0,-6)[c]{\small $i$}}
\put(12,0){\circle*{1,3}}\put(12,0){\makebox(0,-6)[c]{\small $j$}}
\put(19,0){\circle*{1,3}}\put(19,0){\makebox(0,-6)[c]{\small $k$}}
\put(26,0){\circle*{1,3}}\put(26,0){\makebox(0,-6)[c]{\small $l$}}
\red{\qbezier(5,0)(12,-10)(19,0)}
\red{\qbezier(16,2)(17.25, 1)(18.5,0)}
\blue{\qbezier(10.5,0)(18,-10)(25,0)}
\put(15,-6){\makebox(0,-6)[c]{\small Lower crossing of type cpeak}}
%%%%%%%%%%%%%%%%%%%%
\put(43,0){\line(1,0){28}}
\put(47,0){\circle*{1,3}}\put(47,0){\makebox(0,-6)[c]{\small $i$}}
\put(54,0){\circle*{1,3}}\put(54,0){\makebox(0,-6)[c]{\small $j$}}
\put(61,0){\circle*{1,3}}\put(61,0){\makebox(0,-6)[c]{\small $k$}}
\put(68,0){\circle*{1,3}}\put(68,0){\makebox(0,-6)[c]{\small $l$}}
\red{\qbezier(47,0)(54,-10)(61,0)}
\blue{\qbezier(53,0)(60,-10)(67,0)}
\red{\qbezier(59.5,0)(60.75,-1)(62,-2)}
\put(57,-6){\makebox(0,-6)[c]{\small Lower crossing of type cdfall}}
\end{picture}
\\[3.5cm]
\begin{picture}(30,15)(145, 10)
\setlength{\unitlength}{1.5mm}
\linethickness{.5mm}
\put(2,0){\line(1,0){28}}
\put(5,0){\circle*{1,3}}\put(5,0){\makebox(0,-5)[c]{\small $i$}}
\put(12,0){\circle*{1,3}}\put(12,0){\makebox(0,-5)[c]{\small $j$}}
\put(19,0){\circle*{1,3}}\put(19,0){\makebox(0,-5)[c]{\small $k$}}
\put(26,0){\circle*{1,3}}\put(26,0){\makebox(0,-5)[c]{\small $l$}}
\red{\qbezier(5,0)(15.5,10)(26,0)}
\blue{\qbezier(11,0)(14.5,5)(18,0)}
\blue{\qbezier(10.25,0)(11.5,-1)(12.75,-2)}
\put(15,-6){\makebox(0,-6)[c]{\small Upper nesting of type cval}}
%%%%%%%%%%%%%%%%%%%%
\put(43,0){\line(1,0){28}}
\put(47,0){\circle*{1,3}}\put(47,0){\makebox(0,-6)[c]{\small $i$}}
\put(54,0){\circle*{1,3}}\put(54,0){\makebox(0,-6)[c]{\small $j$}}
\put(61,0){\circle*{1,3}}\put(61,0){\makebox(0,-6)[c]{\small $k$}}
\put(68,0){\circle*{1,3}}\put(68,0){\makebox(0,-6)[c]{\small $l$}}
\red{\qbezier(47,0)(57,10)(68,0)}
\blue{\qbezier(53,0)(56.5,5)(60.5,0)}
\blue{\qbezier(50,2)(51.25,1)(52.5,0)}
\put(57,-6){\makebox(0,-6)[c]{\small Upper nesting of type cdrise}}
\end{picture}
\\[3.5cm]
\begin{picture}(30,15)(145, 10)
\setlength{\unitlength}{1.5mm}
\linethickness{.5mm}
\put(2,0){\line(1,0){28}}
\put(5,0){\circle*{1,3}}\put(5,0){\makebox(0,-5)[c]{\small $i$}}
\put(12,0){\circle*{1,3}}\put(12,0){\makebox(0,-5)[c]{\small $j$}}
\put(19,0){\circle*{1,3}}\put(19,0){\makebox(0,-5)[c]{\small $k$}}
\put(26,0){\circle*{1,3}}\put(26,0){\makebox(0,-5)[c]{\small $l$}}
\red{\qbezier(5,0)(15.5,-10)(26,0)}
\blue{\qbezier(11,0)(14.5,-5)(18,0)}
\blue{\qbezier(14.75,2)(16,1)(17.25,0)}
\put(15,-6){\makebox(0,-6)[c]{\small Lower nesting of type cpeak}}
%%%%%%%%%%%%%%%%%%%%
\put(43,0){\line(1,0){28}}
\put(47,0){\circle*{1,3}}\put(47,0){\makebox(0,-5)[c]{\small $i$}}
\put(54,0){\circle*{1,3}}\put(54,0){\makebox(0,-5)[c]{\small $j$}}
\put(61,0){\circle*{1,3}}\put(61,0){\makebox(0,-5)[c]{\small $k$}}
\put(68,0){\circle*{1,3}}\put(68,0){\makebox(0,-5)[c]{\small $l$}}
\red{\qbezier(47,0)(57,-10)(68,0)}
\blue{\qbezier(53,0)(56.5,-5)(60.5,0)}
\blue{\qbezier(59.5,0)(60.75, -1)(62,-2)}
\put(57,-6){\makebox(0,-6)[c]{\small Lower nesting of type cdfall}}
\end{picture}
\vspace*{3cm}
\caption{
   Refined categories of crossing and nesting.
 \label{fig.refined_crossnest}
}
\end{figure}

In fact, a central role in our work will be played
(just as in \cite{Sokal-Zeng_masterpoly})
by a refinement of these statistics:
rather than counting the {\em total}\/ numbers of quadruplets
$i < j < k < l$ that form upper (resp.~lower) crossings or nestings,
we will count the number of upper (resp.~lower) crossings or nestings
that use a particular vertex $j$ (resp.~$k$)
in second (resp.~third) position.
More precisely, we define the
\textbfit{index-refined crossing and nesting statistics}
\begin{subeqnarray}
   \ucross(j,\sigma)
   & = &
   \#\{ i<j<k<l \colon\: k = \sigma(i) \hbox{ and } l = \sigma(j) \}
         \\[2mm]
   \unest(j,\sigma)
   & = &
   \#\{ i<j<k<l \colon\: k = \sigma(j) \hbox{ and } l = \sigma(i) \}
      \\[2mm]
   \lcross(k,\sigma)
   & = &
   \#\{ i<j<k<l \colon\: i = \sigma(k) \hbox{ and } j = \sigma(l) \}
         \\[2mm]
   \lnest(k,\sigma)
   & = &
   \#\{ i<j<k<l \colon\: i = \sigma(l) \hbox{ and } j = \sigma(k) \}
 \label{def.ucrossnestjk}
\end{subeqnarray}
%
% {\bf But maybe we want to define these with second and third positions
%    reversed, i.e.
% \begin{subeqnarray}
%    \ucross^\star(k,\sigma)
%    & = &
%    \#\{ i<j<k<l \colon\: k = \sigma(i) \hbox{ and } l = \sigma(j) \}
%          \\[2mm]
%    \unest^\star(k,\sigma)
%    & = &
%    \#\{ i<j<k<l \colon\: k = \sigma(j) \hbox{ and } l = \sigma(i) \}
%       \\[2mm]
%    \lcross^\star(j,\sigma)
%    & = &
%    \#\{ i<j<k<l \colon\: i = \sigma(k) \hbox{ and } j = \sigma(l) \}
%          \\[2mm]
%    \lnest^\star(j,\sigma)
%    & = &
%    \#\{ i<j<k<l \colon\: i = \sigma(l) \hbox{ and } j = \sigma(k) \}
%  \label{def.ucrossnestjk.star}
% \end{subeqnarray}
% This would not only remove the $\sigma^{-1}$ in \reff{eq.xii.nestings},
% but even more importantly it would allow the subsequent equations
% to be interpreted in terms of $\ucross^\star$ and $\lcross^\star$.}
%
Note that $\ucross(j,\sigma)$ and $\unest(j,\sigma)$ can be nonzero
only when $j$ is an excedance
(that is, a cycle valley or a cycle double rise),
while $\lcross(k,\sigma)$ and $\lnest(k,\sigma)$ can be nonzero
only when $k$ is an anti-excedance
(that is, a cycle peak or a cycle double fall).
In a D-permutation, this means that
$\ucross(j,\sigma)$ and $\unest(j,\sigma)$ can be nonzero
only when $j$ is odd and not a fixed point,
while $\lcross(k,\sigma)$ and $\lnest(k,\sigma)$ can be nonzero
only when $k$ is even and not a fixed point.

When $j$ is a fixed point, we also define the analogous quantity
for pseudo-nestings:
\be
   \psnest(j,\sigma)
   \;\eqdef\;
   \# \{i < j \colon\:  \sigma(i) > j \}
   \;=\;
   \# \{i > j \colon\:  \sigma(i) < j \}
   \;.
 \label{def.psnestj}
\ee
(Here the two expressions are equal because $\sigma$ is a bijection
 from $[1,j) \cup (j,n]$ to itself.)
In \cite[eq.~(2.20)]{Sokal-Zeng_masterpoly}
this quantity was called the {\em level}\/ of the fixed point $j$
and was denoted $\lev(j,\sigma)$.

For some purposes we will also require
%% (in Sections~\ref{subsec.second.master} and \ref{sec.proofs.2})
a variant of \reff{def.ucrossnestjk}
in which the roles of second and third position are interchanged:
\begin{subeqnarray}
   \ucross'(k,\sigma)
   & = &
   \#\{ i<j<k<l \colon\: k = \sigma(i) \hbox{ and } l = \sigma(j) \}
         \\[2mm]
   \unest'(k,\sigma)
   & = &
   \#\{ i<j<k<l \colon\: k = \sigma(j) \hbox{ and } l = \sigma(i) \}
      \\[2mm]
   \lcross'(j,\sigma)
   & = &
   \#\{ i<j<k<l \colon\: i = \sigma(k) \hbox{ and } j = \sigma(l) \}
         \\[2mm]
   \lnest'(j,\sigma)
   & = &
   \#\{ i<j<k<l \colon\: i = \sigma(l) \hbox{ and } j = \sigma(k) \}
 \label{def.ucrossnestjk.prime}
\end{subeqnarray}
We remark that since nestings join the vertices
in second and third positions, we have
\begin{subeqnarray}
   \unest'(k,\sigma) & = & \unest(\sinv(k),\sigma)
      \\[2mm]
   \lnest'(j,\sigma) & = & \lnest(\sinv(j),\sigma)
 \label{eq.nestprime}
\end{subeqnarray}
Note that $\ucross'(k,\sigma)$ and $\unest'(k,\sigma)$ can be nonzero
only when $\sinv(k)$ is an excedance
(that is, when $k$ is a cycle peak or a cycle double rise),
while $\lcross'(j,\sigma)$ and $\lnest'(j,\sigma)$ can be nonzero
only when $\sinv(j)$ is an anti-excedance
(that is, $j$ is a cycle valley or a cycle double fall).
In a D-permutation, this means that
$\ucross'(k,\sigma)$ and $\unest'(k,\sigma)$ can be nonzero
only when $\sinv(k)$ is odd and not a fixed point,
while $\lcross'(j,\sigma)$ and $\lnest'(j,\sigma)$ can be nonzero
only when $\sinv(j)$ is even and not a fixed point.
We call \reff{def.ucrossnestjk.prime} the
\textbfit{variant index-refined crossing and nesting statistics}.

\section{First T-fraction and its generalizations}   \label{sec.first}

In this section we state our first T-fraction for D-permutations,
in three increasingly more general versions.
The first and most basic version (Theorem~\ref{thm.Tfrac.first})
is a T-fraction in 12~variables that enumerates D-permutations
with respect to the parity-refined record-and-cycle classification;
it includes many previously known results as special cases.
%% {\bf Is this true???  Or only continued fractions for the
%%    Genocchi and median Genocchi numbers???}
The second version (Theorem~\ref{thm.Tfrac.first.pq})
is a $(p,q)$-generalization of the first one:
it is a T-fraction in 22~variables that enumerates D-permutations
with respect to the parity-refined record-and-cycle classification
together with four pairs of $(p,q)$-variables
counting the refined categories of crossing and nesting
(Figure~\ref{fig.refined_crossnest})
as well as two variables corresponding to pseudo-nestings of fixed points.
Finally, our third version
--- what we call the ``first master T-fraction''
 (Theorem~\ref{thm.Tfrac.first.master}) ---
is a T-fraction in six infinite families of indeterminates
that generalizes the preceding two
by employing the index-refined crossing and nesting statistics
defined in \reff{def.ucrossnestjk}.

The plan of this section is as follows:
We begin (Section~\ref{subsec.rar}) by establishing
some easy but important preliminary results concerning
record-antirecord fixed points.
Then, in Sections~\ref{subsec.first.statements}--\ref{subsec.first.master}
we state the three versions of our first T-fraction
and note some of their corollaries.
Finally, in Section~\ref{subsec.first.variant}
we state variant forms of the three T-fractions.
All these results will be proven in Section~\ref{sec.proofs.1}.

\subsection{Preliminaries on record-antirecords}  \label{subsec.rar}

An important role in our study of D-permutations
will be played by {\em record-antirecords}\/,
i.e.\ indices~$i$ that are both a record and an antirecord.
It is easy to see that, in any permutation,
a record-antirecord must be a fixed point.
More precisely:

\begin{lemma}
   \label{lemma.recantirec}
Consider a permutation $\sigma \in \Sym_N$.
An index $i \in [N]$ is a record-antirecord if and only if
$\sigma$ maps each of the sets
$\{1,\ldots,i-1\}$, $\{i\}$ and $\{i+1,\ldots,N\}$
onto itself.
\end{lemma}

\proof
If $i$ is a record, then
$\sigma$ maps $\{1,\ldots,i-1\}$ injectively into $\{1,\ldots,\sigma(i)-1\}$,
so that $\sigma(i) \ge i$.
If $i$ is an antirecord, then
$\sigma$ maps $\{i+1,\ldots,N\}$ injectively into $\{\sigma(i)+1,\ldots,N\}$,
so that $\sigma(i) \le i$.
So if $i$ is a record-antirecord,
then $i$ must be a fixed point
and $\sigma$ must map both $\{1,\ldots,i-1\}$ and $\{i+1,\ldots,N\}$
bijectively onto themselves.

The converse is obvious.
\qed

For D-permutations, record-antirecords can occur only in pairs:

% \begin{lemma}
%    \label{lemma.recantirec}
% Let $\sigma \in \dperm_{2n}$ and let $i \in [n]$.
% Then $2i-1$ is a record-antirecord if and only if $2i$ is a record-antirecord.
% This occurs if and only if $\sigma$ maps each of the sets
% $\{1,\ldots,2i-2\}$, $\{2i-1\}$, $\{2i\}$ and $\{2i+1,\ldots,2n\}$
% bijectively onto itself.
% \end{lemma}

\begin{lemma}
   \label{lemma.recantirec_D}
Let $\sigma \in \dperm_{2n}$ and let $i \in [n]$.
Then the following are equivalent:
\begin{itemize}
   \item[(a)]  $2i-1$ is a record-antirecord.
   \item[(b)]  $2i$ is a record-antirecord.
   \item[(c)]  $\sigma$ maps each of the sets
      $\{1,\ldots,2i-1\}$ and $\{2i,\ldots,2n\}$ onto itself.
   \item[(d)]  $\sigma$ maps each of the sets
      $\{1,\ldots,2i-2\}$, $\{2i-1\}$, $\{2i\}$ and $\{2i+1,\ldots,2n\}$
      onto itself.
\end{itemize}
\end{lemma}

\proof
By Lemma~\ref{lemma.recantirec}, (d)$\iff$``(a) and (b)'',
and ``(a) or (b)''$\implies$(c).
%% and (d)$\implies$(c) is obvious.
On the other hand, if $\sigma$ is a D-permutation,
then $\sigma(2i-1) \ge 2i-1$ and $\sigma(2i) \le 2i$,
so (c)$\implies$(d).
\qed

Let us say that a permutation is \textbfit{pure}
if it has no record-antirecords.
We write $\dperm^{\rm pure}_{2n}$
for the set of pure D-permutations of $[2n]$.
From Lemma~\ref{lemma.recantirec} we have
%% https://tex.stackexchange.com/questions/191544/how-to-write-inclusion-symbols-diagonally
\be
\begin{array}{rcccl}
      &     &  \dperm^{\rm e}_{2n}  &     &                \\
      & \rotatebox[origin=c]{45}{$\subseteq$} &      &
        \rotatebox[origin=c]{-45}{$\subseteq$} &           \\
   \dcycle_{2n}  \;\subseteq\;
   \dperm^{\rm eo}_{2n} = \dperm^{\rm e}_{2n} \cap \dperm^{\rm o}_{2n}
         &    &    &   & \dperm^{\rm e}_{2n} \cup \dperm^{\rm o}_{2n}  \\
      & \rotatebox[origin=c]{-45}{$\subseteq$} &      &
        \rotatebox[origin=c]{45}{$\subseteq$} &           \\
      &     &  \dperm^{\rm o}_{2n}  &     &              
\end{array}
\;\subseteq\; \dperm^{\rm pure}_{2n}
\;\subseteq\; \dperm_{2n}
 \;.
\ee
These inclusions are strict for $n \ge 2$:
\bigskip
\begin{table}[H]
\centering
\begin{tabular}{c|rrrrrr}
    & $\dcycle_{2n}$ &
      $\dperm^{\rm eo}_{2n}$ &
      $\dperm^{\rm e/o}_{2n}$ &
      $\dperm^{\rm e}_{2n} \cup \dperm^{\rm o}_{2n}$ &
      $\dperm^{\rm pure}_{2n}$ &
      $\dperm_{2n}$  \\[1mm]
%%\cline{2-6}
$n$ & $g_{n-1}$ & $h_n$ & $g_n$ & $2g_n - h_n$ & $h_{n+1}^\flat$ &  $h_{n+1}$ \\
\hline
\hline
0 &      0 &       1 &      1 &      1\hspace*{3mm} &      1 &      1 \\
1 &      1 &       1 &      1 &      1\hspace*{3mm} &      1 &      2 \\
2 &      1 &       2 &      3 &      4\hspace*{3mm} &      5 &      8 \\
3 &      3 &       8 &     17 &     26\hspace*{3mm} &     41 &     56 \\
4 &     17 &      56 &    155 &    254\hspace*{3mm} &    493 &    608 \\
5 &    155 &     608 &   2073 &   3538\hspace*{3mm} &   8161 &   9440 \\
6 &   2073 &    9440 &  38227 &  67014\hspace*{3mm} & 178469 & 198272 \\
\end{tabular}
\end{table}
\bigskip
\noindent
For instance, the permutation 4231 belongs to
$\dperm^{\rm pure}_{4} \setminus (\dperm^{\rm e}_{4} \cup \dperm^{\rm o}_{4})$:
it has both even and odd fixed points, but they are not record-antirecords.
Likewise, the permutation 513462 belongs to
$\dperm^{\rm pure}_{6} \setminus (\dperm^{\rm e}_{6} \cup \dperm^{\rm o}_{6})$:
it has an odd-even pair of fixed points 34, but they are not record-antirecords.
We will show later that $|\dperm^{\rm pure}_{2n}| = h_{n+1}^\flat$;
we recall that these numbers were defined by the S-fraction
\reff{eq.mediangenocchiflat.Sfrac}/\reff{eq.mediangenocchiflat.Sfrac.weights}.

\subsection{First T-fraction}
   \label{subsec.first.statements}

\subsubsection{Main theorem and its corollaries}

We now introduce a polynomial in 12 variables that enumerates D-permutations
according to the parity-refined record-and-cycle classification:
\begin{eqnarray}
   & &
   P_n(x_1,x_2,y_1,y_2,u_1,u_2,v_1,v_2,\we,\wo,\ze,\zo)
   \;=\;
       \nonumber \\[4mm]
   & & \qquad\qquad
   \sum_{\sigma \in \dperm_{2n}}
   x_1^{\eareccpeak(\sigma)} x_2^{\eareccdfall(\sigma)}
   y_1^{\ereccval(\sigma)} y_2^{\ereccdrise(\sigma)}
   \:\times
       \qquad\qquad
       \nonumber \\[-1mm]
   & & \qquad\qquad\qquad\:
   u_1^{\nrcpeak(\sigma)} u_2^{\nrcdfall(\sigma)}
   v_1^{\nrcval(\sigma)} v_2^{\nrcdrise(\sigma)}
   \:\times
       \qquad\qquad
       \nonumber \\[3mm]
   & & \qquad\qquad\qquad\:
   \we^{\evennrfix(\sigma)} \wo^{\oddnrfix(\sigma)}
   \ze^{\evenrar(\sigma)} \zo^{\oddrar(\sigma)}
   \;.
 \label{def.Pn}
\end{eqnarray}
Thus, the variables $x_1$ and $u_1$ are associated to cycle peaks,
$y_1$ and $v_1$ to cycle valleys,
$x_2$ and $u_2$ to cycle double falls,
$y_2$ and $v_2$ to cycle double rises,
$\we$ and $\wo$ to neither-record-antirecord fixed points,
and $\ze$ and $\zo$ to record-antirecord fixed points.
We remark that \reff{def.Pn} is the same as the polynomial
introduced in \cite[eq.~(2.19)]{Sokal-Zeng_masterpoly},
but restricted to D-permutations
and refined to record the parity of fixed points.
Since in a D-permutation each even (resp.~odd) index $i \in [2n]$
belongs to exactly one of the 6 types
mentioned at the end of Section~\ref{subsec.statistics.1},
it follows that the polynomial $P_n$ is
homogeneous of degree~$n$ in $x_1,x_2,u_1,u_2,\we,\ze$
and also homogeneous of degree~$n$ in $y_1,y_2,v_1,v_2,\wo,\zo$.

The polynomials \reff{def.Pn} have a beautiful T-fraction:

\begin{theorem}[First T-fraction for D-permutations]
   \label{thm.Tfrac.first}
The ordinary generating function of the polynomials \reff{def.Pn} has the
T-type continued fraction
\begin{eqnarray}
   & & \hspace*{-12mm}
   \sum_{n=0}^\infty
   P_n(x_1,x_2,y_1,y_2,u_1,u_2,v_1,v_2,\we,\wo,\ze,\zo) \: t^n
   \;=\;
       \nonumber \\
   & &
   \cfrac{1}{1 - \ze \zo  \,t - \cfrac{x_1 y_1  \,t}{1 -  \cfrac{(x_2\!+\!\we)(y_2\!+\!\wo) \,t}{1 - \cfrac{(x_1\!+\!u_1)(y_1\!+\!v_1)  \,t}{1 - \cfrac{(x_2\!+\!u_2\!+\!\we)(y_2\!+\!v_2\!+\!\wo)  \,t}{1 - \cfrac{(x_1\!+\!2u_1)(y_1\!+\!2v_1)  \,t}{1 - \cfrac{(x_2\!+\!2u_2\!+\!\we)(y_2\!+\!2v_2\!+\!\wo)  \,t}{1 - \cdots}}}}}}}
       \nonumber \\[1mm]
   \label{eq.thm.Tfrac.first}
\end{eqnarray}
with coefficients
\begin{subeqnarray}
   \alpha_{2k-1} & = & [x_1 + (k-1) u_1] \: [y_1 + (k-1) v_1]
        \\[1mm]
   \alpha_{2k}   & = & [x_2 + (k-1) u_2 + \we] \: [y_2 + (k-1) v_2 + \wo]
        \\[1mm]
   \delta_1  & = &   \ze \zo   \\[1mm]
   \delta_n  & = &   0    \qquad\hbox{for $n \ge 2$}
   \label{eq.thm.Tfrac.first.weights}
\end{subeqnarray}
\end{theorem}

\noindent
We will prove Theorem~\ref{thm.Tfrac.first} in Section~\ref{sec.proofs.1}.

Note that each of the coefficients $\alpha_i$ and $\delta_i$
is homogeneous of degree~$1$ in $x_1,x_2,$ $u_1,u_2,\we,\ze$
and also homogeneous of degree~$1$ in $y_1,y_2,v_1,v_2,\wo,\zo$.
This reflects the homogeneities of the $P_n$.

Note also that the involution of $\dperm_{2n}$ defined by
$\sigma \mapsto R \circ \sigma \circ R$
where $R(i) = 2n+1-i$ is the reversal map,
interchanges
$x_1 \leftrightarrow y_1$,
$x_2 \leftrightarrow y_2$,
$u_1 \leftrightarrow v_1$,
$u_2 \leftrightarrow v_2$,
$\ze \leftrightarrow \zo$,
$\we \leftrightarrow \wo$.
The T-fraction \reff{eq.thm.Tfrac.first}/\reff{eq.thm.Tfrac.first.weights}
is invariant under these simultaneous interchanges.

\bigskip

The T-fraction \reff{eq.thm.Tfrac.first}/\reff{eq.thm.Tfrac.first.weights}
has numerous interesting specializations:

\medskip

1) With all variables equal to 1,
it gives the T-fraction
\reff{eq.mediangenocchi.Tfrac}/\reff{eq.mediangenocchi.Tfrac.weights}
for the once-shifted median Genocchi numbers,
and confirms that $|\dperm_{2n}| = h_{n+1}$.

\medskip

2) With $\ze=0$ and/or $\zo=0$ and all other variables equal to 1,
it gives the S-fraction
\reff{eq.mediangenocchiflat.Sfrac}/\reff{eq.mediangenocchiflat.Sfrac.weights}
and shows that
$|\dperm^{\rm pure}_{2n}| = h_{n+1}^\flat$.
More generally, with $\ze=0$ and/or $\zo=0$
and the other variables retained, it gives an S-fraction for pure D-permutations
according to the parity-refined record-and-cycle classification.
Note that since record-antirecords occur in pairs
by Lemma~\ref{lemma.recantirec_D},
setting either $\ze=0$ or $\zo=0$ suffices to suppress them;
this explains why these variables occur in the T-fraction
only as a product $\ze \zo$.

\medskip

3) With $\ze=\we=0$ (resp.\ $\zo=\wo=0$) and all other variables equal to 1,
it gives the S-fraction
\reff{eq.genocchi.Sfrac}/\reff{eq.genocchi.Sfrac.weights}
for the Genocchi numbers,
and confirms that $|\dperm^{\rm e}_{2n}| = |\dperm^{\rm o}_{2n}| = g_n$.
More generally, with $\ze=\we=0$ (resp.\ $\zo=\wo=0$)
and the other variables retained,
it gives an S-fraction for D-e-semiderangements (resp.\ D-o-semiderangements)
according to the parity-refined record-and-cycle classification.

\medskip

4) With $\ze=\zo=\we=\wo=0$ and all other variables equal to 1,
it gives the S-fraction
\reff{eq.mediangenocchi.Sfrac}/\reff{eq.mediangenocchi.Sfrac.weights}
for the median Genocchi numbers,
and confirms that $|\dperm^{\rm eo}_{2n}| = h_n$.
More generally, with $\ze=\zo=\we=\wo=0$ and the other variables retained,
it gives an S-fraction for D-derangements
according to the parity-refined record-and-cycle classification;
this S-fraction is of precisely the form \reff{eq.quadratic_family}
that was proposed in the Introduction.

\medskip

5) If we specialize to $x_1 = u_1$, $y_1 = v_1$, $x_2 = u_2$, $y_2 = v_2$
--- that is, renounce the counting of records ---
then the coefficients (\ref{eq.thm.Tfrac.first.weights}a,b) simplify to
\begin{subeqnarray}
   \alpha_{2k-1} & = &  k^2 \, x_1 y_1
        \\[1mm]
   \alpha_{2k}   & = &  (k x_2 + \we) \, (k y_2 + \wo)
   \label{eq.thm.Tfrac.first.weights.bis}
\end{subeqnarray}

\bigskip

Using Proposition~\ref{prop.augment.Tfrac},
we can alternatively rewrite the T-fraction \reff{eq.thm.Tfrac.first}
as an S-fraction for the ogf of an ``augmented'' sequence:

\begin{corollary}[First S-fraction for augmented D-permutations]
   \label{cor.Sfrac.first.augmented}
The ordinary generating function of the ``augmented'' sequence
of polynomials \reff{def.Pn} has the S-type continued fraction
\begin{eqnarray}
   & & \hspace*{-12mm}
   1 \:+\: \ze \zo t
   \sum_{n=0}^\infty
   P_n(x_1,x_2,y_1,y_2,u_1,u_2,v_1,v_2,\we,\wo,\ze,\zo) \: t^n
   \;=\;
       \nonumber \\
   & &
   \cfrac{1}{1 - \cfrac{\ze \zo  \,t}{1 - \cfrac{x_1 y_1  \,t}{1 -  \cfrac{(x_2\!+\!\we)(y_2\!+\!\wo) \,t}{1 - \cfrac{(x_1\!+\!u_1)(y_1\!+\!v_1)  \,t}{1 - \cfrac{(x_2\!+\!u_2\!+\!\we)(y_2\!+\!v_2\!+\!\wo)  \,t}{1 - \cfrac{(x_1\!+\!2u_1)(y_1\!+\!2v_1)  \,t}{1 - \cfrac{(x_2\!+\!2u_2\!+\!\we)(y_2\!+\!2v_2\!+\!\wo)  \,t}{1 - \cdots}}}}}}}}
       \nonumber \\[1mm]
   \label{eq.cor.Sfrac.first.augmented}
\end{eqnarray}
with coefficients
\begin{subeqnarray}
   \alpha_1  & = &  \ze \zo \\[1mm]
   \alpha_{2k-1}  & = & [x_2 + (k-2) u_2 + \we] \: [y_2 + (k-2) v_2 + \wo]
       \quad\hbox{for $k \ge 2$}
        \\[1mm]
   \alpha_{2k} & = & [x_1 + (k-1) u_1] \: [y_1 + (k-1) v_1]
   \label{eq.cor.Sfrac.first.augmented.weights}
\end{subeqnarray}
\end{corollary}

With all variables equal to 1,
this gives (since $|\dperm_{2n}| = h_{n+1}$) the the S-fraction
\reff{eq.mediangenocchi.Sfrac}/\reff{eq.mediangenocchi.Sfrac.weights}
for the sequence $(h_n)_{n \ge 0}$.

\bigskip

Also, using Proposition~\ref{prop.augment.Tfrac} in the other direction,
we can rewrite the T-fraction \reff{eq.thm.Tfrac.first}
specialized to $\ze = 0$ and/or $\zo = 0$
as a T-fraction for the ogf of the ``restricted'' sequence:

\begin{corollary}[First T-fraction for restricted pure D-permutations]
   \label{cor.Tfrac.first.restricted}
The ordinary generating function of the ``restricted'' sequence
of polynomials \reff{def.Pn} specialized to $\ze = 0$ and/or $\zo = 0$
has the T-type continued fraction
\begin{eqnarray}
   & & \hspace*{-12mm}
   \sum_{n=0}^\infty
   P_{n+1}(x_1,x_2,y_1,y_2,u_1,u_2,v_1,v_2,\we,\wo,0,\zo) \: t^n
   \;=\;
       \nonumber \\
   & & \hspace*{-12mm}
   \sum_{n=0}^\infty
   P_{n+1}(x_1,x_2,y_1,y_2,u_1,u_2,v_1,v_2,\we,\wo,\ze,0) \: t^n
   \;=\;
       \nonumber \\
   & &
   \cfrac{x_1 y_1}{1 - x_1 y_1  \,t -  \cfrac{(x_2\!+\!\we)(y_2\!+\!\wo) \,t}{1 - \cfrac{(x_1\!+\!u_1)(y_1\!+\!v_1)  \,t}{1 - \cfrac{(x_2\!+\!u_2\!+\!\we)(y_2\!+\!v_2\!+\!\wo)  \,t}{1 - \cfrac{(x_1\!+\!2u_1)(y_1\!+\!2v_1)  \,t}{1 - \cfrac{(x_2\!+\!2u_2\!+\!\we)(y_2\!+\!2v_2\!+\!\wo)  \,t}{1 - \cdots}}}}}}
       \nonumber \\[1mm]
   \label{eq.cor.Tfrac.first.restricted}
\end{eqnarray}
with coefficients
\begin{subeqnarray}
   \alpha_{2k-1}  & = & [x_2 + (k-1) u_2 + \we] \: [y_2 + (k-1) v_2 + \wo]
        \\[1mm]
   \alpha_{2k} & = & [x_1 + k\,u_1] \: [y_1 + k\,v_1]
        \\[1mm]
   \delta_1  & = &  x_1 y_1 \\[1mm]
   \delta_n  & = &   0    \qquad\hbox{for $n \ge 2$}
   \label{eq.cor.Tfrac.first.restricted.weights}
\end{subeqnarray}
\end{corollary}

With all variables equal to 1,
this gives (since $|\dperm^{\rm pure}_{2n}| = h^\flat_{n+1}$)
the T-fraction
\reff{eq.mediangenocchiflat.Tfrac}/\reff{eq.mediangenocchiflat.Tfrac.weights}
for the sequence $(h^\flat_{n+2})_{n \ge 0}$.
With $\we = 0$ (or $\wo = 0$) and all other variables equal to 1,
this gives the T-fraction
\reff{eq.genocchi.Tfrac}/\reff{eq.genocchi.Tfrac.weights}
for the once-shifted Genocchi numbers.
And finally,
with $\we = \wo = 0$ and all other variables equal to 1,
this gives the T-fraction
\reff{eq.mediangenocchi.Tfrac}/\reff{eq.mediangenocchi.Tfrac.weights}
for the once-shifted median Genocchi numbers.

\subsubsection{Alternative T-fraction, and specializations leading to S-fraction}

T-fractions --- unlike S-fractions and J-fractions --- are not unique;
indeed, any formal power series can be expressed as a T-fraction in
uncountably many different ways.\footnote{
   See footnote~\ref{footnote_Tfrac_nonunique} below.
}
Of course, most of these T-fractions will be too complicated
to be of any use; but sometimes they can be of interest.
We would now like to exhibit, for the ordinary generating function
of the polynomials \reff{def.Pn}, a T-fraction that is alternative
to the one shown in Theorem~\ref{thm.Tfrac.first}:

\begin{proposition}[Alternative T-fraction for D-permutations]
   \label{prop.Tfrac.first.alternative}
The ordinary generating function of the polynomials \reff{def.Pn}
has the alternative T-type continued fraction
\begin{eqnarray}
   & & \hspace*{-6mm}
    \sum_{n=0}^{\infty}
          P_n(x_1,x_2,y_1,y_2,u_1,u_2,v_1,v_2,\we,\wo,\ze,\zo) \: t^n
   \;=\;
      \nonumber \\
   & & \hspace*{-4mm}
\Scale[0.72]{
  \cfrac{1}{1 \!-\! [(x_1 \!-\! x_2 \!-\! \we) y_1 \!+\! \ze \zo] t \!-\! \cfrac{(x_2\!+\!\we)y_1 t}{1 \!-\!  \cfrac{x_1 (y_2 \!+\! \wo) t}{1 \!-\! [(x_1 \!-\! x_2 \!-\! \we)(y_1 \!-\! y_2 \!+\! v_1 \!-\!\wo) \!+\! (u_1 \!-\! u_2)(y_1 \!+\! v_1)] t \!-\! \cfrac{(x_2\!+\!u_2\!+\!\we)(y_1 \!+\! v_1) t}{1 \!-\! \cfrac{(x_1\!+\!u_1)(y_2 \!+\! v_2 \!+\! \wo) t}{1 \!-\! \cdots}}}}}
}
   \nonumber \\
\label{eq.prop.Tfrac.first.alternative}
\end{eqnarray}
with coefficients
\begin{subeqnarray}
    \alpha'_{2k-1} & = & [x_2 + (k-1)u_2 + \we] \, [y_1 + (k-1) v_1] \\[1mm]
    \alpha'_{2k}   & = & [x_1 + (k-1) u_1] \, [y_2 + (k-1) v_2 + \wo] \\[1mm]
    \delta'_1      & = & (x_1 - x_2 - \we) y_1 \,+\, \ze \zo  \\[1mm]
    \delta'_{2k+1} & = &
        (x_1 - x_2 - \we) \, [y_1 - y_2 + kv_1 - (k-1)v_2 - \wo]
        \;+  \nonumber \\
     & & \quad
         (u_1 - u_2) \, [ky_1 - (k-1)y_2 + k^2 v_1 - (k-1)^2 v_2 - (k-1)\wo]
          \quad\hbox{for $k \ge 1$}  \nonumber \\  \\
    \delta'_{2k}   & = & 0
 \label{eq.prop.Tfrac.first.alternative.weights}
\end{subeqnarray}
\end{proposition}

\proof
It suffices to verify that the T-fractions
\reff{eq.thm.Tfrac.first}/\reff{eq.thm.Tfrac.first.weights}
and
\reff{eq.prop.Tfrac.first.alternative}/\reff{eq.prop.Tfrac.first.alternative.weights}
contract by Proposition~\ref{prop.contraction_even.Ttype}
to the same J-fraction.
Note that the products forming $\alpha_{2k-1}$ and $\alpha_{2k}$
in \reff{eq.thm.Tfrac.first.weights} have been interchanged
in \reff{eq.prop.Tfrac.first.alternative.weights}:
this guarantees that
$\alpha_{2k-1} \alpha_{2k} = \alpha'_{2k-1} \alpha'_{2k}$.
It is then straightforward to verify that the coefficients $\bdelta'$
have been chosen so that
$\alpha_1 + \delta_1 = \alpha'_1 + \delta'_1$
and
$\alpha_{2k} + \alpha_{2k+1} =
 \alpha'_{2k} + \alpha'_{2k+1} + \delta'_{2k+1}$ for all $k \ge 1$.
\qed

{\bf Remark.} There is also an analogous T-fraction in which the roles
of odd and even indices are reversed.
\myendremark

\medskip

Of course, the T-fraction
\reff{eq.prop.Tfrac.first.alternative}/\reff{eq.prop.Tfrac.first.alternative.weights}
is rather ugly, because of the very complicated expression for $\delta_{2k+1}$.
But it simplifies greatly if we specialize to
$y_1 = y_2 = v_1 = v_2 = \wo = \zo$:
that~is, all the weights associated to odd indices are equal,
and hence might as well be set to~1.
In this case we have:

\begin{corollary}[Alternative T-fraction for D-permutations, specialized]
   \label{cor.Tfrac.first.alternative}
The ordinary generating function of the polynomials \reff{def.Pn}
specialized to $y_1 = y_2 = v_1 = v_2 = \wo = \zo = 1$
has the alternative T-type continued fraction
\begin{eqnarray}
   & & \hspace*{-12mm}
      \sum_{n=0}^{\infty} P_n(x_1,x_2,1,1,u_1,u_2,1,1,\we,1,\ze,1) \: t^n
   \;=\;
      \nonumber \\
   & &
  \cfrac{1}{1 - (x_1 - x_2 + \ze - \we)t - \cfrac{(x_2+\we) t}{1 -  \cfrac{2x_1 t}{1 - 2(u_1 - u_2)t - \cfrac{2(x_2+u_2+\we) t}{1 - \cfrac{3(x_1+u_1) t}{1 - \cdots}}}}}
   \qquad
\label{eq.cor.Tfrac.first.alternative}
\end{eqnarray}
with coefficients
\begin{subeqnarray}
    \alpha'_{2k-1} & = & k \, [x_2 + (k-1)u_2 + \we] \\[1mm]
    \alpha'_{2k}   & = & (k+1) \, [x_1 + (k-1) u_1]  \\[1mm]
    \delta'_1      & = & (x_1 - x_2) \,+\, (\ze - \we)  \\[1mm]
    \delta'_{2k+1} & = & (k+1) (u_1 - u_2)  \quad\hbox{for $k \ge 1$}  \\[1mm]
    \delta'_{2k}   & = & 0
 \label{eq.cor.Tfrac.first.alternative.weights}
\end{subeqnarray}
\end{corollary}

And if we further specialize to $x_1 = x_2$, $u_1 = u_2$, $\ze = \we$,
we obtain an S-fraction:

\begin{corollary}[S-fraction for specialized D-permutations]
   \label{cor.Sfrac.first.specialized}
The ordinary generating function of the polynomials $P^\star_n(x,u,\we)$
defined by
\be
   P^\star_n(x,u,\we)
   \;=\;
   \sum_{\sigma\in \dperm_{2n}}
      x^{\earec(\sigma)} u^{\nrcpeak(\sigma)+\nrcdfall(\sigma)}
      \we^{\evenfix(\sigma)}
 \label{eq.def.Pnstar}
\ee
has the S-type continued fraction
\be
\sum_{n=0}^{\infty} P^\star_n(x,u,\we)\: t^n
   \;=\;
  \cfrac{1}{1 - \cfrac{(x+\we) t}{1 -  \cfrac{2x t}{1 - \cfrac{2(x+u+\we) t}{1 - \cfrac{3(x+u) t}{1 - \cdots}}}}}
\label{eq.Sfrac.Pnstar}
\ee
with coefficients
\begin{subeqnarray}
    \alpha_{2k-1} & = & k \, [x + (k-1)u + \we] \\[1mm]
    \alpha_{2k}   & = & (k+1) \, [x + (k-1) u]
 \label{eq.Sfrac.Pnstar.weights}
\end{subeqnarray}
\end{corollary}

\medskip

{\bf Remarks.}
1.  Specializing to $x = u = \we = 1$ gives the S-fraction
\reff{eq.mediangenocchi.Sfrac2}/\reff{eq.mediangenocchi.Sfrac2.weights}
for the once-shifted median Genocchi numbers.

2.  Specializing to $x = u = 1$, $\we = 0$ gives the S-fraction
\reff{eq.genocchi.Sfrac}/\reff{eq.genocchi.Sfrac.weights}
for the Genocchi numbers.

3.  There is also a slightly weaker specialization
of Corollary~\ref{cor.Tfrac.first.alternative} that leads to an S-fraction,
namely, $x_2 = x_1 + \ze - \we$ and $u_1 = u_2$:
it leads to the S-fraction with
$\alpha_{2k-1} = {k \, [x_1 + (k-1)u_1 + \ze]}$,
$\alpha_{2k} = {(k+1) \, [x_1 + (k-1) u_1]}$.
Curiously this is independent of $\we$:
the direct dependence on $\we$ and the dependence via $x_2$
apparently cancel.

4.  See also the generalization in Corollary~\ref{cor.Sfrac.cyc},
as well as Open Problem~\ref{problem.Dperm.S}.
\myendremark

% {\bf BISHAL: I added this, 5~Dec 2022.}
% %%% NO!!!!! THIS IS WEAKER THAN ONE GETS DIRECTLY FROM THE FIRST T-FRACTION
% Going back to Proposition~\ref{prop.Tfrac.first.alternative},
% there is another specialization that leads to an S-fraction,
% namely, $x_1 = x_2$, $u_1 = u_2$ and $\we = \ze = 0$,
% which corresponds to D-e-semiderangements:
% 
% \begin{corollary}[S-fraction for specialized D-e-semiderangements}
%    \label{cor.Sfrac.first.specialized.Desemi}
% The ordinary generating function of the polynomials
% $P^\circ_n(x,u,y_1,y_2,v_1,v_2,\wo)$ defined by
% \begin{eqnarray}
%    & &
%    P^\circ_n(x,u,y_1,y_2,v_1,v_2,\wo)
%    \;=\;
%        \nonumber \\[4mm]
%    & & \qquad\qquad
%    \sum_{\sigma\in \dperm_{2n}^{\rm e}}
%       x^{\earec(\sigma)} u^{\nrcpeak(\sigma)+\nrcdfall(\sigma)}
%    \:\times
%        \qquad\qquad
%        \nonumber \\[-1mm]
%    & & \qquad\qquad\qquad\:
%       y_1^{\ereccval(\sigma)} y_2^{\ereccdrise(\sigma)}
%       v_1^{\nrcval(\sigma)} v_2^{\nrcdrise(\sigma)}
%       \wo^{\oddnrfix(\sigma)}
%  \label{eq.def.Pncirc}
% \end{eqnarray}
% has the S-type continued fraction
% \be
% \sum_{n=0}^{\infty} P^\circ_n(x,u,y_1,y_2,v_1,v_2,\wo) \:  t^n
%    \;=\;
%   \cfrac{1}{1 - \cfrac{x y_1 t}{1 -  \cfrac{x(y_2+\wo) t}{1 - \cfrac{(x+u)(y_1+v_1) t}{1 - \cfrac{(x+u)(y_2+v_2+\wo) t}{1 - \cdots}}}}}
% \label{eq.Sfrac.Pncirc}
% \ee
% with coefficients
% \begin{subeqnarray}
%     \alpha_{2k-1} & = & [x + (k-1)u] \, [y_1 + (k-1)v_1] \\[1mm]
%     \alpha_{2k}   & = & [x + (k-1) u] \, [y_2 + (k-1)v_2 + \wo]
%  \label{eq.Sfrac.Pncirc.weights}
% \end{subeqnarray}
% \end{corollary}

\subsection[$p,q$-generalizations of the first T-fraction]{$\bm{p,q}$-generalizations of the first T-fraction}
   \label{subsec.first.pq}

We can extend Theorem~\ref{thm.Tfrac.first}
by introducing a $p,q$-generalization.
Recall that for integer $n \ge 0$ we define
\be
   [n]_{p,q}
   \;=\;
   {p^n - q^n \over p-q}
   \;=\;
   \sum\limits_{j=0}^{n-1} p^j q^{n-1-j}
\ee
where $p$ and $q$ are indeterminates;
it is a homogeneous polynomial of degree $n-1$ in $p$ and $q$,
which is symmetric in $p$ and $q$.
In particular, $[0]_{p,q} = 0$ and $[1]_{p,q} = 1$;
and for $n \ge 1$ we have the recurrence
\be
   [n]_{p,q}
   \;=\;
   p \, [n-1]_{p,q} \,+\, q^{n-1}
   \;=\;
   q \, [n-1]_{p,q} \,+\, p^{n-1}
   \;.
 \label{eq.recurrence.npq}
\ee
If $p=1$, then $[n]_{1,q}$ is the well-known $q$-integer
\be
   [n]_q
   \;=\;  [n]_{1,q}
   \;=\; {1 - q^n \over 1-q}
   \;=\;  \begin{cases}
               0  & \textrm{if $n=0$}  \\
               1+q+q^2+\ldots+q^{n-1}  & \textrm{if $n \ge 1$}
          \end{cases}
\ee
If $p=0$, then
\be
   [n]_{0,q}
   \;=\;  \begin{cases}
               0  & \textrm{if $n=0$}  \\
               q^{n-1}  & \textrm{if $n \ge 1$}
          \end{cases}
\ee

The statistics on permutations corresponding to the variables $p$ and $q$
will be crossings and nestings,
as defined in Section~\ref{subsec.statistics.2}.
More precisely, we define the following polynomial in 22 variables
that generalizes \reff{def.Pn} by including
four pairs of $(p,q)$-variables
corresponding to the four refined types of crossings and nestings
illustrated in Figure~\ref{fig.refined_crossnest},
as well as two variables corresponding to pseudo-nestings of fixed points:
\begin{eqnarray}
   & &
   \hspace*{-14mm}
   P_n(x_1,x_2,y_1,y_2,u_1,u_2,v_1,v_2,\we,\wo,\ze,\zo,p_{-1},p_{-2},p_{+1},p_{+2},q_{-1},q_{-2},q_{+1},q_{+2},\se,\so)
   \;=\;
   \hspace*{-1cm}
       \nonumber \\[4mm]
   & & \qquad\qquad
   \sum_{\sigma \in \dperm_{2n}}
   x_1^{\eareccpeak(\sigma)} x_2^{\eareccdfall(\sigma)}
   y_1^{\ereccval(\sigma)} y_2^{\ereccdrise(\sigma)}
   \:\times
       \qquad\qquad
       \nonumber \\[-1mm]
   & & \qquad\qquad\qquad\:
   u_1^{\nrcpeak(\sigma)} u_2^{\nrcdfall(\sigma)}
   v_1^{\nrcval(\sigma)} v_2^{\nrcdrise(\sigma)}
   \:\times
       \qquad\qquad
       \nonumber \\[3mm]
   & & \qquad\qquad\qquad\:
   \we^{\evennrfix(\sigma)} \wo^{\oddnrfix(\sigma)}
   \ze^{\evenrar(\sigma)} \zo^{\oddrar(\sigma)}
   \:\times
       \qquad\qquad
       \nonumber \\[3mm]
   & & \qquad\qquad\qquad\:
   p_{-1}^{\lcrosscpeak(\sigma)}
   p_{-2}^{\lcrosscdfall(\sigma)}
   p_{+1}^{\ucrosscval(\sigma)}
   p_{+2}^{\ucrosscdrise(\sigma)}
          \:\times
       \qquad\qquad
       \nonumber \\[3mm]
   & & \qquad\qquad\qquad\:
   q_{-1}^{\lnestcpeak(\sigma)}
   q_{-2}^{\lnestcdfall(\sigma)}
   q_{+1}^{\unestcval(\sigma)}
   q_{+2}^{\unestcdrise(\sigma)}
          \:\times
       \qquad\qquad
       \nonumber \\[3mm]
   & & \qquad\qquad\qquad\:
   \se^{\epsnest(\sigma)}
   \so^{\opsnest(\sigma)}
 \label{def.Pn.pq}
\end{eqnarray}
where
\be
   \epsnest(\sigma)
   \;=\;
   \sum_{i \in \Evenfix} \psnest(i,\sigma)
   \;,\qquad
   \opsnest(\sigma)
   \;=\;
   \sum_{i \in \Oddfix} \psnest(i,\sigma)
   \;.
\ee
We remark that \reff{def.Pn.pq} is essentially the same as the polynomial
introduced in \cite[eq.~(2.51)]{Sokal-Zeng_masterpoly},
but restricted to D-permutations
and refined to record the parity of fixed points.\footnote{
   The polynomial in \cite[eq.~(2.51)]{Sokal-Zeng_masterpoly}
   also included a more refined stratification of fixed points
   by level (= psnest) as defined in \reff{def.psnestj} above.
   That refined stratification is omitted here for simplicity ---
   instead we include only the simple factors $s^\psnest$ ---
   but it is included in \reff{def.Qn.firstmaster} below.
}

We then have the following $p,q$-generalization of
Theorem~\ref{thm.Tfrac.first}:

\begin{theorem}[First T-fraction for D-permutations, $p,q$-generalization]
   \label{thm.Tfrac.first.pq}
The ordinary generating function of the polynomials \reff{def.Pn.pq} has the
T-type continued fraction
\begin{eqnarray}
   & & \hspace*{-7mm}
\Scale[0.92]{
   \sum\limits_{n=0}^\infty
   P_n(x_1,x_2,y_1,y_2,u_1,u_2,v_1,v_2,\we,\wo,\ze,\zo,p_{-1},p_{-2},p_{+1},p_{+2},q_{-1},q_{-2},q_{+1},q_{+2},\se,\so)  \: t^n
   \;=\;
}
       \nonumber \\[2mm]
   & & \hspace*{-3mm}
\Scale[0.87]{
   \cfrac{1}{1 - \ze \zo  \,t - \cfrac{x_1 y_1  \,t}{1 -  \cfrac{(x_2\!+\!\se\we)(y_2\!+\!\so\wo) \,t}{1 - \cfrac{(p_{-1}x_1\!+\!q_{-1}u_1)(p_{+1}y_1\!+\!q_{+1}v_1)  \,t}{1 - \cfrac{(p_{-2}x_2\!+\!q_{-2}u_2\!+\!\se^2\we)(p_{+2}y_2\!+\!q_{+2}v_2\!+\!\so^2\wo)  \,t}{1 - \cfrac{(p_{-1}^2 x_1\!+\! q_{-1} [2]_{p_{-1},q_{-1}}u_1)(p_{+1}^2 y_1\!+\! q_{+1} [2]_{p_{+1},q_{+1}}v_1)  \,t}{1 - \cfrac{(p_{-2}^2 x_2\!+\! q_{-2} [2]_{p_{-2},q_{-2}} u_2\!+\!\se^3\we)(p_{+2}^2 y_2\!+\! q_{+2} [2]_{p_{+2},q_{+2}}v_2\!+\!\so^3\wo)  \,t}{1 - \cdots}}}}}}}
}
       \nonumber \\[1mm]
   \label{eq.thm.Tfrac.first.pq}
\end{eqnarray}
with coefficients
\begin{subeqnarray}
   \alpha_{2k-1} & = &
      \bigl( p_{-1}^{k-1} x_1 + q_{-1} [k-1]_{p_{-1},q_{-1}} u_1 \bigr) \:
      \bigl( p_{+1}^{k-1} y_1 + q_{+1} [k-1]_{p_{+1},q_{+1}} v_1 \bigr)
        \\[2mm]
   \alpha_{2k}   & = &
      \bigl( p_{-2}^{k-1} x_2 + q_{-2} [k-1]_{p_{-2},q_{-2}} u_2 + \se^k \we \bigr) \:
      \bigl( p_{+2}^{k-1} y_2 + q_{+2} [k-1]_{p_{+2},q_{+2}} v_2 + \so^k \wo \bigr)
        \nonumber \\ \\
   \delta_1  & = &   \ze \zo   \\[1mm]
   \delta_n  & = &   0    \qquad\hbox{for $n \ge 2$}
   \label{eq.thm.Tfrac.first.weights.pq}
\end{subeqnarray}
\end{theorem}

\noindent
We will prove Theorem~\ref{thm.Tfrac.first.pq} in Section~\ref{sec.proofs.1}.
Of course we reobtain Theorem~\ref{thm.Tfrac.first} by
making the specialization
$p_{-1} = p_{-2} = p_{+1} = p_{+2} = q_{-1} = q_{-2} = q_{+1} = q_{+2} =
 \se = \so = 1$.

\bigskip

{\bf Remarks.}
1.  If we specialize to $x_1 = u_1$, $y_1 = v_1$, $x_2 = u_2$, $y_2 = v_2$
--- that is, renounce the counting of records ---
then the coefficients \reff{eq.thm.Tfrac.first.weights.pq} simplify to
\begin{subeqnarray}
   \alpha_{2k-1} & = &
      [k]_{p_{-1},q_{-1}} \, [k]_{p_{+1},q_{+1}} \, x_1 y_1
        \\[1mm]
   \alpha_{2k}   & = &
      \bigl( [k]_{p_{-2},q_{-2}} x_2 + \se^k \we \bigr) \:
      \bigl( [k]_{p_{+2},q_{+2}} y_2 + \so^k \wo \bigr)  \\[1mm]
   \delta_1  & = &   \ze \zo   \\[1mm]
   \delta_n  & = &   0    \qquad\hbox{for $n \ge 2$}
   \label{eq.thm.Tfrac.first.weights.pq.bis}
\end{subeqnarray}

2.  If we further specialize to
$x_1 = 1$, $x_2 = y_1 = y_2 = q$,
$p_{-1} = p_{-2} = p_{+1} = p_{+2} = q$
and $q_{-1} = q_{-2} = q_{+1} = q_{+2} = \se = \so = q^2$
and recall \cite[Proposition~2.24]{Sokal-Zeng_masterpoly} that
the number of inversions ($\inv$) of a permutation satisfies
\be
   \inv
   \;=\;
   \cval + \cdrise + \cdfall + \ucross + \lcross
                    + 2(\unest + \lnest + \psnest)
   \;,
\ee
we obtain a T-fraction for D-permutations according to
the number of even and odd fixed points and number of inversions,
with coefficients
$\alpha_{2k-1} = q^{2k-1} [k]_q^2$,
$\alpha_{2k} = q^{2k} ([k]_q + q^k \we) ([k]_q + q^k \wo)$,
$\delta_1 = \ze \zo$.

3. And if we further specialize to $\wo = \zo  = 0$,
we reobtain the S-fraction for D-o-semiderangements
according to the number of fixed points and number of inversions
\cite[Th\'eor\`eme~1.2]{Randrianarivony_97},
with coefficients
$\alpha_{2k-1} = q^{2k-1} [k]_q^2$,
$\alpha_{2k} = q^{2k} [k]_q ([k]_q + q^k \we)$.
\myendremark

\bigskip

We refrain from pursuing the $p,q$-generalizations
of Proposition~\ref{prop.Tfrac.first.alternative}
and Corollary~\ref{cor.Tfrac.first.alternative},
and limit ourselves to giving the $p,q$-generalization
of Corollary~\ref{cor.Sfrac.first.specialized}.
We need to make the specializations
$y_1 = y_2 = v_1 = v_2 = \wo = \zo = p_{+1} = p_{+2} = q_{+1} = q_{+2} = \so = 1$
--- that~is, set all the weights associated to odd indices to 1 ---
and also specialize
$x_1 = x_2$, $u_1 = u_2$, $\ze = \we$, $p_{-1} = p_{-2}$,
$q_{-1} = q_{-2}$ and $\se = 1$:

\begin{corollary}[S-fraction for specialized D-permutations, $p,q$-generalization]
   \label{cor.Sfrac.first.specialized.pq}
The ordinary generating function of the polynomials
$P^\star_n(x,u,\we,p_-,q_-)$ defined by
\be
   P^\star_n(x,u,\we,p_-,q_-)
   \;=\;
   \sum_{\sigma\in \dperm_{2n}}
      x^{\earec(\sigma)} u^{\nrcpeak(\sigma)+\nrcdfall(\sigma)}
      \we^{\evenfix(\sigma)} p_-^{\lcross(\sigma)} q_-^{\lnest(\sigma)}
 \label{eq.def.Pnstar.pq}
\ee
has the S-type continued fraction
\be
\sum_{n=0}^{\infty} P^\star_n(x,u,\we)\: t^n
   \;=\;
  \cfrac{1}{1 - \cfrac{(x+\we) t}{1 -  \cfrac{2x t}{1 - \cfrac{2(p_- x+ q_- u+\we) t}{1 - \cfrac{3(p_- x+ q_- u) t}{1 - \cdots}}}}}
\label{eq.Sfrac.Pnstar.pq}
\ee
with coefficients
\begin{subeqnarray}
    \alpha_{2k-1} & = &
       k \, \bigl( p_{-}^{k-1} x + q_{-} [k-1]_{p_{-},q_{-}} u + \we \bigr)
          \\[1mm]
    \alpha_{2k}   & = &
      (k+1) \, \bigl( p_{-}^{k-1} x + q_{-} [k-1]_{p_{-},q_{-}} u \bigr)
 \label{eq.Sfrac.Pnstar.weights.pq}
\end{subeqnarray}
\end{corollary}

\proof
It suffices to verify that the T-fraction
\reff{eq.thm.Tfrac.first.pq}/\reff{eq.thm.Tfrac.first.weights.pq}
with the given specializations
and the S-fraction
\reff{eq.Sfrac.Pnstar.pq}/\reff{eq.Sfrac.Pnstar.weights.pq}
contract by Proposition~\ref{prop.contraction_even.Ttype}
to the same J-fraction.
\qed

\subsection{First master T-fraction}   \label{subsec.first.master}

In fact, we can go farther,
and introduce a polynomial in six infinite families of indeterminates
$\bsfa = (\sfa_{\ell,\ell'})_{\ell,\ell' \ge 0}$,
$\bsfb = (\sfb_{\ell,\ell'})_{\ell,\ell' \ge 0}$,
$\bsfc = (\sfc_{\ell,\ell'})_{\ell,\ell' \ge 0}$,
$\bsfd = (\sfd_{\ell,\ell'})_{\ell,\ell' \ge 0}$,
$\bsfe = (\sfe_\ell)_{\ell \ge 0}$,
$\bsff = (\sff_\ell)_{\ell \ge 0}$
that will have a nice T-fraction
and that will include the polynomials \reff{def.Pn} and \reff{def.Pn.pq}
as specializations.

Using the index-refined crossing and nesting statistics
defined in \reff{def.ucrossnestjk},
we define the polynomial $Q_n(\bsfa,\bsfb,\bsfc,\bsfd,\bsfe,\bsff)$ by
\begin{eqnarray}
   & & \hspace*{-10mm}
   Q_n(\bsfa,\bsfb,\bsfc,\bsfd,\bsfe,\bsff)
   \;=\;
       \nonumber \\[4mm]
   & &
   \sum_{\sigma \in \dperm_{2n}}
   \;\:
   \prod\limits_{i \in \Cval(\sigma)}  \! \sfa_{\ucross(i,\sigma),\,\unest(i,\sigma)}
   \prod\limits_{i \in \Cpeak(\sigma)} \!\!  \sfb_{\lcross(i,\sigma),\,\lnest(i,\sigma)}
       \:\times
       \qquad\qquad
       \nonumber \\[1mm]
   & & \qquad\;\;\,
   \prod\limits_{i \in \Cdfall(\sigma)} \!\!  \sfc_{\lcross(i,\sigma),\,\lnest(i,\sigma)}
   \;
   \prod\limits_{i \in \Cdrise(\sigma)} \!\!  \sfd_{\ucross(i,\sigma),\,\unest(i,\sigma)}
       \:\times
       \qquad\qquad
       \nonumber \\[1mm]
   & & \qquad\;\:
   \prod\limits_{i \in \Evenfix(\sigma)} \!\!\! \sfe_{\psnest(i,\sigma)}
   \;
   \prod\limits_{i \in \Oddfix(\sigma)}  \!\!\! \sff_{\psnest(i,\sigma)}
   \;.
   \quad
 \label{def.Qn.firstmaster}
\end{eqnarray}
where $\Cval(\sigma) = \{ i\colon \sinv(i) > i < \sigma(i) \}$
and likewise for the others.
We remark that \reff{def.Qn.firstmaster} is the same as the polynomial
introduced in \cite[eq.~(2.77)]{Sokal-Zeng_masterpoly},
but restricted to D-permutations
and refined to record the parity of fixed points.

The polynomials \reff{def.Qn.firstmaster} then have a beautiful T-fraction:

\begin{theorem}[First master T-fraction for D-permutations]
   \label{thm.Tfrac.first.master}
The ordinary generating function of the polynomials
$Q_n(\bsfa,\bsfb,\bsfc,\bsfd,\bsfe,\bsff)$
has the T-type continued fraction
\be
   \sum_{n=0}^\infty Q_n(\bsfa,\bsfb,\bsfc,\bsfd,\bsfe,\bsff) \: t^n
   \;=\;
\Scale[0.95]{
   \cfrac{1}{1 - \sfe_0 \sff_0 t - \cfrac{\sfa_{00} \sfb_{00} t}{1 - \cfrac{(\sfc_{00} + \sfe_1)(\sfd_{00} + \sff_1) t}{1 - \cfrac{(\sfa_{01} + \sfa_{10})(\sfb_{01} + \sfb_{10}) t}{1 - \cfrac{(\sfc_{01} + \sfc_{10} + \sfe_2)(\sfd_{01} + \sfd_{10} + \sff_2)t}{1 - \cdots}}}}}
}
   \label{eq.thm.Tfrac.first.master}
\ee
with coefficients
\begin{subeqnarray}
   \alpha_{2k-1}
   & = &
   \biggl( \sum_{\xi=0}^{k-1} \sfa_{k-1-\xi,\xi} \biggr)
   \biggl( \sum_{\xi=0}^{k-1} \sfb_{k-1-\xi,\xi} \biggr)
        \\[2mm]
   \alpha_{2k}
   & = &
   \biggl( \sfe_k \,+\, \sum_{\xi=0}^{k-1} \sfc_{k-1-\xi,\xi} \biggr)
   \biggl( \sff_k \,+\, \sum_{\xi=0}^{k-1} \sfd_{k-1-\xi,\xi} \biggr)
        \\[2mm]
   \delta_1  & = &  \sfe_0 \sff_0 \\[2mm]
   \delta_n  & = &   0    \qquad\hbox{for $n \ge 2$}
 \label{def.weights.Tfrac.first.master}
\end{subeqnarray}
\end{theorem}

\noindent
We will prove this theorem in Section~\ref{sec.proofs.1};
it is our ``master'' version of the first T-fraction.
It implies Theorems~\ref{thm.Tfrac.first} and \ref{thm.Tfrac.first.pq}
by straightforward specializations.

\subsection{Variant forms of the first T-fractions}
   \label{subsec.first.variant}

Our first T-fractions
(Theorems~\ref{thm.Tfrac.first}, \ref{thm.Tfrac.first.pq}
 and \ref{thm.Tfrac.first.master})
have variant forms
in which we use the variant index-refined crossing and nesting statistics
\reff{def.ucrossnestjk.prime}
in place of the original index-refined crossing and nesting statistics
\reff{def.ucrossnestjk}.

It is convenient to start with the master T-fraction.
Introducing indeterminates $\bsfa,\bsfb,\bsfc,\bsfd,\bsfe,\bsff$ as before,
we define the variant polynomial $Q'_n(\bsfa,\bsfb,\bsfc,\bsfd,\bsfe,\bsff)$ by
\begin{eqnarray}
   & & \hspace*{-10mm}
   Q'_n(\bsfa,\bsfb,\bsfc,\bsfd,\bsfe,\bsff)
   \;=\;
       \nonumber \\[4mm]
   & &
   \sum_{\sigma \in \dperm_{2n}}
   \;\:
   \prod\limits_{i \in \Cval(\sigma)}  \! \sfa_{\lcross'(i,\sigma),\,\lnest'(i,\sigma)}
   \prod\limits_{i \in \Cpeak(\sigma)} \!\!  \sfb_{\ucross'(i,\sigma),\,\unest'(i,\sigma)}
       \:\times
       \qquad\qquad
       \nonumber \\[1mm]
   & & \qquad\;\;\,
   \prod\limits_{i \in \Cdfall(\sigma)} \!\!  \sfc_{\lcross'(i,\sigma),\,\lnest'(i,\sigma)}
   \;
   \prod\limits_{i \in \Cdrise(\sigma)} \!\!  \sfd_{\ucross'(i,\sigma),\,\unest'(i,\sigma)}
       \:\times
       \qquad\qquad
       \nonumber \\[1mm]
   & & \qquad\;\:
   \prod\limits_{i \in \Evenfix(\sigma)} \!\!\! \sfe_{\psnest(i,\sigma)}
   \;
   \prod\limits_{i \in \Oddfix(\sigma)}  \!\!\! \sff_{\psnest(i,\sigma)}
   \;.
   \quad
 \label{def.Qn.firstmaster.variant}
\end{eqnarray}
Note that for cycle valleys and cycle peaks,
the ``u'' and ``l'' have here been interchanged relative to
\reff{def.Qn.firstmaster};
this is in accordance with the explanation after
\reff{def.ucrossnestjk.prime}
about which types of indices can have nonzero values
for the primed crossing and nesting statistics.

We then have the following variant of Theorem~\ref{thm.Tfrac.first.master}:

\begin{theorem}[Variant first master T-fraction for D-permutations]
   \label{thm.Tfrac.first.master.variant}
The ordinary generating function of the polynomials
$Q'_n(\bsfa,\bsfb,\bsfc,\bsfd,\bsfe,\bsff)$
has the same T-type continued fraction
\reff{eq.thm.Tfrac.first.master}/\reff{def.weights.Tfrac.first.master}
as the polynomials $Q_n(\bsfa,\bsfb,\bsfc,\bsfd,\bsfe,\bsff)$.
Therefore
\be
   Q_n(\bsfa,\bsfb,\bsfc,\bsfd,\bsfe,\bsff)
   \;=\;
   Q'_n(\bsfa,\bsfb,\bsfc,\bsfd,\bsfe,\bsff)
   \;.
\ee
\end{theorem}

\noindent
We will prove this theorem in Section~\ref{subsec.proofs.alternative}.

\bigskip

Next we define the variant $p,q$-generalized polynomials $P'_n$:
\begin{eqnarray}
   & &
   \hspace*{-14mm}
   P'_n(x_1,x_2,y_1,y_2,u_1,u_2,v_1,v_2,\we,\wo,\ze,\zo,p_{-1},p_{-2},p_{+1},p_{+2},q_{-1},q_{-2},q_{+1},q_{+2},\se,\so)
   \;=\;
   \hspace*{-1cm}
       \nonumber \\[4mm]
   & & \qquad\qquad
   \sum_{\sigma \in \dperm_{2n}}
   x_1^{\ereccpeak'(\sigma)} x_2^{\eareccdfall'(\sigma)}
   y_1^{\eareccval'(\sigma)} y_2^{\ereccdrise'(\sigma)}
   \:\times
       \qquad\qquad
       \nonumber \\[-1mm]
   & & \qquad\qquad\qquad\:
   u_1^{\nrcpeak'(\sigma)} u_2^{\nrcdfall'(\sigma)}
   v_1^{\nrcval'(\sigma)} v_2^{\nrcdrise'(\sigma)}
   \:\times
       \qquad\qquad
       \nonumber \\[3mm]
   & & \qquad\qquad\qquad\:
   \we^{\evennrfix(\sigma)} \wo^{\oddnrfix(\sigma)}
   \ze^{\evenrar(\sigma)} \zo^{\oddrar(\sigma)}
   \:\times
       \qquad\qquad
       \nonumber \\[3mm]
   & & \qquad\qquad\qquad\:
   p_{-1}^{\ucrosscpeak'(\sigma)}
   p_{-2}^{\lcrosscdfall'(\sigma)}
   p_{+1}^{\lcrosscval'(\sigma)}
   p_{+2}^{\ucrosscdrise'(\sigma)}
          \:\times
       \qquad\qquad
       \nonumber \\[3mm]
   & & \qquad\qquad\qquad\:
   q_{-1}^{\unestcpeak'(\sigma)}
   q_{-2}^{\lnestcdfall'(\sigma)}
   q_{+1}^{\lnestcval'(\sigma)}
   q_{+2}^{\unestcdrise'(\sigma)}
          \:\times
       \qquad\qquad
       \nonumber \\[3mm]
   & & \qquad\qquad\qquad\:
   \se^{\epsnest(\sigma)}
   \so^{\opsnest(\sigma)}
 \label{def.Pn.pq.variant}
\end{eqnarray}
where, for example,
the statistics $\ereccpeak'$ and $\nrcpeak'$ measure
whether a cycle peak $i$ is a record {\em value}\/,
i.e.~whether $\sinv(i)$ is a record position.
That is, we define
\begin{subeqnarray}
   \ereccpeak'(\sigma)
   & = &
   \#\{ i \colon\: \sinv(i) < i > \sigma(i) \hbox{ and }
           \sinv(i) \hbox{ is a record}  \}
      \\[2mm]
   \nrcpeak'(\sigma)
   & = &
   \#\{ i \colon\: \sinv(i) < i > \sigma(i) \hbox{ and }
           \sinv(i) \hbox{ is not a record}  \}
     \qquad
\end{subeqnarray}
and likewise for the others.
(Since a record position is necessarily a weak excedance,
 which here means $\sinv(i) \le i$,
 erec here corresponds to cpeak and cdrise;
 likewise earec corresponds to cval and cdfall.)
Similarly, we define
\begin{subeqnarray}
\ucrosscpeak'(\sigma) & = & \sum_{k\in \Cpeak(\sigma)} \ucross'(k,\sigma) 
   \\[2mm]
\unestcpeak'(\sigma) & = & \sum_{k\in \Cpeak(\sigma)} \unest'(k,\sigma) 
\end{subeqnarray}
where $\ucross'$ and $\unest'$ have been defined in~(\ref{def.ucrossnestjk.prime}a,b),
and likewise for the others.

We then have the following variant of Theorem~\ref{thm.Tfrac.first.pq}:

\begin{theorem}[Variant first T-fraction for D-permutations, $p,q$-generalization]
   \label{thm.Tfrac.first.pq.variant}
The ordinary generating function of the polynomials $P'_n$
has the same T-type continued fraction
\reff{eq.thm.Tfrac.first.pq}/\reff{eq.thm.Tfrac.first.weights.pq}
as the polynomials $P_n$.
Therefore
\begin{eqnarray}
   & &  \hspace*{-7mm}
\Scale[0.96]{
   P_n(x_1,x_2,y_1,y_2,u_1,u_2,v_1,v_2,\we,\wo,\ze,\zo,p_{-1},p_{-2},p_{+1},p_{+2},q_{-1},q_{-2},q_{+1},q_{+2},\se,\so)
}
            \nonumber \\
   & &  \hspace*{-7mm}
\Scale[0.96]{
   =\;
   P'_n(x_1,x_2,y_1,y_2,u_1,u_2,v_1,v_2,\we,\wo,\ze,\zo,p_{-1},p_{-2},p_{+1},p_{+2},q_{-1},q_{-2},q_{+1},q_{+2},\se,\so)
   \;.
}
            \nonumber \\
\end{eqnarray}
\end{theorem}

And finally, we can define variant polynomials
$P'_n(x_1,x_2,y_1,y_2,u_1,u_2,v_1,v_2,\we,\wo,$ $\ze,\zo)$
by specializing
$p_{-1} = p_{-2} = p_{+1} = p_{+2} = q_{-1} = q_{-2} = q_{+1} = q_{+2} =
 \se = \so = 1$,
and thereby obtain a variant version of Theorem~\ref{thm.Tfrac.first}.
We leave the details to the reader.

\bigskip

{\bf Remark.}
A four-variable special case of Theorem~\ref{thm.Tfrac.first.pq.variant}
was found by Randrianarivony and Zeng
\cite[Proposition~10]{Randrianarivony_96b}
for D-o-semiderangements.
Note first that for a D-o-semiderangement,
each fixed point (which is necessarily even)
must be a neither-record-antirecord:
for if a fixed point $i$ were a record or antirecord, then it would be both;
but by Lemma~\ref{lemma.recantirec_D},
$i\!-\!1$ would then also be a record-antirecord fixed point,
which is impossible since $\sigma$ is a D-o-semiderangement.

Let us now look at the statistics defined by Randrianarivony and Zeng
\cite[p.~2]{Randrianarivony_96b}.
Their statistic ``lema'' is a left-to-right maximum (i.e.~a record)
whose {\em value}\/ is even:
thus we can say that $\sinv(i)$ is a record {\em position}\/
and that $i$ is even.
But a record is always a weak excedance, hence $\sinv(i) \le i$;
and $i$ is even, hence $i \ge \sigma(i)$.
Since $i$ is a record, it cannot be a fixed point,
so it must be a cycle peak, and we have ${\rm lema} = \ereccpeak'$.
Similarly, their statistic ``remi'' (resp.~``romi'')
is a right-to-left minimum (i.e.~an antirecord)
whose {\em value}\/ is even (resp.~odd);
by similar reasoning we obtain
${\rm remi} = \eareccdfall'$ and ${\rm romi} = \eareccval'$.
Their polynomial \cite[eq.~(3.3)]{Randrianarivony_96b}
\be
   R_n(x,y,\bar{x},\bar{y})
   \;=\;
   \sum_{\sigma \in \dperm^{\rm o}_{2n}}
      x^{{\rm lema}(\sigma)}
      y^{{\rm romi}(\sigma)}
      \bar{x}^{{\rm fix}(\sigma)}
      \bar{y}^{{\rm remi}(\sigma)}
\ee
thus corresponds to our $P'_n$
specialized to $x_1 = x$, $y_1 = y$, $\we = \bar{x}$, $x_2 = \bar{y}$,
$\wo = \zo = 0$, and all other variables equal to 1.
Their S-fraction \cite[Proposition~10]{Randrianarivony_96b}
is then a special case of Theorem~\ref{thm.Tfrac.first.pq.variant}.
\myendremark

\section{Second T-fraction and its generalizations}   \label{sec.second}

It is natural to want to refine the foregoing polynomials
by keeping track also of the number of cycles (cyc).
Unfortunately, $\cyc$ does not seem to mesh well with
the record classification: even the three-variable polynomials 
\be
   P_n(x,y,\lambda)
   \;=\;
   \sum_{\sigma\in \dperm_{2n}}
       x^{\arec(\sigma)}y^{\erec(\sigma)} \lambda^{\cyc(\sigma)}
 \label{eq.Pn.xylam}
\ee
do not have a J-fraction with polynomial coefficients
(see Appendix~\ref{app.Pnxylam}).
%% {\bf Give first few $P_n$ and first few $\beta,\gamma$
%%    up to the first nonpolynomial one!!!}
Using the contraction formula (Proposition~\ref{prop.contraction_even.Ttype}),
it follows from this that $P_n(x,y,\lambda)$
also cannot have a T-fraction with polynomial coefficients
and $\delta_2 = \delta_4 = \ldots = 0$.\footnote{
   This does not exclude that $P_n(x,y,\lambda)$ might have
   a general T-fraction with polynomial coefficients.
   Indeed, it is easy to see \cite{Sokal_totalpos} that,
   for {\em any}\/ sequence $\ba = (a_n)_{n \ge 0}$ with $a_0 = 1$
   in a commutative ring $R$,
   and {\em any}\/ sequence $\balpha = (\alpha_n)_{n \ge 1}$
   of {\em invertible}\/ elements of $R$,
   there exists a unique sequence $\bdelta = (\delta_n)_{n \ge 1}$ in $R$
   such that the ordinary generating function $\sum_{n=0}^\infty a_n t^n$
   is represented by the T-fraction with coefficients $\balpha$ and $\bdelta$.
   So the polynomials $P_n(x,y,\lambda)$ are in fact represented by
   {\em uncountably many}\/ different T-fractions
   with coefficients $\alpha_n \in \Z \setminus \{0\}$
   (or even $\alpha_n \in \{1,2\}$)
   and $\delta_n \in \Q[x,y,\lambda]$.
   Furthermore, such T-fractions might exist also for certain
   noninvertible $\balpha$, subject to suitable divisibility conditions.
   But it is far from clear whether any of these T-fractions
   are simple enough to be found explicitly.
 \label{footnote_Tfrac_nonunique}
}

Nevertheless, it turns out that cyc {\em almost}\/ meshes with the
complete parity-refined record-and-cycle classification:
it suffices to make one (conjecturally) or two (provably)
specializations in which we partially renounce recording the record status.
We explain this in the next subsection.

\subsection{Second T-fraction}   \label{subsec.second.statements}

\subsubsection{Conjectured T-fraction and main theorem}

We begin by introducing a polynomial in 13~variables
that generalizes \reff{def.Pn}
by keeping track also of the number of cycles (cyc):
\begin{eqnarray}
   & &
   \Scale[0.97]{
   \widehat{P}_n(x_1,x_2,y_1,y_2,u_1,u_2,v_1,v_2,\we,\wo,\ze,\zo,\lambda)}
   \;=\;
       \nonumber \\[4mm]
   & & \qquad\qquad
   \sum_{\sigma \in \dperm_{2n}}
   x_1^{\eareccpeak(\sigma)} x_2^{\eareccdfall(\sigma)}
   y_1^{\ereccval(\sigma)} y_2^{\ereccdrise(\sigma)}
   \:\times
       \qquad\qquad
       \nonumber \\[-1mm]
   & & \qquad\qquad\qquad\:
   u_1^{\nrcpeak(\sigma)} u_2^{\nrcdfall(\sigma)}
   v_1^{\nrcval(\sigma)} v_2^{\nrcdrise(\sigma)}
   \:\times
       \qquad\qquad
       \nonumber \\[3mm]
   & & \qquad\qquad\qquad\:
   \we^{\evennrfix(\sigma)} \wo^{\oddnrfix(\sigma)}
   \ze^{\evenrar(\sigma)} \zo^{\oddrar(\sigma)}
   \, \lambda^{\cyc(\sigma)}
   \;.
 \label{def.Pnhat}
\end{eqnarray}
Of course there is no hope that $\widehat{P}_n$ has a J-fraction
(or a T-fraction with $\delta_2 = \delta_4 = \ldots = 0$)
with polynomial coefficients,
because even the specialization \reff{eq.Pn.xylam} does not have one.
Nevertheless, we find {\em empirically}\/ that we need to make
only one specialization
--- either $u_1 = x_1$ or $v_1 = y_1$ ---
to obtain a good T-fraction.
In other words, it suffices to renounce distinguishing
either the antirecord status of cycle peaks
or the record status of cycle valleys.
For concreteness we show the second of these conjectures:

\begin{conjecture}
   \label{conj.Tfrac.second}
The ordinary generating function of the polynomials \reff{def.Pnhat}
specialized to $v_1 = y_1$ has the T-type continued fraction
\begin{eqnarray}
   & & \hspace*{-12mm}
   \sum_{n=0}^\infty
   \widehat{P}_n(x_1,x_2,y_1,y_2,u_1,u_2,y_1,v_2,\we,\wo,\ze,\zo,\lambda) \: t^n
   \;=\;
       \nonumber \\
   & &
   \cfrac{1}{1 - \lambda^2 \ze \zo  \,t - \cfrac{\lambda x_1 y_1  \,t}{1 -  \cfrac{(x_2\!+\!\lambda\we)(y_2\!+\!\lambda\wo) \,t}{1 - \cfrac{(\lambda+1)(x_1\!+\!u_1) y_1 \,t}{1 - \cfrac{(x_2\!+\!u_2\!+\!\lambda\we)(y_2\!+\!v_2\!+\!\lambda\wo)  \,t}{1 - \cfrac{(\lambda+2)(x_1\!+\!2u_1) y_1 \,t}{1 - \cfrac{(x_2\!+\!2u_2\!+\!\lambda\we)(y_2\!+\!2v_2\!+\!\lambda\wo)  \,t}{1 - \cdots}}}}}}}
       \nonumber \\[1mm]
   \label{eq.conj.Tfrac.second}
\end{eqnarray}
with coefficients
\begin{subeqnarray}
   \alpha_{2k-1} & = & (\lambda + k-1) \: [x_1 + (k-1) u_1] \: y_1
        \\[1mm]
   \alpha_{2k}   & = & [x_2 + (k-1) u_2 + \lambda\we] \: [y_2 + (k-1)v_2 + \lambda\wo]
        \\[1mm]
   \delta_1  & = &   \lambda^2 \ze \zo   \\[1mm]
   \delta_n  & = &   0    \qquad\hbox{for $n \ge 2$}
   \label{eq.conj.Tfrac.second.weights}
\end{subeqnarray}
\end{conjecture}

\noindent
We have tested this conjecture through $n=6$.

Alas, we are unable at present to prove Conjecture~\ref{conj.Tfrac.second};
we are only able to prove the weaker version
in which we make the {\em two}\/ specializations
$v_1 = y_1$ {\em and}\/ $v_2 = y_2$,
i.e.~we renounce distinguishing
the record status of cycle valleys {\em and}\/ cycle double rises.
However, we can do better than this:
namely, we replace the pair $y_2,v_2$
by a pair of new variables $\yhat_2,\vhat_2$ that measure,
not whether a cycle double rise $i$ is a record {\em position}\/,
but rather whether it is a record {\em value}\/,
i.e.~whether $\sinv(i)$ is a record.
That is, we define
\begin{subeqnarray}
   \ereccdrise'(\sigma)
   & = &
   \#\{ i \colon\: \sinv(i) < i < \sigma(i) \hbox{ and }
           \sinv(i) \hbox{ is a record}  \}
      \\[2mm]
   \nrcdrise'(\sigma)
   & = &
   \#\{ i \colon\: \sinv(i) < i < \sigma(i) \hbox{ and }
           \sinv(i) \hbox{ is not a record}  \}
     \qquad
\end{subeqnarray}
[These statistics have already been used in \reff{def.Pn.pq.variant}.]
We then define the modified polynomials
\begin{eqnarray}
   & &
   \widehat{\widehat{P}}_n(x_1,x_2,y_1,\yhat_2,u_1,u_2,v_1,\vhat_2,\we,\wo,\ze,\zo,\lambda)
   \;=\;
       \nonumber \\[4mm]
   & & \qquad\qquad
   \sum_{\sigma \in \dperm_{2n}}
   x_1^{\eareccpeak(\sigma)} x_2^{\eareccdfall(\sigma)}
   y_1^{\ereccval(\sigma)} \yhat_2^{\ereccdrise'(\sigma)}
   \:\times
       \qquad\qquad
       \nonumber \\[-1mm]
   & & \qquad\qquad\qquad\:
   u_1^{\nrcpeak(\sigma)} u_2^{\nrcdfall(\sigma)}
   v_1^{\nrcval(\sigma)} \vhat_2^{\nrcdrise'(\sigma)}
   \:\times
       \qquad\qquad
       \nonumber \\[3mm]
   & & \qquad\qquad\qquad\:
   \we^{\evennrfix(\sigma)} \wo^{\oddnrfix(\sigma)}
   \ze^{\evenrar(\sigma)} \zo^{\oddrar(\sigma)}
   \, \lambda^{\cyc(\sigma)}
   \;.
 \label{def.Pnhathat}
\end{eqnarray}
We are then able to prove:

\begin{theorem}[Second T-fraction for D-permutations]
   \label{thm.Tfrac.second}
The ordinary generating function of the polynomials \reff{def.Pnhathat}
specialized to $v_1 = y_1$ has the T-type continued fraction
\begin{eqnarray}
   & & \hspace*{-12mm}
   \sum_{n=0}^\infty
   \widehat{\widehat{P}}_n(x_1,x_2,y_1,\yhat_2,u_1,u_2,y_1,\vhat_2,\we,\wo,\ze,\zo,\lambda) \: t^n
   \;=\;
       \nonumber \\
   & &
   \cfrac{1}{1 - \lambda^2 \ze \zo  \,t - \cfrac{\lambda x_1 y_1  \,t}{1 -  \cfrac{(x_2\!+\!\lambda\we)(\yhat_2\!+\!\lambda\wo) \,t}{1 - \cfrac{(\lambda+1)(x_1\!+\!u_1) y_1 \,t}{1 - \cfrac{(x_2\!+\!u_2\!+\!\lambda\we)(\yhat_2\!+\!\vhat_2\!+\!\lambda\wo)  \,t}{1 - \cfrac{(\lambda+2)(x_1\!+\!2u_1) y_1 \,t}{1 - \cfrac{(x_2\!+\!2u_2\!+\!\lambda\we)(\yhat_2\!+\!2\vhat_2\!+\!\lambda\wo)  \,t}{1 - \cdots}}}}}}}
   \label{eq.thm.Tfrac.second}
\end{eqnarray}
with coefficients
\begin{subeqnarray}
   \alpha_{2k-1} & = & (\lambda + k-1) \: [x_1 + (k-1) u_1] \: y_1
        \\[1mm]
   \alpha_{2k}   & = & [x_2 + (k-1) u_2 + \lambda\we] \: [\yhat_2 + (k-1)\vhat_2 + \lambda\wo]
        \\[1mm]
   \delta_1  & = &   \lambda^2 \ze \zo   \\[1mm]
   \delta_n  & = &   0    \qquad\hbox{for $n \ge 2$}
   \label{eq.thm.Tfrac.second.weights}
\end{subeqnarray}
\end{theorem}

\noindent
This continued fraction has exactly the same form as
Conjecture~\ref{conj.Tfrac.second},
except that $y_2,v_2$ are replaced by $\yhat_2,\vhat_2$.
We will prove Theorem~\ref{thm.Tfrac.second} in Section~\ref{sec.proofs.2}.

Note that each of the coefficients $\alpha_i$ and $\delta_i$
in \reff{eq.thm.Tfrac.second.weights}
is homogeneous of degree~$1$ in $x_1,x_2,$ $u_1,u_2,\we,\ze$
and also homogeneous of degree~$1$ in $y_1,\yhat_2,\vhat_2,\wo,\zo$.
As before, this reflects the homogeneities of the $\widehat{\widehat{P}}_n$.

\bigskip

{\bf Remark.}
Specializing Theorem~\ref{thm.Tfrac.second} to $\we = \wo = \ze = \zo = 0$
and all other variables except $\lambda$ to 1,
we obtain an S-fraction with
$\alpha_{2k-1} = k(\lambda+k-1)$, $\alpha_{2k} = k^2$
for counting D-derangements by number of cycles
\cite[eq.~(4.4)]{Pan_21}.
\myendremark

\bigskip

It is perhaps worth observing that, in view of Theorem~\ref{thm.Tfrac.second},
Conjecture~\ref{conj.Tfrac.second} is equivalent to
the following assertion about the equidistribution of statistics:

\addtocounter{theorem}{-2}
\begin{conjecture}\hspace*{-3mm}${}^{\bf\prime}$
There exists a bijection $\psi_n \colon \dperm_{2n} \to \dperm_{2n}$
that maps the 12-tuple
\begin{eqnarray*}
   & &  (\eareccpeak, \eareccdfall, \cval, \ereccdrise,
         \nrcpeak, \nrcdfall, \nrcdrise,           \qquad \\
   & &  \;\,\evennrfix, \oddnrfix, \evenrar, \oddrar, \cyc)
\end{eqnarray*}
onto
\begin{eqnarray*}
   & &  (\eareccpeak, \eareccdfall, \cval, \ereccdrise',
         \nrcpeak, \nrcdfall, \nrcdrise',           \qquad \\
   & &  \;\,\evennrfix, \oddnrfix, \evenrar, \oddrar, \cyc)
        \;.
\end{eqnarray*}
\end{conjecture}
\addtocounter{theorem}{1}

%%\bigdash
\subsubsection{Specialization leading to S-fraction}

%% {\bf Maybe also in this case I can first get an alternative T-fraction???}

We refrain from pursuing the cycle-counting generalizations
of Proposition~\ref{prop.Tfrac.first.alternative}
and Corollary~\ref{cor.Tfrac.first.alternative},
and limit ourselves to giving the cycle-counting generalization
of Corollary~\ref{cor.Sfrac.first.specialized}.
If in the T-fraction
\reff{eq.thm.Tfrac.second}/\reff{eq.thm.Tfrac.second.weights}
we specialize $y_1 = \yhat_2 = \vhat_2 = 1 = \zo = \wo = 1$
--- that~is, all the weights associated to odd indices are set to~1 ---
and also specialize $\ze = \we$, $x_1 = x_2$ and $u_1 = u_2$,
the resulting T-fraction can be rewritten as an S-fraction
by noticing that both of them contract
(Proposition~\ref{prop.contraction_even.Ttype})
to the same J-fraction,
and we obtain a generalization of Corollary~\ref{cor.Sfrac.first.specialized}
that includes counting the number of cycles:

\begin{corollary}[S-fraction for specialized D-permutations, counting number of cycles] 
\label{cor.Sfrac.cyc}
The ordinary generating function of the polynomials $P^\star_n(x,u,\we,\lambda)$
defined by
\be
   P^\star_n(x,u,\we,\lambda)
   \;=\;
   \sum_{\sigma\in \dperm_{2n}}
      x^{\earec(\sigma)} u^{\nrcpeak(\sigma)+\nrcdfall(\sigma)}
      \we^{\evenfix(\sigma)} \lambda^{\cyc(\sigma)}
 \label{eq.def.cyccount}
\ee
has the S-type continued fraction
\be
\sum_{n=0}^{\infty} P^\star_n(x,u,\we,\lambda)\: t^n
   \;=\;
  \cfrac{1}{1 - \cfrac{\lambda(x+\lambda\we) t}{1 -  \cfrac{(\lambda+1)x t}{1 - \cfrac{(\lambda+1)(x+u+\lambda\we) t}{1 - \cfrac{(\lambda+2)(x+u) t}{1 - \cdots}}}}}
\label{eq.Sfrac.cyc}
\ee
with coefficients
\begin{subeqnarray}
    \alpha_{2k-1} & = & (\lambda + k-1) \, [x + (k-1)u + \lambda\we] \\[1mm]
    \alpha_{2k}   & = & (\lambda+k) \, [x + (k-1) u]
 \label{labels.Sfrac.cyc}
\end{subeqnarray}
\end{corollary}

\noindent
Of course there is also an analogous S-fraction in which the roles
of odd and even indices are reversed.

\bigskip

{\bf Remarks.}
1. Setting $x=u=1$ and $\we = 0$ in
\reff{eq.Sfrac.cyc}/\reff{labels.Sfrac.cyc},
we obtain the S-fraction for D-semiderangements enumerated by number of cycles,
with coefficients $\alpha_{2k-1} = k(\lambda+k-1)$
and $\alpha_{2k} = k(\lambda+k)$,
which was found a quarter-century ago by
Randrianarivony and Zeng \cite[Corollary~13]{Randrianarivony_96b}.

2. Setting $x=u=\we = 1$ in
\reff{eq.Sfrac.cyc}/\reff{labels.Sfrac.cyc},
we obtain an S-fraction for D-permutations enumerated by number of cycles,
with coefficients $\alpha_{2k-1} = \alpha_{2k} = k(\lambda+k)$.
This S-fraction appeared in \cite[eq.~(14)]{Han_99a}
(in fact in a $q$-generalization),
but with a different combinatorial interpretation.
Further specializing to $\lambda=1$,
we obtain the S-fraction
\reff{eq.mediangenocchi.Sfrac2}/\reff{eq.mediangenocchi.Sfrac2.weights}
for the once-shifted Genocchi medians $h_{n+1}$.
\myendremark

Since the once-shifted Genocchi medians $h_{n+1}$
count D-permutations of $[2n]$,
the preceding remark suggests that one should try to find
an S-fraction for D-permutations that generalizes
\reff{eq.mediangenocchi.Sfrac2}/\reff{eq.mediangenocchi.Sfrac2.weights}:

\begin{openproblem}
   \label{problem.Dperm.S}
\rm
Find the statistics on D-permutations
that would lead to a more general S-fraction of the form
\begin{subeqnarray}
   \alpha_{2k-1}  & = &  [x_1 + (k-1) u_1] \, [x_2 + k u_2] \\
   \alpha_{2k}    & = &  [y_1 + (k-1) v_1] \, [y_2 + k v_2]
 \label{eq.quadratic_family.bis1}
\end{subeqnarray}
or
\begin{subeqnarray}
   \alpha_{2k-1}  & = &  [x_1 + (k-1) u_1] \, [x_2 + (k-1) u_2 + w] \\
   \alpha_{2k}    & = &  [y_1 + (k-1) v_1] \, [y_2 + (k-1) v_2 + w']
 \label{eq.quadratic_family.bis2}
\end{subeqnarray}
\myendremark
\end{openproblem}

%%\bigdash
\subsubsection{Reformulation using cycle valley minima}

We can rephrase Conjecture~\ref{conj.Tfrac.second}
and Theorem~\ref{thm.Tfrac.second}
in a more suggestive form by observing that
each non-singleton cycle contains
precisely one maximum element, which is necessarily a cycle peak,
and precisely one minimum element, which is necessarily a cycle valley.
We therefore define four new statistics:
\begin{itemize}
   \item cycle peak maximum (${\maxpeak}$):  cycle peak that is the largest in its cycle;
   \item cycle peak non-maximum ($\nmaxpeak$):  cycle peak that is not the largest in its cycle;
   \item cycle valley minimum (${\minval}$):  cycle valley that is the smallest in its cycle;
   \item cycle valley non-minimum ($\nminval$):  cycle valley that is not the smallest in its cycle.
\end{itemize}
We now introduce a polynomial that is similar to \reff{def.Pn}
except that the classification of cycle valleys as records or non-records
--- which we have already renounced in Conjecture~\ref{conj.Tfrac.second} ---
is replaced by the classification of cycle valleys as minima or non-minima:
\begin{eqnarray}
   & &
   \widetilde{P}_n(x_1,x_2,\ytilde_1,y_2,u_1,u_2,\vtilde_1,v_2,\we,\wo,\ze,\zo)
   \;=\;
       \nonumber \\[4mm]
   & & \qquad\qquad
   \sum_{\sigma \in \dperm_{2n}}
   x_1^{\eareccpeak(\sigma)} x_2^{\eareccdfall(\sigma)}
   \ytilde_1^{\minval(\sigma)} y_2^{\ereccdrise(\sigma)}
   \:\times
       \qquad\qquad
       \nonumber \\[-1mm]
   & & \qquad\qquad\qquad\:
   u_1^{\nrcpeak(\sigma)} u_2^{\nrcdfall(\sigma)}
   \vtilde_1^{\nminval(\sigma)} v_2^{\nrcdrise(\sigma)}
   \:\times
       \qquad\qquad
       \nonumber \\[3mm]
   & & \qquad\qquad\qquad\:
   \we^{\evennrfix(\sigma)} \wo^{\oddnrfix(\sigma)}
   \ze^{\evenrar(\sigma)} \zo^{\oddrar(\sigma)}
   \;.
 \label{def.Pntilde}
\end{eqnarray}
Then the factor $(\lambda + k-1) y_1$
in Conjecture~\ref{conj.Tfrac.second} and Theorem~\ref{thm.Tfrac.second}
is replaced by $\ytilde_1 + (k-1)\vtilde_1$,
and Conjecture~\ref{conj.Tfrac.second} can be rewritten as:

\addtocounter{theorem}{-3}
\begin{conjecture}\hspace*{-3mm}${}^{\bf\prime\prime}$
The ordinary generating function of the polynomials \reff{def.Pntilde}
has the T-type continued fraction
\begin{eqnarray}
   & & \hspace*{-12mm}
   \sum_{n=0}^\infty
   \widetilde{P}_n(x_1,x_2,\ytilde_1,y_2,u_1,u_2,\vtilde_1,v_2,\we,\wo,\ze,\zo) \: t^n
   \;=\;
       \nonumber \\
   & &
   \cfrac{1}{1 - \ze \zo  \,t - \cfrac{x_1 \ytilde_1  \,t}{1 -  \cfrac{(x_2\!+\!\we)(y_2\!+\!\wo) \,t}{1 - \cfrac{(x_1\!+\!u_1) (\ytilde_1 + \vtilde_1) \,t}{1 - \cfrac{(x_2\!+\!u_2\!+\!\we)(y_2\!+\!v_2\!+\!\wo)  \,t}{1 - \cfrac{(x_1\!+\!2u_1) (\ytilde_1 + 2\vtilde_1) \,t}{1 - \cfrac{(x_2\!+\!2u_2\!+\!\we)(y_2\!+\!2v_2\!+\!\wo)  \,t}{1 - \cdots}}}}}}}
       \nonumber \\[1mm]
   \label{eq.conj.Tfrac.second.bis}
\end{eqnarray}
with coefficients
\begin{subeqnarray}
   \alpha_{2k-1} & = & [x_1 + (k-1) u_1] \: [\ytilde_1 + (k-1) \vtilde_1]
        \\[1mm]
   \alpha_{2k}   & = & [x_2 + (k-1) u_2 + \we] \: [y_2 + (k-1)v_2 + \wo]
        \\[1mm]
   \delta_1  & = &   \ze \zo   \\[1mm]
   \delta_n  & = &   0    \qquad\hbox{for $n \ge 2$}
   \label{eq.conj.Tfrac.second.weights.bis}
\end{subeqnarray}
\end{conjecture}
\addtocounter{theorem}{2}

\noindent
Note that the factors of $\lambda$ multiplying $\we,\wo,\ze,\zo$
in Conjecture~\ref{conj.Tfrac.second}
have now disappeared because we are no longer giving such factors
to singleton cycles.
Even more strikingly, Conjecture~\ref{conj.Tfrac.second}${}'$
now looks identical to Theorem~\ref{thm.Tfrac.first}
except that $y_1,v_1$ have been replaced by
$\ytilde_1,\vtilde_1$.
Moreover, Theorem~\ref{thm.Tfrac.second} can be rewritten in a similar way,
replacing $y_1,v_1$ by $\ytilde_1,\vtilde_1$ as in \reff{def.Pntilde}
{\em and}\/ replacing $y_2,v_2$ by $\yhat_2,\vhat_2$
as in \reff{def.Pnhathat};
for brevity we leave this reformulation to the reader.

\subsubsection{S-fraction for D-cycles}

We can also enumerate D-cycles by extracting the coefficient of $\lambda^1$
in Conjecture~\ref{conj.Tfrac.second} and Theorem~\ref{thm.Tfrac.second}.
We begin by defining the polynomials analogous to \reff{def.Pn}
but restricted to D-cycles (which we recall have no fixed points):
\begin{eqnarray}
   & &
   P^{\dcycle}_n(x_1,x_2,y_1,y_2,u_1,u_2,v_1,v_2)
   \;=\;
       \nonumber \\[4mm]
   & & \qquad\qquad
   \sum_{\sigma \in \dcycle_{2n}}
   x_1^{\eareccpeak(\sigma)} x_2^{\eareccdfall(\sigma)}
   y_1^{\ereccval(\sigma)} y_2^{\ereccdrise(\sigma)}
   \:\times
       \qquad\qquad
       \nonumber \\[-1mm]
   & & \qquad\qquad\qquad\quad\!
   u_1^{\nrcpeak(\sigma)} u_2^{\nrcdfall(\sigma)}
   v_1^{\nrcval(\sigma)} v_2^{\nrcdrise(\sigma)}
   \;.
 \label{def.Pn.dcycle}
\end{eqnarray}
Conjecture~\ref{conj.Tfrac.second} then implies the following:

\begin{conjecture}[Conjectured S-fraction for D-cycles]
   \label{conj.Tfrac.second.dcycle}
The ordinary generating function of the polynomials \reff{def.Pn.dcycle}
specialized to $v_1 = y_1$ has the S-type continued fraction
\be
   \!\!\!
   \sum_{n=0}^\infty
   P^{\dcycle}_{n+1}(x_1,x_2,y_1,y_2,u_1,u_2,y_1,v_2) \: t^n
   \;=\;
   \cfrac{x_1 y_1}{1 - \cfrac{x_2 y_2  \,t}{1 - \cfrac{(x_1\!+\!u_1) y_1 \,t}{1 - \cfrac{(x_2\!+\!u_2\!)(y_2\!+\!v_2)  \,t}{1 - \cfrac{(x_1\!+\!2u_1) 2y_1 \,t}{1 - \cfrac{(x_2\!+\!2u_2)(y_2\!+\!2v_2)  \,t}{1 - \cdots}}}}}}
  \label{eq.conj.Tfrac.second.dcycle}
\ee
with coefficients
\begin{subeqnarray}
   \alpha_{2k-1}   & = & [x_2 + (k-1) u_2] \: [y_2 + (k-1)v_2]
        \\[1mm]
   \alpha_{2k}     & = & [x_1 + k u_1] \: k y_1
   \label{eq.conj.Tfrac.second.weights.dcycle}
\end{subeqnarray}
\end{conjecture}

\noindent
We have tested this conjecture through $n=5$.

Similarly, we can define the polynomials analogous to \reff{def.Pnhathat}
but restricted to D-cycles:
%% https://tex.stackexchange.com/questions/601576/how-to-lower-the-superscript-in-a-math-mode-manually
\begin{eqnarray}
   & &
   \widehat{\widehat{P}}{}^{\dcycle}_n(x_1,x_2,y_1,\yhat_2,u_1,u_2,v_1,\vhat_2)
   \;=\;
       \nonumber \\[4mm]
   & & \qquad\qquad
   \sum_{\sigma \in \dcycle_{2n}}
   x_1^{\eareccpeak(\sigma)} x_2^{\eareccdfall(\sigma)}
   y_1^{\ereccval(\sigma)} \yhat_2^{\ereccdrise'(\sigma)}
   \:\times
       \qquad\qquad
       \nonumber \\[-1mm]
   & & \qquad\qquad\qquad\quad\!
   u_1^{\nrcpeak(\sigma)} u_2^{\nrcdfall(\sigma)}
   v_1^{\nrcval(\sigma)} \vhat_2^{\nrcdrise'(\sigma)}
   \;.
 \label{def.Pn.dcyclehathat}
\end{eqnarray}
Theorem~\ref{thm.Tfrac.second} implies the following:

\begin{corollary}[S-fraction for D-cycles]
   \label{cor.Tfrac.second.dcycle}
The ordinary generating function of the polynomials
\reff{def.Pn.dcyclehathat}
specialized to $v_1 = y_1$ has the S-type continued fraction
\be
   \!\!\!
   \sum_{n=0}^\infty
   \widehat{\widehat{P}}{}^{\dcycle}_{n+1}(x_1,x_2,y_1,\yhat_2,u_1,u_2,y_1,\vhat_2) \: t^n
   \;=\;
   \cfrac{x_1 y_1}{1 - \cfrac{x_2 \yhat_2  \,t}{1 - \cfrac{(x_1\!+\!u_1) y_1 \,t}{1 - \cfrac{(x_2\!+\!u_2) (\yhat_2\!+\!\vhat_2)  \,t}{1 - \cfrac{(x_1\!+\!2u_1) 2y_1 \,t}{1 - \cfrac{(x_2\!+\!2u_2) (\yhat_2\!+\!2\vhat_2)  \,t}{1 - \cdots}}}}}}
 \label{eq.cor.Tfrac.second.dcycle}
\ee
with coefficients
\begin{subeqnarray}
   \alpha_{2k-1}   & = & [x_2 + (k-1) u_2] \: [\yhat_2 + (k-1) \vhat_2]
        \\[1mm]
   \alpha_{2k}     & = & [x_1 + k u_1] \: k y_1
   \label{eq.cor.Tfrac.second.weights.dcycle}
\end{subeqnarray}
\end{corollary}

Specializing Corollary~\ref{cor.Tfrac.second.dcycle}
by setting all variables equal to 1,
we obtain an S-fraction with
$\alpha_{2k-1} = k^2$ and $\alpha_{2k} = k(k+1)$,
which we recognize as the S-fraction
\reff{eq.genocchi.Sfrac}/\reff{eq.genocchi.Sfrac.weights}
for the Genocchi numbers $g_n$.
We have therefore recovered the known fact \cite{Lazar_20,Lazar_22}
that $|\dcycle_{2n+2}| = g_n$,
or equivalently that $|\dcycle_{2n}| = g_{n-1}$.

%% {\bf Does anyone have continued fractions with more variables
%% for D-cycles???}

% \begin{quote}
% \small
% {\bf OLD STUFF --- DO WE NEED ANY OF THIS???}
% 
% Using Corollary~\ref{cor.Sfrac.cyc}, we can enumerate the number of D-cycles. 
% 
% \begin{corollary}\label{cor.Dcycle.count} 
% 	The cardinality of the $\dcycle_{2n}$, the set of D-cycles on $[2n]$, is the Genocchi number $g_{n-1}$.
% \end{corollary}
% 
% \begin{proof} Using Proposition~\ref{prop.augment.Tfrac}, we can rewrite the S-fraction~\reff{eq.Sfrac.cyc} as
% \begin{eqnarray}
% 	1+\cfrac{\lambda(1+\lambda) t}{1- \lambda(1+\lambda) t -\cfrac{(1+\lambda) t}{1 -\cfrac{(1+\lambda)(2+\lambda) t}{1-\cfrac{2(2+\lambda) t}{1-\cdots}}}} \nonumber\\
% \qquad = 1+(\lambda + \lambda^2)t F(\lambda; t)
% 	\label{eq.augment.Sfrac.cyc}
% \end{eqnarray}
% where $F(\lambda;t)$ is defined to be the formal power series that satisfies Equation~\reff{eq.augment.Sfrac.cyc}.
% Thus, the cardinality of $\dcycle_{2n}$ is given by the coefficient of $t^{n-1}$ in $F(0;t)$.
% However, notice that $F(0;t)$ is given by the classical 
% S-fraction~\reff{eq.genocchi.Sfrac} which is S-fraction expansion 
% of the ordinary generating function for Genocchi numbers. \qed
% \end{proof}
% 
% {\bf What else to say????}
% \end{quote}

\subsection[$p,q$-generalizations of the second T-fraction]{$\bm{p,q}$-generalizations of the second T-fraction}

We can also make a $p,q$-generalization of the second T-fraction
in Theorem~\ref{thm.Tfrac.second}.
%% involving the number of cycles ($\cyc$).
Let us define the polynomial in 23~variables
that generalizes \reff{def.Pnhathat} by including
four pairs of $(p,q)$-variables
corresponding to the four refined types of crossings and nestings
illustrated in Figure~\ref{fig.refined_crossnest},
as well as two variables corresponding to pseudo-nestings of fixed points:

\vspace*{-5mm}

\begin{eqnarray}
   & & \hspace*{-8mm}
   \Scale[0.95]{
   \widehat{\widehat{P}}_n(x_1,x_2,y_1,\yhat_2,u_1,u_2,v_1,\vhat_2,\we,\wo,\ze,\zo,p_{-1},p_{-2},p_{+1},\phat_{+2},q_{-1},q_{-2},q_{+1},\qhat_{+2},\se,\so,\lambda)
   \;=\; }
       \nonumber \\[4mm]
   & & \qquad\qquad
   \sum_{\sigma \in \dperm_{2n}}
   x_1^{\eareccpeak(\sigma)} x_2^{\eareccdfall(\sigma)}
   y_1^{\ereccval(\sigma)} \yhat_2^{\ereccdrise'(\sigma)}
   \:\times
       \qquad\qquad\:
       \nonumber \\[-1mm]
   & & \qquad\qquad\qquad\;\,
   u_1^{\nrcpeak(\sigma)} u_2^{\nrcdfall(\sigma)}
   v_1^{\nrcval(\sigma)} \vhat_2^{\nrcdrise'(\sigma)}
   \:\times
       \qquad\qquad\;
       \nonumber \\[3mm]
   & & \qquad\qquad\qquad\;\,
   \we^{\evennrfix(\sigma)} \wo^{\oddnrfix(\sigma)}
   \ze^{\evenrar(\sigma)} \zo^{\oddrar(\sigma)}
   \:\times
       \qquad\qquad\;
       \nonumber \\[3mm]
   & & \qquad\qquad\qquad\;\,
   p_{-1}^{\lcrosscpeak(\sigma)}
   p_{-2}^{\lcrosscdfall(\sigma)}
   p_{+1}^{\ucrosscval(\sigma)}
   \phat_{+2}^{\ucrosscdrise'(\sigma)}
          \:\times
       \qquad\qquad\;
       \nonumber \\[3mm]
   & & \qquad\qquad\qquad\;\,
   q_{-1}^{\lnestcpeak(\sigma)}
   q_{-2}^{\lnestcdfall(\sigma)}
   q_{+1}^{\unestcval(\sigma)}
   \qhat_{+2}^{\unestcdrise'(\sigma)}
          \:\times
       \qquad\qquad
       \nonumber \\[3mm]
   & & \qquad\qquad\qquad\;\,
   \se^{\epsnest(\sigma)}
   \so^{\opsnest(\sigma)}
   \, \lambda^{\cyc(\sigma)}
   \;.
 \label{def.Pn.pq.second}
\end{eqnarray}
This is the same as \reff{def.Pn.pq}, except for the
inclusion of the factor $\lambda^{\cyc(\sigma)}$
and the replacement of $y_2^{\ereccdrise(\sigma)}, v_2^{\nrcdrise(\sigma)},
p_{+2}^{\ucrosscdrise(\sigma)}, q_{+2}^{\unestcdrise(\sigma)}$
by $\yhat_2^{\ereccdrise'(\sigma)}, \vhat_2^{\nrcdrise'(\sigma)}$,
\linebreak
 $\phat_{+2}^{\ucrosscdrise'(\sigma)}, \qhat_{+2}^{\unestcdrise'(\sigma)}$,
respectively, where the statistics
$\ucrosscdrise'$ and
\linebreak
$\unestcdrise'$ are defined as
% 
% \begin{subeqnarray}
%    \ucrosscdrise'(\sigma)
%    & = &
%    \#\{ i<j<k<l: k\in \cdrise(\sigma) \hbox{ and }
%    \\ \nonumber
%           & & \qquad\qquad\qquad \qquad   k = \sigma(i), l = \sigma(j)  \}
%       \\[2mm]
%    \unestcdrise'(\sigma)
%    & = &
%    \#\{ i<j<k<l: k\in \cdrise(\sigma) \hbox{ and }
%    \\ \nonumber
%           & & \qquad\qquad\qquad \qquad   k = \sigma(j), l = \sigma(i)  \}
% \end{subeqnarray}
% i.e.,
\begin{subeqnarray}
\ucrosscdrise'(\sigma) & = & \sum_{k\in \Cdrise(\sigma)} \ucross'(k,\sigma) 
   \\[2mm]
\unestcdrise'(\sigma) & = & \sum_{k\in \Cdrise(\sigma)} \unest'(k,\sigma) 
\end{subeqnarray}
%% Here $\ucross'$ and $\unest'$ have been defined in~(\ref{def.ucrossnestjk.prime}a,b).
[These statistics have already been used in \reff{def.Pn.pq.variant}.]

We refrain from attempting to generalize Conjecture~\ref{conj.Tfrac.second}
and simply limit ourselves to stating the $p,q$-generalization of
Theorem~\ref{thm.Tfrac.second} that we are able to prove.
It turns out that we need to make the specializations
$v_1=y_1$ and $q_{+1} = p_{+1}$.
The result is then the following:

\begin{theorem}[Second T-fraction for D-permutations, $p,q$-generalization]
   \label{thm.Tfrac.second.pq}
The ordinary generating function of the polynomials \reff{def.Pn.pq.second} specialized to $v_1=y_1$ and $q_{+1} = p_{+1}$ has the
T-type continued fraction
\begin{eqnarray}
   & & \hspace*{-7mm}
\Scale[0.89]{
   \sum\limits_{n=0}^\infty
   \widehat{\widehat{P}}_n(x_1,x_2,y_1,\yhat_2,u_1,u_2,y_1,\vhat_2,\we,\wo,\ze,\zo,p_{-1},p_{-2},p_{+1},\phat_{+2},q_{-1},q_{-2},p_{+1},\qhat_{+2},\se,\so,\lambda) \: t^n
   \;=\;
}
       \nonumber \\[2mm]
   & & \hspace*{-3mm}
\Scale[0.82]{
   \cfrac{1}{1 - \lambda^2 \ze \zo  \,t - \cfrac{\lambda p_{+1} x_1 y_1  \,t}{1 -  \cfrac{(x_2\!+\!\lambda\se\we)(\yhat_2\!+\!\lambda\so\wo) \,t}{1 - \cfrac{(\lambda + 1) p_{+1}^2 y_1 (p_{-1}x_1\!+\!q_{-1}u_1)  \,t}{1 - \cfrac{(p_{-2}x_2\!+\!q_{-2}u_2\!+\!\lambda\se^2\we)(\phat_{+2}\yhat_2\!+\!\qhat_{+2}\vhat_2\!+\!\lambda\so^2\wo)  \,t}{1 - \cfrac{(\lambda + 2)p_{+1}^3 y_1(p_{-1}^2 x_1\!+\! q_{-1} [2]_{p_{-1},q_{-1}}u_1)  \,t}{1 - \cfrac{(p_{-2}^2 x_2\!+\! q_{-2} [2]_{p_{-2},q_{-2}} u_2\!+\!\lambda\se^3\we)(\phat_{+2}^2 \yhat_2\!+\! \qhat_{+2} [2]_{\phat_{+2},\qhat_{+2}}\vhat_2\!+\!\lambda\so^3\wo)  \,t}{1 - \cdots}}}}}}}
}
       \nonumber \\[1mm]
   \label{eq.thm.Tfrac.second.pq}
\end{eqnarray}
with coefficients 
\begin{subeqnarray}
   \alpha_{2k-1} & = &
   (\lambda + k-1) \: p_{+1}^{k} y_1 \:
      \bigl( p_{-1}^{k-1} x_1 + q_{-1} [k-1]_{p_{-1},q_{-1}} u_1 \bigr) \:
        \\[1mm]
   \alpha_{2k}   & = &
      \bigl( p_{-2}^{k-1} x_2 + q_{-2} [k-1]_{p_{-2},q_{-2}} u_2 + \lambda\se^k\we \bigr) \:
      \bigl( \phat_{+2}^{k-1} \yhat_2 + \qhat_{+2} [k-1]_{\phat_{+2},\qhat_{+2}} \vhat_2 + \lambda\so^k\wo \bigr)
        \nonumber \\ \\
   \delta_1  & = &   \lambda^2\ze \zo   \\[1mm]
   \delta_n  & = &   0    \qquad\hbox{for $n \ge 2$}
   \label{eq.thm.Tfrac.second.weights.pq}
\end{subeqnarray}
\end{theorem}

\noindent
We will prove Theorem~\ref{thm.Tfrac.second.pq} in Section~\ref{sec.proofs.2}.
Of course we reobtain Theorem~\ref{thm.Tfrac.second} by
making the specialization
$p_{-1} = p_{-2} = p_{+1} = \phat_{+2} = q_{-1} = q_{-2} = \qhat_{+2} =
 \se = \so = 1$.

\bigskip

{\bf Remark.}
If we specialize to $x_1 = u_1$, $x_2 = u_2$, $\yhat_2 = \vhat_2$
--- that is, renounce the counting of records ---
then the coefficients \reff{eq.thm.Tfrac.second.weights.pq} simplify to
\begin{subeqnarray}
   \alpha_{2k-1} & = &
   (\lambda + k-1) \: p_{+1}^{k} y_1 \: [k]_{p_{-1},q_{-1}} x_1
        \\[1mm]
   \alpha_{2k}   & = &
      \bigl( [k]_{p_{-2},q_{-2}} x_2 + \lambda\se^k\we \bigr) \:
      \bigl( [k]_{\phat_{+2},\qhat_{+2}} \yhat_2 + \lambda\so^k\wo \bigr)
        \\[1mm]
   \delta_1  & = &   \ze \zo   \\[1mm]
   \delta_n  & = &   0    \qquad\hbox{for $n \ge 2$}
   \label{eq.thm.Tfrac.second.weights.pq.bis}
\end{subeqnarray}
\myendremark

\bigskip

We can also obtain a $p,q$-generalization of Corollary~\ref{cor.Sfrac.cyc},
or equivalently a cycle-counting generalization of
Corollary~\ref{cor.Sfrac.first.specialized.pq}.
If in the T-fraction
\reff{eq.thm.Tfrac.second.pq}/\reff{eq.thm.Tfrac.second.weights.pq}
we specialize $\zo = \wo = y_1 = \yhat_2 = \vhat_2 = 
   p_{+1} = \phat_{+2} = \qhat_{+2} = \so = 1$
--- that~is, all the weights associated to odd indices are set to 1 ---
and also specialize $\ze = \we$, $x_1 = x_2$, $u_1 = u_2$,
$p_{-1} = p_{-2}$, $q_{-1} = q_{-2}$ and $\se = 1$,
the resulting T-fraction can be rewritten as an S-fraction
by noticing that they contract
(Proposition~\ref{prop.contraction_even.Ttype})
to the same J-fraction:

\begin{corollary}[S-fraction for D-permutations counting number of cycles,
   $p,q$-generalization]
\label{cor.Sfrac.cyc.pq}
The ordinary generating function of the polynomials
$P^\star_n(x,u,\we,p_-,q_-,\lambda)$ defined by
\begin{eqnarray}
   & &
   P^\star_n(x,u,\we,p_-,q_-,\lambda)
   \;=\;
       \nonumber \\[3mm]
   & & \qquad
   \sum_{\sigma\in \dperm_{2n}}
      x^{\earec(\sigma)} u^{\nrcpeak(\sigma)+\nrcdfall(\sigma)}
      \we^{\evenfix(\sigma)}
      p_-^{\lcross(\sigma)} q_-^{\lnest(\sigma)} \,
     \lambda^{\cyc(\sigma)}
     \qquad\qquad
 \label{eq.def.cyccount.pq}
\end{eqnarray}
has the S-type continued fraction
\be
\sum_{n=0}^{\infty} P^\star_n(x,u,\we,p_-,q_-,\lambda)\: t^n
   \;=\;
  \cfrac{1}{1 - \cfrac{\lambda(x+\lambda\we) t}{1 -  \cfrac{(\lambda+1)x t}{1 - \cfrac{(\lambda+1)(p_- x+ q_- u+\lambda\we) t}{1 - \cfrac{(\lambda+2)(p_- x+ q_- u) t}{1 - \cdots}}}}}
\label{eq.Sfrac.cyc.pq}
\ee
with coefficients
\begin{subeqnarray}
    \alpha_{2k-1} & = & (\lambda + k-1) \, 
      \bigl( p_-^{k-1} x + q_- [k-1]_{p_-,q_-} u + \lambda\we \bigr)
            \\[2mm]
    \alpha_{2k}   & = & (\lambda+k) \, 
      \bigl( p_-^{k-1} x + q_- [k-1]_{p_-,q_-} u \bigr)
 \label{labels.Sfrac.cyc.pq}
\end{subeqnarray}
\end{corollary}

\bigskip

Finally, we can enumerate D-cycles by extracting the coefficient of $\lambda^1$
in Theorem~\ref{thm.Tfrac.second.pq}.
We begin by defining the polynomials analogous to \reff{def.Pn.pq.second}
but restricted to D-cycles (which we recall have no fixed points):
\begin{eqnarray}
   & & \hspace*{-8mm}
   \widehat{\widehat{P}}{}^{\dcycle}_n(x_1,x_2,y_1,\yhat_2,u_1,u_2,v_1,\vhat_2,p_{-1},p_{-2},p_{+1},\phat_{+2},q_{-1},q_{-2},q_{+1},\qhat_{+2},\lambda)
   \;=\;
       \nonumber \\[4mm]
   & & \qquad\qquad
   \sum_{\sigma \in \dperm_{2n}}
   x_1^{\eareccpeak(\sigma)} x_2^{\eareccdfall(\sigma)}
   y_1^{\ereccval(\sigma)} \yhat_2^{\ereccdrise'(\sigma)}
   \:\times
       \qquad\qquad
       \nonumber \\[-1mm]
   & & \qquad\qquad\qquad\;\,
   u_1^{\nrcpeak(\sigma)} u_2^{\nrcdfall(\sigma)}
   v_1^{\nrcval(\sigma)} \vhat_2^{\nrcdrise'(\sigma)}
   \:\times
       \qquad\qquad
       \nonumber \\[3mm]
   & & \qquad\qquad\qquad\;\,
   p_{-1}^{\lcrosscpeak(\sigma)}
   p_{-2}^{\lcrosscdfall(\sigma)}
   p_{+1}^{\ucrosscval(\sigma)}
   \phat_{+2}^{\ucrosscdrise'(\sigma)}
          \:\times
       \qquad\qquad
       \nonumber \\[3mm]
   & & \qquad\qquad\qquad\;\,
   q_{-1}^{\lnestcpeak(\sigma)}
   q_{-2}^{\lnestcdfall(\sigma)}
   q_{+1}^{\unestcval(\sigma)}
   \qhat_{+2}^{\unestcdrise'(\sigma)}
   \, \lambda^{\cyc(\sigma)}
   \;.
 \label{def.Pn.dcyclehathat.pq}
\end{eqnarray}
Theorem~\ref{thm.Tfrac.second.pq} then implies the following
generalization of Corollary~\ref{cor.Tfrac.second.dcycle}:

\begin{corollary}[S-fraction for D-cycles, $p,q$-generalization]
   \label{cor.Tfrac.second.dcycle.pq}
The ordinary generating function of the polynomials
\reff{def.Pn.dcyclehathat.pq}
specialized to $v_1 = y_1$ and $q_{+1} = p_{+1}$
has the S-type continued fraction
\begin{eqnarray}
   & &
   \hspace*{-8mm}
   \sum_{n=0}^\infty
   \widehat{\widehat{P}}{}^{\dcycle}_{n+1}(x_1,x_2,y_1,\yhat_2,u_1,u_2,y_1,\vhat_2,p_{-1},p_{-2},p_{+1},\phat_{+2},q_{-1},q_{-2},p_{+1},\qhat_{+2}) \: t^n
   \;=\;
       \nonumber \\[2mm]
   & &
   \cfrac{p_{+1} x_1 y_1}{1 - \cfrac{x_2 \yhat_2  \,t}{1 - \cfrac{p_{+1}^2 (p_{-1}x_1\!+\!q_{-1}u_1) y_1 \,t}{1 - \cfrac{(p_{-2}x_2\!+\!q_{-2}u_2) (\phat_{+2}\yhat_2\!+\!\qhat_{+2}\vhat_2)  \,t}{1 - \cfrac{(p_{-1}^2x_1\!+\!q_{-1} [2]_{p_{-1},q_{-1}}u_1) 2 p_{+1}^3 y_1 \,t}{1 - \cfrac{(p_{-2}^2 x_2\!+\! q_{-2} [2]_{p_{-2},q_{-2}}u_2) (p_{+2}^2 \yhat_2\!+\! q_{+2} [2]_{p_{+2},q_{+2}}\vhat_2)  \,t}{1 - \cdots}}}}}}
  \qquad\qquad
 \label{eq.cor.Tfrac.second.dcycle.pq}
\end{eqnarray}
with coefficients
\begin{subeqnarray}
   \alpha_{2k-1}   & = &
      \bigl( p_{-2}^{k-1} x_2 + q_{-2} [k-1]_{p_{-2},q_{-2}} u_2 \bigr) \:
      \bigl( \phat_{+2}^{k-1} \yhat_2 + \qhat_{+2} [k-1]_{\phat_{+2},\qhat_{+2}} \vhat_2 \bigr)
     \qquad
        \\[1mm]
   \alpha_{2k}     & = & 
      \bigl( p_{-1}^{k} x_1 + q_{-1} [k]_{p_{-1},q_{-1}} u_1 \bigr) \:
      k \, p_{+1}^{k+1} y_1 \:
   \label{eq.cor.Tfrac.second.weights.dcycle.pq}
\end{subeqnarray}
\end{corollary}

\subsection{Second master T-fraction}  \label{subsec.second.master}

As with the first T-fraction, we can go much farther, and obtain a T-fraction
in six infinite families of indeterminates:
$\bsfa = (\sfa_{\ell})_{\ell \ge 0}$,
$\bsfb = (\sfb_{\ell,\ell'})_{\ell,\ell' \ge 0}$,
$\bsfc = (\sfc_{\ell,\ell'})_{\ell,\ell' \ge 0}$,
$\bsfd = (\sfd_{\ell,\ell'})_{\ell,\ell' \ge 0}$,
$\bsfe = (\sfe_\ell)_{\ell \ge 0}$,
$\bsff = (\sff_\ell)_{\ell \ge 0}$;
please note that $\bsfa$ now has one index rather than two. We
now define a sequence of polynomials that will include 
the polynomials \reff{def.Pnhathat} and \reff{def.Pn.pq.second}
as specializations:
\begin{eqnarray}
   & & \hspace*{-10mm}
   \widehat{\widehat{Q}}_n(\bsfa,\bsfb,\bsfc,\bsfd,\bsfe,\bsff, \lambda)
   \;=\;
       \nonumber \\[4mm]
   & &
   \sum_{\sigma \in \dperm_{2n}}
   \;\:
	\lambda^{{\rm cyc}(\sigma)}
        \prod\limits_{i \in \Cval(\sigma)}  \! \sfa_{\ucross(i,\sigma) + \unest(i,\sigma) }
	\prod\limits_{i \in \Cpeak(\sigma)} \!\!  \sfb_{\lcross(i,\sigma),\,\lnest(i,\sigma) }
       \:\times
       \qquad\qquad
       \nonumber \\[1mm]
   & & \qquad\qquad\quad\;\;
   \prod\limits_{i \in \Cdfall(\sigma)} \!\!  \sfc_{\lcross(i,\sigma),\,\lnest(i,\sigma)}
   \;
	\prod\limits_{i \in \Cdrise(\sigma)} \!\!  \sfd_{\ucross'(i,\sigma),\,\unest'(i,\sigma)}
       \:\times
       \qquad\qquad
       \nonumber \\[1mm]
   & & \qquad\qquad\quad\;
   \prod\limits_{i \in \Evenfix(\sigma)} \!\!\! \sfe_{\psnest(i,\sigma)}
   \;
   \prod\limits_{i \in \Oddfix(\sigma)}  \!\!\! \sff_{\psnest(i,\sigma)}
   \;.
   \quad
 \label{def.Qhatn.secondmaster}
\end{eqnarray}
We remark that \reff{def.Qhatn.secondmaster}
is {\em almost}\/ the same as the polynomial
introduced in \cite[eq.~(2.100)]{Sokal-Zeng_masterpoly},
but restricted to D-permutations
and refined to record the parity of fixed points;
the main difference is that the treatment of $\sfd$ is a bit nicer here,
using the statistics $\ucross'$ and $\unest'$.

Note that here, in contrast to the first master T-fraction,
$\widehat{\widehat{Q}}_n$ depends on
$\ucross(i,\sigma)$ and $\unest(i,\sigma)$ only via their sum:
that is the price we have to pay in order to include the statistic $\cyc$.
Furthermore, the indices of $\sfd$ involve $\ucross'$ and $\unest'$
instead of $\ucross$ and $\unest$.

The polynomials \reff{def.Qhatn.secondmaster} then have a nice T-fraction:

\begin{theorem}[Second master T-fraction for D-permutations]
   \label{thm.Tfrac.second.master}
The ordinary generating function of the polynomials
$\widehat{\widehat{Q}}_n(\bsfa,\bsfb,\bsfc,\bsfd,\bsfe,\bsff,\lambda)$
has the T-type continued fraction
\be
   \sum_{n=0}^\infty \widehat{\widehat{Q}}_n(\bsfa,\bsfb,\bsfc,\bsfd,\bsfe,\bsff) \: t^n
   \;=\;
\Scale[0.95]{
	\cfrac{1}{1 - \lambda^2 \sfe_0 \sff_0 t - \cfrac{\lambda \sfa_{0} \sfb_{00} t}{1 - \cfrac{(\sfc_{00} + \lambda\sfe_1)(\sfd_{00} + \lambda\sff_1) t}{1 - \cfrac{(\lambda+1)\sfa_{1}(\sfb_{01} + \sfb_{10}) t}{1 - \cfrac{(\sfc_{01} + \sfc_{10} + \lambda\sfe_2)(\sfd_{01} + \sfd_{10} + \lambda\sff_2)t}{1 - \cdots}}}}}
}
   \label{eq.thm.Tfrac.second.master}
\ee
with coefficients
\begin{subeqnarray}
   \alpha_{2k-1}
   & = &
	(\lambda+k-1) \sfa_{k-1}
   \biggl( \sum_{\xi=0}^{k-1} \sfb_{k-1-\xi,\xi} \biggr)
        \\[2mm]
   \alpha_{2k}
   & = &
   \biggl( \lambda\sfe_k \,+\, \sum_{\xi=0}^{k-1} \sfc_{k-1-\xi,\xi} \biggr)
   \biggl( \lambda\sff_k \,+\, \sum_{\xi=0}^{k-1} \sfd_{k-1-\xi,\xi} \biggr)
        \\[2mm]
   \delta_1  & = &  \lambda^2\sfe_0 \sff_0 \\[2mm]
   \delta_n  & = &   0    \qquad\hbox{for $n \ge 2$}
 \label{def.weights.Tfrac.second.master}
\end{subeqnarray}
\end{theorem}

\noindent
We will prove this theorem in Section~\ref{sec.proofs.2}.
It implies Theorems~\ref{thm.Tfrac.second} and \ref{thm.Tfrac.second.pq}
by straightforward specializations.

\section{Preliminaries for the proofs}   \label{sec.prelimproofs}

Our proofs are based on Flajolet's \cite{Flajolet_80}
combinatorial interpretation of continued fractions
in terms of Dyck and Motzkin paths
and its generalization
\cite{Fusy_15,Oste_15,Josuat-Verges_18,Sokal_totalpos,Elvey-Price-Sokal_wardpoly}
to Schr\"oder paths,
together with some bijections mapping combinatorial objects
(in~particular, D-permutations)
to labeled Dyck, Motzkin or Schr\"oder paths.
We begin by reviewing briefly these two ingredients.

\subsection{Combinatorial interpretation of continued fractions}
   \label{subsec.prelimproofs.1}

Recall that a \textbfit{Motzkin path} of length $n \ge 0$
is a path $\omega = (\omega_0,\ldots,\omega_n)$
in the right quadrant $\N \times \N$,
starting at $\omega_0 = (0,0)$ and ending at $\omega_n = (n,0)$,
whose steps $s_j = \omega_j - \omega_{j-1}$
are $(1,1)$ [``rise'' or ``up step''], $(1,-1)$ [``fall'' or ``down step'']
or $(1,0)$ [``level step''].
We write $h_j$ for the \textbfit{height} of the Motzkin path at abscissa~$j$,
i.e.\ $\omega_j = (j,h_j)$;
note in particular that $h_0 = h_n = 0$.
We write $\scrm_n$ for the set of Motzkin paths of length~$n$,
and $\scrm = \bigcup_{n=0}^\infty \scrm_n$.
A Motzkin path is called a \textbfit{Dyck path} if it has no level steps.
A Dyck path always has even length;
we write $\scrd_{2n}$ for the set of Dyck paths of length~$2n$,
and $\scrd = \bigcup_{n=0}^\infty \scrd_{2n}$.

Let ${\bf a} = (a_i)_{i \ge 0}$, ${\bf b} = (b_i)_{i \ge 1}$
and ${\bf c} = (c_i)_{i \ge 0}$ be indeterminates;
we will work in the ring $\Z[[{\bf a},{\bf b},{\bf c}]]$
of formal power series in these indeterminates.
To each Motzkin path $\omega$ we assign a weight
$W(\omega) \in \Z[{\bf a},{\bf b},{\bf c}]$
that is the product of the weights for the individual steps,
where a rise starting at height~$i$ gets weight~$a_i$,
a~fall starting at height~$i$ gets weight~$b_i$,
and a level step at height~$i$ gets weight~$c_i$.
Flajolet \cite{Flajolet_80} showed that
the generating function of Motzkin paths
can be expressed as a continued fraction:

\begin{theorem}[Flajolet's master theorem]
   \label{thm.flajolet}
We have
\be
   \sum_{\omega \in \scrm}  W(\omega)
   \;=\;
   \cfrac{1}{1 - c_0 - \cfrac{a_0 b_1}{1 - c_1 - \cfrac{a_1 b_2}{1- c_2 - \cfrac{a_2 b_3}{1- \cdots}}}}
 \label{eq.thm.flajolet}
\ee
as an identity in $\Z[[{\bf a},{\bf b},{\bf c}]]$.
\end{theorem}

In particular, if $a_{i-1} b_i = \beta_i t^2$ and $c_i = \gamma_i t$
(note that the parameter $t$ is conjugate to the length of the Motzkin path),
we have
\be
   \sum_{n=0}^\infty t^n \sum_{\omega \in \scrm_n}  W(\omega)
   \;=\;
   \cfrac{1}{1 - \gamma_0 t - \cfrac{\beta_1 t^2}{1 - \gamma_1 t - \cfrac{\beta_2 t^2}{1 - \cdots}}}
   \;\,,
 \label{eq.flajolet.motzkin}
\ee
so that the generating function of Motzkin paths with height-dependent weights
is given by the J-type continued fraction \reff{def.Jtype}.
Similarly, if $a_{i-1} b_i = \alpha_i t$ and $c_i = 0$
(note that $t$ is now conjugate to the semi-length of the Dyck path), we have
\be
   \sum_{n=0}^\infty t^n \sum_{\omega \in \scrd_{2n}}  W(\omega)
   \;=\;
   \cfrac{1}{1 - \cfrac{\alpha_1 t}{1 - \cfrac{\alpha_2 t}{1 - \cdots}}}
   \;\,,
 \label{eq.flajolet.dyck}
\ee
so that the generating function of Dyck paths with height-dependent weights
is given by the S-type continued fraction \reff{def.Stype}.

Let us now show how to handle Schr\"oder paths within this framework.
A \textbfit{Schr\"oder path} of length $2n$ ($n \ge 0$)
is a path $\omega = (\omega_0,\ldots,\omega_{2n})$
in the right quadrant $\N \times \N$,
starting at $\omega_0 = (0,0)$ and ending at $\omega_{2n} = (2n,0)$,
whose steps are $(1,1)$ [``rise'' or ``up step''],
$(1,-1)$ [``fall'' or ``down step'']
or $(2,0)$ [``long level step''].
We write $s_j$ for the step starting at abscissa $j-1$.
If the step $s_j$ is a rise or a fall,
we set $s_j = \omega_j - \omega_{j-1}$ as before.
If the step $s_j$ is a long level step,
we set $s_j = \omega_{j+1} - \omega_{j-1}$ and leave $\omega_j$ undefined;
furthermore, in this case there is no step $s_{j+1}$.
We write $h_j$ for the height of the Schr\"oder path at abscissa~$j$
whenever this is defined, i.e.\ $\omega_j = (j,h_j)$.
Please note that $\omega_{2n} = (2n,0)$ and $h_{2n} = 0$
are always well-defined,
because there cannot be a long level step starting at abscissa $2n-1$.
Note also that a long level step at even (resp.~odd) height
can occur only at an odd-numbered (resp.~even-numbered) step.
We write $\scrs_{2n}$ for the set of Schr\"oder paths of length~$2n$,
and $\scrs = \bigcup_{n=0}^\infty \scrs_{2n}$.

There is an obvious bijection between Schr\"oder paths and Motzkin paths:
namely, every long level step is mapped onto a level step.
If we apply Flajolet's master theorem with
$a_{i-1} b_i = \alpha_i t$ and $c_i = \delta_{i+1} t$
to the resulting Motzkin path
(note that $t$ is now conjugate to the semi-length
 of the underlying Schr\"oder path),
we obtain
\be
   \sum_{n=0}^\infty t^n \sum_{\omega \in \scrs_{2n}}  W(\omega)
   \;=\;
   \cfrac{1}{1 - \delta_1 t - \cfrac{\alpha_1 t}{1 - \delta_2 t - \cfrac{\alpha_2 t}{1 - \cdots}}}
   \;\,,
 \label{eq.flajolet.schroder}
\ee
so that the generating function of Schr\"oder paths
with height-dependent weights
is given by the T-type continued fraction \reff{def.Ttype}.
More precisely, every rise gets a weight~1,
every fall starting at height~$i$ gets a weight $\alpha_i$,
and every long level step at height~$i$ gets a weight $\delta_{i+1}$.
This combinatorial interpretation of T-fractions in terms of Schr\"oder paths
was found recently by several authors
\cite{Fusy_15,Oste_15,Josuat-Verges_18,Sokal_totalpos}.

\subsection{Labeled Dyck, Motzkin and Schr\"oder paths}

Let $\bfscra = (\scra_h)_{h \ge 0}$, $\bfscrb = (\scrb_h)_{h \ge 1}$
and $\bfscrc = (\scrc_h)_{h \ge 0}$ be sequences of finite sets.
An
\textbfit{$(\bfscra,\bfscrb,\bfscrc)$-labeled Motzkin path of length $\bm{n}$}
is a pair $(\omega,\xi)$
where $\omega = (\omega_0,\ldots,\omega_n)$
is a Motzkin path of length $n$,
%% from $\omega_0 = (0,0)$ to $\omega_n = (n,0)$,
and $\xi = (\xi_1,\ldots,\xi_n)$ is a sequence satisfying
\be
   \xi_i  \:\in\:
   \begin{cases}
       \scra(h_{i-1})  & \textrm{if step $i$ is a rise (i.e.\ $h_i = h_{i-1} + 1$)}
              \\[1mm]
       \scrb(h_{i-1})  & \textrm{if step $i$ is a fall (i.e.\ $h_i = h_{i-1} - 1$)}
              \\[1mm]
       \scrc(h_{i-1})  & \textrm{if step $i$ is a level step (i.e.\ $h_i = h_{i-1}$)}
   \end{cases}
 \label{eq.xi.ineq}
\ee
where $h_{i-1}$ (resp.~$h_i$) is the height of the Motzkin path
before (resp.~after) step $i$.
[For typographical clarity
 we have here written $\scra(h)$ as a synonym for $\scra_h$, etc.]
We call $\xi_i$ the \textbfit{label} associated to step $i$.
We call the pair $(\omega,\xi)$
an \textbfit{$(\bfscra,\bfscrb)$-labeled Dyck path}
if $\omega$ is a Dyck path (in this case $\bfscrc$ plays no role).
We denote by $\scrm_n(\bfscra,\bfscrb,\bfscrc)$
the set of $(\bfscra,\bfscrb,\bfscrc)$-labeled Motzkin paths of length $n$,
and by $\scrd_{2n}(\bfscra,\bfscrb)$
the set of $(\bfscra,\bfscrb)$-labeled Dyck paths of length $2n$.

We define a \textbfit{$(\bfscra,\bfscrb,\bfscrc)$-labeled 
Schr\"oder path}
in an analogous way;
now the sets $\scrc_h$ refer to long level steps.
We denote by $\scrs_{2n}(\bfscra,\bfscrb,\bfscrc)$
the set of $(\bfscra,\bfscrb,\bfscrc)$-labeled Schr\"oder paths
of length $2n$.

Let us stress that the sets $\scra_h$, $\scrb_h$ and $\scrc_h$ are allowed
to be empty.
Whenever this happens, the path $\omega$ is forbidden to take a step
of the specified kind starting at the specified height.

%% We shall also make use of multicolored Motzkin paths.
%% An {\em $\ell$-colored Motzkin path}\/ is simply a Motzkin path
%% in which each level step has been given a ``color''
%% from the set $\{1,2,\ldots,\ell\}$.
%% In other words, we distinguish $\ell$ different types of level steps.
%% An {\em $(\bfscra,\bfscrb,\bfscrc^{(1)},\ldots,\bfscrc^{(\ell)})$-labeled
%%  $\ell$-colored Motzkin path of length $n$}\/
%% is then defined in the obvious way,
%% where we use the sequence $\bfscrc^{(j)}$ to bound
%% the label $\xi_i$ when step $i$ is a level step of type $j$.
%% We denote by $\scrm_n(\bfscra,\bfscrb,\bfscrc^{(1)},\ldots,\bfscrc^{(\ell)})$
%% the set of $(\bfscra,\bfscrb,\bfscrc^{(1)},\ldots,\bfscrc^{(\ell)})$-labeled
%% $\ell$-colored Motzkin paths of length $n$.

\bigskip

{\bf Remark.}  What we have called an
$(\bfscra,\bfscrb,\bfscrc)$-labeled Motzkin path
is (up to small changes in notation)
called a {\em path diagramme}\/ by Flajolet \cite[p.~136]{Flajolet_80}
and a {\em history}\/ by Viennot \cite[p.~II-9]{Viennot_83}.
Often the label sets $\scra_h, \scrb_h, \scrc_h$ are intervals of integers,
e.g.\ $\scra_h = \{1,\ldots,A_h\}$ or $\{0,\ldots,A_h\}$;
in this case the triplet $({\bf A},{\bf B},{\bf C})$
of sequences of maximum values is called a {\em possibility function}\/.
On the other hand, it is sometimes useful to employ labels
that are {\em pairs}\/ of integers
(e.g.\ \cite[Section~6.2]{Sokal-Zeng_masterpoly}
 and Section~\ref{sec.proofs.2} below).
It therefore seems preferable to state the general theory
without any specific assumption about the nature of the label sets.
\myendremark

\bigskip

Following Flajolet \cite[Proposition~7A]{Flajolet_80},
we can state a ``master J-fraction'' for
$(\bfscra,\bfscrb,\bfscrc)$-labeled Motzkin paths.
Let ${\bf a} = (a_{h,\xi})_{h \ge 0 ,\; \xi \in \scra(h)}$,
${\bf b} = (b_{h,\xi})_{h \ge 1 ,\; \xi \in \scrb(h)}$
and ${\bf c} = (c_{h,\xi})_{h \ge 0 ,\; \xi \in \scrc(h)}$
be indeterminates;
we give an $(\bfscra,\bfscrb,\bfscrc)$-labeled Motzkin path $(\omega,\xi)$
a weight $W(\omega,\xi)$
that is the product of the weights for the individual steps,
where a rise starting at height~$h$ with label $\xi$ gets weight~$a_{h,\xi}$,
a~fall starting at height~$h$ with label $\xi$ gets weight~$b_{h,\xi}$,
and a level step at height~$h$ with label $\xi$ gets weight~$c_{h,\xi}$.
Then:

\begin{theorem}[Flajolet's master theorem for labeled Motzkin paths]
   \label{thm.flajolet_master_labeled_Motzkin}
We have
\be
   \sum_{n=0}^\infty t^n
   \!\!
   \sum_{(\omega,\xi) \in \scrm_n(\bfscra,\bfscrb,\bfscrc)} \!\!\!  W(\omega,\xi)
   \;=\;
   \cfrac{1}{1 - c_0 t - \cfrac{a_0 b_1 t^2}{1 - c_1 t - \cfrac{a_1 b_2 t^2}{1- c_2 t - \cfrac{a_2 b_3 t^2}{1- \cdots}}}}
\ee
as an identity in $\Z[{\bf a},{\bf b},{\bf c}][[t]]$, where
\be
   a_h  \;=\;  \sum_{\xi \in \scra(h)} a_{h,\xi}
   \;,\qquad
   b_h  \;=\;  \sum_{\xi \in \scrb(h)} b_{h,\xi}
   \;,\qquad
   c_h  \;=\;  \sum_{\xi \in \scrc(h)} c_{h,\xi}
   \;.
 \label{def.weights.akbkck}
\ee
\end{theorem}

\noindent
This is an immediate consequence of Theorem~\ref{thm.flajolet}
together with the definitions.
% There is obviously also a similar theorem for
% $(\bfscra,\bfscrb,\bfscrc^{(1)},\ldots,\bfscrc^{(\ell)})$-labeled
% $\ell$-colored Motzkin paths,
% in which $c_k$ involves a sum over the colors of the level steps.

By specializing to ${\bf c} = \bzero$ and replacing $t^2$ by $t$,
we obtain the corresponding theorem
for $(\bfscra,\bfscrb)$-labeled Dyck paths:

\begin{corollary}[Flajolet's master theorem for labeled Dyck paths]
   \label{cor.flajolet_master_labeled_Dyck}
We have
\be
   \sum_{n=0}^\infty t^n
   \!\!
   \sum_{(\omega,\xi) \in \scrd_{2n}(\bfscra,\bfscrb)} \!\!\!  W(\omega,\xi)
   \;=\;
   \cfrac{1}{1 - \cfrac{a_0 b_1 t}{1 - \cfrac{a_1 b_2 t}{1- \cfrac{a_2 b_3 t}{1- \cdots}}}}
\ee
as an identity in $\Z[{\bf a},{\bf b}][[t]]$, where
$a_h$ and $b_h$ are defined by \reff{def.weights.akbkck}.
\end{corollary}

Similarly, for labeled Schr\"oder paths we have:

\begin{theorem}[Flajolet's master theorem for labeled Schr\"oder paths]
   \label{thm.flajolet_master_labeled_Schroder}
We have
\be
   \sum_{n=0}^\infty t^n
   \!\!
   \sum_{(\omega,\xi) \in \scrs_{2n}(\bfscra,\bfscrb,\bfscrc)} \!\!\!  W(\omega,\xi)
   \;=\;
   \cfrac{1}{1 - c_0 t - \cfrac{a_0 b_1 t}{1 - c_1 t - \cfrac{a_1 b_2 t}{1- c_2 t - \cfrac{a_2 b_3 t}{1- \cdots}}}}
\ee
as an identity in $\Z[{\bf a},{\bf b},{\bf c}][[t]]$, where
$a_h, b_h, c_h$ are defined by \reff{def.weights.akbkck},
with $c_{h,\xi}$ now referring to long level steps.
\end{theorem}

\section{First T-fraction: Proof of Theorems~\ref{thm.Tfrac.first},
   \ref{thm.Tfrac.first.pq}, \ref{thm.Tfrac.first.master},
   \ref{thm.Tfrac.first.master.variant} and \ref{thm.Tfrac.first.pq.variant}}
       \label{sec.proofs.1}

In this section we prove the first master T-fraction
(Theorem~\ref{thm.Tfrac.first.master})
by a bijection from D-permutations to labeled Schr\"oder paths.
Our construction combines ideas of
Randrianarivony \cite{Randrianarivony_97}
with a variant \cite[Section~6.1]{Sokal-Zeng_masterpoly}
of Foata--Zeilberger \cite{Foata_90},
together with some new ingredients.
After proving Theorem~\ref{thm.Tfrac.first.master},
we deduce Theorems~\ref{thm.Tfrac.first} and \ref{thm.Tfrac.first.pq}
by specialization.
Then, in Section~\ref{subsec.proofs.alternative},
we prove the variant T-fractions
of Theorems~\ref{thm.Tfrac.first.master.variant} 
and \ref{thm.Tfrac.first.pq.variant}.

%% \bigskip

Let us define an \textbfit{almost-Dyck path} of length $2n$
to be a path $\omega = (\omega_0,\ldots,\omega_{2n})$
in the right half-plane $\N \times \Z$,
starting at $\omega_0 = (0,0)$ and ending at $\omega_{2n} = (2n,0)$,
using the steps $(1,1)$ and $(1,-1)$,
that stays always at height $\ge -1$.
Thus, an almost-Dyck path is like a Dyck path {\em except that}\/
a down step from height 0 to height $-1$ is allowed;
note, however, that it must be immediately followed by
an up step back to height 0.
Each non-Dyck part of the path is therefore of the form
$(h_{2i-2},h_{2i-1},h_{2i}) = (0,-1,0)$.
We write $\scrd^\sharp_{2n}$ for the set of almost-Dyck paths of length $2n$.

Next let us define a \textbfit{0-Schr\"oder path}
to be a Schr\"oder path in which long level steps, if any,
occur only at height 0.
We write $\scrs^0_{2n}$ for the set of 0-Schr\"oder paths of length $2n$.
There is an obvious bijection $\psi\colon \scrd^\sharp_{2n} \to \scrs^0_{2n}$
from almost-Dyck paths to 0-Schr\"oder paths:
namely, we replace each down-up pair starting and ending at height 0
with a long level step at height 0.

In this section we will construct a bijection from 
D-permutations of $[2n]$ onto
labeled 0-Schr\"oder paths of length $2n$, as follows:
We first define the path by constructing an almost-Dyck path $\omega$
and then transforming it into a 0-Schr\"oder path $\omegahat = \psi(\omega)$.
Then we define the labels $\xi_i$, which will lie in the sets
\begin{subeqnarray}
   \scra_h        & = &  \{0,\ldots, \lceil h/2 \rceil \}  \qquad\quad\;\;\hbox{for $h \ge 0$}  \\
   \scrb_h        & = &  \{0,\ldots, \lceil (h-1)/2 \rceil \}    \quad\hbox{for $h \ge 1$}  \\
   \scrc_0        & = &  \{0\}       \\
   \scrc_h        & = &  \emptyset   \qquad\qquad\qquad\qquad\quad\hbox{for $h \ge 1$}
 \label{def.abc}
\end{subeqnarray}
We also interpret our crossing, nesting and record statistics
in terms of the heights and labels.
Next we prove that the map $\sigma \mapsto (\omegahat,\xi)$
really is a bijection from the set $\dperm_{2n}$ of D-permutations of $[2n]$
onto the set $\scrs_{2n}(\bfscra,\bfscrb,\bfscrc)$
of $(\bfscra,\bfscrb,\bfscrc)$-labeled Schr\"oder paths of length $2n$.
%% After that, {\bf or before???}
%% we will translate the various statistics from $\dperm_{2n}$
%% to our labeled Schr\"oder paths.
Finally, we sum over the labels $\xi$ to obtain the weight $W(\omegahat)$
associated to a Schr\"oder path $\omegahat$,
which upon applying \reff{eq.flajolet.schroder}
will yield Theorem~\ref{thm.Tfrac.first.master}.

%\bigskip

\subsection{Step 1: Definition of the almost-Dyck path}

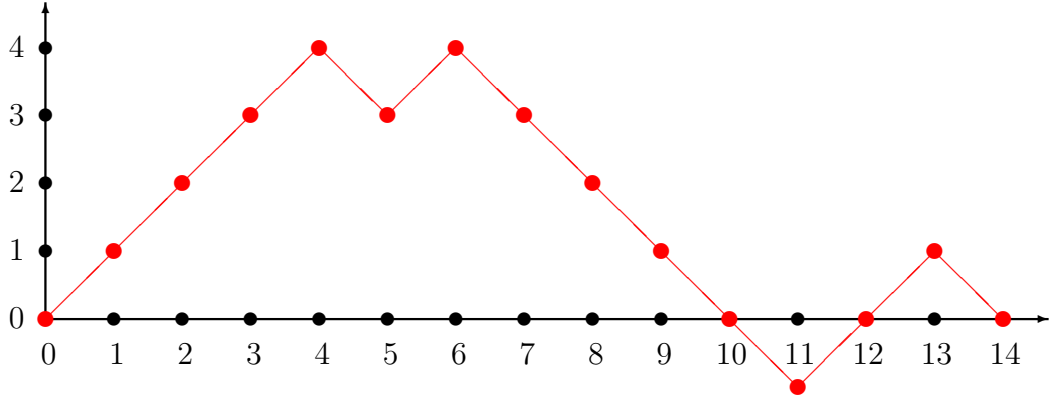
\begin{figure}[t]
\vspace*{1.5cm}
\hspace*{2cm}
\centering
\begin{picture}(50,120)(200,-40)
\setlength{\unitlength}{1.2cm}
\linethickness{.3mm}
\put(0,0){\vector(1,0){11}}\put(0,0){\vector(0,1){3.5}}
\put(-0.05,-0.5){$0$}
\put(0.7,-0.5){$1$}\put(0.75,0){\circle*{.15}}
\put(1.45,-0.5){$2$}\put(1.5,0){\circle*{.15}}
\put(2.2,-0.5){$3$}\put(2.25,0){\circle*{.15}}
\put(2.95,-0.5){$4$}\put(3,0){\circle*{.15}}
\put(3.7,-0.5){$5$}\put(3.75,0){\circle*{.15}}
\put(4.45,-0.5){$6$}\put(4.5,0){\circle*{.15}}
\put(5.2,-0.5){$7$}\put(5.25,0){\circle*{.15}}
\put(5.95,-0.5){$8$}\put(6,0){\circle*{.15}}
\put(6.7,-0.5){$9$}\put(6.75,0){\circle*{.15}}
\put(7.35,-0.5){$10$}\put(7.5,0){\circle*{.15}}
\put(8.1,-0.5){$11$}\put(8.25,0){\circle*{.15}}
\put(8.85,-0.5){$12$}\put(9,0){\circle*{.15}}
\put(9.6,-0.5){$13$}\put(9.75,0){\circle*{.15}}
\put(10.35,-0.5){$14$}\put(10.5,0){\circle*{.15}}
%\put(3.35,-1.2){(b)}
%%%%%%%%%%%%%%%%%
\put(-0.4,-0.1){$0$}
\put(-0.4,0.65){$1$}\put(0,0.75){\circle*{.15}}
\put(-0.4,1.4){$2$}\put(0,1.5){\circle*{.15}}
\put(-0.4,2.15){$3$}\put(0,2.25){\circle*{.15}}
\put(-0.4,2.9){$4$}\put(0,3){\circle*{.15}}
%%%%%%%%%%%%%%%%%%
\put(0,0){\red{\line(1,1){0.75}}}
\put(0.75,0.75){\red{\line(1,1){0.75}}}
\put(1.5,1.5){\red{\line(1,1){0.75}}}
\put(2.25,2.25){\red{\line(1,1){0.75}}}
\put(3,3){\red{\line(1,-1){0.75}}}
\put(3.75,2.25){\red{\line(1,1){0.75}}}
\put(4.5,3){\red{\line(1,-1){0.75}}}
\put(5.25,2.25){\red{\line(1,-1){0.75}}}
\put(6,1.5){\red{\line(1,-1){0.75}}}
\put(6.75,0.75){\red{\line(1,-1){0.75}}}
\put(7.5,0){\red{\line(1,-1){0.75}}}
\put(8.25,-0.75){\red{\line(1,1){0.75}}}
\put(9,0){\red{\line(1,1){0.75}}}
\put(9.75,0.75){\red{\line(1,-1){0.75}}}
%%%%%%%%
\put(0,0){\red{\circle*{.17}}}
\put(0.75,0.75){\red{\circle*{.17}}}
\put(1.5,1.5){\red{\circle*{.17}}}
\put(2.25,2.25){\red{\circle*{.17}}}
\put(3,3){\red{\circle*{.17}}}
\put(3.75,2.25){\red{\circle*{.17}}}
\put(4.5,3){\red{\circle*{.17}}}
\put(5.25,2.25){\red{\circle*{.17}}}
\put(6,1.5){\red{\circle*{.17}}}
\put(6.75,0.75){\red{\circle*{.17}}}
\put(7.5,0){\red{\circle*{.17}}}
\put(8.25,-0.75){\red{\circle*{.17}}}
\put(9,0){\red{\circle*{.17}}}
\put(9.75,0.75){\red{\circle*{.17}}}
\put(10.5,0){\red{\circle*{.17}}}
%\put(2,2){\red{\circle*{.17}}}
%\put(3,2){\red{\circle*{.17}}}
%\put(4,2){\red{\circle*{.17}}}
%\put(5,1){\red{\circle*{.17}}}
%\put(6,1){\red{\circle*{.17}}}
%\put(7,0){\red{\circle*{.17}}}
\end{picture}
\vspace*{-5mm}
\caption{Almost-Dyck path $\omega$ corresponding to the D-permutation
 $\sigma = 7\, 1\, 9\, 2\, 5\, 4\, 8\, 6\, 10\, 3\, 11\, 12\, 14\, 13
   \;=\; (1,7,8,6,4,2)\,(3,9,10)\,(5)\,(11)\,(12)\,(13,14)$.
 This is the same D-permutation that was illustrated in
 Figure~\ref{fig.pictorial}.
}
\vspace*{7mm}
   \label{fig.almostdyck}
\end{figure}

Given a D-permutation $\sigma \in \dperm_{2n}$,
we define a path $\omega = (\omega_0,\ldots,\omega_{2n})$
starting at $\omega_0 = (0,0)$,
%% and ending at $\omega_{2n} = (2n,0)$,
with steps $s_1,\ldots,s_{2n}$ as follows:
\begin{itemize}
   \item If $\sinv(i)$ is even, then $s_i$ is a rise.
      (Note that in this case we must have $\sinv(i) \ge i$,
       by definition of D-permutation.)
   \item If $\sinv(i)$ is odd, then $s_i$ is a fall.
      (Note that in this case we must have $\sinv(i) \le i$,
       by definition of D-permutation.)
\end{itemize}
(See Figure~\ref{fig.almostdyck} for an example.)
An alternative way of saying this is:
\begin{itemize}
   \item If $\sinv(i) > i$, then $s_i$ is a rise.
      (In this case $\sinv(i)$ must be even.)
   \item If $\sinv(i) < i$, then $s_i$ is a fall.
      (In this case $\sinv(i)$ must be odd.)
   \item If $i$ is a fixed point, then $s_i$ is a rise if $i$ is even,
       and a fall if $i$ is odd.
\end{itemize}
Yet another alternative way of saying this is:
\begin{itemize}
   \item If $i$ is a cycle valley, cycle double fall or even fixed point,
      then $s_i$ is a rise.
   \item If $i$ is a cycle peak, cycle double rise or odd fixed point,
      then $s_i$ is a fall.
\end{itemize}
Of course we need to prove that this $\omega$ is indeed an almost-Dyck path,
i.e.\ that all the heights $h_i$ are $\ge -1$ and that $h_{2n} = 0$.
We will do this by obtaining a precise interpretation of the heights $h_i$.

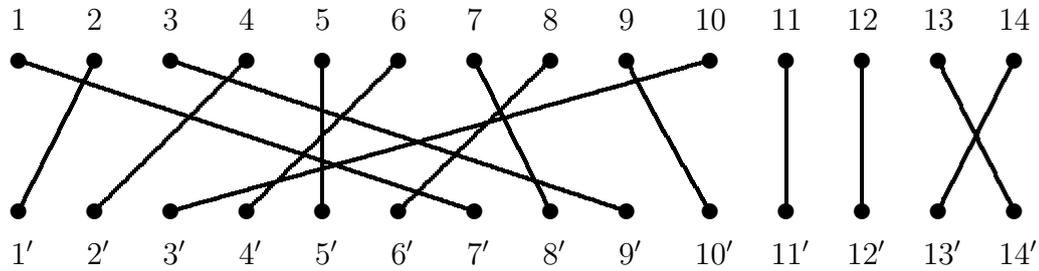
\begin{figure}[t]
\vspace*{3cm}
\hspace*{-3mm}
\begin{picture}(20,3)(-35,-20)
\setlength{\unitlength}{10mm}
\linethickness{.5mm}
\put(-0.1,2.4){$1$}\put(0,2){\circle*{0.2}}
\put(0.9,2.4){$2$}\put(1,2){\circle*{0.2}}
\put(1.9,2.4){$3$}\put(2,2){\circle*{0.2}}
\put(2.9,2.4){$4$}\put(3,2){\circle*{0.2}}
\put(3.9,2.4){$5$}\put(4,2){\circle*{0.2}}
\put(4.9,2.4){$6$}\put(5,2){\circle*{0.2}}
\put(5.9,2.4){$7$}\put(6,2){\circle*{0.2}}
\put(6.9,2.4){$8$}\put(7,2){\circle*{0.2}}
\put(7.9,2.4){$9$}\put(8,2){\circle*{0.2}}
\put(8.9,2.4){$10$}\put(9.1,2){\circle*{0.2}}
\put(9.9,2.4){$11$}\put(10.1,2){\circle*{0.2}}
\put(10.9,2.4){$12$}\put(11.1,2){\circle*{0.2}}
\put(11.9,2.4){$13$}\put(12.1,2){\circle*{0.2}}
\put(12.9,2.4){$14$}\put(13.1,2){\circle*{0.2}}
%%%
\put(-0.1,-0.7){$1'$}\put(0,0){\circle*{0.2}}
\put(0.9,-0.7){$2'$}\put(1,0){\circle*{0.2}}
\put(1.9,-0.7){$3'$}\put(2,0){\circle*{0.2}}
\put(2.9,-0.7){$4'$}\put(3,0){\circle*{0.2}}
\put(3.9,-0.7){$5'$}\put(4,0){\circle*{0.2}}
\put(4.9,-0.7){$6'$}\put(5,0){\circle*{0.2}}
\put(5.9,-0.7){$7'$}\put(6,0){\circle*{0.2}}
\put(6.9,-0.7){$8'$}\put(7,0){\circle*{0.2}}
\put(7.9,-0.7){$9'$}\put(8,0){\circle*{0.2}}
\put(8.9,-0.7){$10'$}\put(9.1,0){\circle*{0.2}}
\put(9.9,-0.7){$11'$}\put(10.1,0){\circle*{0.2}}
\put(10.9,-0.7){$12'$}\put(11.1,0){\circle*{0.2}}
\put(11.9,-0.7){$13'$}\put(12.1,0){\circle*{0.2}}
\put(12.9,-0.7){$14'$}\put(13.1,0){\circle*{0.2}}
%%%%%%%
\qbezier(0,2)(3,1)(6,0)
\qbezier(1,2)(0.5,1)(0,0)
\qbezier(2,2)(5,1)(8,0)
\qbezier(3,2)(2,1)(1,0)
\qbezier(4,2)(4,1)(4,0)
\qbezier(5,2)(4,1)(3,0)
\qbezier(6,2)(6.5,1)(7,0)
\qbezier(7,2)(6,1)(5,0)
\qbezier(8,2)(8.5,1)(9.1,0)
\qbezier(9.1,2)(5.55,1)(2,0)
\qbezier(10.1,2)(10.1,1)(10.1,0)
\qbezier(11.1,2)(11.1,1)(11.1,0)
\qbezier(12.1,2)(12.6,1)(13.1,0)
\qbezier(13.1,2)(12.6,1)(12.1,0)
\end{picture}
\caption{
   Bipartite digraph representing the permutation
   $\sigma = 7\, 1\, 9\, 2\, 5\, 4\, 8\, 6\, 10\, 3\, 11\, 12\, 14\, 13
     \;=\; (1,7,8,6,4,2)\,(3,9,10)\,(5)\,(11)\,(12)\,(13,14)$.
   Arrows run from the top row to the bottom row and are suppressed
   for clarity.
 \label{fig.bipartite}
 \vspace*{7mm}
}
\end{figure}

In what follows, it will be convenient to represent
a permutation $\sigma \in \Sym_N$ by a bipartite digraph
$\Gamma = \Gamma(\sigma)$
in which the top row of vertices is labeled $1,\ldots,N$
and the bottom row $1',\ldots,N'$,
and we draw an arrow from $i$ to $j'$ in case $\sigma(i) = j$
(see Figure~\ref{fig.bipartite}).
For $k \in [N]$, we denote by $\Gamma_k$
the induced subgraph of $\Gamma$ on the vertex set
$\{1,\ldots,k\} \cup \{1',\ldots,k'\}$.
Thus, the edges of $\Gamma_k$ are arrows $i \to j'$
drawn whenever $\sigma(i) = j$ and $i \le k$ and $j \le k$.
We say that a vertex of $\Gamma_k$ is \textbfit{free}
if no arrow is incident on it.
We write
\begin{subeqnarray}
   f_k
   & \eqdef &
   \hbox{\# of free vertices in the top row of $\Gamma_k$}
      \\[2mm]
   & = &
   \#\{ i \le k \colon\: \sigma(i) > k \}
   \;.
 \label{eq.fk.1}
\end{subeqnarray}
Of course, we also have
\begin{subeqnarray}
   f_k
   & = &
   \hbox{\# of free vertices in the bottom row of $\Gamma_k$}
      \\[2mm]
   & = &
   \#\{ j \le k \colon\: \sigma^{-1}(j) > k \}
   \;.
 \label{eq.fk.2}
\end{subeqnarray}
The total number of free vertices in $\Gamma_k$ is therefore $2f_k$.
Note that $f_k = 0$ if and only if $\sigma$ maps $\{1,\ldots,k\}$ onto itself.
Note also that
\be
   f_k - f_{k-1}
   \;=\;
   \begin{cases}
       +1     & \textrm{if $k$ is a cycle valley}  \\
       -1     & \textrm{if $k$ is a cycle peak}  \\
       0      & \textrm{if $k$ is a cycle double rise,
                                    cycle double fall, or fixed point}
   \end{cases}
 \label{eq.fk.fk1}
\ee

{\bf Remark.}
The sequence $(f_0,\ldots,f_{2n})$ is a Motzkin path;
in fact, it is precisely the Motzkin path associated to
the permutation $\sigma$ by the Foata--Zeilberger (or Biane) bijection.
To see this, compare \reff{eq.fk.1}/\reff{eq.fk.2}
with \cite[eq.~(6.4)]{Sokal-Zeng_masterpoly};
or equivalently, compare \reff{eq.fk.fk1}
with \cite[definition of steps $s_i$ preceding (6.2)]{Sokal-Zeng_masterpoly}.
\myendremark

\medskip

We can now give the promised interpretation of the heights:

\begin{lemma}[Interpretation of the heights]
   \label{lemma.heights}
For $k \in [2n]$ we have
\be
   h_k
   \;=\;
   \begin{cases}
       2f_k - 1   & \textrm{if $k$ is odd}  \\
       2f_k       & \textrm{if $k$ is even}
   \end{cases}
 \label{eq.lemma.heights}
\ee
In particular, $h_k \ge -1$ and $h_{2n} = 0$,
so that $\omega$ is an almost-Dyck path.

Furthermore, we have $(h_{2i-2}, h_{2i-1}, h_{2i}) = (0, -1, 0)$
if and only if $2i-1$ and $2i$ are record-antirecords.
\end{lemma}

\noindent
Please note that, by \reff{eq.lemma.heights},
the parity of $h_k$ equals the parity of $k$;
this reflects the fact that $\omega$ is an almost-Dyck path.
Note also that \reff{eq.lemma.heights} can be rewritten as
\be
   f_k  \;=\; \left\lceil \dfrac{h_k}{2} \right\rceil
     \;=\;
   \begin{cases}
       (h_k+1)/2   & \textrm{if $k$ is odd}  \\[1mm]
       h_k/2        & \textrm{if $k$ is even}
   \end{cases}
 \label{eq.lemma.heights.bis}
\ee
And recall, finally, from Lemma~\ref{lemma.recantirec_D}
that record-antirecords come in pairs:
$2i-1$ is a record-antirecord if and only if $2i$ is a record-antirecord.

\proofof{Lemma~\ref{lemma.heights}}
We have
\begin{subeqnarray}
   h_k
   & = &
   \#\{j\leq k \colon\: \sinv(j) \text{ is even}\}
      \:-\: \#\{j\leq k \colon\: \sinv(j) \text{ is odd}\}
         \\[3mm]
   & = &
   2\, \#\{j\leq k \colon\: \sinv(j) \text{ is even}\} \:-\: k
         \\[3mm]
   & = &
   2\, \bigl[ \#\{j\leq k\colon\: \sinv(j)>k\}  \:+\:
              \#\{j\leq k \colon\: k\geq\sinv(j)>j\}
         \nonumber \\[1mm]
   &  & \qquad +\: \#\{j\leq k \colon\: \sinv(j) = j \text{ is even}\} \bigr]
             \;-\; k
         \\[3mm]
   & = &
   2\, \bigl[ \#\{j\leq k: \sinv(j)>k\}  \:+\:
              \#\{i\leq k \colon\: i>\sigma(i)\}
         \nonumber \\[1mm]
   &  & \qquad +\: \#\{i\leq k \colon\: \sigma(i)=i \text{ is even}\} \bigr]
             \;-\; k
         \\[3mm]
   & = &
   2\, \bigl[ \#\{j\leq k \colon\: \sinv(j)>k\} \:+\:
             \#\{i\leq k \colon\: i\text{ is even}\} \bigr]
             \;-\; k
         \\[3mm]
   & = &
   2f_k \:+\: 2\left\lceil \frac{k-1}{2} \right\rceil \:-\: k 
   \;,
\end{subeqnarray}
which is \reff{eq.lemma.heights}.

Furthermore, $h_k = -1$ occurs when and only when
$k$ is odd (say, $k=2i-1$) and $f_k = 0$.
The latter statement means that $\sigma$ maps
$\{1,\ldots,2i-1\}$ onto itself.
By Lemma~\ref{lemma.recantirec_D}(c)$\implies$(a,b)
this means that $2i-1$ and $2i$ are record-antirecords.
Conversely, if $2i-1$ is a record-antirecord,
then by Lemma~\ref{lemma.recantirec_D}(a)$\implies$(d)
we have $f_{2i-2} = f_{2i-1} = f_{2i} = 0$,
so that $(h_{2i-2}, h_{2i-1}, h_{2i}) = (0, -1, 0)$.
\qed

\smallskip

{\bf Remarks.}
1. Randrianarivony \cite[Section~6]{Randrianarivony_97}
considered the special case of this construction
in which $\sigma$ is a D-o-semiderangement
(in his terminology, a ``Genocchi permutation'');
in this case $\omega$ is a Dyck path.
Our proof of Lemma~\ref{lemma.heights} is a very slight modification of his,
designed to allow for fixed points of both parities.
See also Han and Zeng \cite[pp.~126--127]{Han_99a}
for the interpretation in terms of bipartite graphs and free vertices.

2. The number of almost-Dyck paths of length $2n$ is $C_{n+1}$,
where $C_n = \frac{1}{n+1} \binom{2n}{n}$ is the $n$th Catalan number:
it suffices to observe that an almost-Dyck path of semilength $n$
can be converted to a Dyck path of semilength $n+1$
by adding a rise at the beginning and a fall at the end,
and conversely.
Equivalently, the number of 0-Schr\"oder paths of length $2n$ is $C_{n+1}$:
this follows from \reff{eq.flajolet.schroder}
with $\delta_1 = 1$, $\delta_n = 0$ for $n \ge 2$,
and $\alpha_n = 1$ for $n \ge 1$
together with the identity
\be
   {1 \over 1 - t - tC(t)}  \;=\;  {C(t) - 1 \over t}
\ee
for the Catalan generating function
$C(t) = \sum\limits_{n=0}^\infty C_n \, t^n = (1 - \sqrt{1-4t})/(2t)$.
\myendremark

%\medskip

\subsection[Step 2: Definition of the labels $\xi_i$]{Step 2: Definition of the labels $\bm{\xi_i}$} 

%% This is Bishal's "version 2".
We now define
\be
   \xi_i
   \;=\;
   \begin{cases}
       \#\{j \colon\: \sigma(j) < \sigma(i) \leq i <j\}
           & \textrm{if $i$ is even}
              \\[2mm]
     \#\{j \colon\: j<i \leq \sigma(i)<\sigma(j) \}
           & \textrm{if $i$ is odd}
   \end{cases}
 \label{def2.xi}
\ee
Note that the middle inequalities in this definition hold automatically:
a D-permutation {\em always}\/ has $\sigma(i) \leq i$ if $i$ is even,
and $i \leq \sigma(i)$ if $i$ is odd.
The definition \reff{def2.xi} can be written equivalently as
\be
   \xi_i
   \;=\;
   \begin{cases}
       \#\{2l > 2k \colon\: \sigma(2l) < \sigma(2k) \}
           & \textrm{if $i=2k$}
              \\[2mm]
     \#\{2l-1 < 2k-1 \colon\: \sigma(2l-1) > \sigma(2k-1) \}
           & \textrm{if $i=2k-1$}
   \end{cases}
 \label{def2.xi.bis}
\ee
since $\sigma(j) < j$ implies that $j$ is even,
and $j < \sigma(j)$ implies that $j$ is odd.

It is worth remarking that the labels $\xi_i$ defined in \reff{def2.xi}
are the same as those in the variant Foata--Zeilberger bijection
\cite[eq.~(6.5)]{Sokal-Zeng_masterpoly}
whenever $i$ is not a fixed point.\footnote{
   With the only difference that in \cite{Sokal-Zeng_masterpoly}
   the labels were defined to start at~1,
   whereas here they start at~0.
}

The label $\xi_i$ has a simple interpretation in terms of the
index-refined nesting statistics defined in
\reff{def.ucrossnestjk}/\reff{def.psnestj}:

\begin{figure}[t]
\centering
\begin{picture}(30,15)(160, 10)
\setlength{\unitlength}{4mm}
\linethickness{.5mm}
\put(2,0){\line(1,0){27}}
\put(5,0){\circle*{1,3}}\put(5,0){\makebox(0,3)[c]{\small $\sigma(j)$}}
\put(12,0){\circle*{1,3}}\put(12,0){\makebox(0,3)[c]{\small $\sigma(i)$}}
\put(19,0){\circle*{1,3}}\put(19,0){\makebox(0,3)[c]{\small $i$}}
\put(26,0){\circle*{1,3}}\put(26,0){\makebox(0,3)[c]{\small $j$}}
\red{\qbezier(5,0)(15.5,-8)(26,0)}
\blue{\qbezier(11.5,0)(15.0,-4)(18.5,0)}
\put(15,-3){\makebox(0,-6)[c]{\small When $i$ is even and $\neq \sigma(i)$}}
%%%%%%%%%%%%%%%%%%%%
\end{picture}
\\[5cm]
\begin{picture}(30,15)(640, 10)
\setlength{\unitlength}{4mm}
\linethickness{.5mm}
\put(44,0){\line(1,0){27}}
\put(47,0){\circle*{1,3}}\put(47,0){\makebox(0,-3)[c]{\small $j$}}
\put(54,0){\circle*{1,3}}\put(54,0){\makebox(0,-3)[c]{\small $i$}}
\put(61,0){\circle*{1,3}}\put(61,0){\makebox(0,-3)[c]{\small $\sigma(i)$}}
\put(68,0){\circle*{1,3}}\put(68,0){\makebox(0,-3)[c]{\small $\sigma(j)$}}
\red{\qbezier(47,0)(57.5,8)(68,0)}
\blue{\qbezier(53.5,0)(57.0,4)(60.5,0)}
\put(57,-1){\makebox(0,-6)[c]{\small When $i$ is odd and $\neq \sigma(i)$}}
\end{picture}
\\[2.5cm]
\caption{
   Nestings involved in the definition of the label $\xi_i$.
 \label{fig.nestings.labels}
}
\vspace*{2.5cm}
\end{figure}

\begin{lemma}[Nesting statistics]
   \label{lemma.xii.nesting}
We have
\be
   \xi_i
   \;=\;
   \begin{cases}
        \lnest(i,\sigma)
           & \hbox{\rm if $i$ is even and $\neq \sigma(i)$}
             \;\, \hbox{\rm [equivalently, $i > \sigma(i)$]}
              \\[0.5mm]
        \unest(i,\sigma)
           & \hbox{\rm if $i$ is odd and $\neq \sigma(i)$}
             \;\, \hbox{\rm [equivalently, $i < \sigma(i)$]}
              \\[0.5mm]
        \psnest(i,\sigma)
           & \hbox{\rm if $i = \sigma(i)$  (that is, $i$ is a fixed point)}
   \end{cases}
 \label{eq.xii.nestings}
\ee
\end{lemma}

\noindent
See Figure~\ref{fig.nestings.labels}.
Note also that if $i$ is a fixed point, we have
\be
   \xi_i  \;=\;  f_i  \;=\; \left\lceil \dfrac{h_i}{2} \right\rceil
     \;=\;
   \begin{cases}
       h_i/2        & \textrm{if $i$ is even} \\[1mm]
       (h_i+1)/2   & \textrm{if $i$ is odd}
   \end{cases}
 \label{eq.xii.fp}
\ee

Now we shall show, as required by \reff{def.abc},
that the following inequalities hold true:

\begin{lemma}[Inequalities satisfied by the labels]
   \label{lemma.ineqs.labels}
We have
\begin{subeqnarray}
   & 0 \;\le\; \xi_i  \;\le\; \left\lceil \dfrac{h_i-1}{2} \right\rceil \;=\;
                      \left\lceil\dfrac{h_{i-1}}{2} \right\rceil
       & \textrm{if $\sinv(i) $ is even (i.e., $s_i$ is a rise)}
    \qquad
 \slabel{eq.xi.ineqs.a}
      \\[3mm]
   & 0 \;\le\; \xi_i  \;\le\; \left\lceil \dfrac{h_i}{2} \right\rceil  \;=\;
                      \left\lceil\dfrac{h_{i-1}-1}{2} \right\rceil
       & \textrm{if $\sinv(i) $ is odd (i.e., $s_i$ is a fall)}
 \slabel{eq.xi.ineqs.b}
 \label{eq.xi.ineqs}
\end{subeqnarray}
\end{lemma}

\noindent
We remark that a fall starting at height $h_{i-1} = 0$,
or a rise starting at height $h_{i-1} = -1$,
always gets the label $\xi_i  = 0$.
When we pass from the almost-Dyck path $\omega$
to the 0-Schr\"oder path $\omegahat = \psi(\omega)$,
these labels become the label $\xi = 0$ for the long level step at height~0.

To prove the inequalities \reff{eq.xi.ineqs}, we will interpret
$\left\lceil \dfrac{h_i-1}{2} \right\rceil - \xi_i$ when $s_i$ is a rise,
and $\left\lceil \dfrac{h_i}{2} \right\rceil - \xi_i$ when $s_i$ is a fall,
in terms of crossing statistics, as follows:

\begin{lemma}[Crossing statistics]
   \label{lemma.xii.crossing}
\nopagebreak
\quad\hfill
\vspace*{-1mm}
\begin{itemize}
   \item[(a)]  If $s_i$ a rise and $i$ (hence also $h_i$) is odd, then
\be
   \left\lceil \dfrac{h_i - 1}{2} \right\rceil - \xi_i
   \;=\;
   \ucross(i,\sigma)
   \;.
 \label{eq.lemma.xii.crossing.a}
\ee
   \item[(b)]  If $s_i$ a rise and $i$ (hence also $h_i$) is even, then
\begin{subeqnarray}
   \left\lceil \dfrac{h_i - 1}{2} \right\rceil - \xi_i
   & = &
   \lcross(i,\sigma) \,+\, {\rm I}[\sigma(i) \neq i]
       \\[-1mm]
   & = &
   \lcross(i,\sigma) \,+\, {\rm I}[\hbox{\rm $i$ is a cycle double fall}]
   \;.
 \label{eq.lemma.xii.crossing.b}
\end{subeqnarray}
   \item[(c)]  If $s_i$ a fall and $i$ (hence also $h_i$) is odd, then
\begin{subeqnarray}
   \left\lceil \dfrac{h_i}{2} \right\rceil - \xi_i
   & = &
   \ucross(i,\sigma) \,+\, {\rm I}[\sigma(i) \neq i]
       \\[-1mm]
   & = &
   \ucross(i,\sigma) \,+\, {\rm I}[\hbox{\rm $i$ is a cycle double rise}]
   \;.
 \label{eq.lemma.xii.crossing.c}
\end{subeqnarray}
   \item[(d)]  If $s_i$ a fall and $i$ (hence also $h_i$) is even, then
\be
   \left\lceil \dfrac{h_i}{2} \right\rceil - \xi_i
   \;=\;
   \lcross(i,\sigma)
   \;.
 \label{eq.lemma.xii.crossing.d}
\ee
\end{itemize}
(Here ${\rm I}[\hbox{\sl proposition}] = 1$ if {\sl proposition} is true,
and 0 if it is false.)
\end{lemma}

Since the right-hand sides of
\reff{eq.lemma.xii.crossing.a}--\reff{eq.lemma.xii.crossing.d}
are manifestly nonnegative,
Lemma~\ref{lemma.ineqs.labels} is an immediate consequence of
Lemma~\ref{lemma.xii.crossing}.

\begin{figure}[t]
\centering
\begin{picture}(30,15)(160, 10)
\setlength{\unitlength}{4mm}
\linethickness{.5mm}
\put(2,0){\line(1,0){27}}
\put(5,0){\circle*{1,3}}\put(5,0){\makebox(0,3)[c]{\small $\sigma(i)$}}
\put(12,0){\circle*{1,3}}\put(12,0){\makebox(0,3)[c]{\small $\sigma(j)$}}
\put(19,0){\circle*{1,3}}\put(19,0){\makebox(0,3)[c]{\small $i$}}
\put(26,0){\circle*{1,3}}\put(26,0){\makebox(0,3)[c]{\small $j$}}
\red{\qbezier(5,0)(12,-6)(19,0)}
\blue{\qbezier(11.5,0)(18.5,-6)(25.5,0)}
\put(15,-2){\makebox(0,-6)[c]{\small (b,d) When $i$ is even and $\neq \sigma(i)$}}
% %%%%%%%%%%%%%%%%%%%%
\end{picture}
\\[4.3cm]
\begin{picture}(30,15)(640, 10)
\setlength{\unitlength}{4mm}
\linethickness{.5mm}
\put(44,0){\line(1,0){27}}
\put(47,0){\circle*{1,3}}\put(47,0){\makebox(0,-3)[c]{\small $j$}}
\put(54,0){\circle*{1,3}}\put(54,0){\makebox(0,-3)[c]{\small $i$}}
\put(61,0){\circle*{1,3}}\put(61,0){\makebox(0,-3)[c]{\small $\sigma(j)$}}
\put(68,0){\circle*{1,3}}\put(68,0){\makebox(0,-3)[c]{\small $\sigma(i)$}}
\red{\qbezier(47,0)(54,6)(61,0)}
\blue{\qbezier(53.5,0)(60.5,6)(67.5,0)}
\put(57,-0.5){\makebox(0,-6)[c]{\small (a,c) When $i$ is odd and $\neq \sigma(i)$}}
\end{picture}
\\[2.5cm]
\caption{
   Crossings involved in the identities for the label $\xi_i$.
 \label{fig.xi.cross}
}
\vspace*{2cm}
\end{figure}

\proofof{Lemma~\ref{lemma.xii.crossing}}
(a,b) If $s_i$ is a rise, then $\sinv(i)$ is even and $i \le \sinv(i)$.
    We now consider separately the cases of $i$ odd and $i$ even.

\medskip

(a) If $i$ is odd, then $h_i$ is odd;
      moreover, since $\sinv(i)$ is even,
      $i$ cannot be a fixed point,
      and we have the strict inequalities $i < \sinv(i)$
      and $i < \sigma(i)$.
      Then
\begin{subeqnarray}
   \left\lceil \dfrac{h_i - 1}{2} \right\rceil - \xi_i
   & = &
   \dfrac{h_i - 1}{2} \,-\, \xi_i
       \\[2mm]
   & = &
   f_i \,-\, 1 \,-\, \xi_i
       \\[2mm]
   & = &
   \#\{ j \le i \colon\: \sigma(j) > i \}   \,-\,  1
       \,-\, \#\{ j \colon\: j < i \le \sigma(i) < \sigma(j) \}
       \qquad\qquad \\[2mm]
   & = &
   \#\{ j < i \colon\: \sigma(j) > i \}
       \,-\, \#\{ j \colon\: j < i < \sigma(i) < \sigma(j) \}
       \qquad\qquad \\[2mm]
   & = &
   \#\{ j \colon\:  j < i < \sigma(j) < \sigma(i) \}
       \qquad\qquad \\[2mm]
   & = &
   \ucross(i,\sigma)
   \;.
\end{subeqnarray}
See Figure~\ref{fig.xi.cross}(a,c).

\medskip

(b) If $i$ is even, then $h_i$ is even, and $\sigma(i) \le i$.
So either $i$ is a fixed point (hence $\sigma(i) = i = \sinv(i)$)
or else $\sigma(i) < i < \sinv(i)$.
Then
\begin{subeqnarray}
   \left\lceil \dfrac{h_i - 1}{2} \right\rceil - \xi_i
   & = &
   \dfrac{h_i}{2} \,-\, \xi_i
       \\[2mm]
   & = &
   f_i \,-\, \xi_i
       \\[2mm]
   & = &
   \#\{ j \le i \colon\: \sinv(j) > i \}   
       \,-\, \#\{j \colon\: \sigma(j) < \sigma(i) \leq i <j \}
       \qquad\qquad \\[2mm]
   & = &
   \#\{ j < i \colon\: \sinv(j) > i \}  \,+\, {\rm I}[\sinv(i) > i]
       \nonumber \\
   & & \hspace*{4cm}
       \,-\, \#\{j \colon\: \sigma(j) < \sigma(i) \leq i <j \}
       \qquad\qquad \\[2mm]
   & = &
   \#\{ j \colon\: \sigma(j) < i < j \} \,+\, {\rm I}[\sigma(i) \neq i]
       \nonumber \\
   & & \hspace*{3.7cm}
       \,-\, \#\{j \colon\: \sigma(j) < \sigma(i) \leq i <j \}
       \qquad\qquad \\[2mm]
   & = &
   \#\{ j \colon\: \sigma(i) < \sigma(j) < i < j \}
       \,+\, {\rm I}[\sigma(i) \neq i]
       \qquad\qquad \\[2mm]
   & = &
   \lcross(i,\sigma)
       \,+\, {\rm I}[\sigma(i) \neq i]
   \;.
 \label{eq.labeldiff.case.b}
\end{subeqnarray}
See Figure~\ref{fig.xi.cross}(b,d).
Note that the identity \reff{eq.labeldiff.case.b} holds also when $h_i = 0$,
i.e.\ when the step $s_i$ is a rise from height $h_{i-1} = -1$:
in this case $i$ is a record-antirecord fixed point
and we have $f_i = \xi_i = 0$,
so that both sides of \reff{eq.labeldiff.case.b} are zero.

\bigskip

(c,d) If $s_i$ is a fall, then $\sinv(i)$ is odd
      and $\sinv(i) \le i$.
    We again consider separately the cases of $i$ odd and $i$ even.

\medskip

(c) If $i$ is odd, then $h_i$ is odd, and $i \le \sigma(i)$.
So either $i$ is a fixed point (hence $\sinv(i) = i = \sigma(i)$)
or else $\sinv(i) < i < \sigma(i)$.
Then
\begin{subeqnarray}
   \left\lceil \dfrac{h_i}{2} \right\rceil - \xi_i
   & = &
   \dfrac{h_i + 1}{2} \,-\, \xi_i
       \\[2mm]
   & = &
   f_i \,-\, \xi_i
       \\[2mm]
   & = &
   \#\{ j \le i \colon\: \sigma(j) > i \}   
       \,-\, \#\{j \colon\: j < i \le \sigma(i) < \sigma(j) \}
       \qquad\qquad \\[2mm]
   & = &
   \#\{ j < i \colon\: \sigma(j) > i \}   \,+\, {\rm I}[\sigma(i) \neq i]
       \nonumber \\
   & & \hspace*{3.6cm}
       \,-\, \#\{j \colon\: j < i < \sigma(i) < \sigma(j) \}
       \qquad\qquad \\[2mm]
   & = &
   \#\{ j \colon\: j < i < \sigma(j) < \sigma(i) \}
       \,+\, {\rm I}[\sigma(i) \neq i]
       \\[2mm]
   & = &
   \ucross(i,\sigma) \,+\, {\rm I}[\sigma(i) \neq i]
   \;.
 \label{eq.labeldiff.case.c}
\end{subeqnarray}
See again Figure~\ref{fig.xi.cross}(a,c).
%% {\bf How should we interpret the added term?????}
Note that the identity \reff{eq.labeldiff.case.c} holds also when $h_i = -1$,
i.e.\ when the step $s_i$ is a fall from height $h_{i-1} = 0$;
in this case $i$ is a record-antirecord fixed point
and we have $f_i = \xi_i = 0$,
so that both sides of \reff{eq.labeldiff.case.c} are zero.

\medskip

(d) If $i$ is even, then $h_i$ is even;
      moreover, since $\sinv(i)$ is odd,
      $i$ cannot be a fixed point,
      and we have the strict inequalities $\sinv(i) < i$
      and $\sigma(i) < i$.
      Then
\begin{subeqnarray}
   \left\lceil \dfrac{h_i}{2} \right\rceil - \xi_i
   & = &
   \dfrac{h_i}{2} \,-\, \xi_i
       \\[2mm]
   & = &
   f_i \,-\, \xi_i
       \\[2mm]
   & = &
   \#\{ j \le i \colon\: \sinv(j) > i \}   
       \,-\,  \#\{j \colon\: \sigma(j) < \sigma(i) \leq i <j\}
       \qquad\qquad \\[2mm]
   & = &
   \#\{ j < i \colon\: \sinv(j) > i \}   
       \,-\,  \#\{j \colon\: \sigma(j) < \sigma(i) < i <j\}
       \qquad\qquad \\[2mm]
   & = &
   \#\{ j \colon\: \sigma(j) < i < j \}   
       \,-\,  \#\{j \colon\: \sigma(j) < \sigma(i) < i <j\}
       \qquad\qquad \\[2mm]
   & = &
   \#\{ j \colon\: \sigma(i) < \sigma(j) < i < j \}
       \qquad\qquad \\[2mm]
   & = &
   \lcross(i,\sigma)
   \;.
\end{subeqnarray}
See again Figure~\ref{fig.xi.cross}(b,d).
\qed

\smallskip

We now consider the four possible combinations of $s_i$ (rise or fall)
and parity of $h_i$ (odd or even),
and determine in each case the cycle classification of the index~$i$.
By definition $s_i$ tells us the parity of $\sinv(i)$,
while the parity of $h_i$ equals the parity of $i$.
So these two pieces of information tell us what was recorded in
\reff{eq.parities.1}/\reff{eq.parities.2}:
\begin{itemize}
   \item $\sinv(i)$ even and $i$ odd $\iff$ $i$ is a cycle valley
   \item $\sinv(i)$ even and $i$ even $\iff$ $i$ is either
            a cycle double fall or an even fixed point
   \item $\sinv(i)$ odd and $i$ odd $\iff$ $i$ is either
            a cycle double rise or an odd fixed point
   \item $\sinv(i)$ odd and $i$ even $\iff$ $i$ is a cycle peak
\end{itemize}
So we need only disambiguate the fixed points from the
cycle double falls/rises in the middle two cases;
we will see that in these cases $i$ is a fixed point
if and only if $\xi_i$ takes its maximum allowed value.
More precisely:

\begin{lemma}[Cycle classification]
   \label{lemma.xi.maxval}
\nopagebreak
\quad\hfill
\vspace*{-1mm}
\nopagebreak
\begin{itemize}
   \item[(a)]  If $s_i$ a rise and $h_i$ is odd (hence $h_{i-1}$ is even),
      then $i$ is a cycle valley.
      %% (that is, $\sinv(i)$ an anti-excedance and $i$ is an excedance).
   \item[(b)]  If $s_i$ a rise and $h_i$ is even (hence $h_{i-1}$ is odd),
      then $i$ is an even fixed point
      in case $\xi_i = \left\lceil \dfrac{h_i - 1}{2} \right\rceil$
      ($= h_i/2 = f_i)$;
      otherwise it is a cycle double fall.
      %% (that is, $\sinv(i)$ and $i$ are anti-excedances).
   \item[(c)]  If $s_i$ is a fall and $h_i$ is odd (hence $h_{i-1}$ is even),
      then $i$ is an odd fixed point
      in case $\xi_i = \left\lceil \dfrac{h_i}{2} \right\rceil$
      ($= (h_i +1)/2 = f_i$);
      otherwise it is a cycle double rise.
      %% (that is, $\sinv(i)$ and $i$ are excedances).
   \item[(d)]  If $s_i$ is a fall and $h_i$ is even (hence $h_{i-1}$ is odd),
      then $i$ is a cycle peak.
      %% (that is, $\sinv(i)$ an excedance and $i$ is an anti-excedance).
\end{itemize}
\end{lemma}

\proof
(a,d) follow immediately from \reff{eq.parities.1}/\reff{eq.parities.2}.

(b) If $s_i$ a rise and $h_i$ is even (hence $i$ is even),
then by \reff{eq.lemma.xii.crossing.b} we have
\be
   \left\lceil \dfrac{h_i - 1}{2} \right\rceil - \xi_i
   \;=\;
   \lcross(i,\sigma) \,+\, {\rm I}[\sigma(i) \neq i]
   \;,
\ee
so $i$ is a fixed point if and only if
$\xi_i = \left\lceil \dfrac{h_i - 1}{2} \right\rceil$.
Otherwise $i$ is a cycle double fall.

(c) If $s_i$ is a fall and $h_i$ is odd (hence $i$ is odd),
by \reff{eq.lemma.xii.crossing.c} we have
\be
   \left\lceil \dfrac{h_i}{2} \right\rceil - \xi_i
   \;=\;
   \ucross(i,\sigma) \,+\, {\rm I}[\sigma(i) \neq i]
   \;,
\ee
so $i$ is a fixed point if and only if
$\xi_i = \left\lceil \dfrac{h_i}{2} \right\rceil$.
Otherwise $i$ is a cycle double rise.
\qed

{\bf Remark.}
We already saw in \reff{eq.xii.fp} that if $i$ is a fixed point,
then $\xi_i$ takes the value specified in (b) or (c),
which is also the maximum allowed value according to \reff{eq.xi.ineqs}.
Now we see the converse.
\myendremark

At the other extreme, it is easy to see that the index $i$
is a record or antirecord
if and only if $\xi_i$ takes its minimum allowed value (namely, zero):

\begin{lemma}[Record statistics]
   \label{lemma.record}
\nopagebreak
\quad\hfill
\vspace*{-1mm}
\begin{itemize}
   \item[(a)]  If $i$ is odd,
      then the index $i$ is a record if and only if $\xi_i = 0$.
   \item[(b)]  If $i$ is odd,
      then the index $i$ is an antirecord if and only if
      $h_i = -1$ and $\xi_i = 0$,
      in which case $i$ is a record-antirecord fixed point.
   \item[(c)]  If $i$ is even,
      then the index $i$ is an antirecord if and only if $\xi_i = 0$.
   \item[(d)]  If $i$ is even,
      then the index $i$ is a record if and only if
      $h_i = 0$ and $\xi_i = 0$,
      in which case $i$ is a record-antirecord fixed point.
\end{itemize}
\end{lemma}

\proof
(a,c) This is an immediate consequence of the definition \reff{def2.xi}.

(b) Every antirecord is a weak anti-excedance,
so an odd index in a D-permutation can be an antirecord
only if it is a fixed point,
in which case it is a record-antirecord fixed point.
This happens if and only if $h_i = -1$ and $\xi_i = 0$.

(d) Similar to (b), using the fact that every record is a weak excedance.
\qed

% \begin{quote}
% \small
% {\bf OLD STUFF:}
% \begin{itemize}
% 
% \item When $i = 2k$ and $s_{2k}$ is a rise, we get that $\sinv(2k)$ is even.
% 
% \begin{subeqnarray}
% \left\lceil \dfrac{h_{2k}-1}{2} \right\rceil - \xi_{2k} &=& \dfrac{h_{2k}}{2} - \xi_{2k} \\&=& |\{j\leq 2k \colon\: \sinv(j)>2k\}| \\
% &&- |\{\sinv(2k)<2l\colon\: 2k>\sigma(2l)\}|\\
% &=& |\{j \colon\: j\leq 2k < \sinv(j) \leq \sinv(2k)\}|.
% \end{subeqnarray}
% 
% 
% \item When $i = 2k-1$ and $s_{2k-1}$ is a fall, we get that $\sinv(2k-1)$ is odd.
% 
% \begin{subeqnarray}
% \left\lceil \dfrac{h_{2k-1}}{2} \right\rceil - \xi_{2k-1}
%  &=& \dfrac{h_{2k-1}+1}{2} - \xi_{2k-1} \\&=& |\{j\leq 2k-1 \colon\: \sigma(j)>2k-1\}| \\
% &&- |\{2l-1<\sinv(2k-1)\colon\: \sigma(2l-1)>2k-1\}|\\
% &=& |\{j\colon\: \sinv(2k-1)\leq j \leq 2k-1 < \sigma(j)\}|.
% \end{subeqnarray}
% 
% 
% \item When $i = 2k$ and $s_{2k}$ is a fall, we get that $\sinv(2k)$ is odd.
% 
% \begin{subeqnarray}
% \left\lceil \dfrac{h_{2k}}{2} \right\rceil - \xi_{2k}
%  &=& \dfrac{h_{2k}}{2} - \xi_{2k} \\&=& |\{j\leq 2k \colon\: \sigma(j)>2k\}| \\
% &&- |\{2l-1<\sinv(2k)\colon\: \sigma(2l-1)>2k\}|\\
% &=& |\{j\colon\: \sinv(2k)\leq j \leq 2k < \sigma(j)\}|.
% \end{subeqnarray}
% 
% \end{itemize}
% \end{quote}
% 
% 
% Since the quantities \textbf{!!!refer!!!} are manifestly nonnegative,
% the inequalities \textbf{!!!refer!!!} are satisfied.
% 
% \textbf{Interpret in terms of bipartite matchings!!!}

%\bigskip

\subsection{Step 3: Proof of bijection}   \label{subsec.proofs.bijection}

We prove that the map $\sigma \mapsto (\omega, \xi)$ is a bijection
by explicitly describing the inverse map.
That is, we let $\omega$ be any almost-Dyck path of length $2n$
and let $\xi$ be any set of labels satisfying the inequalities
\reff{eq.xi.ineqs}, and we show how to reconstruct the unique D-permutation
$\sigma$ that gives rise to $(\omega, \xi)$ by the foregoing construction.

First, some preliminaries: Given a D-permutation $\sigma\in \dperm_{2n}$
we can define four subsets of $[2n]$:
\begin{subeqnarray}
F & = & \{2,4,\ldots,2n\} \; = \;  \hbox{even positions} \\
F' & = & \{i \colon \sinv(i) \text{ is even} \} \; = \; \{\sigma(2),\sigma(4),\ldots, \sigma(2n)\}\\
G & = & \{1,3,\ldots, 2n-1\} \; = \; \hbox{odd positions}  \\
G' & = & \{i \colon \sinv(i) \text{ is odd} \} \; = \; \{\sigma(1),\sigma(3),\ldots, \sigma(2n-1)\}
\end{subeqnarray}
Note that $F'$ (resp. $G'$) are the positions of the rises (resp. falls) in the
almost-Dyck path $\omega$.

Let us observe that
\begin{subeqnarray}
F\cap F'    &=& \hbox{cycle double falls and even fixed points}\\
G\cap G'    &=& \hbox{cycle double rises and odd fixed points}\\
F\cap G'    &=& \hbox{cycle peaks}\\
F'\cap G    &=& \hbox{cycle valleys} \\
F \cap G    &=& \emptyset \\
F' \cap G'  &=& \emptyset
\end{subeqnarray}

Let us also recall the notion of an {\em inversion table}\/:
Let $S$ be a totally ordered set of cardinality $k$,
and let $\bm{x} = (x_1,\ldots,x_k)$ be a permutation of $S$
(i.e., a word in which each element of $S$ occurs exactly once);
then the (left-to-right) inversion table corresponding to $\bm{x}$
is the sequence $\bm{p} = (p_1,\ldots,p_k)$ of nonnegative integers
defined by $p_\alpha = \#\{\beta < \alpha \colon\: x_\beta > x_\alpha \}$.
Note that $0 \le p_\alpha \le \alpha-1$ for all $\alpha \in [k]$,
so there are exactly $k!$ possible inversion tables.
Given the inversion table $\bm{p}$,
we can reconstruct the sequence $\bm{x}$
by working from right to left, as follows:
There are $p_k$ elements of $S$ larger than $x_k$,
so $x_k$ must be the $(p_k+1)$th largest element of $S$.
Then there are $p_{k-1}$ elements of $S \setminus \{x_k\}$
larger than $x_{k-1}$,
so $x_{k-1}$ must be the $(p_{k-1}+1)$th largest element
of $S \setminus \{x_k\}$.
And so forth. [Analogously, the right-to-left inversion table corresponding to 
$\bm{x}$ is the sequence $\bm{q} = (q_1,\ldots,q_k)$ of nonnegative integers
defined by $q_\alpha = \#\{\beta > \alpha \colon\: x_\beta < x_\alpha \}$.]

With these preliminaries out of the way, we can now describe the map 
$(\omega,\xi)\mapsto \sigma$. Given the almost-Dyck path $\omega$, 
we can immediately reconstruct the sets $F, F',G, G'$. 
We now use the labels $\xi$ to reconstruct the maps
$\sigma \restrict F \colon\: F \to F'$ and
$\sigma \restrict G \colon\: G \to G'$ as follows: 
The even subword $\sigma(2)\sigma(4)\cdots\sigma(2n)$ is a listing of $F'$
whose right-to-left inversion table is given by $q_\alpha = \xi_{2\alpha}$;
this is the content of~(\ref{def2.xi.bis}a).
Similarly, the odd subword $\sigma(1)\sigma(3)\cdots\sigma(2n-1)$ 
is a listing of $G'$
whose left-to-right inversion table is given by $p_\alpha = \xi_{2\alpha-1}$;
this is the content of~ (\ref{def2.xi.bis}b).
See Figure~\ref{fig.bijection.first} for an example.

\begin{figure}[p]
%% Bishal's FZ_v1.txt, edited

\vspace*{-4mm}

\begin{center}
\[
\begin{array}{r}
F = \{2i\;|\;1 \le i \le 7\}  \\[3mm] 
F' = \{\sigma(2i)\; |\; 1\le i \le 7\}   \\[3mm] 
\text{Right-to-left inversion table:} \; \xi_{2i}   \\  
\end{array}
\;
\begin{pmatrix}
\; 2 & 4 & 6 & 8 & 10 & 12 & 14 \; \\[3mm]
\; 1 & 2 & 4 & 6 & 3 & 12 & 13 \; \\[3mm]
\; 0 & 0 & 1 & 1 & 0 & 0 & 0 \; 
\end{pmatrix} 
\]

\[
\begin{array}{r}
G = \{2i-1\;|\;1 \le i \le 7\}  \\[3mm] 
G' = \{\sigma(2i-1)\; |\; 1\le i \le 7\}   \\[3mm] 
\text{Left-to-right inversion table:} \; \xi_{2i-1}   \\  
\end{array}
\;
\begin{pmatrix}
\; 1 & 3 & 5 & 7 & 9 & 11 & 13 \; \\[3mm]
\; 7 & 9 & 5 & 8 & 10 & 11 & 14 \; \\[3mm]
\; 0 & 0 & 2 & 1 & 0 & 0 & 0 \; 
\end{pmatrix} 
\]

\vspace*{1cm}

\begin{tabular}{c| c| c|c}
$2i$ & $\xi_{2i}$ & $F'\backslash\{\sigma(2),\ldots,\sigma(2i-2)\}$ & $\sigma(2i)$\\
\hline
 $2$ & $0$ & $\{\underline{1}, 2,3,4,6,12,13\}$ & $1$ \\[2mm]
$4$ & $0$ & $\{\underline{2},3,4,6,12,13\}$ & $2$ \\[2mm]
$6$ & $1$ & $\{3, \underline{4},6,12,13\}$ & $4$ \\[2mm]
$8$ & $1$ & $\{3, \underline{6},12,13\}$ & $6$ \\[2mm]
$10$ & $0$ & $\{\underline{3}, 12,13\}$ & $3$ \\[2mm]
$12$ & $0$ & $\{\underline{12},13\}$ & $12$ \\[2mm]
$14$ & $0$ & $\{ \underline{13}\}$ & $13$ \\[2mm]
\end{tabular}

\vspace{6mm}

\begin{tabular}{c| c| c|c}
$2i-1$ & $\xi_{2i-1}$ & $G'\backslash\{\sigma(13),\ldots,\sigma(2i+1)\}$ & $\sigma(2i-1)$\\
\hline
 $13$ & $0$ & $\{5,7,8,9,10,11,\underline{14}\}$ & $14$ \\[2mm]
$11$ & $0$ & $\{5,7,8,9,10,\underline{11}\}$ & $11$ \\[2mm]
$9$ & $0$ & $\{5,7,8,9,\underline{10}\}$ & $10$ \\[2mm]
$7$ & $1$ & $\{5,7,\underline{8},9\}$ & $8$ \\[3mm]
$5$ & $2$ & $\{\underline{5}, 7, 9\}$ & $5$ \\[3mm]
$3$ & $0$ & $\{ 7, \underline{9}\}$ & $9$ \\[3mm]
$1$ & $0$ & $\{ \underline{7}\}$ & $7$ \\[3mm]
\end{tabular}
\end{center}

\caption{
   Reconstruction of the permutation
   $\sigma = 7\, 1\, 9\, 2\, 5\, 4\, 8\, 6\, 10\, 3\, 11\, 12\, 14\, 13
     \;=\; (1,7,8,6,4,2)\,(3,9,10)\,(5)\,(11)\,(12)\,(13,14)$
   from its almost-Dyck path $\omega$ and labels $\xi$.
   The value $\sigma(2i)$ is chosen so that it has $\xi_{2i}$ entries
   to its left in the remaining subset of $F'$ in increasing order;
   the value $\sigma(2i-1)$ is chosen so that it has $\xi_{2i-1}$ entries
   to its right in the remaining subset of $G'$ in increasing order.
}
   \label{fig.bijection.first}
\end{figure}

The only thing that remains to be shown is that the $\sigma$ thus constructed
is indeed a D-permutation. For this, we need to show that the following
inequalities hold:
\begin{subeqnarray}
	\sigma(2k) & \leq & 2k  \\[1mm]
	\sigma(2k-1) & \geq & 2k-1
\end{subeqnarray}

Let us do a double counting of the number of rises occuring after the step $s_{2k}$. As $2n-2k$ is the total number of steps after $2k$, the number of rises after $s_{2k}$ is $ (n-k) - h_{2k}/2$.
Thus, we have that
\begin{subeqnarray}
	(n-k) \,-\, \frac{h_{2k}}{2}
	& = & \#\{i>2k \colon\, \sinv(i) \hbox{ is even}\}   \\
	& \geq & \#\{i > 2k \colon\, \sinv(i)> 2k \text{ and }
           \sinv(i) \hbox{ is even}\}  \\[2mm]
	& = &\#\{2l > 2k \colon\, \sigma(2l) > 2k \}   \\[2mm]
	& = & (n-k) \,-\, \#\{2l > 2k \colon\, \sigma(2l) \leq 2k \}
   \label{eq.even.arg1}
\end{subeqnarray}
and hence,
\be
   \frac{h_{2k}}{2} \;\:\leq\;\: \#\{2l > 2k \colon\, \sigma(2l) \leq 2k \}  \;.
\label{eq.even.arg1a}
\ee
If $k$ is such that $\sigma(2k) > 2k$, then \reff{eq.even.arg1a} becomes strict
[because $i=\sigma(2k)$ contributes to (\ref{eq.even.arg1}a)
 but not to (\ref{eq.even.arg1}b)].
Furthermore, if $\sigma(2k) > 2k$, we also have
\be
   \#\{2l > 2k \colon\, \sigma(2l) \leq 2k \}
   \;\:\le\;\:
   \#\{2l > 2k \colon\, \sigma(2l) < \sigma(2k) \}
   \;\:=\;\:
   \xi_{2k}
  \label{eq.even.arg2}
\ee
by \reff{def2.xi.bis}.
In this situation,
equations~\reff{eq.even.arg1a} and~\reff{eq.even.arg2} together give
\be
  \xi_{2k}  \;>\;  \frac{h_{2k}}{2} 
            \;=\;  \left\lceil \frac{h_{2k}}{2} \right\rceil 
            \;=\; \left\lceil \frac{h_{2k} -1}{2} \right\rceil
  \;,
\ee
where the equalities occur because $h_{2k}$ is even.
This contradicts~(\ref{eq.xi.ineqs}a,b),
and proves that $\sigma(2k) \le 2k$.

%\textbf{Argument for odd case 2k-1!!!}
%
%Notice that there are $2k-2$ steps before the $(2k-1){\rm th}$ step, and hence, the number of falls before $s_{2k-1}$ is $(k-1)-h_{2k-2}/2$. Thus, we have that
%\begin{subeqnarray}
%	(k-1)-\frac{h_{2k-2}}{2}
%	& = & \#\{i<2k-1: \sinv(i) \hbox{ is odd}\}\\
%	& \geq & \#\{i<2k-1: \sinv(i)< 2k-1 \text{ and } \sinv(i) \hbox{ is odd }\}\\
%	& = & \#\{2l-1 < 2k-1: \sigma(2l-1)<2k-1 \}\\
%	& = & (k-1) - \#\{2l-1 < 2k-1: \sigma(2l-1) \geq 2k-1 \}.
%	\label{eq.odd.arg1}
%\end{subeqnarray}
%If we assume $\sigma(2k-1)<2k-1$, then the inequality in~(\ref{eq.odd.arg1}b)~becomes strict. Furthermore, we also have
%\begin{subeqnarray}
%	& & \#\{2l-1 < 2k-1: \sigma(2l-1) \geq 2k-1 \}\\
%	& \leq & \#\{2l-1 < 2k-1: \sigma(2l-1) > \sigma(2k-1) \} = \xi_{2k-1}.
%	\label{eq.odd.arg2}
%\end{subeqnarray}
%Thus, equations~\reff{eq.odd.arg1} and~\reff{eq.odd.arg2} together give
%\be
%\xi_{2k-1} > \frac{h_{2k-2}}{2} = \left\lceil \frac{h_{2k-2}}{2} \right\rceil = \left\lceil \frac{h_{2k-2} -1}{2} \right\rceil.
%\ee
%The equalities here occur due to $h_{2k-2}$ being even. This immediately gives us a contradiction.

The proof that $\sigma(2k-1)\geq 2k-1$ uses a similar double-counting argument
for the number of falls before $s_{2k-1}$.
\qed

\subsection{Step 4: Computation of the weights}
   \label{subsec.proof.computation}

We can now compute the weights associated to the
0-Schr\"oder path $\omegahat$
in Theorem~\ref{thm.flajolet_master_labeled_Schroder},
which we recall are
$a_{h,\xi}$ for a rise starting at height $h$ with label $\xi$,
$b_{h,\xi}$ for a fall starting at height $h$ with label $\xi$, and
$c_{h,\xi}$ for a long level step at height $h$ with label $\xi$.
(Of course, in the present case we have long level steps only at height~0.)
We do this by putting together the information
collected in Lemmas~\ref{lemma.xii.nesting}--\ref{lemma.record}:
\begin{itemize}
   \item[(a)] Rise from height $h_{i-1} = 2k$ to height $h_i = 2k+1$
      (hence $i$ odd):
\begin{itemize}
   \item By Lemma~\ref{lemma.ineqs.labels},
      the label satisfies $0 \le \xi_i \le k$.
   \item By Lemma~\ref{lemma.xi.maxval}(a), this is a cycle valley.
   \item By Lemma~\ref{lemma.xii.nesting}, $\unest(i,\sigma) = \xi_i$.
   \item By Lemma~\ref{lemma.xii.crossing}(a), $\ucross(i,\sigma) = k - \xi_i$.
\end{itemize}
Therefore, from \reff{def.Qn.firstmaster}, the weight for this step is
\be
   a_{2k,\xi}  \;=\;  \sfa_{k-\xi,\xi}
   \;.
\ee
   \item[(b)] Rise from height $h_{i-1} = 2k-1$ to height $h_i = 2k$
      (hence $i$ even):
\begin{itemize}
   \item By Lemma~\ref{lemma.ineqs.labels},
      the label satisfies $0 \le \xi_i \le k$.
   \item By Lemma~\ref{lemma.xi.maxval}(b), this is a cycle double fall
      if $0 \le \xi_i < k$, and a fixed point if $\xi_i = k$.
   \item By Lemma~\ref{lemma.xii.nesting},
      $\xi_i = \displaystyle
       \begin{cases}
           \lnest(i,\sigma)  &  \textrm{if $i$ is not a fixed point} \\
           \psnest(i,\sigma)  &  \textrm{if $i$ is a fixed point}
       \end{cases}
      $
   \item By Lemma~\ref{lemma.xii.crossing}(b), $\lcross(i,\sigma) = k-1-\xi_i$
      when $i$ is not a fixed point.
\end{itemize}
Therefore, from \reff{def.Qn.firstmaster}, the weight for this step is
\be
   a_{2k-1,\xi}  \;=\;
   \begin{cases}
      \sfc_{k-1-\xi,\xi}   & \textrm{if $0 \le \xi < k$}  \\[1mm]
      \sfe_k               & \textrm{if $\xi = k$}
   \end{cases}
\ee
   \item[(c)] Fall from height $h_{i-1} = 2k$ to height $h_i = 2k-1$
      (hence $i$ odd):
\begin{itemize}
   \item By Lemma~\ref{lemma.ineqs.labels},
      the label satisfies $0 \le \xi_i \le k$.
   \item By Lemma~\ref{lemma.xi.maxval}(c), this is a cycle double rise
      if $0 \le \xi_i < k$, and a fixed point if $\xi_i = k$.
   \item By Lemma~\ref{lemma.xii.nesting},
      $\xi_i = \displaystyle
       \begin{cases}
           \unest(i,\sigma)  &  \textrm{if $i$ is not a fixed point} \\
           \psnest(i,\sigma)  &  \textrm{if $i$ is a fixed point}
       \end{cases}
      $
   \item By Lemma~\ref{lemma.xii.crossing}(c), $\ucross(i,\sigma) = k-1-\xi_i$
      when $i$ is not a fixed point.
\end{itemize}
Therefore, from \reff{def.Qn.firstmaster}, the weight for this step is
\be
   b_{2k,\xi}  \;=\;
   \begin{cases}
      \sfd_{k-1-\xi,\xi}   & \textrm{if $0 \le \xi < k$}  \\[1mm]
      \sff_k               & \textrm{if $\xi = k$}
   \end{cases}
\ee
   \item[(d)] Fall from height $h_{i-1} = 2k+1$ to height $h_i = 2k$
      (hence $i$ even):
\begin{itemize}
   \item By Lemma~\ref{lemma.ineqs.labels},
      the label satisfies $0 \le \xi_i \le k$.
   \item By Lemma~\ref{lemma.xi.maxval}(d), this is a cycle peak.
   \item By Lemma~\ref{lemma.xii.nesting}, $\lnest(i,\sigma) = \xi_i$.
   \item By Lemma~\ref{lemma.xii.crossing}(d), $\lcross(i,\sigma) = k - \xi_i$.
\end{itemize}
Therefore, from \reff{def.Qn.firstmaster}, the weight for this step is
\be
   b_{2k+1,\xi}  \;=\;  \sfb_{k-\xi,\xi}
   \;.
\ee
   \item[(e)] Long level step at height~$0$: \\[2mm]
This corresponds in the almost-Dyck path $\omega$
to a fall from height~$0$ to height~$-1$,
followed by a rise from height~$-1$ to height~$0$.
Applying case~(d) with $k=0$ and $\xi=0$,
followed by case~(b) with $k=0$ and $\xi=0$,
we obtain a weight
\be
   c_{0,0}  \;=\;  \sfe_0 \sff_0
   \;.
\ee
\end{itemize}

Putting this all together
in Theorem~\ref{thm.flajolet_master_labeled_Schroder},
we obtain a T-fraction with
\begin{eqnarray}
   \alpha_{2k-1}
   & = &
   (\hbox{rise from $2k-2$ to $2k-1$})
   \,\times\,
   (\hbox{fall from $2k-1$ to $2k-2$})
       \nonumber \\[2mm]
   & = &
   \biggl( \sum_{\xi=0}^{k-1} \sfa_{k-1-\xi,\xi} \biggr)
   \biggl( \sum_{\xi=0}^{k-1} \sfb_{k-1-\xi,\xi} \biggr)
       \label{alpha.2k-1.master}  \\[4mm]
   \alpha_{2k}
   & = &
   (\hbox{rise from $2k-1$ to $2k$})
   \,\times\,
   (\hbox{fall from $2k$ to $2k-1$})
       \nonumber \\[2mm]
   & = &
   \biggl( \sfe_k \,+\, \sum_{\xi=0}^{k-1} \sfc_{k-1-\xi,\xi} \biggr)
   \biggl( \sff_k \,+\, \sum_{\xi=0}^{k-1} \sfd_{k-1-\xi,\xi} \biggr)
       \label{alpha.2k.master}  \\[2mm]
   \delta_1  & = &  \sfe_0 \sff_0
       \label{delta1.master} \\[1mm]
   \delta_n  & = &   0    \qquad\hbox{for $n \ge 2$}
       \label{deltan.master}
\end{eqnarray}
This completes the proof of Theorem~\ref{thm.Tfrac.first.master}.
\qed

We can now deduce Theorem~\ref{thm.Tfrac.first.pq} as a corollary:

\proofof{Theorem~\ref{thm.Tfrac.first.pq}}
Comparing \reff{def.Pn.pq} with \reff{def.Qn.firstmaster}
and using Lemma~\ref{lemma.record},
we see that the needed weights in \reff{def.Qn.firstmaster} are
\begin{eqnarray}
   \sfa_{k-1-\xi,\xi}
   & = &
   p_{+1}^{k-1-\xi} q_{+1}^\xi  \,\times\,
   \begin{cases}
      y_1  & \textrm{if $\xi = 0$}   \\
      v_1  & \textrm{if $1 \le \xi \le k-1$}
   \end{cases}
        \label{eq.proof.weights.sfa}  \\[2mm]
   \sfb_{k-1-\xi,\xi}
   & = &
   p_{-1}^{k-1-\xi} q_{-1}^\xi  \,\times\,
   \begin{cases}
      x_1  & \textrm{if $\xi = 0$}   \\
      u_1  & \textrm{if $1 \le \xi \le k-1$}
   \end{cases}
        \\[2mm]
   \sfc_{k-1-\xi,\xi}
   & = &
   p_{-2}^{k-1-\xi} q_{-2}^\xi  \,\times\,
   \begin{cases}
      x_2  & \textrm{if $\xi = 0$}   \\
      u_2  & \textrm{if $1 \le \xi \le k-1$}
   \end{cases}
        \\[2mm]
   \sfd_{k-1-\xi,\xi}
   & = &
   p_{+2}^{k-1-\xi} q_{+2}^\xi  \,\times\,
   \begin{cases}
      y_2  & \textrm{if $\xi = 0$}   \\
      v_2  & \textrm{if $1 \le \xi \le k-1$}
   \end{cases}
        \\[2mm]
   \sfe_k
   & = &
   \begin{cases}
      \ze  & \textrm{if $k = 0$}   \\[1mm]
      \se^k \we  & \textrm{if $k \ge 1$}
   \end{cases}
        \\[2mm]
   \sff_k
   & = &
   \begin{cases}
      \zo  & \textrm{if $k = 0$}   \\[1mm]
      \so^k \wo  & \textrm{if $k \ge 1$}
   \end{cases}
        \label{eq.proof.weights.sff}
\end{eqnarray}
Inserting these into \reff{alpha.2k-1.master}--\reff{deltan.master}
yields the continued-fraction coefficients \reff{eq.thm.Tfrac.first.weights.pq}.
\qed\hspace*{-5mm}

\proofof{Theorem~\ref{thm.Tfrac.first}}
Specialize Theorem~\ref{thm.Tfrac.first.pq} to
$p_{-1} = p_{-2} = p_{+1} = p_{+2} = q_{-1} = q_{-2} = q_{+1} = q_{+2} =
 \se = \so = 1$.
\qed

%\bigskip

\subsection[An alternative label $\widehat{\xi}_i$: Proof of Theorems
   \ref{thm.Tfrac.first.master.variant} and \ref{thm.Tfrac.first.pq.variant}]%
{An alternative label $\bm{\widehat{\xi}_i}$: Proof of Theorems
   \ref{thm.Tfrac.first.master.variant} and \ref{thm.Tfrac.first.pq.variant}}
         \label{subsec.proofs.alternative}

As mentioned earlier,
Randrianarivony \cite[Section~6]{Randrianarivony_97}
employed a very similar construction
in the special case where $\sigma$ is a D-o-semiderangement.
Our definition of the almost-Dyck path $\omega$
is essentially the same as his Dyck path,
modified slightly to allow for fixed points of both parities.
However, he used a very different definition of the labels,
namely \cite[eq.~(6.2)]{Randrianarivony_97}
\begin{subeqnarray}
   \widehat{\xi}_i
   & = &
   \begin{cases}
       \#\{2l > 2k \colon\: \sigma(2l) < \sigma(2k) \}
           & \textrm{if $i = \sigma(2k)$}
              \\[2mm]
     \#\{2l-1< 2k-1 \colon\: \sigma(2l-1) > \sigma(2k-1) \}
           & \textrm{if $i= \sigma(2k-1) $}
   \end{cases}
   \qquad
 \slabel{def.xi.a}
       \\[3mm]
   & = &
   \begin{cases}
       \#\{j \colon\:  j < i \le \sinv(i) < \sinv(j) \}
           & \textrm{if $\sinv(i)$ is even}
              \\[2mm]
       \#\{j \colon\:  \sinv(j) < \sinv(i) \le i < j \}
           & \textrm{if $\sinv(i)$ is odd}
   \end{cases}
 \slabel{def.xi.b}
 \label{def.xi}
\end{subeqnarray}
% This can equivalently be written as
% \be
%    \xi_i
%    \;=\;
%    \begin{cases}
%        \#\{j > 2k \colon\: \sigma(j) < \sigma(2k) \}
%            & \textrm{if $i = \sigma(2k)$}
%               \\[2mm]
%      \#\{j < 2k-1 \colon\: \sigma(j) > \sigma(2k-1) \}
%            & \textrm{if $i= \sigma(2k-1) $}
%    \end{cases}
%  \label{def.xi.2}
% \ee
% since in the first case, if $j > 2k$ is odd,
% then $\sigma(j) \ge j > 2k \ge \sigma(2k)$ by definition of D-permutation,
% and in the second case, if $j < 2k-1$ is even,
% then $\sigma(j) \le j < 2k-1 \le \sigma(2k-1)$ for the same reason.
%
We would now like to show how the alternative label $\widehat{\xi}_i$
can be use to prove the variant forms of our T-fractions
(Theorems~\ref{thm.Tfrac.first.master.variant}
 and \ref{thm.Tfrac.first.pq.variant}).

% We began by trying to use these labels to prove our T-fractions,
% but we found that they lead to somewhat different
% --- and in our opinion less natural ---
% interpretations in terms of crossing and nesting statistics.
% More precisely, for these labels
% we require the modified crossing and nesting statistics
% defined in \reff{def.ucrossnestjk.prime}:  we have

Just as the labels $\xi_i$ are related to the
index-refined crossing and nesting statistics \reff{def.ucrossnestjk},
so the alternative labels $\widehat{\xi}_i$ are related to the variant
index-refined crossing and nesting statistics \reff{def.ucrossnestjk.prime}.
For the nesting statistics this is immediate from the
definition \reff{def.xi.b}:

\begin{lemma}[Nesting statistics for the alternative labels]
   \label{lemma.xii.nesting.alternative}
We have
\be
   \widehat{\xi}_i
   \;=\;
   \begin{cases}
        \lnest'(i,\sigma)
           & \hbox{\rm if $\sinv(i) $ is even (i.e., $s_i$ is a rise)
                       and $\neq i$}
              \\[1mm]
        \unest'(i,\sigma)
           & \hbox{\rm if $\sinv(i) $ is odd (i.e., $s_i$ is a fall)
                       and $\neq i$}
              \\[1mm]
        \psnest(i,\sigma)
           & \hbox{\rm if $\sinv(i) = i$ (i.e., $i$ is a fixed point)}
   \end{cases}
 \label{eq.xii.nestings.version1}
\ee
\end{lemma}

%% together with the identities for crossing statistics
%% stated in Lemma~\ref{lemma.xihat.crossing} below.

% \begin{itemize}
%    \item[(a)]  If $s_i$ a rise, then
% \begin{subeqnarray}
%    \left\lceil \dfrac{h_i - 1}{2} \right\rceil - \widehat{\xi}_i
%    & = &
%    \lcross'(i,\sigma)
%        \,+\, {\rm I}[\hbox{\rm $i$ is even and } \sigma(i) \neq i]
%      \\[-1mm]
%    & = &
%    \lcross'(i,\sigma)
%        \,+\, {\rm I}[\hbox{\rm $i$ is a cycle double fall}]
%    \;.
%  \label{eq.lemma.xii.crossing.1.version1}
% \end{subeqnarray}
%    \item[(b)]  If $s_i$ a fall, then
% \begin{subeqnarray}
%    \left\lceil \dfrac{h_i}{2} \right\rceil - \widehat{\xi}_i
%    & = &
%    \ucross'(i,\sigma)
%        \,+\, {\rm I}[\hbox{\rm $i$ is odd and } \sigma(i) \neq i]
%      \\[-1mm]
%    & = &
%    \ucross'(i,\sigma)
%        \,+\, {\rm I}[\hbox{\rm $i$ is a cycle double rise}]
%    \;.
%  \label{eq.lemma.xii.crossing.2.version1}
% \end{subeqnarray}
% \end{itemize}
%% where $\lcross^\star$ and $\ucross^\star$ were defined in
%% \reff{def.ulcross.star}.

%% We therefore sought new labels that would be better adapted to
%% the statistics we have used, and we came up with \reff{def2.xi}.
%% Only later did we realize that these are exactly the
%% variant Foata--Zeilberger labels whenever $i$ is not a fixed point.

\bigskip

{\bf Remarks.}
1.  By comparing \reff{def.xi.b} with \reff{def2.xi}, we see that
\be
   \widehat{\xi}_i  \;=\;  \xi_{\sinv(i)}
   \;.
 \label{eq.compare.v1.v2}
\ee
This explains \reff{eq.xii.nestings.version1}:
combine \reff{eq.xii.nestings} with \reff{eq.nestprime}.

2.  The distinction between $\xi_i$ and $\widehat{\xi}_i$
is also related to the distinction between two different notions
of ``inversion table''.
In the approach used here and in
\cite[Section~6.1, Step~3]{Sokal-Zeng_masterpoly},
the label $\xi_i$ is the number of inversions
associated to the {\em position}\/ $i$.
By contrast, in \cite[p.~88]{Randrianarivony_97} and \cite[p.~51]{Foata_90},
the label $\widehat{\xi}_i$ is the number of inversions
associated to the {\em value}\/ $i$,
i.e.~to the position $\sinv(i)$.
Thus, whenever $i$ is not a fixed point,
$\widehat{\xi}_i$ corresponds to
the original Foata--Zeilberger \cite{Foata_90} labels,
while $\xi_i$ corresponds to
the modified labels used in \cite[Section~6.1]{Sokal-Zeng_masterpoly}.
\myendremark

% Since these alternative labels $\widehat{\xi}_i$
% will turn out to play a central role in our proof of
% the {\em second}\/ T-fraction
% (Sections~\ref{sec.second} and \ref{sec.proofs.2}),
% we prove now two key lemmas:
% that the labels $\widehat{\xi}_i$ satisfy the same inequalities
% as the labels $\xi_i$;
% and the identities analogous to Lemma~\ref{lemma.xii.crossing}
% that relate the labels $\widehat{\xi}_i$ to crossing statistics.
%% and that the identities
%% \reff{eq.lemma.xii.crossing.1.version1}/%
%% \reff{eq.lemma.xii.crossing.2.version1} hold.

The alternative labels $\widehat{\xi}_i$ satisfy the same inequalities
as the original labels $\xi_i$:

\begin{lemma}[Inequalities satisfied by the alternative labels]
   \label{lemma.ineqs.labels.xihat}
We have
\begin{subeqnarray}
   & 0 \;\le\; \widehat{\xi}_i  \;\le\; \left\lceil \dfrac{h_i-1}{2} \right\rceil \;=\;
                      \left\lceil\dfrac{h_{i-1}}{2} \right\rceil
       & \textrm{if $\sinv(i) $ is even (i.e., $s_i$ is a rise)}
    \qquad
 \slabel{eq.xihat.ineqs.a}
      \\[3mm]
   & 0 \;\le\; \widehat{\xi}_i  \;\le\; \left\lceil \dfrac{h_i}{2} \right\rceil  \;=\;
                      \left\lceil\dfrac{h_{i-1}-1}{2} \right\rceil
       & \textrm{if $\sinv(i) $ is odd (i.e., $s_i$ is a fall)}
 \slabel{eq.xihat.ineqs.b}
 \label{eq.xihat.ineqs}
\end{subeqnarray}
\end{lemma}

Lemma~\ref{lemma.ineqs.labels.xihat} will be an immediate consequence
of the following identities:

\begin{lemma}[Crossing statistics for the alternative labels]
   \label{lemma.xihat.crossing}
\nopagebreak
\quad\hfill
\vspace*{-1mm}
\begin{itemize}
   \item[(a)]  If $s_i$ a rise (i.e.~$\sinv(i)$ is even), then
\begin{subeqnarray}
   \left\lceil \dfrac{h_i - 1}{2} \right\rceil - \widehat{\xi}_i
   & = &
   \lcross'(i,\sigma)
       \,+\, {\rm I}[\hbox{\rm $i$ is even and } \sigma(i) \neq i]
     \\[-1mm]
   & = &
   \lcross'(i,\sigma)
       \,+\, {\rm I}[\hbox{\rm $i$ is a cycle double fall}]
   \;.
 \label{eq.lemma.xihat.crossing.1}
\end{subeqnarray}
   \item[(b)]  If $s_i$ a fall (i.e.~$\sinv(i)$ is odd), then
\begin{subeqnarray}
   \left\lceil \dfrac{h_i}{2} \right\rceil - \widehat{\xi}_i
   & = &
   \ucross'(i,\sigma)
       \,+\, {\rm I}[\hbox{\rm $i$ is odd and } \sigma(i) \neq i]
     \\[-1mm]
   & = &
   \ucross'(i,\sigma)
       \,+\, {\rm I}[\hbox{\rm $i$ is a cycle double rise}]
   \;.
 \label{eq.lemma.xihat.crossing.2}
\end{subeqnarray}
\end{itemize}
\end{lemma}

\begin{figure}[t]
\centering
\begin{picture}(30,15)(160, 10)
\setlength{\unitlength}{4mm}
\linethickness{.5mm}
\put(2,0){\line(1,0){27}}
\put(5,0){\circle*{1,3}}\put(5,0){\makebox(0,3)[c]{\small $j$}}
\put(12,0){\circle*{1,3}}\put(12,0){\makebox(0,3)[c]{\small $i=\sigma(2k)$}}
\put(19,0){\circle*{1,3}}\put(19,0){\makebox(0,3)[c]{\small $\sinv(j)$}}
\put(26,0){\circle*{1,3}}\put(26,0){\makebox(0,3)[c]{\small $\sinv(i)=2k$}}
\red{\qbezier(5,0)(12,-6)(19,0)}
\blue{\qbezier(11.5,0)(18.5,-6)(25.5,0)}
\put(15,-2){\makebox(0,-6)[c]{\small When $\sigma^{-1}(i)$ is even}}
% %%%%%%%%%%%%%%%%%%%%
\end{picture}
\\[4.3cm]
\begin{picture}(30,15)(640, 10)
\setlength{\unitlength}{4mm}
\linethickness{.5mm}
\put(44,0){\line(1,0){27}}
\put(47,0){\circle*{1,3}}\put(47,0){\makebox(0,-3)[c]{\small $\sinv(i)=2k\!-\!1$}}
\put(54,0){\circle*{1,3}}\put(54,0){\makebox(0,-3)[c]{\small $j$}}
\put(61,0){\circle*{1,3}}\put(61,0){\makebox(0,-3)[c]{\small $i=\sigma(2k\!-\!1)$}}
\put(68,0){\circle*{1,3}}\put(68,0){\makebox(0,-3)[c]{\small $\sigma(j)$}}
\red{\qbezier(47,0)(54,6)(61,0)}
\blue{\qbezier(53.5,0)(60.5,6)(67.5,0)}
\put(57,-0.5){\makebox(0,-6)[c]{\small When $\sigma^{-1}(i)$ is odd}}
\end{picture}
\\[2.5cm]
\caption{
   Crossings involved in the inequalities for the label $\widehat{\xi}_i$.
 \label{fig.xihat.cross}
}
\vspace*{2cm}
\end{figure}
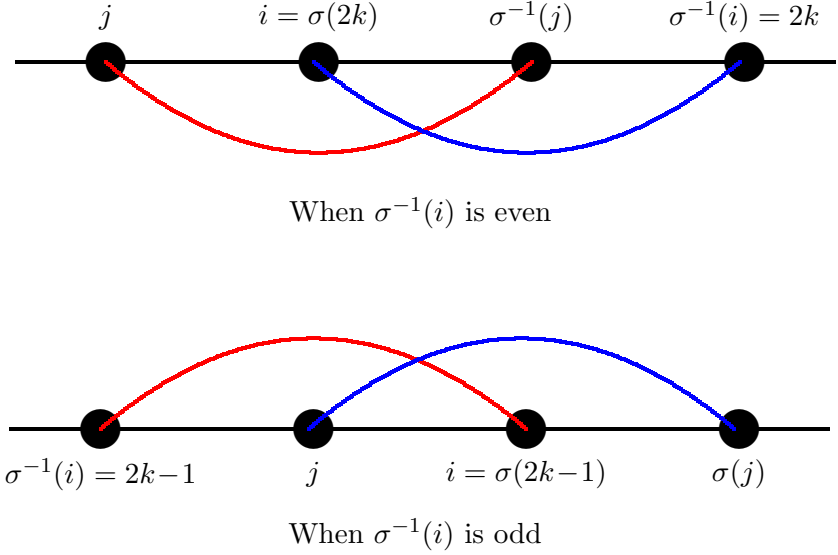

\proof
(a) If $s_i$ is a rise, then $\sinv(i)$ is even, and $i \le \sinv(i)$.
    We now consider separately the cases of $i$ odd and $i$ even.
\begin{itemize}
   \item[(i)] If $i$ is odd, then $h_i$ is odd;
      moreover, we have the strict inequality $i < \sinv(i)$
      (hence $i$ is not a fixed point).  Then
\begin{subeqnarray}
   \hspace*{-6mm}
   \left\lceil \dfrac{h_i - 1}{2} \right\rceil - \widehat{\xi}_i
   & = &
   \dfrac{h_i - 1}{2} \,-\, \widehat{\xi}_i
       \\
   & = &
   f_i \,-\, 1 \,-\, \widehat{\xi}_i
       \\[2mm]
%    & = &
%    \#\{ j \le i \colon\: \sinv(j) > i \}   \,-\,  1
%        \,-\, \#\{ 2l > \sinv(i) \colon\: \sigma(2l) < i \}
%        \qquad\qquad \\[2mm]
   & = &
   \#\{ j \le i \colon\: \sinv(j) > i \}   \,-\,  1
          \nonumber \\
   & & \hspace*{4cm}
       \,-\, \#\{ j < i \colon\: \sinv(j) > \sinv(i) \}
       \qquad\qquad \\[2mm]
   & = &
   \#\{ j < i \colon\: \sinv(j) > i \}
       \,-\, \#\{ j < i \colon\: \sinv(j) > \sinv(i) \}
       \qquad\qquad \\[2mm]
   & = &
   \#\{ j \colon\:  j < i < \sinv(j) \le \sinv(i) \}
       \\[2mm]
   & = &
   \#\{ j \colon\:  j < i < \sinv(j) < \sinv(i) \}
         \nonumber \\
   &  &  \qquad\qquad
       \hbox{since $j \neq i$ implies $\sinv(j) \neq \sinv(i)$}
       \\[2mm]
   & = &
   \lcross'(i,\sigma)
   \;.
 \label{eq.labeldiff.case.ai}
\end{subeqnarray}
See Figure~\ref{fig.xihat.cross}(a).
   \item[(ii)] If $i$ is even, then $h_i$ is even, and
\begin{subeqnarray}
   \hspace*{-6mm}
   \left\lceil \dfrac{h_i - 1}{2} \right\rceil - \widehat{\xi}_i
   & = &
   \dfrac{h_i}{2} \,-\, \widehat{\xi}_i
       \\
   & = &
   f_i \,-\, \widehat{\xi}_i
       \\[2mm]
%    & = &
%    \#\{ j \le i \colon\: \sinv(j) > i \}   
%        \,-\, \#\{ 2l > \sinv(i) \colon\: \sigma(2l) < i \}
%        \qquad\qquad \\[2mm]
   & = &
   \#\{ j \le i \colon\: \sinv(j) > i \}  
       \,-\, \#\{ j < i \colon\: \sinv(j) > \sinv(i) \}
       \qquad\qquad \\[2mm]
   & = &
   \#\{ j < i \colon\: \sinv(j) > i \}
       \,+\, {\rm I}[\sigma(i) \neq i]
          \nonumber \\
   &  & \hspace*{4cm}
       \,-\, \#\{ j < i \colon\: \sinv(j) > \sinv(i) \}
       \\[2mm]
   & = &
   \#\{ j \colon\: j < i < \sinv(j) \le \sinv(i) \}
           \,+\, {\rm I}[\sigma(i) \neq i]
       \qquad\qquad \\[2mm]
   & = &
   \#\{ j \colon\: j < i < \sinv(j) < \sinv(i) \}
           \,+\, {\rm I}[\sigma(i) \neq i]
        \nonumber \\
   & & \qquad\qquad
       \hbox{since $j \neq i$ implies $\sinv(j) \neq \sinv(i)$}
       \qquad \\[2mm]
   & = &
   \lcross'(i,\sigma) \,+\, {\rm I}[\sigma(i) \neq i]
   \;.
 \label{eq.labeldiff.case.aii}
\end{subeqnarray}
See again Figure~\ref{fig.xihat.cross}(a).
%% {\bf How should we interpret the added term?????}
Note that the identity \reff{eq.labeldiff.case.aii} holds also when $h_i = 0$,
i.e.\ when the step $s_i$ is a rise from height $h_{i-1} = -1$;
in this case $i$ is a record-antirecord fixed point
and we have $\widehat{\xi}_i = 0$,
so that both sides of \reff{eq.labeldiff.case.aii} are zero.
\end{itemize}
Combining \reff{eq.labeldiff.case.ai} and \reff{eq.labeldiff.case.aii}
yields \reff{eq.lemma.xihat.crossing.1}.

\medskip

(b) If $s_i$ is a fall, then $\sinv(i)$ is odd, and $\sinv(i) \le i$.
    We again consider separately the cases of $i$ odd and $i$ even.
\begin{itemize}
   \item[(i)] If $i$ is odd, then $h_i$ is odd, and
\begin{subeqnarray}
   \left\lceil \dfrac{h_i}{2} \right\rceil - \widehat{\xi}_i
   & = &
   \dfrac{h_i + 1}{2} \,-\, \widehat{\xi}_i
       \\
   & = &
   f_i \,-\, \widehat{\xi}_i
       \\[2mm]
   & = &
   \#\{ j \le i \colon\: \sigma(j) > i \}   
       \,-\, \#\{ j \colon\: \sinv(j) < \sinv(i) \le i < j \}
       \\[2mm]
   & = &
   \#\{ j \le i \colon\: \sigma(j) > i \}   
       \,-\, \#\{ j < \sinv(i) \colon\: \sigma(j) > i \}
       \qquad\qquad \\[2mm]
   & = &
   \#\{ j \colon\: \sinv(i) \le j \le i < \sigma(j) \}
       \\[2mm]
   & = &
   \#\{ j \colon\: \sinv(i) < j \le i < \sigma(j) \}
       \nonumber \\
   & & \qquad
       \hbox{since $i \neq \sigma(j)$ implies $j \neq \sinv(i)$}
       \\[2mm]
   & = &
   \#\{ j \colon\: \sinv(i) < j < i < \sigma(j) \}
      \,+\, {\rm I}[\sinv(i) < i < \sigma(i)]
       \qquad\qquad \\[2mm]
   & = &
   \ucross'(i,\sigma) \,+\, {\rm I}[\sigma(i) \neq i]
   \;.
 \label{eq.labeldiff.case.bi}
\end{subeqnarray}
See Figure~\ref{fig.xihat.cross}(b).
%% {\bf How should we interpret the added term?????}
Note that the identity \reff{eq.labeldiff.case.bi} holds also when $h_i = -1$,
i.e.\ when the step $s_i$ is a fall from height $h_{i-1} = 0$;
in this case $i$ is a record-antirecord fixed point
and we have $\widehat{\xi}_i = 0$,
so that both sides of \reff{eq.labeldiff.case.bi} are zero.
   \item[(ii)] If $i$ is even, then $h_i$ is even;
      moreover, we have the strict inequality $\sinv(i) < i$
      (hence $i$ is not a fixed point).  Then
\begin{subeqnarray}
   \left\lceil \dfrac{h_i}{2} \right\rceil - \widehat{\xi}_i
   & = &
   \dfrac{h_i}{2} \,-\, \widehat{\xi}_i
       \\
   & = &
   f_i \,-\, \widehat{\xi}_i
       \\[2mm]
   & = &
   \#\{ j \le i \colon\: \sigma(j) > i \}   
       \,-\, \#\{ j < \sinv(i) \colon\: \sigma(j) > i \}
       \qquad\qquad \\[2mm]
   & = &
   \#\{ j \colon\: \sinv(i) \le j \le i < \sigma(j) \}
       \\[2mm]
   & = &
   \#\{ j \colon\: \sinv(i) < j \le i < \sigma(j) \}
       \nonumber \\
   & & \qquad
       \hbox{since $i \neq \sigma(j)$ implies $j \neq \sinv(i) = \sinv(i)$}
       \\[2mm]
   & = &
   \#\{ j \colon\: \sinv(i) < j < i < \sigma(j) \}
       \quad\hbox{since $i \ge \sigma(i)$}
       \qquad\qquad \\[2mm]
   & = &
   \ucross'(i,\sigma)
   \;.
 \label{eq.labeldiff.case.bii}
\end{subeqnarray}
See again Figure~\ref{fig.xihat.cross}(b).
\end{itemize}
Combining \reff{eq.labeldiff.case.bi} and \reff{eq.labeldiff.case.bii}
yields \reff{eq.lemma.xihat.crossing.2}.
\qed

\begin{lemma}[Cycle classification for the alternative labels]
   \label{lemma.xi.maxval.alternative}
\nopagebreak
\quad\hfill
\vspace*{-1mm}
\nopagebreak
\begin{itemize}
   \item[(a)]  If $s_i$ a rise and $h_i$ is odd (hence $h_{i-1}$ is even),
      then $i$ is a cycle valley.
      %% (that is, $\sinv(i)$ an anti-excedance and $i$ is an excedance).
   \item[(b)]  If $s_i$ a rise and $h_i$ is even (hence $h_{i-1}$ is odd),
      then $i$ is an even fixed point
      in case $\widehat{\xi}_i = \left\lceil \dfrac{h_i - 1}{2} \right\rceil$
      ($= h_i/2 = f_i)$;
      otherwise it is a cycle double fall.
      %% (that is, $\sinv(i)$ and $i$ are anti-excedances).
   \item[(c)]  If $s_i$ is a fall and $h_i$ is odd (hence $h_{i-1}$ is even),
      then $i$ is an odd fixed point
      in case $\widehat{\xi}_i = \left\lceil \dfrac{h_i}{2} \right\rceil$
      ($= (h_i +1)/2 = f_i$);
      otherwise it is a cycle double rise.
      %% (that is, $\sinv(i)$ and $i$ are excedances).
   \item[(d)]  If $s_i$ is a fall and $h_i$ is even (hence $h_{i-1}$ is odd),
      then $i$ is a cycle peak.
      %% (that is, $\sinv(i)$ an excedance and $i$ is an anti-excedance).
\end{itemize}
\end{lemma}

\proof
Completely analogous to the proof of Lemma~\ref{lemma.xi.maxval}:
just use Lemma~\ref{lemma.xihat.crossing} in place of
Lemma~\ref{lemma.xii.crossing}.
\qed

\begin{lemma}[Record statistics for the alternative labels]
   \label{lemma.record.alternative}
\nopagebreak
\quad\hfill
\vspace*{-1mm}
\begin{itemize}
   \item[(a)]  If $\sinv(i)$ is odd,
      then the index $\sinv(i)$ is a record if and only if
      $\widehat{\xi}_i = 0$.
   \item[(b)]  If $\sinv(i)$ is even,
      then the index $\sinv(i)$ is an antirecord if and only if
      $\widehat{\xi}_i = 0$.
\end{itemize}
\end{lemma}

\proof
This is an immediate consequence of the definition \reff{def.xi.b}.
\qed

\bigskip

{\bf Proof of bijection.}
The proof is similar to that presented in
Section~\ref{subsec.proofs.bijection},
but using a value-based rather than position-based notion
of inversion table.
Recall that if $S = \{s_1 < s_2 < \ldots < s_k \}$
is a totally ordered set of cardinality $k$,
and $\bm{x} = (x_1,\ldots,x_k)$ is a permutation of $S$,
then the (left-to-right) (position-based) inversion table
corresponding to $\bm{x}$
is the sequence $\bm{p} = (p_1,\ldots,p_k)$ of nonnegative integers
defined by $p_\alpha = \#\{\beta < \alpha \colon\: x_\beta > x_\alpha \}$.
We now define the (left-to-right) {\em value-based}\/ inversion table
$\bm{p}'$ by $p'_{x_i} = p_i$;
note that $\bm{p}'$ is a map from $S$ to $\{0,\ldots,k-1\}$,
such that $p'_{x_i}$ is the number of entries to the left of $x_i$
(in the word $\bm{x}$) that are larger than $x_i$.
In particular, $0 \le p'_{s_i} \le k-i$.
Given the value-based inversion table $\bm{p}'$,
we can reconstruct the sequence $\bm{x}$
by working from largest to smallest value,
as follows \cite[section~5.1.1]{Knuth_98}:
We start from an empty sequence, and insert $s_k$.
Then we insert $s_{k-1}$ so that the resulting word has $p'_{s_{k-1}}$
entries to its left.
Next we insert $s_{k-2}$ so that the resulting word has $p'_{s_{k-2}}$
entries to its left, and so on.
[The right-to-left value-based inversion table $\bm{q}'$
 is defined analogously, and the reconstruction proceeds from
 smallest to largest.]

We now recall the definitions
\begin{subeqnarray}
F & = & \{2,4,\ldots,2n\} \; = \;  \hbox{even positions} \\
F' & = & \{i \colon \sinv(i) \text{ is even} \} \; = \; \{\sigma(2),\sigma(4),\ldots, \sigma(2n)\}\\
G & = & \{1,3,\ldots, 2n-1\} \; = \; \hbox{odd positions}  \\
G' & = & \{i \colon \sinv(i) \text{ is odd} \} \; = \; \{\sigma(1),\sigma(3),\ldots, \sigma(2n-1)\}
\end{subeqnarray}
Note that $F'$ (resp. $G'$) are the positions of the rises (resp. falls)
in the almost-Dyck path $\omega$.

We can now describe the map
$(\omega,\widehat{\xi})\mapsto \sigma$.
Given the almost-Dyck path $\omega$,
we can immediately reconstruct the sets $F, F',G, G'$.
We now use the labels $\widehat{\xi}$ to reconstruct the maps
$\sigma \restrict F \colon\: F \to F'$ and
$\sigma \restrict G \colon\: G \to G'$ as follows:
The even subword $\sigma(2)\sigma(4)\cdots\sigma(2n)$ is a listing of $F'$
whose right-to-left value-based inversion table
is given by $q'_i = \widehat{\xi}_i$ for all $i \in F'$;
this is the content of \reff{def.xi.a}.
Similarly, the odd subword $\sigma(1)\sigma(3)\cdots\sigma(2n-1)$
is a listing of $G'$
whose left-to-right value-based inversion table
is given by $p'_i = \widehat{\xi}_i$ for all $i \in G'$;
this again is the content of \reff{def.xi.a}.
See Figure~\ref{fig.bijection.first_v2} for an example.

\begin{figure}[p]
%% Bishal's FZ_v2.txt, edited

\vspace*{-4mm}

\begin{center}
\[
\begin{array}{r}
F = \{2i\;|\;1 \le i \le 7\}  \\[3mm] 
F' = \{\sigma(2i)\; |\; 1\le i \le 7\}   \\[3mm] 
\text{Right-to-left inversion table:} \; \widehat{\xi}_{\sigma(2i)}   \\  
\end{array}
\;
\begin{pmatrix}
\; 2 & 4 & 6 & 8 & 10 & 12 & 14 \; \\[3mm]
\; 1 & 2 & 4 & 6 & 3 & 12 & 13 \; \\[3mm]
\; 0 & 0 & 1 & 1 & 0 & 0 & 0 \; 
\end{pmatrix} 
\]

\[
\begin{array}{r}
G = \{2i-1\;|\;1 \le i \le 7\}  \\[3mm] 
G' = \{\sigma(2i-1)\; |\; 1\le i \le 7\}   \\[1em] 
\text{Left-to-right inversion table:} \; \widehat{\xi}_{\sigma(2i-1)}   \\  
\end{array}
\;
\begin{pmatrix}
\; 1 & 3 & 5 & 7 & 9 & 11 & 13 \; \\[3mm]
\; 7 & 9 & 5 & 8 & 10 & 11 & 14 \; \\[3mm]
\; 0 & 0 & 2 & 1 & 0 & 0 & 0 \; 
\end{pmatrix} 
\]

\vspace*{1cm}

\begin{tabular}{c| c| l}
$f_i$ & $\widehat{\xi}_{f_i}$ & Partial subword of $F'$ \\
\hline
$1$ & $0$ & $\underline{1}$ \\[2mm]
$2$ & $0$ & $1 \; \underline{2}$  \\[2mm]
$3$ & $0$ & $1 \; 2 \; \underline{3}$ \\[2mm]
$4$ & $1$ & $1 \; 2 \; \underline{4} \; 3$ \\[2mm]
$6$ & $1$ & $1 \; 2 \; 4\; \underline{6} \; 3$ \\[2mm]
$12$ & $0$ & $1 \; 2 \; 4\; 6 \; 3 \; \underline{12}$ \\[2mm]
$13$ & $0$ & $1 \; 2 \; 4\; 6 \; 3 \; 12 \; \underline{13}$
\end{tabular}

\vspace{5mm}

\begin{tabular}{c| c| l}
$g_i$ & $\widehat{\xi}_{g_i}$ & Partial subword of $G'$ \\
\hline
$14$ & $0$ & $\underline{14}$ \\[2mm]
$11$ & $0$ & $\underline{11} \; 14 $ \\[2mm]
$10$ & $0$ & $\underline{10} \; 11 \; 14 $ \\[2mm]
$9$ & $0$ & $\underline{9} \; 10 \; 11 \; 14 $ \\[2mm]
$8$ & $1$ & $9 \; \underline{8} \; 10 \; 11 \; 14 $ \\[2mm]
$7$ & $0$ & $\underline{7} \; 9 \;8 \; 10 \; 11 \; 14 $ \\[2mm]
$5$ & $2$ & $7\; 9\;  \underline{5} \; 8 \; 10 \; 11 \; 14 $
\end{tabular}
\end{center}

\caption{
   Reconstruction of the permutation
   $\sigma = 7\, 1\, 9\, 2\, 5\, 4\, 8\, 6\, 10\, 3\, 11\, 12\, 14\, 13
     \;=\; (1,7,8,6,4,2)\,(3,9,10)\,(5)\,(11)\,(12)\,(13,14)$
   from its almost-Dyck path $\omega$ and labels $\widehat{\xi}$.
   The value $f_i$ is inserted so that it has
   $\widehat{\xi}_{f_i}$ entries to its right in the partial subword;
   the value $g_i$ is inserted so that it has
   $\widehat{\xi}_{g_i}$ entries to its left in the partial subword.
}
   \label{fig.bijection.first_v2}
\end{figure}

The only thing that remains to be shown is that the $\sigma$
thus constructed is indeed a D-permutation.
For this, we need to show that the following inequalities hold:
\begin{subeqnarray}
\sigma(2k-1) &\geq &  2k-1
   \slabel{eq.ineqodd} \\[1mm]
\sigma(2k)   &\leq & 2k
   \slabel{eq.ineqeven}
\end{subeqnarray}

 We begin with a lemma:
 
 \begin{lemma}
    \label{lemma.bij}
 \nopagebreak
 \quad\hfill
 \vspace*{-1mm}
 \begin{itemize}
    \item[(a)]
 Let $f_1< f_2 < \ldots <f_n$ be the elements of $F'$ in increasing order.
 Then $f_j \leq 2j$ for every $j\in [n]$.
 \item[(b)]
 Let $g_1< g_2 < \ldots < g_n$ be the elements of $G'$ in increasing order.
 Then $g_j \geq 2j-1$ for every $j\in [n]$.
 \end{itemize}
 \end{lemma}
 
\proof
 (a) Notice that the number of rises among the steps $s_1,s_2,\ldots,s_i$
 is $(i+h_i)/2$ (and the number of falls is $(i-h_i)/2$).
 As there are exactly $j$ rises among the steps $s_1,s_2,\ldots,s_{f_j}$,
 we have
 \be
    \frac{f_j + h_{f_j}}{2} \;=\;  j
  \label{eq.evenineqpi0}
 \ee
 and hence
\be
    f_j \;=\; 2j - h_{f_j} \;\leq\; 2j
  \label{eq.oddineqpi}
 \ee
 since $h_{f_j}\geq 0$.
 
\medskip
 
% (b) Notice that the number of falls among the steps $s_1,s_2,\ldots,s_i$
% is $(i-h_i)/2$ (and the number of rises is $(i+h_i)/2$).
% As there are exactly $j$ falls among the steps $s_1,s_2,\ldots,s_{g_j}$,
% we have
% \be
%    \frac{g_j - h_{g_j}}{2} \;=\;  j
%  \label{eq.oddineqpi0}
% \ee
% and hence
% \be
%    g_j \;=\; 2j + h_{g_j} \;\geq\; 2j - 1
%  \label{eq.oddineqpi}
% \ee
% since $h_{g_j}\geq -1$.
 
 The proof of (b) is similar, using the fact that $h_{g_j}\geq -1$.
 \qed
 
% Now consider any index $k \in [n]$.
% Let $j'$ be the index for which $g_{j'} = \sigma(2k-1)$.
% From the definition of left-to-right inversion table,
% we know that there are $\widehat{\xi}_{\sigma(2k-1)}$
% elements to the left of $\sigma(2k-1)$
% in the word $\sigma(1)\sigma(3)\cdots\sigma(2n-1)$
% which are larger than $\sigma(2k-1)$.
% On the other hand, there are $k-1$ elements
% to the left of $\sigma(2k-1)$
% in the word $\sigma(1)\sigma(3)\cdots\sigma(2n-1)$;
% and there are $j' - 1$ elements in $G'$ that are smaller than $g_{j'}$.
% Therefore, there are at least $(k-1) - (j'-1) = k - j'$ elements
% to the left of $\sigma(2k-1)$
% in the word $\sigma(1)\sigma(3)\cdots\sigma(2n-1)$
% that are larger than $g_{j'} = \sigma(2k-1)$.
% Therefore,
% \begin{subeqnarray}
%    \widehat{\xi}_{\sigma(2k-1)}  & \ge &  k - j'
%           \\[1mm]
%    & = &  k \,-\, {g_{j'} - h_{g_{j'}} \over 2}
%                \qquad\hbox{[by \reff{eq.oddineqpi0}]}   \\[1mm]
%    & = &  {h_{\sigma(2k-1)} \over 2}  \,+\, {2k - \sigma(2k-1) \over 2}
%    \;.
% \end{subeqnarray}
% On the other hand, from \reff{eq.xihat.ineqs.b} we know that
% \be
%    \widehat{\xi}_{\sigma(2k-1)}
%    \;\le\;
%    \left\lceil \dfrac{h_{\sigma(2k-1)}}{2} \right\rceil 
%    \;\le\;
%    \dfrac{h_{\sigma(2k-1)} + 1}{2}
%    \;.
% \ee
% Combining these two inequalities, we conclude that $\sigma(2k-1) \ge 2k-1$.

 Now consider any index $k \in [n]$.
 Let $j'$ be the index for which $f_{j'} = \sigma(2k)$.
 From the definition of right-to-left inversion table,
 we know that there are $\widehat{\xi}_{\sigma(2k)}$
 elements to the right of $\sigma(2k)$
 in the word $\sigma(2)\sigma(4)\cdots\sigma(2n)$
 which are smaller than $\sigma(2k)$.
 On the other hand, there are $n-k$ elements
 to the right of $\sigma(2k)$
 in the word $\sigma(2)\sigma(4)\cdots\sigma(2n)$;
 and there are $n-j'$ elements in $F'$ that are larger than $f_{j'}$.
 Therefore, there are at least $(n-k) - (n-j') = j'-k$ elements
 to the right of $\sigma(2k)$
 in the word $\sigma(2)\sigma(4)\cdots\sigma(2n)$
 that are smaller than $f_{j'} = \sigma(2k)$.
 Therefore,
 \begin{subeqnarray}
    \widehat{\xi}_{\sigma(2k)}  & \ge &  j'-k
           \\[1mm]
    & = &   {f_{j'} + h_{f_{j'}} \over 2} \,-\, k
                \qquad\hbox{[by \reff{eq.evenineqpi0}]}   \\[1mm]
    & = &  {h_{\sigma(2k)} \over 2} +  { \sigma(2k) -2k \over 2} 
    \;.
 \end{subeqnarray}
 On the other hand, from \reff{eq.xihat.ineqs.a} we know that
 \be
    \widehat{\xi}_{\sigma(2k)}
    \;\le\;
    \left\lceil \dfrac{h_{\sigma(2k)}-1}{2} \right\rceil 
    \;\le\;
    \dfrac{h_{\sigma(2k)} }{2}
    \;.
 \ee
 Combining these two inequalities, we conclude that $\sigma(2k) \le 2k$.
 
 The proof that $\sigma(2k-1) \ge 2k-1$ is similar,
 using \reff{eq.xihat.ineqs.b}.
\qed

\medskip

{\bf Remark.}  This proof of bijection is very close in spirit to that
of Randrianarivony \cite[pp.~89--90]{Randrianarivony_97}.
\myendremark

\bigskip

\proofof{Theorem~\ref{thm.Tfrac.first.master.variant}}
The computation of the weights is completely analogous to what was done in
Section~\ref{subsec.proof.computation},
but using Lemmas~\ref{lemma.xii.nesting.alternative}--\ref{lemma.record.alternative}
in place of Lemmas~\ref{lemma.xii.nesting}--\ref{lemma.record}.
We leave the details to the reader:
the upshot is that for cycle valleys and cycle peaks,
``u'' and ``l'' are interchanged compared to
Section~\ref{subsec.proof.computation},
and all statistics are primed.
This is exactly what we have in \reff{def.Qn.firstmaster.variant}.
It therefore completes the proof of
Theorem~\ref{thm.Tfrac.first.master.variant}.
\qed

\proofof{Theorem~\ref{thm.Tfrac.first.pq.variant}}
Comparing \reff{def.Pn.pq.variant} with \reff{def.Qn.firstmaster.variant}
and using Lemma~\ref{lemma.record.alternative},
we see that the needed weights in \reff{def.Qn.firstmaster.variant} are
the same as given in \reff{eq.proof.weights.sfa}--\reff{eq.proof.weights.sff}.
Inserting these into Theorem~\ref{thm.Tfrac.first.master.variant}
gives Theorem~\ref{thm.Tfrac.first.pq.variant}.
\qed

\section{Second T-fraction: Proof of Theorems~\ref{thm.Tfrac.second},
    \ref{thm.Tfrac.second.pq} and \ref{thm.Tfrac.second.master}}
       \label{sec.proofs.2}

In this section we prove the second master T-fraction
(Theorem~\ref{thm.Tfrac.second.master})
by a bijection from D-permutations to labeled Schr\"oder paths.
Our construction combines ideas of
Randrianarivony \cite{Randrianarivony_97}
and Biane \cite{Biane_93}
together with some new ingredients.
After proving Theorem~\ref{thm.Tfrac.second.master},
we deduce Theorems~\ref{thm.Tfrac.second} and \ref{thm.Tfrac.second.pq}
by specialization.

Here we need to construct a bijection that will allow us to count the number
of cycles (cyc), which is a global variable. To do this, we employ
a modification of the Biane \cite{Biane_93} bijection,
just as in Section~\ref{sec.proofs.1} we employed a modification of the
Foata--Zeilberger \cite{Foata_90} bijection.
Our bijection maps $\dperm_{2n}$ to the set of
\textbfit{$\bm{(\bfscra,\bfscrb,\bfscrc)}$-labeled {0-Schr\"oder} paths of
length $\bm{2n}$},
where each label $\xi_i$ is a {\em pair}\/ of nonnegative integers
$\xi_i = (\xi'_i, \xi''_i)$ as follows:

\begin{subeqnarray}
   \scra_h & = &  \{0 \}\times \{0\}  \qquad\qquad\qquad\qquad\;\hbox{for $h \ge 0$ and $h$ even}
    \slabel{def.abc.2ndTfrac.a}  \\[1mm]
   \scra_h & = &  \{0\}\times \{0,\ldots, \lceil h/2 \rceil \}  \qquad\quad\;\;\hbox{for $h \ge 0$ and $h$ odd}
    \slabel{def.abc.2ndTfrac.b}  \\[1mm]
   \scrb_h & = &  \{0,\ldots, \lceil (h-1)/2 \rceil \}\times \{0\} \quad\hbox{for $h \ge 1$ and $h$ even}
    \slabel{def.abc.2ndTfrac.c}  \\[1mm]
   \scrb_h & = & \{0,\ldots, \lceil (h-1)/2 \rceil \}^2 \qquad\quad\,\,\hbox{for $h \ge 1$ and $h$ odd} 
    \slabel{def.abc.2ndTfrac.d}  \\[1mm]
   \scrc_0  & = &  \{0\}\times \{0\}   \\[1mm]
   \scrc_h  & = &  \emptyset   \qquad\qquad\qquad\qquad\qquad\quad\;\;\,\hbox{for $h \ge 1$}
 \label{def.abc.2ndTfrac}
\end{subeqnarray}
or equivalently
\begin{subeqnarray}
   \scra_h & = &  \{0 \}\times \{0\}  \qquad\qquad\;\hbox{for $h \ge 0$ and $h=2k$}
    \slabel{def.abc.2ndTfrac.a.bis}  \\[1mm]
   \scra_h & = &  \{0\}\times \{0,\ldots, k \}  \quad\;\:\hbox{for $h \ge 0$ and $h=2k-1$}
    \slabel{def.abc.2ndTfrac.b.bis}  \\[1mm]
   \scrb_h & = &  \{0,\ldots, k \}\times \{0\} \quad\;\:\hbox{for $h \ge 1$ and $h=2k$}
    \slabel{def.abc.2ndTfrac.c.bis}  \\[1mm]
   \scrb_h & = & \{0,\ldots, k \}^2 \qquad\qquad\!\hbox{for $h \ge 1$ and $h=2k+1$} 
    \slabel{def.abc.2ndTfrac.d.bis}  \\[1mm]
   \scrc_0  & = &  \{0\}\times \{0\}   \\[1mm]
   \scrc_h  & = &  \emptyset   \qquad\qquad\qquad\quad\;\;\,\hbox{for $h \ge 1$}
 \label{def.abc.2ndTfrac.bis}
\end{subeqnarray}

Our presentation of this bijection will follow the same steps
as in Section~\ref{sec.proofs.1}.

%\bigskip

\subsection{Step 1: Definition of the almost-Dyck path}

The almost-Dyck path $\omega$ 
associated to a D-permutation $\sigma \in \dperm_{2n}$
is identical to the one employed in Section~\ref{sec.proofs.1}. That is:
\begin{itemize}
   \item If $\sinv(i)$ is even, then $s_i$ is a rise.
   \item If $\sinv(i)$ is odd, then $s_i$ is a fall.
\end{itemize}
The interpretation of the heights $h_i$ is thus exactly as in
Lemma~\ref{lemma.heights}.
We then define the 0-Schr\"oder path $\omegahat = \psi(\omega)$ as before.

%\bigskip

\subsection[Step 2: Definition of the labels $\xi_i = (\xi'_i,\xi''_i)$]{Step 2: Definition of the labels $\bm{\xi_i = (\xi'_i,\xi''_i)}$}

We define the labels $\xi_i = (\xi'_i,\xi''_i)$  as follows:
\begin{eqnarray}
   \xi_i'
   & = &
   \begin{cases}
     0     & \textrm{if $\sinv(i)$ is even} \\[1mm]
     \#\{j \colon\:  \sinv(j) < \sinv(i) \leq i < j \}
           & \textrm{if $\sinv(i)$ is odd} 
   \end{cases}
 \label{def2.xi.2ndTfrac.1}  \\[5mm]
   \xi_i''
   & = &
   \begin{cases}
     0     & \textrm{if $i$ is odd}  \\[1mm]
     \#\{j \colon\: \sigma(j) < \sigma(i)\leq i < j \}
           & \textrm{if $i$ is even}
   \end{cases}
 \label{def2.xi.2ndTfrac.2}
\end{eqnarray}
These labels $\xi'_i,\xi''_i$ needed for the proof of the second T-fraction
are related to the labels $\xi_i,\widehat{\xi}_i$
defined in \reff{def2.xi} and \reff{def.xi}
and employed in Section~\ref{sec.proofs.1}
for the proof of the first T-fraction, as follows:
\begin{subeqnarray}
   & &
   \sinv(i) \hbox{ is odd}  \iff
	i \in \Cpeak \cup \Cdrise \cup \Oddfix \colon\quad
   \nonumber \\
   & & \hspace*{3cm}
	\xi'_i(\textrm{Second}) \:=\: \widehat{\xi}_i(\textrm{First})
                                \:=\:  \xi_{\sigma^{-1}(i)}(\textrm{First})
     \hspace*{1cm}
     \slabel{eq.1st.2nd.AA} \\[3mm]
   & &
   i \hbox{ is even}  \iff
        i \in \Cpeak \cup \Cdfall \cup \Evenfix \colon\quad
   \nonumber \\
   & & \hspace*{3cm}
	\xi''_i(\textrm{Second})  =  \xi_i(\textrm{First})
     \slabel{eq.1st.2nd.BB}
 \label{eq.1st.2nd}
\end{subeqnarray}
Note that \reff{eq.1st.2nd.AA} refers to the cases
(\ref{def.abc.2ndTfrac}c,d) where $\xi'_i$ is nontrivial
[$\sinv(i)$~odd means that $s_i$ is a fall],
while \reff{eq.1st.2nd.BB} refers to the cases
(\ref{def.abc.2ndTfrac}b,d) where $\xi''_i$ is nontrivial
[$i$~even $\iff$ $h_i$~even $\iff$ $h = h_{i-1}$~odd].
Note in particular that cycle peaks belong to both cases
[cf.~(\ref{def.abc.2ndTfrac}d)],
while cycle valleys belong to neither
[cf.~(\ref{def.abc.2ndTfrac}a)].

% {\bf SAVE THIS FOR LATER!!!}
% This interpretation, along with Lemma~\ref{lemma.ineqs.labels}
% and Equations~\reff{eq.lemma.xii.crossing.1.version1}/\reff{eq.lemma.xii.crossing.2.version1},
% shows in particular that the containment relations
% \reff{eq.xi.ineq}/\reff{def.abc.2ndTfrac} are satisfied.
% {\bf CHECK THIS AND STATE AS LEMMA!!!}

It is worth remarking that these labels $\xi_i = (\xi'_i,\xi''_i)$
are the same as those in the Biane bijection
\cite[Section~6.2, Step~2]{Sokal-Zeng_masterpoly}
whenever $i$ is not a fixed point.
Compare also \reff{eq.1st.2nd} to \cite[eq.~(6.24)]{Sokal-Zeng_masterpoly}.

% \begin{quote}
% \small
% {\bf KEEP FOR LATER REMARK???}
% 
% Let us first make an observation before introducing the labels.
% In the cycle classification in Lemma~\ref{lemma.xi.maxval},
% the labels are only used in (b) and (c) 
% to separate even and odd fixed points from cycle double falls 
% and cycle double rises, respectively.
% The rest of the proof depended entirely of the definition of the
% almost-Dyck path.
% \end{quote}

We can give these labels a nice interpretation
by using (as in Section~\ref{sec.proofs.1})
the representation of a permutation $\sigma \in \Sym_N$
by a bipartite digraph $\Gamma = \Gamma(\sigma)$
in which the top row of vertices is labeled $1,\ldots,N$
and the bottom row $1',\ldots,N'$,
and we draw an arrow from $i$ to $j'$ in case $\sigma(i) = j$.
Recall that, for $k \in [N]$, we denote by $\Gamma_k$
the induced subgraph of $\Gamma$ on the vertex set
$\{1,\ldots,k\} \cup \{1',\ldots,k'\}$.
We can consider the ``history''
$\emptyset = \Gamma_0 \subset \Gamma_1 \subset \Gamma_2 \subset \ldots \subset
 \Gamma_N = \Gamma$
as a process of building up the permutation $\sigma$
by successively considering the status of indices $1,2,\ldots,N$.

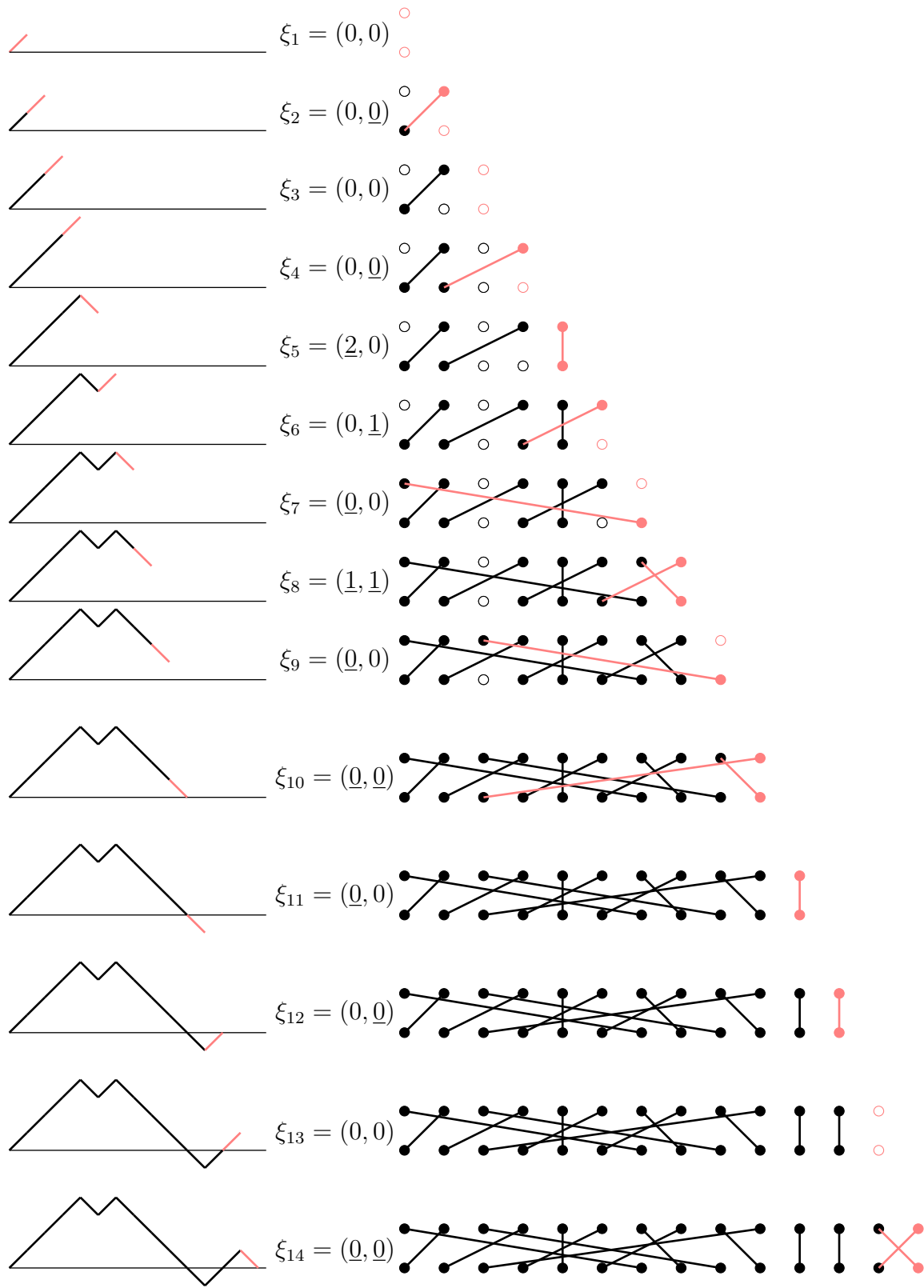
\begin{figure}[p]
\vspace*{-8mm}
\begin{tikzpicture}[scale = 0.325]

%\draw[gray,dotted] (0,0) grid (46,65);

%1
\foreach \i in {1,...,6}
{
\pgfmathsetmacro{\a}{(\i - 1) * 6 +1};
\pgfmathsetmacro{\b}{(\i - 1) * 6 +1.9}; 
\pgfmathsetmacro{\c}{(\i - 1) * 6 +3};

\draw[rounded corners=1, color=black, line width=1] (0,\a)-- (0.9,\b);
\draw[rounded corners=1, color=black, line width=0.5] (0,\a)-- (13,\a);

\draw[black,fill=black] (20,\a) circle (.25);
\draw[black,fill=black] (20,\c) circle (.25);

}
\foreach \i in {14,...,14}
{
\pgfmathsetmacro{\a}{(\i - 6) * 4 + 31};
\pgfmathsetmacro{\b}{(\i - 6) * 4 +0.9+31 };
\pgfmathsetmacro{\c}{(\i - 6) * 4 +31+2};
\draw[rounded corners=1, color=red!50!white, line width=1] (0,\a)-- (0.9,\b);
\draw[rounded corners=1, color=black, line width=0.5] (0,\a)-- (13,\a);

\draw[red!50!white,fill=white] (20,\a) circle (.25);
\draw[red!50!white,fill=white] (20,\c) circle (.25);

\node at (16.5,\b) {$\xi_1 = (0,0)$}; 
}

\foreach \i in {7,...,13}
{
\pgfmathsetmacro{\a}{(\i - 6) * 4 + 31};
\pgfmathsetmacro{\b}{(\i - 6) * 4 +0.9+31 };
\pgfmathsetmacro{\c}{(\i - 6) * 4 +31+2};
\draw[rounded corners=1, color=black, line width=1] (0,\a)-- (0.9,\b);
\draw[rounded corners=1, color=black, line width=0.5] (0,\a)-- (13,\a);

\draw[black,fill=black] (20,\a) circle (.25);
\draw[black,fill=black] (20,\c) circle (.25);
}

%2
\foreach \i in {1,...,6}
{
\pgfmathsetmacro{\a}{(\i - 1) * 6 +1.9};
\pgfmathsetmacro{\b}{(\i - 1) * 6 +2.8 };
\draw[rounded corners=1, color=black, line width=1] (0.9,\a)-- (1.8,\b);

\pgfmathsetmacro{\c}{(\i - 1) * 6 +1};
\pgfmathsetmacro{\d}{(\i - 1) * 6 +3};
\draw[black,fill=black] (22,\c) circle (.25);
\draw[black,fill=black] (22,\d) circle (.25);
\draw[rounded corners=1, color=black, line width=1] (22,\d)-- (20,\c);
}
\foreach \i in {13,...,13}
{
\pgfmathsetmacro{\a}{(\i - 6) * 4 + 31 + 0.9};
\pgfmathsetmacro{\b}{(\i - 6) * 4 +31+ 1.8 };
\pgfmathsetmacro{\c}{(\i - 6) * 4 +31+2};
\draw[rounded corners=1, color=red!50!white, line width=1] (0.9,\a)-- (1.8,\b);

\pgfmathsetmacro{\c}{(\i - 6) * 4 +31};
\pgfmathsetmacro{\d}{(\i - 6) * 4 +31 +2};
\draw[red!50!white,fill=white] (22,\c) circle (.25);
\draw[red!50!white,fill=red!50!white] (22,\d) circle (.25);
\draw[rounded corners=1, color=red!50!white, line width=1] (22,\d)-- (20,\c);

\node at (16.5,\a) {$\xi_2 = (0,\underline{0})$};
}

\draw[black, fill=white] (20,61) circle (.25);

\foreach \i in {7,...,12}
{
\pgfmathsetmacro{\a}{(\i - 6) * 4 + 31 + 0.9};
\pgfmathsetmacro{\b}{(\i - 6) * 4 +31+ 1.8 };
\pgfmathsetmacro{\c}{(\i - 6) * 4 +31+2};
\draw[rounded corners=1, color=black, line width=1] (0.9,\a)-- (1.8,\b);

\pgfmathsetmacro{\c}{(\i - 6) * 4 +31};
\pgfmathsetmacro{\d}{(\i - 6) * 4 +31 +2};
\draw[black,fill=black] (22,\c) circle (.25);
\draw[black,fill=black] (22,\d) circle (.25);
\draw[rounded corners=1, color=black, line width=1] (22,\d)-- (20,\c);
}

%3
\foreach \i in {1,...,6}
{
\pgfmathsetmacro{\a}{(\i - 1) * 6 +2.8};
\pgfmathsetmacro{\b}{(\i - 1) * 6 +3.7 };
\draw[rounded corners=1, color=black, line width=1] (1.8,\a)-- (2.7,\b);

\pgfmathsetmacro{\c}{(\i - 1) * 6 +1};
\pgfmathsetmacro{\d}{(\i - 1) * 6 +3};
\draw[black,fill=black] (24,\c) circle (.25);
\draw[black,fill=black] (24,\d) circle (.25);
}
\foreach \i in {7,...,11}
{
\pgfmathsetmacro{\a}{(\i - 6) * 4 + 31 + 1.8};
\pgfmathsetmacro{\b}{(\i - 6) * 4 +31+ 2.7 };
\draw[rounded corners=1, color=black, line width=1] (1.8,\a)-- (2.7,\b);

\pgfmathsetmacro{\c}{(\i - 6) * 4 +31};
\pgfmathsetmacro{\d}{(\i - 6) * 4 +31 +2};
\draw[black,fill=black] (24,\c) circle (.25);
\draw[black,fill=black] (24,\d) circle (.25);
}

\foreach \i in {12,...,12}
{
\pgfmathsetmacro{\a}{(\i - 6) * 4 + 31 + 1.8};
\pgfmathsetmacro{\b}{(\i - 6) * 4 +31+ 2.7 };
\draw[rounded corners=1, color=red!50!white, line width=1] (1.8,\a)-- (2.7,\b);

\pgfmathsetmacro{\c}{(\i - 6) * 4 +31};
\pgfmathsetmacro{\d}{(\i - 6) * 4 +31 +2};
\draw[red!50!white,fill=white] (24,\c) circle (.25);
\draw[red!50!white,fill=white] (24,\d) circle (.25);

\pgfmathsetmacro{\e}{(\i - 6) * 4 +31+1};
\node at (16.5,\e) {$\xi_3 = (0,0)$}; 
}

\draw[black,fill=white] (20,57) circle (.25);
\draw[black,fill=white] (22,55) circle (.25);

%4
\foreach \i in {1,...,6}
{
\pgfmathsetmacro{\a}{(\i - 1) * 6 +3.7};
\pgfmathsetmacro{\b}{(\i - 1) * 6 +4.6 };
\draw[rounded corners=1, color=black, line width=1] (2.7,\a)-- (3.6,\b);

\pgfmathsetmacro{\c}{(\i - 1) * 6 +1};
\pgfmathsetmacro{\d}{(\i - 1) * 6 +3};
\draw[black,fill=black] (26,\c) circle (.25);
\draw[black,fill=black] (26,\d) circle (.25);
\draw[rounded corners=1, color=black, line width=1] (26,\d)-- (22,\c);
}
\foreach \i in {7,...,10}
{
\pgfmathsetmacro{\a}{(\i - 6) * 4 + 31 + 2.7};
\pgfmathsetmacro{\b}{(\i - 6) * 4 +31+ 3.6 };
\draw[rounded corners=1, color=black, line width=1] (2.7,\a)-- (3.6,\b);

\pgfmathsetmacro{\c}{(\i - 6) * 4 +31};
\pgfmathsetmacro{\d}{(\i - 6) * 4 +31 +2};
\draw[black,fill=black] (26,\c) circle (.25);
\draw[black,fill=black] (26,\d) circle (.25);
\draw[rounded corners=1, color=black, line width=1] (26,\d)-- (22,\c);
}

\foreach \i in {11,...,11}
{
\pgfmathsetmacro{\a}{(\i - 6) * 4 + 31 + 2.7};
\pgfmathsetmacro{\b}{(\i - 6) * 4 +31+ 3.6 };
\draw[rounded corners=1, color=red!50!white, line width=1] (2.7,\a)-- (3.6,\b);

\pgfmathsetmacro{\c}{(\i - 6) * 4 +31};
\pgfmathsetmacro{\d}{(\i - 6) * 4 +31 +2};
\draw[red!50!white,fill=white] (26,\c) circle (.25);
\draw[red!50!white,fill=red!50!white] (26,\d) circle (.25);
\draw[rounded corners=1, color=red!50!white, line width=1] (26,\d)-- (22,\c);

\pgfmathsetmacro{\e}{(\i - 6) * 4 +31+1};
\node at (16.5,\e) {$\xi_4 = (0,\underline{0})$}; 
}
\draw[black,fill=white] (24,53) circle (.25);
\draw[black,fill=white] (24,51) circle (.25);
\draw[black,fill=white] (20,53) circle (.25);

%5
\foreach \i in {1,...,6}
{
\pgfmathsetmacro{\a}{(\i - 1) * 6 +4.6};
\pgfmathsetmacro{\b}{(\i - 1) * 6 +3.7 };
\draw[rounded corners=1, color=black, line width=1] (3.6,\a)-- (4.5,\b);

\pgfmathsetmacro{\c}{(\i - 1) * 6 +1};
\pgfmathsetmacro{\d}{(\i - 1) * 6 +3};
\draw[black,fill=black] (28,\c) circle (.25);
\draw[black,fill=black] (28,\d) circle (.25);
\draw[rounded corners=1, color=black, line width=1] (28,\c)-- (28,\d);
}
\foreach \i in {7,...,9}
{
\pgfmathsetmacro{\a}{(\i - 6) * 4 + 31 + 3.6};
\pgfmathsetmacro{\b}{(\i - 6) * 4 +31+ 2.7 };
\draw[rounded corners=1, color=black, line width=1] (3.6,\a)-- (4.5,\b);

\pgfmathsetmacro{\c}{(\i - 6) * 4 +31};
\pgfmathsetmacro{\d}{(\i - 6) * 4 +31 +2};
\draw[black,fill=black] (28,\c) circle (.25);
\draw[black,fill=black] (28,\d) circle (.25);
\draw[rounded corners=1, color=black, line width=1] (28,\c)-- (28,\d);
}

\foreach \i in {10,...,10}
{
\pgfmathsetmacro{\a}{(\i - 6) * 4 + 31 + 3.6};
\pgfmathsetmacro{\b}{(\i - 6) * 4 +31+ 2.7 };
\draw[rounded corners=1, color=red!50!white, line width=1] (3.6,\a)-- (4.5,\b);

\pgfmathsetmacro{\c}{(\i - 6) * 4 +31};
\pgfmathsetmacro{\d}{(\i - 6) * 4 +31 +2};
\draw[red!50!white,fill=red!50!white] (28,\c) circle (.25);
\draw[red!50!white,fill=red!50!white] (28,\d) circle (.25);
\draw[rounded corners=1, color=red!50!white, line width=1] (28,\c)-- (28,\d);

\pgfmathsetmacro{\e}{(\i - 6) * 4 +31+1};
\node at (16.5,\e) {$\xi_5 = (\underline{2}, 0)$}; 

}

\draw[black,fill=white] (26,47) circle (.25);
\draw[black,fill=white] (24,49) circle (.25);
\draw[black,fill=white] (24,47) circle (.25);
\draw[black,fill=white] (20,49) circle (.25);

%6
\foreach \i in {1,...,6}
{
\pgfmathsetmacro{\a}{(\i - 1) * 6 +3.7};
\pgfmathsetmacro{\b}{(\i - 1) * 6 +4.6 };
\draw[rounded corners=1, color=black, line width=1] (4.5,\a)-- (5.4,\b);

\pgfmathsetmacro{\c}{(\i - 1) * 6 +1};
\pgfmathsetmacro{\d}{(\i - 1) * 6 +3};
\draw[black,fill=black] (30,\c) circle (.25);
\draw[black,fill=black] (30,\d) circle (.25);
\draw[rounded corners=1, color=black, line width=1] (30,\d)-- (26,\c);
}
\foreach \i in {7,...,8}
{
\pgfmathsetmacro{\a}{(\i - 6) * 4 + 31 + 2.7};
\pgfmathsetmacro{\b}{(\i - 6) * 4 +31+ 3.6 };
\draw[rounded corners=1, color=black, line width=1] (4.5,\a)-- (5.4,\b);

\pgfmathsetmacro{\c}{(\i - 6) * 4 +31};
\pgfmathsetmacro{\d}{(\i - 6) * 4 +31 +2};
\draw[black,fill=black] (30,\c) circle (.25);
\draw[black,fill=black] (30,\d) circle (.25);
\draw[rounded corners=1, color=black, line width=1] (30,\d)-- (26,\c);
}

\foreach \i in {9,...,9}
{
\pgfmathsetmacro{\a}{(\i - 6) * 4 + 31 + 2.7};
\pgfmathsetmacro{\b}{(\i - 6) * 4 +31+ 3.6 };
\draw[rounded corners=1, color=red!50!white, line width=1] (4.5,\a)-- (5.4,\b);

\pgfmathsetmacro{\c}{(\i - 6) * 4 +31};
\pgfmathsetmacro{\d}{(\i - 6) * 4 +31 +2};
\draw[red!50!white,fill=white] (30,\c) circle (.25);
\draw[red!50!white,fill=red!50!white] (30,\d) circle (.25);
\draw[rounded corners=1, color=red!50!white, line width=1] (30,\d)-- (26,\c);

\pgfmathsetmacro{\e}{(\i - 6) * 4 +31+1};
\node at (16.5,\e) {$\xi_6 = (0,\underline{1})$};

}

\draw[black,fill=white] (24,45) circle (.25);
\draw[black,fill=white] (24,43) circle (.25);
\draw[black,fill=white] (20,45) circle (.25);

%7
\foreach \i in {1,...,6}
{
\pgfmathsetmacro{\a}{(\i - 1) * 6 +4.6};
\pgfmathsetmacro{\b}{(\i - 1) * 6 +3.7 };
\draw[rounded corners=1, color=black, line width=1] (5.4,\a)-- (6.3,\b);

\pgfmathsetmacro{\c}{(\i - 1) * 6 +1};
\pgfmathsetmacro{\d}{(\i - 1) * 6 +3};
\draw[black,fill=black] (32,\c) circle (.25);
\draw[black,fill=black] (32,\d) circle (.25);
\draw[rounded corners=1, color=black, line width=1] (32,\c)-- (20,\d);
}
\foreach \i in {7,...,7}
{
\pgfmathsetmacro{\a}{(\i - 6) * 4 + 31 + 3.6};
\pgfmathsetmacro{\b}{(\i - 6) * 4 +31+ 2.7 };
\draw[rounded corners=1, color=black, line width=1] (5.4,\a)-- (6.3,\b);

\pgfmathsetmacro{\c}{(\i - 6) * 4 +31};
\pgfmathsetmacro{\d}{(\i - 6) * 4 +31 +2};
\draw[black,fill=black] (32,\c) circle (.25);
\draw[black,fill=black] (32,\d) circle (.25);
\draw[rounded corners=1, color=black, line width=1] (32,\c)-- (20,\d);
}

\foreach \i in {8,...,8}
{
\pgfmathsetmacro{\a}{(\i - 6) * 4 + 31 + 3.6};
\pgfmathsetmacro{\b}{(\i - 6) * 4 +31+ 2.7 };
\draw[rounded corners=1, color=red!50!white, line width=1] (5.4,\a)-- (6.3,\b);

\pgfmathsetmacro{\c}{(\i - 6) * 4 +31};
\pgfmathsetmacro{\d}{(\i - 6) * 4 +31 +2};
\draw[red!50!white,fill=red!50!white] (32,\c) circle (.25);
\draw[red!50!white,fill=white] (32,\d) circle (.25);
\draw[rounded corners=1, color=red!50!white, line width=1] (32,\c)-- (20,\d);

\pgfmathsetmacro{\e}{(\i - 6) * 4 +31+1};
\node at (16.5,\e) {$\xi_7 = (\underline{0},0)$}; 

}

\draw[black,fill=white] (24,41) circle (.25);
\draw[black,fill=white] (24,39) circle (.25);
\draw[black,fill=white] (30,39) circle (.25);

%8
\foreach \i in {1,...,6}
{
\pgfmathsetmacro{\a}{(\i - 1) * 6 +3.7};
\pgfmathsetmacro{\b}{(\i - 1) * 6 +2.8 };
\draw[rounded corners=1, color=black, line width=1] (6.3,\a)-- (7.2,\b);

\pgfmathsetmacro{\c}{(\i - 1) * 6 +1};
\pgfmathsetmacro{\d}{(\i - 1) * 6 +3};
\draw[black,fill=black] (34,\c) circle (.25);
\draw[black,fill=black] (34,\d) circle (.25);
\draw[rounded corners=1, color=black, line width=1] (34,\d)-- (30,\c);
\draw[rounded corners=1, color=black, line width=1] (34,\c)-- (32,\d);
}
\foreach \i in {7,...,7}
{
\pgfmathsetmacro{\a}{(\i - 6) * 4 + 31 + 2.7};
\pgfmathsetmacro{\b}{(\i - 6) * 4 +31 +1.8 };
\draw[rounded corners=1, color=red!50!white, line width=1] (6.3,\a)-- (7.2,\b);

\pgfmathsetmacro{\c}{(\i - 6) * 4 +31};
\pgfmathsetmacro{\d}{(\i - 6) * 4 +31 +2};
\draw[red!50!white,fill=red!50!white] (34,\c) circle (.25);
\draw[red!50!white,fill=red!50!white] (34,\d) circle (.25);
\draw[rounded corners=1, color=red!50!white, line width=1] (34,\d)-- (30,\c);
\draw[rounded corners=1, color=red!50!white, line width=1] (34,\c)-- (32,\d);

\pgfmathsetmacro{\e}{(\i - 6) * 4 +31 +1};
\node at (16.5,\e) {$\xi_8 = (\underline{1},\underline{1})$}; 

}
\draw[black,fill=white] (24,37) circle (.25);
\draw[black,fill=white] (24,35) circle (.25);

%9
\foreach \i in {6,...,6}
{
\pgfmathsetmacro{\a}{(\i - 1) * 6 +2.8};
\pgfmathsetmacro{\b}{(\i - 1) * 6 +1.9 };
\draw[rounded corners=1, color=red!50!white, line width=1] (7.2,\a)-- (8.1,\b);

\pgfmathsetmacro{\c}{(\i - 1) * 6 +1};
\pgfmathsetmacro{\d}{(\i - 1) * 6 +3};
\draw[red!50!white,fill=red!50!white] (36,\c) circle (.25);
\draw[red!50!white,fill=white] (36,\d) circle (.25);
\draw[rounded corners=1, color=red!50!white, line width=1] (36,\c)-- (24,\d);

\pgfmathsetmacro{\e}{(\i - 1) * 6 +2};
\node at (16.5,\e) {$\xi_9 = (\underline{0},0)$}; 

}

\foreach \i in {1,...,5}
{
\pgfmathsetmacro{\a}{(\i - 1) * 6 +2.8};
\pgfmathsetmacro{\b}{(\i - 1) * 6 +1.9 };
\draw[rounded corners=1, color=black, line width=1] (7.2,\a)-- (8.1,\b);

\pgfmathsetmacro{\c}{(\i - 1) * 6 +1};
\pgfmathsetmacro{\d}{(\i - 1) * 6 +3};
\draw[black,fill=black] (36,\c) circle (.25);
\draw[black,fill=black] (36,\d) circle (.25);
\draw[rounded corners=1, color=black, line width=1] (36,\c)-- (24,\d);
}
\draw[black,fill=white] (24,31) circle (.25);

%10
\foreach \i in {5,...,5}
{
\pgfmathsetmacro{\a}{(\i - 1) * 6 +1.9 };
\pgfmathsetmacro{\b}{(\i - 1) * 6 +1};
\draw[rounded corners=1, color=red!50!white, line width=1] (8.1,\a)-- (9,\b);

\pgfmathsetmacro{\c}{(\i - 1) * 6 +1};
\pgfmathsetmacro{\d}{(\i - 1) * 6 +3};
\draw[red!50!white,fill=red!50!white] (38,\c) circle (.25);
\draw[red!50!white,fill=red!50!white] (38,\d) circle (.25);
\draw[rounded corners=1, color=red!50!white, line width=1] (38,\c)-- (36,\d);
\draw[rounded corners=1, color=red!50!white, line width=1] (38,\d)-- (24,\c);

\pgfmathsetmacro{\e}{(\i - 1) * 6 +2};
\node at (16.5,\e) {$\xi_{10} = (\underline{0},\underline{0})$}; 

}

\foreach \i in {1,...,4}
{
\pgfmathsetmacro{\a}{(\i - 1) * 6 +1.9 };
\pgfmathsetmacro{\b}{(\i - 1) * 6 + 1 };
\draw[rounded corners=1, color=black, line width=1] (8.1,\a)-- (9,\b);

\pgfmathsetmacro{\c}{(\i - 1) * 6 +1};
\pgfmathsetmacro{\d}{(\i - 1) * 6 +3};
\draw[black,fill=black] (38,\c) circle (.25);
\draw[black,fill=black] (38,\d) circle (.25);
\draw[rounded corners=1, color=black, line width=1] (38,\c)-- (36,\d);
\draw[rounded corners=1, color=black, line width=1] (38,\d)-- (24,\c);
}

%11
\foreach \i in {4,...,4}
{
\pgfmathsetmacro{\a}{(\i - 1) * 6 +1};
\pgfmathsetmacro{\b}{(\i - 1) * 6 +0.1 };
\draw[rounded corners=1, color=red!50!white, line width=1] (9,\a)-- (9.9,\b);

\pgfmathsetmacro{\c}{(\i - 1) * 6 +1};
\pgfmathsetmacro{\d}{(\i - 1) * 6 +3};
\draw[red!50!white,fill=red!50!white] (40,\c) circle (.25);
\draw[red!50!white,fill=red!50!white] (40,\d) circle (.25);
\draw[rounded corners=1, color=red!50!white, line width=1] (40,\d)-- (40,\c);

\pgfmathsetmacro{\e}{(\i - 1) * 6 +2};
\node at (16.5,\e) {$\xi_{11} = (\underline{0},0)$}; 
}

\foreach \i in {1,...,3}
{
\pgfmathsetmacro{\a}{(\i - 1) * 6 +1};
\pgfmathsetmacro{\b}{(\i - 1) * 6 +0.1 };
\draw[rounded corners=1, color=black, line width=1] (9,\a)-- (9.9,\b);

\pgfmathsetmacro{\c}{(\i - 1) * 6 +1};
\pgfmathsetmacro{\d}{(\i - 1) * 6 +3};
\draw[black,fill=black] (40,\c) circle (.25);
\draw[black,fill=black] (40,\d) circle (.25);
\draw[rounded corners=1, color=black, line width=1] (40,\d)-- (40,\c);
}

%12
\foreach \i in {3,...,3}
{
\pgfmathsetmacro{\a}{(\i - 1) * 6+0.1};
\pgfmathsetmacro{\b}{(\i - 1) * 6 +1 };
\draw[rounded corners=1, color=red!50!white, line width=1] (9.9,\a)-- (10.8,\b);

\pgfmathsetmacro{\c}{(\i - 1) * 6 +1};
\pgfmathsetmacro{\d}{(\i - 1) * 6 +3};
\draw[red!50!white,fill=red!50!white] (42,\c) circle (.25);
\draw[red!50!white,fill=red!50!white] (42,\d) circle (.25);
\draw[rounded corners=1, color=red!50!white, line width=1] (42,\d)-- (42,\c);

\pgfmathsetmacro{\e}{(\i - 1) * 6 +2};
\node at (16.5,\e) {$\xi_{12} = (0,\underline{0})$}; 
}

\foreach \i in {1,...,2}
{
\pgfmathsetmacro{\a}{(\i - 1) * 6 +0.1};
\pgfmathsetmacro{\b}{(\i - 1) * 6 +1 };
\draw[rounded corners=1, color=black, line width=1] (9.9,\a)-- (10.8,\b);

\pgfmathsetmacro{\c}{(\i - 1) * 6 +1};
\pgfmathsetmacro{\d}{(\i - 1) * 6 +3};
\draw[black,fill=black] (42,\c) circle (.25);
\draw[black,fill=black] (42,\d) circle (.25);
\draw[rounded corners=1, color=black, line width=1] (42,\d)-- (42,\c);
}

%13
\foreach \i in {2,...,2}
{
\pgfmathsetmacro{\a}{(\i - 1) * 6 +1};
\pgfmathsetmacro{\b}{(\i - 1) * 6 +1.9 };
\draw[rounded corners=1, color=red!50!white, line width=1] (10.8,\a)-- (11.7,\b);

\pgfmathsetmacro{\c}{(\i - 1) * 6 +1};
\pgfmathsetmacro{\d}{(\i - 1) * 6 +3};
\draw[red!50!white,fill=white] (44,\c) circle (.25);
\draw[red!50!white,fill=white] (44,\d) circle (.25);

\node at (16.5,\b) {$\xi_{13} = (0,0)$}; 
}

\foreach \i in {1,...,1}
{
\pgfmathsetmacro{\a}{(\i - 1) * 6 +1};
\pgfmathsetmacro{\b}{(\i - 1) * 6 +1.9 };
\draw[rounded corners=1, color=black, line width=1] (10.8,\a)-- (11.7,\b);

\pgfmathsetmacro{\c}{(\i - 1) * 6 +1};
\pgfmathsetmacro{\d}{(\i - 1) * 6 +3};
\draw[black,fill=black] (44,\c) circle (.25);
\draw[black,fill=black] (44,\d) circle (.25);

}

%14
\foreach \i in {1,...,1}
{
\pgfmathsetmacro{\a}{(\i - 1) * 6 +1.9};
\pgfmathsetmacro{\b}{(\i - 1) * 6 +1 };
\draw[rounded corners=1, color=red!50!white, line width=1] (11.7,\a)-- (12.6,\b);

\pgfmathsetmacro{\c}{(\i - 1) * 6 +1};
\pgfmathsetmacro{\d}{(\i - 1) * 6 +3};
\draw[red!50!white,fill=red!50!white] (46,\c) circle (.25);
\draw[red!50!white,fill=red!50!white] (46,\d) circle (.25);
\draw[rounded corners=1, color=red!50!white, line width=1] (46,\c)-- (44,\d);
\draw[rounded corners=1, color=red!50!white, line width=1] (46,\d)-- (44,\c);

\node at (16.5,\a) {$\xi_{14} = (\underline{0},\underline{0})$}; 
}

\end{tikzpicture}
\caption{
   History of the permutation
   $\sigma = 7\, 1\, 9\, 2\, 5\, 4\, 8\, 6\, 10\, 3\, 11\, 12\, 14\, 13 
    = (1,5,7,8,2)\,(3,6,4)\,(9)\,(10)\,(11,14,12)\,(13) \in \dperm_{14}$
   for the second bijection.
   Each stage of the construction shows:
   the partial almost-Dyck path $(\omega_0,\ldots,\omega_i)$,
   with the most recent step in pink;
   the label $\xi_i = (\xi'_i,\xi''_i)$,
   with entries $\xi'_i$ or $\xi''_i$ corresponding to new edges
   in the bipartite digraph being underlined;
   and the bipartite digraph $\Gamma_i$,
   with the new vertices and edges shown in pink.
   Note that the final almost-Dyck path
   is that of Figure~\ref{fig.almostdyck},
   and the final bipartite digraph is that of Figure~\ref{fig.bipartite}.
}
   \label{fig.biane}
\end{figure}

Since we have here a D-permutation $\sigma \in \dperm_{2n}$,
we will have vertices $1,\ldots,2n$ and $1',\ldots,2n'$.
First recall from \reff{eq.fk.1}/\reff{eq.fk.2} that $f_{i-1}$ is
the number of free vertices in the top row of $\Gamma_{i-1}$,
and also the number of free vertices in the bottom row of $\Gamma_{i-1}$;
and recall from \reff{eq.lemma.heights.bis} that
$f_{i-1} = \lceil h_{i-1}/2 \rceil$.
We index the free vertices on each row of $\Gamma_{i-1}$ starting from~0:
the indices are thus $0,\ldots,f_{i-1} - 1$.
We then start from the digraph $\Gamma_{i-1}$
and look at what happens at stage $i$ (see Figure~\ref{fig.biane}):
\begin{itemize}
   \item If $i$ is a cycle valley, then at stage $i$ we add no arrows.
       Since no choices are being made at this stage, we set
       $\xi'_i = \xi''_i = 0$.
   \item If $i$ is a cycle double fall or an even fixed point,
	then at stage $i$ we add an arrow from $i$ on the top row
	to an unconnected dot $j'$ on the bottom row,
       where $j = \sigma(i) \leq i$;
       then $\xi''_i$ is the index (left-to-right, counting from~0)
       of the unconnected dot $j'$
       among all the unconnected dots on the bottom row
       (together with the new dot $i'$) ---
       that is the content of \reff{def2.xi.2ndTfrac.2}.
       Note that $j=i$ (an even fixed point) corresponds to
       $\xi''_i = f_{i-1}$.

       Since no unconnected dot on the top row was touched,
       we set $\xi'_i = 0$.
   \item Similarly, if $i$ is a cycle double rise or an odd fixed point,
	we add an arrow from an unconnected dot $j$ on the top row
	to $i'$ on the bottom row, where $j = \sigma^{-1}(i) \leq i$;
        then $\xi'_i$ is the index (left-to-right, counting from~0)
        of the unconnected dot $j$
	among all the unconnected dots on the top row
        (together with the new dot $i$) ---
        that is the content of \reff{def2.xi.2ndTfrac.1}.
        Note that $j=i$ (an odd fixed point) corresponds to
        $\xi'_i = f_{i-1}$.

        Since no unconnected dot on the bottom row was touched,
        we set $\xi''_i = 0$.
   \item If $i$ is a cycle peak, then we add two arrows:
	from $i$ on the top row to the unconnected dot $j'$ on the bottom row,
       where $j = \sigma(i) < i$;
       and also from the unconnected dot $k$ on the top row
       to $i'$ on the bottom row, where $k = \sigma^{-1}(i) < i$.
       Then $\xi'_i$ (resp.\ $\xi''_i$) is the index of $k$ (resp.\ $j'$)
       among the unconnected dots on the top (resp.\ bottom) row ---
       that is the content of
       \reff{def2.xi.2ndTfrac.1}/\reff{def2.xi.2ndTfrac.2}.
\end{itemize}
%% The indices $\xi'$ and $\xi''$, wherever used, start at $0$.
All this is closely analogous to what was done in
\cite[Section~6.2, Step~2]{Sokal-Zeng_masterpoly},
but with fixed points treated differently.

\medskip

\begin{samepage}
\begin{lemma}[Nesting statistics for the second labels]
   \label{lemma.xi.nesting.second}
\nopagebreak
\quad\hfill
\vspace*{-1mm}
\begin{itemize}
   \item[(a)]  If $\sinv(i)$ is odd (i.e., $s_i$ is a fall), then
\be
   \xi'_i
   \;=\;
   \begin{cases}
      \unest'(i,\sigma)   & \hbox{\rm if $i$ is not a fixed point} \\[1mm]
      \psnest(i,\sigma)   & \hbox{\rm if $i$ is a fixed point}
   \end{cases}
 \label{eq.lemma.xi.nesting.second.a}
\ee
   \item[(b)]  If $i$ is even, then
\be
   \xi''_i
   \;=\;
   \begin{cases}
        \lnest(i,\sigma) & \hbox{\rm if $i$ is not a fixed point}
              \\[1mm]
        \psnest(i,\sigma) & \hbox{\rm if $i$ is a fixed point}
   \end{cases}
 \label{eq.lemma.xi.nesting.second.b}
\ee
\end{itemize}
\end{lemma}
\end{samepage}

\proof
(a) follows from \reff{eq.1st.2nd.AA} and \reff{eq.xii.nestings.version1}.

(b) follows from \reff{eq.1st.2nd.BB} and \reff{eq.xii.nestings}.
\qed

\begin{lemma}[Inequalities satisfied by the second labels]
   \label{lemma.ineqs.labels.second}
\nopagebreak
\quad\hfill
\vspace*{-1mm}
\begin{itemize}
   \item[(a)]  If $\sinv(i)$ is odd (i.e., $s_i$ is a fall), then
\be
   0 \;\le\; \xi'_i  \;\le\;
                      \left\lceil\dfrac{h_i}{2} \right\rceil
                      \;=\;
                      \left\lceil\dfrac{h_{i-1} -1}{2} \right\rceil
 \label{eq.xisecond.ineqs.a}
\ee
   \item[(b)]  If $i$ is even, then
\be
   0 \;\le\; \xi''_i  \;\le\; \dfrac{h_{i}}{2} 
      \;=\;
      \begin{cases}
          \left\lceil\dfrac{h_{i-1}}{2} \right\rceil
              & \hbox{\rm if $\sinv(i)$ is even (i.e., $s_i$ is a rise)} \\[5mm]
          \dfrac{h_{i-1}-1}{2}
              & \hbox{\rm if $\sinv(i)$ is odd (i.e., $s_i$ is a fall)}
      \end{cases}
      \quad
 \label{eq.xisecond.ineqs.b}
\ee
\end{itemize}
\end{lemma}

Lemma~\ref{lemma.ineqs.labels.second} will be an immediate consequence
of the following identities:

\begin{lemma}[Crossing statistics for the second labels]
   \label{lemma.xi.crossing.second}
\nopagebreak
\quad\hfill
\vspace*{-1mm}
\begin{itemize}
   \item[(a)]  If $\sinv(i)$ is odd (i.e., $s_i$ is a fall), then
\begin{subeqnarray}
   \left\lceil \dfrac{h_i}{2} \right\rceil \,-\, \xi'_i
   & = &
   \ucross'(i,\sigma)
       \,+\, {\rm I}[\hbox{\rm $i$ is odd and } \sigma(i) \neq i]
     \\[-1mm]
   & = &
   \ucross'(i,\sigma)
       \,+\, {\rm I}[\hbox{\rm $i$ is a cycle double rise}]
   \;.
 \label{eq.lemma.xi.crossing.second.a}
\end{subeqnarray}
   \item[(b)]  If $i$ is even, then
\begin{subeqnarray}
   \dfrac{h_i}{2} \,-\, \xi''_i
   & = &
   \lcross(i,\sigma)
       \,+\, {\rm I}[\sinv(i) \hbox{\rm is even and } \sigma(i) \neq i]
     \\
   & = &
   \lcross(i,\sigma)
       \,+\, {\rm I}[\hbox{\rm $i$ is a cycle double fall}]
   \;.
 \label{eq.lemma.xi.crossing.second.b}
\end{subeqnarray}
\end{itemize}
\end{lemma}

\proof
(a)  Combining \reff{eq.1st.2nd.AA} with \reff{eq.lemma.xihat.crossing.2}
yields \reff{eq.lemma.xi.crossing.second.a}.

(b) If $i$ is even, then $h_i$ is also even.
Combining \reff{eq.1st.2nd.BB} with
\reff{eq.lemma.xii.crossing.b}/\reff{eq.lemma.xii.crossing.d}
yields \reff{eq.lemma.xi.crossing.second.b}.
\qed

\medskip

% Similarly to the cycle classification in Lemma~\ref{lemma.xi.maxval}
% for our first bijection, we now consider the four possible combinations
% of $s_i$ (rise or fall) and parity of $h_i$ (odd or even):
% {\bf CHECK THIS!!!!}

We now consider the four possible combinations of $s_i$ (rise or fall)
and parity of $h_i$ (odd or even),
and determine in each case the cycle classification of the index~$i$.
Exactly as with the first bijection,
$s_i$ tells us the parity of $\sinv(i)$,
while the parity of $h_i$ equals the parity of $i$.
So these two pieces of information again tell us what was recorded in
\reff{eq.parities.1}/\reff{eq.parities.2}:
\begin{itemize}
   \item $\sinv(i)$ even and $i$ odd $\iff$ $i$ is a cycle valley
   \item $\sinv(i)$ even and $i$ even $\iff$ $i$ is either
            a cycle double fall or an even fixed point
   \item $\sinv(i)$ odd and $i$ odd $\iff$ $i$ is either
            a cycle double rise or an odd fixed point
   \item $\sinv(i)$ odd and $i$ even $\iff$ $i$ is a cycle peak
\end{itemize}
So, once again, we need only disambiguate the fixed points from the
cycle double falls/rises in the middle two cases.
We have:

\begin{lemma}[Cycle classification for second bijection]
   \label{lemma.xi.maxval.2nd}
\nopagebreak
\quad\hfill
\vspace*{-1mm}
\nopagebreak
\begin{itemize}
   \item[(a)]  If $s_i$ a rise and $h_i$ is odd (hence $h_{i-1}$ is even),
      then $i$ is a cycle valley.
   \item[(b)]  If $s_i$ a rise and $h_i$ is even (hence $h_{i-1}$ is odd),
      then $i$ is an even fixed point
      in case $\xi_i'' = \left\lceil \dfrac{h_i - 1}{2} \right\rceil$
      ($= h_i/2 = f_i$);
      otherwise it is a cycle double fall.
   \item[(c)]  If $s_i$ is a fall and $h_i$ is odd (hence $h_{i-1}$ is even),
      then $i$ is an odd fixed point
      in case $\xi_i' = \left\lceil \dfrac{h_i}{2} \right\rceil$
      ($= (h_i +1)/2 = f_i$);
      otherwise it is a cycle double rise.
   \item[(d)]  If $s_i$ is a fall and $h_i$ is even (hence $h_{i-1}$ is odd),
      then $i$ is a cycle peak.
\end{itemize}
\end{lemma}

\noindent
We remark that (b) also includes the case in which $(h_{i-1}, h_i)=(-1,0)$,
and (c) also includes the case in which $(h_{i-1}, h_i)=(0,-1)$.

\proof
(a,d) follow immediately from \reff{eq.parities.1}/\reff{eq.parities.2}.

(b) From Lemma~\ref{lemma.xi.crossing.second}(b), we have
\be
   \dfrac{h_i}{2} \,-\, \xi''_i
   \;\;
   \begin{cases}
      \ge 1   & \textrm{if $i$ is a cycle double fall}  \\[1mm]
      = 0     & \textrm{if $i$ is an even fixed point}
   \end{cases}
\ee

(c) Lemma~\ref{lemma.xi.crossing.second}(a), we have
\be
   \left\lceil \dfrac{h_i}{2} \right\rceil \,-\, \xi'_i
   \;\;
   \begin{cases}
      \ge 1   & \textrm{if $i$ is a cycle double rise}  \\[1mm]
      = 0     & \textrm{if $i$ is an odd fixed point}
   \end{cases}
\ee
\qed

\begin{lemma}[Record statistics for the second bijection]
   \label{lemma.record.second}
\nopagebreak
\quad\hfill
\vspace*{-1mm}
\begin{itemize}
   \item[(a)]  If $\sinv(i)$ is odd,
      then the index $\sinv(i)$ is a record if and only if $\xi'_i = 0$.
   \item[(b)]  If $i$ is even,
      then the index $i$ is an antirecord if and only if $\xi''_i = 0$.
\end{itemize}
\end{lemma}

\proof
This is an immediate consequence of the definitions
\reff{def2.xi.2ndTfrac.1} and \reff{def2.xi.2ndTfrac.2}.
\qed

\subsection{Step 3: Proof of bijection}

We prove that the map $\sigma \mapsto (\omega, \xi)$ is a bijection
by explicitly describing the inverse map.
That is, we let $\omega$ be any almost-Dyck path of length $2n$
and let $\xi$ be any set of labels satisfying the inequalities
\reff{def.abc.2ndTfrac}/\reff{def.abc.2ndTfrac.bis},
and we show how to reconstruct the unique D-permutation $\sigma$
that gives rise to $(\omega, \xi)$ by the foregoing construction.

In fact, the interpretation as a bipartite digraph
shows how to build the digraph,
and hence reconstruct the D-permutation $\sigma$,
by successively reading the steps $s_i$ and labels $\xi_i$.
Specifically, at stage $i$ one starts from the digraph $\Gamma_{i-1}$
and proceeds as follows (see again Figure~\ref{fig.biane}):
\begin{itemize}
  \item[(i)]
  If $s_i$ is a rise from height $h_{i-1} = 2k$ to height $h_i = 2k+1$
	(by Lemma~\ref{lemma.xi.maxval.2nd}(a) this corresponds to
         $i$ being a cycle valley),
      then we add no arrows.
  \item[(ii)]
  If $s_i$ is a rise from height $h_{i-1} = 2k-1$ to height $h_i = 2k$
           and $\xi_i = (0,k)$
        [by Lemma~\ref{lemma.xi.maxval.2nd}(b) this corresponds to
	$i$ being an even fixed point],
        we add an arrow from $i$ to $i'$.
  \item[(iii)]
  If $s_i$ is a rise from height $h_{i-1} = 2k-1$ to height $h_i = 2k$
           with $\xi_i = (0,m)$ with $0 \le m < k$
	(by Lemma~\ref{lemma.xi.maxval.2nd}(b) this corresponds to
         $i$ being a cycle double fall),
	we add an arrow from $i$ on the top row
	to the $m$th free vertex (counting from~0) on the bottom row,
        call it $j'$.
        Of course $j < i$, because these are the only vertices
        visible in the digraph $\Gamma_{i-1}$.
  \item[(iv)]
  If $s_i$ is a fall from height $h_{i-1} = 2k$ to height $h_i = 2k-1$
           and $\xi_i = (k,0)$
        [by Lemma~\ref{lemma.xi.maxval.2nd}(c) this corresponds to
        $i$ being an odd fixed point],
        we add an arrow from $i$ to $i'$.
  \item[(v)]
  If $s_i$ is a fall from height $h_{i-1} = 2k$ to height $h_i = 2k-1$
           and $\xi_i = (l,0)$ with $0 \le l < k$
        (by Lemma~\ref{lemma.xi.maxval.2nd}(c) this corresponds to
         $i$ being a cycle double rise),
        we add an arrow from the $l$th free vertex (counting from~0)
        on the top row (which is of course $< i$)
        to $i'$ on the bottom row.
  \item[(vi)]
  If $s_i$ is a fall from height $h_{i-1} = 2k+1$ to height $h_i = 2k$
           and $\xi = (l,m)$ with $0 \le l,m \le k$
	(by Lemma~\ref{lemma.xi.maxval.2nd}(d) this corresponds to
         $i$ being a cycle peak),
	we add two arrows:
	one going from $i$ on the top row to the 
	$m$th free vertex (counting from~0) on the bottom row
        (which is of course $< i$);
	and the other going from the $l$th free vertex (counting from~0)
        on the top row to $i'$ on the bottom row.
\end{itemize}
%% The indices $l$ and $m$, wherever used, start at $0$.

\medskip

Clearly, once a vertex has become the source or sink of an arrow,
it plays no further role in the construction and in particular receives no
further arrows.
Moreover, since $h_{2n} = 0$, at the end of the construction
there are no unconnected vertices.
The final result of the construction thus corresponds to a bijection
between $\{1,\ldots, 2n\}$ and $\{1',\ldots,2n'\}$,
or in other words to permutation $\sigma \in \Sym_{2n}$.

Note now that a free vertex $i$ in the top row can be created
only in situations (i) and (v), and in this case $i$ is always odd.
This means that a free vertex on the top row is always odd.
Therefore, when in situations (v) and (vi) we connect a vertex $r < i$
on the top row to $i'$ on the bottom row, $r$ is always odd.
On the other hand, we have already seen that in situations (iii) and (vi)
when we connect the vertex $i$ on the top row to a vertex $j'$ with $j<i$
on the bottom row, $i$ is always even.
These two facts together show that $\sigma$ is a D-permutation.

%% {\bf Mention similarities with de Medecis-Viennot exposition
%% of Foata-Zeilberger????}

%\bigskip

\subsection{Step 4: Translation of the statistics}

\medskip

We begin by compiling the interpretations of heights and labels
in terms of crossing and nesting statistics:

\begin{lemma}[Crossing and nesting statistics]  \label{lem.cross.nest.stat}
We have
\vspace*{-2mm}
\begin{itemize}
    \item[(a)] When $i$ is even and $h_i = 2k$,
\begin{subeqnarray}
   \xi_i'' &=& \lnest(i,\sigma)
         \qquad\qquad\; \hbox{\rm if } i\in \Cpeak \cup \Cdfall \\[1mm]
   k-\xi_i'' &=& \lcross(i,\sigma)
         \qquad\qquad \hbox{\rm if } i\in  \Cpeak \\[1mm]
   k-1-\xi_i'' &=& \lcross(i,\sigma)
         \qquad\qquad \hbox{\rm if } i\in \Cdfall  \\[1mm] 
   \xi''_i \,=\, k &=& \psnest(i,\sigma)
         \qquad\qquad\!\! \hbox{\rm if } i\in \Evenfix \\[1mm]
   \xi'_i  & = &  \unest'(i,\sigma)
         \qquad\qquad\! \hbox{\rm if } i\in  \Cpeak \\[1mm]
   k-\xi'_i  & = &  \ucross'(i,\sigma)
         \qquad\qquad\!\!\! \hbox{\rm if } i\in  \Cpeak
  \label{eq.cross.nest.stat.even}
\end{subeqnarray}
    \item[(b)] When $i$ is odd and $h_{i-1} = 2k$,
\begin{subeqnarray}
   \xi_i' &=& \unest'(i,\sigma)
        \qquad\qquad\! \hbox{\rm if } i\in \Cdrise \\[1mm]
   k-1-\xi_i' &=& \ucross'(i,\sigma)
        \qquad\qquad\!\!\! \hbox{\rm if } i\in \Cdrise \\[1mm]
    k &=& \ucross(i,\sigma) +\unest(i,\sigma)
        \qquad \hbox{\rm if } i\in \Cval \\[1mm]
    \xi'_i  \,=\, k &=& \psnest(i,\sigma)
        \qquad\qquad\!\! \hbox{\rm if } i\in \Oddfix
  \label{eq.cross.nest.stat.odd}
\end{subeqnarray}
\end{itemize}
\end{lemma}

\proof
(a) When $i$ is even, $i$ can either be a cycle peak, a cycle double fall
or an even fixed point. 
% Then (\ref{eq.cross.nest.stat.even}a) follows from
% Lemma~\ref{lemma.xi.nesting.second}(b);
% (\ref{eq.cross.nest.stat.even}b,c) follow from
% Lemma~\ref{lemma.xi.crossing.second}(b);
% (\ref{eq.cross.nest.stat.even}d) follows from
% Lemmas~\ref{lemma.xi.nesting.second}(b)
% and \ref{lemma.xi.maxval.2nd}(b);
% and (\ref{eq.cross.nest.stat.even}e) follows from
% Lemma~\ref{lemma.xi.crossing.second}(a).
%
\begin{itemize}
   \item When $i$ is a cycle double fall, $\sinv(i)$ is even,
so $s_i$ is a rise from height $h_{i-1} = 2k-1$ to height $h_i = 2k$.
Then (\ref{eq.cross.nest.stat.even}a) follows from
Lemma~\ref{lemma.xi.nesting.second}(b),
and (\ref{eq.cross.nest.stat.even}c) follows from
Lemma~\ref{lemma.xi.crossing.second}(b).
   \item When $i$ is a cycle peak, $\sinv(i)$ is odd,
so $s_i$ is a fall from height $h_{i-1} = 2k+1$ to height $h_i = 2k$.
Then (\ref{eq.cross.nest.stat.even}a) follows from
Lemma~\ref{lemma.xi.nesting.second}(b),
(\ref{eq.cross.nest.stat.even}b) follows from
Lemma~\ref{lemma.xi.crossing.second}(b),
(\ref{eq.cross.nest.stat.even}e) follows from
Lemma~\ref{lemma.xi.nesting.second}(a),
and (\ref{eq.cross.nest.stat.even}f) follows from
Lemma~\ref{lemma.xi.crossing.second}(a).
   \item When $i$ is an even fixed point, $\sinv(i)$ is even,
so $s_i$ is a rise from height $h_{i-1} = 2k-1$ to height $h_i = 2k$.
Then (\ref{eq.cross.nest.stat.even}d) follows from
Lemmas~\ref{lemma.xi.nesting.second}(b)
and \ref{lemma.xi.maxval.2nd}(b).
\end{itemize} 

% \begin{itemize}
% \item When $i$ is a cycle peak or a cycle double fall, 
% we use Equation~\reff{eq.1st.2nd.BB} and 
% Lemma~\ref{lemma.xii.nesting} to obtain $\xi_i'' = \lnest(i,\sigma)$.
% 
% \item When $i$ is a cycle double fall, $s_i$ is a rise from height
% $h_{i-1} = 2k-1$ to height $h_i = 2k$.
% Using Equation~\reff{eq.1st.2nd.BB} and Lemma~\ref{lemma.xii.crossing}(b)
% we get $k-1-\xi_i'' = \lcross(i,\sigma)$.
% 
% \item When $i$ is a cycle peak, $s_i$ is a fall from height
% $h_{i-1} = 2k+1$ to height $h_i = 2k$.
% Using Equation~\reff{eq.1st.2nd.BB} and Lemma~\ref{lemma.xii.crossing}(d)
% we get $k-\xi_i'' = \lcross(i,\sigma)$.
% 
% \item When $i$ is an even fixed point, we use
% Equations~\reff{eq.1st.2nd.BB} and~\reff{eq.xii.fp} and 
% Lemma~\ref{lemma.xii.nesting} to obtain $k = \psnest(i,\sigma)$.
% \end{itemize}

(b) When $i$ is odd, $i$ can either be a cycle valley, a cycle double rise
or an odd fixed point.
\begin{itemize}
   \item When $i$ is a cycle double rise, $\sinv(i)$ is odd,
so $s_i$ is a fall from height $h_{i-1} = 2k$ to height $h_i = 2k-1$.
Then (\ref{eq.cross.nest.stat.odd}a) follows from
Lemma~\ref{lemma.xi.nesting.second}(a),
and (\ref{eq.cross.nest.stat.odd}b) follows from
Lemma~\ref{lemma.xi.crossing.second}(a).
   \item When $i$ is a cycle valley, $\sinv(i)$ is even,
so $s_i$ is a rise from height $h_{i-1} = 2k$ to height $h_i = 2k+1$.
As the almost-Dyck path is the same as in the first bijection,
we can use Lemmas~\ref{lemma.xii.nesting} and \ref{lemma.xii.crossing}(a)
to obtain $k = \ucross(i,\sigma) +\unest(i,\sigma)$,
which is \ref{eq.cross.nest.stat.odd}c).
   \item When $i$ is an odd fixed point, $\sinv(i)$ is odd,
so $s_i$ is a fall from height $h_{i-1} = 2k$ to height $h_i = 2k-1$.
Then (\ref{eq.cross.nest.stat.odd}d) follows from
Lemmas~\ref{lemma.xi.nesting.second}(a)
and \ref{lemma.xi.maxval.2nd}(c).
\end{itemize} 
\qed

{\bf Remarks.}
1.  When $i$ is a cycle valley, $\sinv(i)$ is even and $i$ is odd,
so Lemmas~\ref{lemma.xi.nesting.second} and \ref{lemma.xi.crossing.second}
tell us nothing about ucross and unest;
and the labels $\xi'_i = \xi''_i = 0$ carry no information.
That is why in this case we learn only about the sum $\ucross + \unest$,
which can be deduced from the heights alone.

2. Our treatment of cycle double rises is slightly nicer than
that achieved in \cite[Lemma~6.4]{Sokal-Zeng_masterpoly},
thanks to the introduction of $\ucross'$ and $\unest'$.

3. In (\ref{eq.cross.nest.stat.even}c) or (\ref{eq.cross.nest.stat.odd}b),
$k=0$ is impossible.  These would correspond, respectively,
to a rise from height $-1$ or a fall from height $0$,
which by Lemma~\ref{lemma.heights} occur only when
$i$ is a record-antirecord fixed point,
not when it is a cycle double fall or cycle double rise.
\myendremark

Finally, and most crucially, we come to the counting of cycles (cyc).
We use the term {\em cycle closer}\/ to denote the largest element
in a non-singleton cycle.
(This is the same as the ``cycle peak maximum'' defined in
 Section~\ref{subsec.second.statements}.)
Obviously every non-singleton cycle has precisely one cycle closer.
A cycle closer is always a cycle peak, but not conversely.
So we need to know how many of the cycle peaks are the cycle closers.
The answer is as follows:

\begin{lemma}[Counting of cycles]
   \label{lemma.cycles.2}
Fix $i=2k$ such that $k \in [n]$,
and fix $(s_1,\ldots,s_{i-1})$ and $(\xi_1,\ldots,\xi_{i-1})$.
Consider all permutations $\sigma \in \dperm_{2n}$
that have those given values for the first $i-1$ steps and labels
and for which $i$ is a cycle peak.
Then:
\begin{itemize}
   \item[(a)] The value of $\xi_i = (\xi'_i,\xi''_i)$ completely determines
      whether $i$ is a cycle closer or not.
   \item[(b)] For each value $\xi'_i \in [0, (h_{i-1}-1)/2]$
	   there is precisely one value
           $\xi''_i \in$ \linebreak $[0, (h_{i-1}-1)/2]$
      that makes $i$ a~cycle closer, and conversely.
\end{itemize}
\end{lemma}

The proof is similar to the proof of \cite[Lemma~6.5]{Sokal-Zeng_masterpoly}.

\proof
We use once again the bipartite digraph
%% employed in the proof of Lemma~\ref{lemma.heights},
of Figure~\ref{fig.bipartite},
and let us also draw a vertical dotted line
(with an upwards arrow) to connect each pair $j' \to j$.
Now consider the restriction of this digraph to
the vertex set $\{1,\ldots,{i-1},1',\ldots,{(i-1)'}\}$:
as discussed in Step~3, this restriction can be reconstructed
from the steps $(s_1,\ldots,s_{i-1})$
and the labels $(\xi_1,\ldots,\xi_{i-1})$.
The connected components of this restriction are of two types:
complete directed cycles and directed open chains;
they correspond to cycles of $\sigma$
whose cycle closers are, respectively, $\le i-1$ and $> i-1$.
Each directed open chain runs from an unconnected dot on the bottom row
to an unconnected dot on the top row.

Now suppose that $i$ is a cycle peak.
Then at stage $i$ we add two arrows:
from $i$ on the top row to an unconnected dot $j'$ on the bottom row;
and also from an unconnected dot $k$ on the top row
to $i'$ on the bottom row.
Here $\xi'_i$ (resp.\ $\xi''_i$) is the index of $k$ (resp.\ $j'$)
among the unconnected dots on the top (resp.\ bottom) row.

Now the point is simply this: $i$ is a cycle closer if and only if
$j'$ and $k$ belong to the {\em same}\/ directed open chain
(with $j'$ being its starting point and $k$ being its ending point).
So for each value $\xi'_i \in [0, (h_{i-1}-1)/2]$
there is precisely one value $\xi''_i \in [0, (h_{i-1}-1)/2]$
that makes $i$ a cycle closer, and conversely.
\qed

% {\bf The proof is exactly the same except for the intervals.
%   Should the details be omitted and a note on Laguerre digraphs
%   be provided instead?????}

%\bigskip

\subsection{Step 5: Computation of the weights}
We can now compute the weights associated to the
0-Schr\"oder path $\omegahat$
in Theorem~\ref{thm.flajolet_master_labeled_Schroder},
which we recall are
$a_{h,\xi}$ for a rise starting at height $h$ with label $\xi$,
$b_{h,\xi}$ for a fall starting at height $h$ with label $\xi$, and
$c_{h,\xi}$ for a long level step at height $h$ with label $\xi$.
(Of course, in the present case we have long level steps only at height~0.)
We do this by putting together the information
collected in Lemmas~\ref{lemma.ineqs.labels.second},
\ref{lemma.xi.maxval.2nd}, \ref{lem.cross.nest.stat} and \ref{lemma.cycles.2}:
\begin{itemize}
   \item[(a)] Rise from height $h_{i-1} = 2k$ to height $h_i = 2k+1$
      (hence $i$ odd):
\begin{itemize}
   \item By (\ref{def.abc.2ndTfrac.bis}a), the label is $\xi_i = (0,0)$.
   \item By Lemma~\ref{lemma.xi.maxval.2nd}(a), this is a cycle valley.
   \item By (\ref{eq.cross.nest.stat.odd}c),
	   $k = \ucross(i,\sigma) + \unest(i,\sigma)$.	   
\end{itemize}
Therefore, from \reff{def.Qhatn.secondmaster}, the weight for this step is
\be
   a_{2k,(0,0)}  \;=\;  \sfa_{k}
   \;.
\ee
   \item[(b)] Rise from height $h_{i-1} = 2k-1$ to height $h_i = 2k$
      (hence $i$ even):
\begin{itemize}
   \item By (\ref{def.abc.2ndTfrac.bis}b)
      and Lemma~\ref{lemma.ineqs.labels.second}(b),
      the label is $\xi_i = (0,\xi''_i)$ with $0 \le \xi''_i \le k$.
   \item By Lemma~\ref{lemma.xi.maxval.2nd}(b), this is a cycle double fall
	   if $0\le \xi_i''<k$, and an even fixed point if $\xi_i'' = k$.
   \item By (\ref{eq.cross.nest.stat.even}a,d),
\be
   \xi_i'' 
   \;=\;
   \begin{cases}
       \lnest(i,\sigma)  &  \textrm{if $i$ is a cycle double fall} \\
       \psnest(i,\sigma) &  \textrm{if $i$ is an even fixed point}
   \end{cases}
\ee
   \item By (\ref{eq.cross.nest.stat.even}c),
      $\lcross(i,\sigma) = k-1-\xi_i''$ when $i$ is a cycle double fall.
\end{itemize}
Therefore, from \reff{def.Qhatn.secondmaster}, the weight for this step is
\be
    a_{2k-1,(0,\xi'')}  \;=\;
   \begin{cases}
      \sfc_{k-1-\xi'',\xi''}   & \textrm{if $0 \le \xi'' < k$}  \\[1mm]
      \lambda \sfe_k               & \textrm{if $\xi'' = k$}
   \end{cases}
\ee
   \item[(c)] Fall from height $h_{i-1} = 2k$ to height $h_i = 2k-1$
      (hence $i$ odd):
\begin{itemize}
   \item By (\ref{def.abc.2ndTfrac.bis}c)
      and Lemma~\ref{lemma.ineqs.labels.second}(a),
      the label is $\xi_i = (\xi'_i,0)$ with $0 \le \xi'_i \le k$.
   \item By Lemma~\ref{lemma.xi.maxval.2nd}(c), this is a cycle double rise
           if $0\le \xi_i'<k$, and an odd fixed point if $\xi_i' = k$.
   \item By (\ref{eq.cross.nest.stat.odd}a,d),
\be
   \xi_i' 
   \;=\;
   \begin{cases}
      \unest'(i,\sigma)  &  \textrm{if $i$ is a cycle double rise} \\
      \psnest(i,\sigma)  &  \textrm{if $i$ is an odd fixed point}
   \end{cases}
\ee
   \item By (\ref{eq.cross.nest.stat.odd}b),
       $\ucross'(i,\sigma) = k-1 - \xi_i'$ when $i$ is a cycle double rise.
\end{itemize}
Therefore, from \reff{def.Qhatn.secondmaster}, the weight for this step is
\be
   b_{2k,(\xi',0)}  \;=\;
   \begin{cases}
      \sfd_{k-1-\xi',\xi'}   & \textrm{if $0 \le \xi' < k$}  \\[1mm]
      \lambda \sff_k         & \textrm{if $\xi' = k$}
   \end{cases}
\ee
   \item[(d)] Fall from height $h_{i-1} = 2k+1$ to height $h_i = 2k$
      (hence $i$ even):
\begin{itemize}
   \item By (\ref{def.abc.2ndTfrac.bis}d)
      and Lemma~\ref{lemma.ineqs.labels.second}(a,b),
      the label is $\xi_i = (\xi'_i,\xi''_i)$ with $0 \le \xi'_i \le k$
      and $0 \le \xi''_i \le k$.
   \item By Lemma~\ref{lemma.xi.maxval.2nd}(d), this is a cycle peak.
   \item By (\ref{eq.cross.nest.stat.even}a), $\xi_i'' = \lnest(i,\sigma) $.
   \item By (\ref{eq.cross.nest.stat.even}b), $k-\xi_i'' = \lcross(i,\sigma)$.
\end{itemize}
For each choice of $\xi_i''\in[0,k]$
there are $k+1$ possible choices of $\xi_i'$,
of which one closes a cycle and the rest don't; this is the content of 
Lemma~\ref{lemma.cycles.2}. 
Therefore, from \reff{def.Qhatn.secondmaster}, the total weight of all such steps is
\be
   b_{2k+1}
   \;\eqdef\;
   \sum_{\xi',\xi''} b_{2k+1,(\xi',\xi'')}
   \;=\;
  (\lambda + k) \sum_{\xi'' = 0}^{k}\sfb_{k-\xi'',\xi''}
   \;.
\ee
   \item[(e)] Long level step at height~$0$: \\[2mm]
This corresponds in the almost-Dyck path $\omega$
to a fall from height~$0$ to height~$-1$,
followed by a rise from height~$-1$ to height~$0$.
Applying case~(c) with $k=0$ and $\xi=0$,
followed by case~(b) with $k=0$ and $\xi=0$,
we obtain a weight
\be
   c_{0,0}  \;=\;  \lambda^2 \sfe_0 \sff_0
   \;.
\ee
\end{itemize}

Putting this all together
in Theorem~\ref{thm.flajolet_master_labeled_Schroder},
we obtain a T-fraction with
\begin{eqnarray}
   \alpha_{2k-1}
   & = &
   (\hbox{rise from $2k-2$ to $2k-1$})
   \,\times\,
   (\hbox{fall from $2k-1$ to $2k-2$})
       \nonumber \\[2mm]
   & = &
   \sfa_{k-1} \:
   (\lambda+k-1)
   \biggl( \sum_{\xi=0}^{k-1} \sfb_{k-1-\xi,\xi} \biggr)
       \label{alpha.2k-1.master.2nd}  \\[4mm]
   \alpha_{2k}
   & = &
   (\hbox{rise from $2k-1$ to $2k$})
   \,\times\,
   (\hbox{fall from $2k$ to $2k-1$})
       \nonumber \\[2mm]
   & = &
   \biggl( \lambda\sfe_k \,+\, \sum_{\xi=0}^{k-1} \sfc_{k-1-\xi,\xi} \biggr)
   \biggl( \lambda\sff_k \,+\, \sum_{\xi=0}^{k-1} \sfd_{k-1-\xi,\xi} \biggr)
       \label{alpha.2k.master.2nd}  \\[2mm]
   \delta_1  & = &  \lambda^2\sfe_0 \sff_0
       \label{delta1.master.2nd} \\[1mm]
   \delta_n  & = &   0    \qquad\hbox{for $n \ge 2$}
       \label{deltan.master.2nd}
\end{eqnarray}
This completes the proof of Theorem~\ref{thm.Tfrac.second.master}.
\qed

\proofof{Theorem~\ref{thm.Tfrac.second.pq}}
Comparing \reff{def.Pn.pq.second} with \reff{def.Qhatn.secondmaster}
and using Lemma~\ref{lemma.record.second},
and recalling that we are making the specializations
$v_1 = y_1$ and $q_{+1} = p_{+1}$,
we see that the needed weights in \reff{def.Qhatn.secondmaster} are
\begin{eqnarray}
   \sfa_{k}
   & = &
   p_{+1}^{k} \, y_1
        \\[2mm]
   \sfb_{k-1-\xi,\xi}
   & = &
   p_{-1}^{k-1-\xi} q_{-1}^\xi  \,\times\,
   \begin{cases}
      x_1  & \textrm{if $\xi = 0$}   \\
      u_1  & \textrm{if $1 \le \xi \le k-1$}
   \end{cases}
        \\[2mm]
   \sfc_{k-1-\xi,\xi}
   & = &
   p_{-2}^{k-1-\xi} q_{-2}^\xi  \,\times\,
   \begin{cases}
      x_2  & \textrm{if $\xi = 0$}   \\
      u_2  & \textrm{if $1 \le \xi \le k-1$}
   \end{cases}
        \\[2mm]
   \sfd_{k-1-\xi,\xi}
   & = &
   \phat_{+2}^{k-1-\xi} \qhat_{+2}^\xi  \,\times\,
   \begin{cases}
      \yhat_2  & \textrm{if $\xi = 0$}   \\
      \vhat_2  & \textrm{if $1 \le \xi \le k-1$}
   \end{cases}
        \\[2mm]
   \sfe_k
   & = &
   \begin{cases}
      \ze  & \textrm{if $k = 0$}   \\[1mm]
      \se^k \we  & \textrm{if $k \ge 1$}
   \end{cases}
        \\[2mm]
   \sff_k
   & = &
   \begin{cases}
      \zo  & \textrm{if $k = 0$}   \\[1mm]
      \so^k \wo  & \textrm{if $k \ge 1$}
   \end{cases}
\end{eqnarray}
Inserting these into \reff{alpha.2k-1.master.2nd}--\reff{deltan.master.2nd}
yields the continued-fraction coefficients
\reff{eq.thm.Tfrac.second.weights.pq}.
\qed\hspace*{-5mm}

{\bf Remarks.}
1.  We needed to make the specialization $v_1 = y_1$
because Lemma~\ref{lemma.record.second} tells us nothing about
the record status of cycle valleys,
for which $\sinv(i)$ is even and $i$ is odd.
Similarly, we needed to make the specialization $q_{+1} = p_{+1}$
because, for cycle valleys, Lemma~\ref{lem.cross.nest.stat}(b)
does not tell us about $\ucross$ and $\unest$ individually,
but only their sum.

2. There is a variant master polynomial that we could have treated:
instead of weighting cycle peaks using
$\sfb_{\lcross(i,\sigma),\,\lnest(i,\sigma)}$
as in \reff{def.Qhatn.secondmaster},
we could instead weight them as
$\sfb_{\ucross'(i,\sigma),\,\unest'(i,\sigma)}$.
Then we would use (\ref{eq.cross.nest.stat.even}e,f)
instead of (\ref{eq.cross.nest.stat.even}a,b);
the resulting T-fraction would be the same.
\myendremark

\proofof{Theorem~\ref{thm.Tfrac.second}}
Specialize Theorem~\ref{thm.Tfrac.second.pq} to
$p_{-1} = p_{-2} = p_{+1} = \phat_{+2} = q_{-1} = q_{-2} = \qhat_{+2} =
 \se = \so = 1$.
\qed

% \subsection{Proof of a conjecture of Randrianarivony and Zeng}
%   \label{subsec.RZ}
% 
% {\bf Need to define all the statistics (lema etc should \emph{not} include
%    fixed points!!) and define the polynomials $R_n(x,y,\bar{x},\bar{y})$
%    and $G_n(x,y,\bar{x},\bar{y})$ and then prove continued fractions
%    for both.  The equality $G_n = R_n$ will be an immediate corollary.
% 
% The first Proposition (for $R_n$) will be our recounting of what is
% essentially R--Z's proof, using the variant first bijection
% with labels $\widehat{\xi}_i$.
% The second Proposition (for $G_n$), however, will need our second bijection,
% which was not available to R--Z.}

\section{Some final remarks}  \label{sec.final}

When we began this work, we envisioned it as analogous to,
though probably more complicated than, the study of the ``linear family''
[cf.~\reff{eq.linear_family}]
that was undertaken in \cite{Sokal-Zeng_masterpoly}:
our goal was to make a corresponding study of the ``quadratic family''
[cf.~\reff{eq.quadratic_family}].
What has surprised us is that the final results,
as well as the associated methods of proof,
turned out to be, not merely {\em analogous}\/
to those of \cite{Sokal-Zeng_masterpoly},
but in fact very closely {\em parallel}\/
--- much more closely parallel than we expected.

The key objects of \cite{Sokal-Zeng_masterpoly} were permutations of $[n]$;
the key results were J-fractions;
and the proofs involved bijections from permutations of $[n]$
to labeled Motzkin paths of length $n$,
using the Foata--Zeilberger bijection for the first J-fraction
(which did not involve the counting of cycles)
and the Biane bijection for the second J-fraction
(which involved the counting of cycles).
The two bijections have the same paths but different labels.

By contrast, the key objects of the present paper are D-permutations of $[2n]$;
the key results are 0-T-fractions
(that is, T-fractions that have $\delta_i = 0$ for all $i \ge 2$);
and the proofs involve bijections from D-permutations of $[2n]$
to labeled 0-Schr\"oder paths
(or equivalently, labeled almost-Dyck paths) of length $2n$.
Once again, the bijections for the first and second continued fractions
have the same paths but different labels.
Since these 0-Schr\"oder paths are very different from
the Motzkin paths of \cite{Sokal-Zeng_masterpoly},
the definitions of the paths must also be very different, and indeed they are.
The surprise was that the definitions of the {\em labels}\/
in our constructions turned out to be almost identical
to those employed in \cite{Sokal-Zeng_masterpoly}
(where ``almost'' means that fixed points are treated differently):
once again, we use the Foata--Zeilberger labels for the first T-fraction
(which does not involve the counting of cycles)
and the Biane labels for the second T-fraction
(which involves the counting of cycles).

The final results also turned out to be amazingly similar.
Compare, for instance, Theorem~\ref{thm.Tfrac.first.pq}
with \cite[Theorem~2.7]{Sokal-Zeng_masterpoly}:
our coefficient $\alpha_{2k-1}$ is identical
to the coefficient $\beta_k$ in \cite[eq.~(2.53c)]{Sokal-Zeng_masterpoly};
our coefficient $\alpha_{2k}$ includes the same ingredients
as the coefficient $\gamma_k$ in \cite[eq.~(2.53b)]{Sokal-Zeng_masterpoly},
but combined as a product rather than a sum.
An analogous comparison holds between our first master T-fraction
(Theorem~\ref{thm.Tfrac.first.master})
and the first master J-fraction in the earlier work
\cite[Theorem~2.9]{Sokal-Zeng_masterpoly}.
Furthermore, an analogous comparison
holds between our conjectured second T-fraction
(Conjecture~\ref{conj.Tfrac.second})
and the conjectured second J-fraction found in the earlier work
\cite[Conjecture~2.3]{Sokal-Zeng_masterpoly},
and between our proved second T-fraction
(Theorem~\ref{thm.Tfrac.second})
when specialized to $\vhat_2 = \yhat_2$
and the proved second J-fraction found in the earlier work
\cite[Theorem~2.4]{Sokal-Zeng_masterpoly}.
Finally, an almost analogous comparison holds between
our second master T-fraction (Theorem~\ref{thm.Tfrac.second.master})
and the second master J-fraction in the earlier work
\cite[Theorem~2.14]{Sokal-Zeng_masterpoly},
the only difference being that the treatment of $\sfd$
looks a bit more natural in the present work.

All this suggests to us that D-permutations are,
among all the combinatorial models of the
Genocchi and median Genocchi numbers,
a particularly well-behaved one,
which is closely analogous to ordinary permutations.
It would be interesting to know whether similar continued fractions,
in a large (or infinite) number of independent indeterminates,
can be found for some of the other models
of the Genocchi and median Genocchi numbers.
Our approach in \cite{Sokal-Zeng_masterpoly} and the present paper
has been to introduce a natural classification of indices
into mutually exclusive categories
--- here the parity-refined record-and-cycle classification
(Section~\ref{subsec.statistics.1}) ---
and to build a homogeneous multivariate generating polynomial
implementing this classification.
It would be interesting to find analogous natural classifications
for the other combinatorial models.

\section*{Acknowledgments}

We wish to thank Be\'ata B\'enyi, Sheila Sundaram and Jiang Zeng
for helpful conversations and/or correspondence.

The first author (B.D.)\ also wishes to thank the organizers of
the Graduate Student Combinatorics Conference (GSCC 2022),
the Graduate Online Combinatorics Colloquium (GOCC),
the combinatorics seminar series at the 
Institute of Mathematical Sciences (Chennai, India),
and the Heilbronn workshop on
Positivity Problems Associated to Permutation Patterns (Lancaster, UK)
for providing the opportunity to present preliminary versions of this work.

This research was supported in part by the
U.K.~Engineering and Physical Sciences Research Council grant EP/N025636/1,
and by a teaching assistantship from the Department of Mathematics,
University College London.

\appendix
\section{J-fraction for the polynomials \reff{eq.Pn.xylam}}  \label{app.Pnxylam}

% From Bishal's file ugly_Jfrac.tex

The first few polynomials
\be
   P_n(x,y,\lambda)
   \;=\;
   \sum_{\sigma\in \dperm_{2n}}
       x^{\arec(\sigma)}y^{\erec(\sigma)} \lambda^{\cyc(\sigma)}
 \label{eq.Pn.xylam.bis}
\ee
are too complicated to print even at $n=3$.
But we can give the first few J-fraction coefficients:  they are
\begin{subeqnarray}
\gamma_0 &=& \lambda x (\lambda x + y)  \\[1mm]
\beta_1  &=& \lambda xy (\lambda + x) (\lambda + y) \\[1mm]
\gamma_1 &=& (1 + \lambda) (\lambda + x + y + x y) \\[1mm]
\beta_2  &=& \lambda^3 (1 + x y) + \lambda^2 (2 + x + x^2 + y + 4 x y + y^2) + 
 \nonumber \\
 && \lambda (1 + 3 x + x^2 + 3 y + 4 x y + x^2 y + y^2 + x y^2 + x^2 y^2) +
 \nonumber \\
 && 2 (x + y + x^2 y + x y^2)
\end{subeqnarray}
followed by
\be
\gamma_2 \;=\; \dfrac{N}{D}
\ee
where
%{\bf $N$ has a factor $1+\lambda$ that needs to be pulled out!!!
%  And collect factors as polynomials in $\lambda$!!!}
\begin{eqnarray}
\nonumber
N &=& (1+\lambda) \, \left[ \lambda^4\,(1 + x y) \; + \; \lambda^3\, (5 + x + 2 x^2 + y + 10 x y + 2 y^2 + x^2 y^2)  \right.\\
&& + \; \lambda^2 \, (7 + 8 x + 8 x^2 + 8 y + 26 x y + 4 x^2 y + 2 x^3 y + 8 y^2 + 
 4 x y^2 + 7 x^2 y^2 + 2 x y^3) \nonumber\\
   && +\; \lambda\, (3 + 15 x + 8 x^2 + 2 x^3 + 15 y + 22 x y + 18 x^2 y + 5 x^3 y + 
 8 y^2 + 18 x y^2   \nonumber\\
   &&  \qquad + 13 x^2 y^2 +  x^3 y^2 + 2 y^3 + 5 x y^3 + x^2 y^3) \nonumber\\
   && + \; (8 x + 2 x^2 + 2 x^3 + 8 y + 5 x y + 18 x^2 y + x^3 y + 2 y^2 + 
 18 x y^2 + 4 x^2 y^2 + 4 x^3 y^2  \nonumber\\
   &&  \quad \;\left. + 2 y^3 + x y^3 + 4 x^2 y^3 + 
 x^3 y^3)\right] \nonumber\\
\end{eqnarray}
and
\begin{eqnarray}
 D &=& \lambda^3\, (1 + x y) \; + \; \lambda^2\,(2 + x + x^2 + y + 4 x y + y^2) \nonumber
        \\
   && + \;\lambda \,(1 + 3 x + x^2 + 3 y + 4 x y + x^2 y + y^2 + x y^2 + x^2 y^2) \; + \; 2 (x + y) (1 + x y)\nonumber\\
   &&
\end{eqnarray}
It can then be shown that
\begin{itemize}
\item[(a)] $\gamma_2$ is not a polynomial in $x$ (when $y$ and $\lambda$ are given fixed real values) unless  $\lambda \in \{-2,-1,+1\}$ or $y \in  \{-1,+1\}$.

\item[(b)] $\gamma_2$ is not a polynomial in $y$ unless $\lambda \in \{-2,-1,+1\}$ or $x \in  \{-1,+1\}$.

\item[(c)] $\gamma_2$ is not a polynomial in $\lambda$ unless $x \in \{-1,+1\}$ or $y \in  \{-1,+1\}$ or $x=y=0$.
\end{itemize}

We state the polynomials obtained from the foregoing
specializations of $\gamma_2$:

\begin{center}
\begin{tabular}{c|c}
Specialization & $\gamma_2$\\
\hline
$\lambda=+1$ & $2 (2 + x) (2 + y)$ \\
$\lambda=-1$ & $0$ \\
$\lambda=-2$ & $-1 + x + y + x y/2$ \\
$x=+1$ & $(2 + \lambda) (3 + \lambda + 2 y)$ \\
$x=-1$ & $(1 + \lambda) (4 + \lambda - 2 y)$\\
$y=+1$ & $(2 + \lambda) (3 + \lambda + 2 x)$ \\
$y=-1$ & $(1 + \lambda) (4 + \lambda - 2 x)$\\
$x=y=0$ & $(1 + \lambda) (3 + \lambda)$ \\
\end{tabular}
\end{center}

These give rise to continued fractions as follows:

\medskip

$\bm{\lambda=+1}$.
By Theorem~\ref{thm.Tfrac.first}
specialized to $x_1 = x_2 = \ze = \zo = x$, $y_1 = y_2 = y$,
$u_1 = u_2 = v_1 = v_2 = \we = \wo = 1$,
we obtain a T-fraction with
$\alpha_{2k-1} = (x+k-1)(y+k-1)$, $\alpha_{2k} = (x+k)(y+k)$,
$\delta_1 = x^2$
and hence by contraction (Proposition~\ref{prop.contraction_even.Ttype})
a J-fraction with $\gamma_0 = x(x+y)$,
$\gamma_n = 2(x+n)(y+n)$ for $n \ge 1$,
$\beta_n = (x+n-1)(x+n)(y+n-1)(y+n)$.

\medskip

$\bm{x=+1}$.
By the reversal map $i \mapsto 2n+1-i$,
the weight $x^{\arec(\sigma)}y^{\erec(\sigma)}$
is equivalent to $x^{\rec(\sigma)}y^{\earec(\sigma)}$.
Now apply Theorem~\ref{thm.Tfrac.second}
specialized to $x_1 = x_2 = y$,
$y_1 = \yhat_2 = u_1 = u_2 = v_1 = \vhat_2 = \we = \wo = \ze = \zo = 1$:
we obtain a T-fraction with
$\alpha_{2k-1} = (\lambda+k-1)(y+k-1)$,
$\alpha_{2k} = (\lambda+k)(y+k-1+\lambda)$,
$\delta_1 = \lambda^2$.
By contraction this yields a J-fraction with
$\gamma_0 = \lambda (\lambda+y)$,
$\gamma_n = (\lambda+n) (\lambda+2n-1+2y)$ for $n \ge 1$,
$\beta_n = (\lambda+n-1)(\lambda+n)(y+n-1)(y+n-1+\lambda)$.

\medskip

$\bm{y=+1}$.
By Theorem~\ref{thm.Tfrac.second}
specialized to $x_1 = x_2 = \ze = \zo = x$,
$y_1 = \yhat_2 = u_1 = u_2 = v_1 = \vhat_2 = \we = \wo = 1$,
we obtain a T-fraction with
$\alpha_{2k-1} = {(\lambda+k-1)(x+k-1)}$,
$\alpha_{2k} = (\lambda+k)(x+k-1+\lambda)$,
$\delta_1 = \lambda^2 x^2$.
By contraction this yields a J-fraction with
$\gamma_0 = \lambda x (1 + \lambda x)$,
$\gamma_n = (\lambda+n) (\lambda+2n-1+2x)$ for $n \ge 1$,
$\beta_n = (\lambda+n-1)(\lambda+n)(x+n-1)(x+n-1+\lambda)$.

\medskip

$\bm{\lambda=-1}$.
We {\em conjecture}\/ a J-fraction with
$\gamma_0 = x(x-y)$, $\gamma_n  = 0$ for $n \ge 1$,
$\beta_n = -xy (1-x)(1-y)$.
We have verified this through $\gamma_2$ and $\beta_3$.

\medskip

$\bm{\lambda=-2}$, $\bm{x=-1}$, $\bm{y=-1}$.
These are polynomial through $\gamma_2$ and $\beta_3$,
but we are unable to guess the general formula,
and indeed we do not know whether they will continue to be polynomial
at higher orders.

\medskip

$\bm{x=0}$.
For all $n \ge 1$ and all $\sigma \in \Sym_n$,
the index $n$ is an antirecord.
Therefore, setting $x=0$ suppresses all permutations for $n \ge 1$,
and we have $P_n(0,y,\lambda) = \delta_{n0}$ (Kronecker delta).
The value of $\gamma_2$ stated in the table above is completely irrelevant,
because $\beta_1 = 0$.

\medskip

$\bm{y=0}$.
It is not difficult to show that
the only permutation $\sigma \in \Sym_n$
with no exclusive records is the identity permutation.
(The index~1 is always a record, so if there are no exclusive records
 it must be a fixed point;
 then the index~2 will be a record, and so forth.)
Therefore $P_n(x,0,\lambda) = \lambda^{2n} x^{2n}$.
This gives an S-fraction with $\alpha_1 = \lambda^2 x^2$
and $\alpha_n = 0$ for $n \ge 2$.
Equivalently, it gives a J-fraction with $\gamma_0 = \lambda^2 x^2$
and all other coefficients zero.
The value of $\gamma_2$ stated in the table above is again
completely irrelevant, because $\beta_1 = 0$.

%% \appendix
%% \section{Combinatorial models for the Genocchi and median Genocchi numbers}
%%    \label{app.models}


\addcontentsline{toc}{section}{Bibliography}
\begin{thebibliography}{199}

% \bibitem{Aigner_99}  M. Aigner, Catalan-like numbers and determinants,
%    J. Combin. Theory A {\bf 87}, 33--51 (1999).
% 
% \bibitem{Aigner_01a}  M. Aigner, Catalan and other numbers: a recurrent theme,
%    in: {\em Algebraic Combinatorics and Computer Science}\/,
%    edited by H.~Crapo and D.~Senato
%    (Springer-Verlag Italia, Milan, 2001), pp.~347--390.
% 
% \bibitem{Aigner_01b}  M. Aigner, Lattice paths and determinants,
%    in:  {\em Computational Discrete Mathematics: Advanced Lectures}\/,
%    edited by H.~Alt,
%    Lecture Notes in Computer Science \#2122)
%    (Springer-Verlag, Berlin, 2001), pp.~1--12.
% 
% \bibitem{Aigner_14}  M. Aigner and G.M. Ziegler,
%    {\em Proofs from THE BOOK}\/, 5th ed.\ 
%    (Springer-Verlag, Berlin--Heidelberg, 2014).
%    %% {\bf Change citations to forthcoming 6th ed., 2018!!!!!}
% 
% \bibitem{Akhiezer_65}  N.I. Akhiezer, {\em The Classical Moment Problem
%    and Some Related Questions in Analysis}\/,
%    translated by N.~Kemmer
%    (Hafner, New York, 1965).
% 
% \bibitem{Albenque_12}  M. Albenque and J. Bouttier,
%    Constellations and multicontinued fractions: Application to Eulerian
%    triangulations,
%    in {\em 24th International Conference on Formal Power Series and
%     Algebraic Combinatorics (FPSAC 2012)}\/,
%    Discrete Mathematics \& Theoretical Computer Science Proceedings
%    (Nancy, France, 2012), pp.~805--816.
% 
% % \bibitem{Anderson_18}  D. Anderson, E.S. Egge, M. Riehl, L. Ryan,
% %    R. Steinke and Y. Vaughan,
% %    Pattern avoiding linear extensions of rectangular posets,
% %    J. Comb. {\bf 9}, 185--220 (2018).
% %    %% arXiv:1605.06825 [math.CO]
% 
% \bibitem{Andrews_99}  G.E. Andrews, R. Askey and R. Roy,
%    {\em Special Functions}\/ (Cambridge University Press, Cambridge, 1999).
% 
% \bibitem{Armstrong_09}  D. Armstrong, Generalized noncrossing partitions
%    and combinatorics of Coxeter groups,
%    Mem. Amer. Math. Soc. {\bf 202}, no.~949 (2009).
% 
% \bibitem{Arques_84}  D. Arqu\`es and J. Fran\c{c}on,
%    Arbres bien \'etiquet\'es et fractions multicontinues,
%    in {\em Ninth Colloquium on Trees in Algebra and Programming}\/,
%    edited by B.~Courcelle (Cambridge University Press, Cambridge, 1984),
%    pp.~51--61.
% 
% \bibitem{Askey_85}  R. Askey, Appendix to
%    ``Chapter 12 of Ramanujan's second notebook: continued fractions'',
%    Rocky Mountain J. Math. {\bf 15}, 311--318 (1985).
% 
% \bibitem{Asner_70}  B.A. Asner, Jr., On the total nonnegativity of the
%    Hurwitz matrix, SIAM J. Appl. Math. {\bf 18}, 407--414 (1970).
% 
% \bibitem{Aval_08}  J.-C. Aval, Multivariate Fuss-Catalan numbers,
%    Discrete Math. {\bf 308}, 4660--4669 (2008).
% 
% \bibitem{Bailey_35}  W.N. Bailey, {\em Generalized Hypergeometric Series}\/
%    (Cambridge University Press, Cambridge, 1935).
%    Reprinted by Stechert-Hafner, New York, 1964.
% 
% %% \bibitem{Bailey_37}  W.N. Bailey, Associated hypergeometric series,
% %%    Quart. J. Math. (Oxford) {\bf 8}, 115--118 (1937).
% 
% \bibitem{Banderier_02}  C. Banderier and P. Flajolet,
%    Basic analytic combinatorics of directed lattice paths,
%    Theoret. Comput. Sci. {\bf 281}, 37--80 (2002).
% 
% \bibitem{Barbero_15}  J.F. Barbero G., J. Salas and E.J.S. Villase\~nor,
%    Generalized Stirling permutations and forests: higher-order Eulerian
%    and Ward numbers, Electron. J. Combin. {\bf 22}, no.~3, paper 3.37 (2015).

\bibitem{Barsky_81}  D. Barsky, 
   Congruences pour les nombres de Genocchi de 2e esp\`ece
   [extrait d'un travail en commun avec Dominique Dumont],
   Groupe d'\'etude d'Analyse ultram\'etrique, 8e ann\'ee (1980/81),
   Expos\'e no.~34, 13~pp.
   % http://www.numdam.org/item/GAU_1979-1981__7-8__A15_0

% \bibitem{Barry_16}  P. Barry, {\em Riordan Arrays: A Primer}\/
%    (Logic Press, County Kildare, Ireland, 2016).
% 
% \bibitem{Bergeron_12}  F. Bergeron, Combinatorics of $r$-Dyck paths,
%    $r$-parking functions, and the $r$-Tamari lattices,
%    arXiv:1202.6269 [math.CO].
% 
% \bibitem{Bergeron_92}  F. Bergeron, P. Flajolet and B. Salvy,
%    Varieties of increasing trees,
%    in {\em CAAP '92}\/, edited by J.-C.~Raoult,
%    Lecture Notes in Computer Science \#581
%    (Springer-Verlag, Berlin, 1992), pp.~24--48.
% 
% \bibitem{Bergeron_98}  F. Bergeron, G. Labelle and P. Leroux,
%       {\em Combinatorial Species and Tree-Like Structures}\/
%       (Cambridge University Press, Cambridge--New York, 1998).
% 
% \bibitem{Berndt_85}  B.C. Berndt, R.L. Lamphere and B.M. Wilson,
%    Chapter 12 of Ramanujan's second notebook: Continued fractions,
%    Rocky Mountain J. Math. {\bf 15}, 235--310 (1985).

\bibitem{Biane_93}  P. Biane, Permutations suivant le type d'exc\'edance
   et le nombre d'inversions et interpr\'etation combinatoire d'une fraction
   continue de Heine,
   European J. Combin. {\bf 14}, 277--284 (1993).

% \bibitem{Borodin_17}  A. Borodin and G. Olshanski,
%    {\em Representations of the Infinite Symmetric Group}\/
%    (Cambridge University Press, Cambridge, 2017).
% 
% \bibitem{Brenti_89}  F. Brenti,
%    Unimodal, log-concave and P\'olya frequency sequences in combinatorics,
%    Mem. Amer. Math. Soc. {\bf 81}, no.~413 (1989).
% 
% \bibitem{Brenti_95}  F. Brenti, Combinatorics and total positivity,
%    J. Combin. Theory A {\bf 71}, 175--218 (1995).
% 
% \bibitem{Brenti_96}  F. Brenti, The applications of total positivity
%    to combinatorics, and conversely,
%    in {\em Total Positivity and its Applications}\/,
%    edited by M.~Gasca and C.A.~Micchelli
%    (Kluwer, Dordrecht, 1996), pp.~451--473.
% 
% \bibitem{Brenti_98}  F. Brenti, Hilbert polynomials in combinatorics,
%    J. Algebraic Combin. {\bf 7}, 127--156 (1998).
% 
% \bibitem{Brumfiel_79}  G.W. Brumfiel,
%    {\em Partially Ordered Rings and Semi-Algebraic Geometry}\/,
%    London Mathematical Society Lecture Note Series \#37
%    (Cambridge University Press, Cambridge--New York, 1979).

\bibitem{Burstein_06}  A. Burstein, S. Elizalde and T. Mansour,
   Restricted Dumont permutations, Dyck paths, and noncrossing partitions,
   Disc. Math. {\bf 306}, 2851--2869 (2006).

\bibitem{Burstein_21}  A. Burstein and O. Jones,
   Enumeration of Dumont permutations avoiding certain four-letter patterns,
   Discrete Math. Theor. Comput. Sci. {\bf 22}, no.~2, paper \#7 (2021).

% % \bibitem{Bytev_10}  V.V. Bytev, M.Yu. Kalmykov and B.A. Kniehl,
% %    Differential reduction of generalized hypergeometric functions
% %    from Feynman diagrams: one-variable case,
% %    Nuclear Phys. B {\bf 836}, 129--170 (2010).
% 
% \bibitem{Cameron_16}  N.T. Cameron and J.E. McLeod,
%    Returns and hills on generalized Dyck paths,
%    J. Integer Seq. {\bf 19}, article 16.6.1 (2016), 28~pp. 
% 

% \bibitem{Carlitz_72}  L. Carlitz, A conjecture concerning Genocchi numbers,
%    K. Norske Vidensk. Selsk. Skr. (Trondheim) {\bf 9}, 1--4 (1972).

% \bibitem{Carlitz_73}  L. Carlitz, Eulerian numbers and operators,
%    Collectanea Math. {\bf 24}, 175--200 (1973).
% 
% \bibitem{Chang_16}  X.-K. Chang, X.-B. Hu, H. Lei and Y.-N. Yeh,
%    Combinatorial proofs of addition formulas,
%    Electron. J. Combin. {\bf 23}, no.~1, \#P1.8 (2016).
% 
% \bibitem{Chen_15a}  X. Chen, H. Liang and Y. Wang,
%    Total positivity of Riordan arrays,
%    European J. Combin. {\bf 46}, 68--74 (2015).
% 
% \bibitem{Chen_15b}  X. Chen, H. Liang and Y. Wang,
%    Total positivity of recursive matrices,
%    Lin. Alg. Appl. {\bf 471}, 383--393 (2015).
%    %% See letter_from_wang.Jan_4_15
% 
% \bibitem{Cigler_87}  J. Cigler, Some remarks on Catalan families,
%    European J. Combin. {\bf 8}, 261--267 (1987).
% 
% \bibitem{Curtis_98}  E.B. Curtis, D. Ingerman and J.A. Morrow,
%    Circular planar graphs and resistor networks,
%    Lin. Alg. Appl. {\bf 283}, 115--150 (1998).
% 
% \bibitem{Cuyt_08}  A. Cuyt, V.B. Petersen, B. Verdonk, H. Waadeland
%    and W.B. Jones,
%    {\em Handbook of Continued Fractions for Special Functions}\/
%    (Springer-Verlag, New York, 2008).

% \bibitem{DeMedicis_93}  A. de M\'edicis and P. Leroux,
% A unified combinatorial approach for $q$- (and $p,q$-) Stirling numbers,
% J. Statist. Plann. Inference {\bf 34}, 89--105 (1993);
% erratum {\bf 35}, 267 (1993).

% \bibitem{Deutsch_05}  E. Deutsch, L. Ferrari and S. Rinaldi,
%    Production matrices,  Adv. Appl. Math. {\bf 34}, 101--122 (2005).
% 
% \bibitem{Deutsch_09}  E. Deutsch, L. Ferrari and S. Rinaldi,
%    Production matrices and Riordan arrays,
%    Ann. Comb. {\bf 13}, 65--85 (2009).

% \bibitem{Domaratzki_04}  M. Domaratzki, 
%    Combinatorial interpretations of a generalization of the Genocchi numbers,
%    J. Integer Seq. {\bf 7}, no.~3, article 04.3.6 (2004).

% \bibitem{Duchon_00}  P. Duchon, On the enumeration and generation
%    of generalized Dyck words,
%    Discrete Math. {\bf 225}, 121--135 (2000).

% \bibitem{Dumont_72}  D. Dumont,
%    Sur une conjecture de Gandhi concernant les nombres de Genocchi,
%    Discrete Math. {\bf 1}, 321--327 (1972).

\bibitem{Dumont_74}  D. Dumont,
   Interpr\'etations combinatoires des nombres de Genocchi,
   Duke Math. J. {\bf 41}, 305--318 (1974).

% \bibitem{Dumont_80}  D. Dumont, Une g\'en\'eralisation trivari\'ee
%    sym\'etrique des nombres eul\'eriens,
%    J. Combin. Theory A {\bf 28}, 307--320 (1980).

\bibitem{Dumont_86}  D. Dumont, Pics de cycle et d\'eriv\'ees partielles,
   S\'eminaire Lotharingien de Combinatoire {\bf 13}, article B13a (1986).

\bibitem{Dumont_95}  D. Dumont, Further triangles of Seidel--Arnold type
   and continued fractions related to Euler and Springer numbers,
   Adv. Appl. Math. {\bf 16}, 275--296 (1995).

% \bibitem{Dumont_95c}  D. Dumont, Conjectures sur des sym\'etries ternaires
%   li\'ees aux nombres de Genocchi,
%   Discrete Math. {\bf 139}, 469--472 (1995).

% \bibitem{Dumont_76}  D. Dumont and D. Foata,
%    Une propri\'et\'e de sym\'etrie des nombres de Genocchi,
%    Bull. Soc. Math. France {\bf 104}, 433--451 (1976).

% \bibitem{Dumont_88}  D. Dumont and G. Kreweras,
%    Sur le d\'eveloppement d'une fraction continue li\'ee \`a la s\'erie
%    hyperg\'eom\'etrique et son interpr\'etation en termes de records et
%    anti-records dans les permutations,
%    European J. Combin. {\bf 9}, 27--32 (1988).

\bibitem{Dumont_94}  D. Dumont and A. Randrianarivony,
   D\'erangements et nombres de Genocchi,
   Discrete Math. {\bf 132}, 37--49 (1994).

\bibitem{Dumont_95b}  D. Dumont and A. Randrianarivony,
   Sur une extension des nombres de Genocchi,
   European J. Combin. {\bf 16}, 147--151 (1995).

\bibitem{Dumont_80b}  D. Dumont and G. Viennot,
   A combinatorial interpretation of the Seidel generation of Genocchi numbers,
   Ann. Discrete Math. {\bf 6}, 77--87 (1980).

\bibitem{Dumont_94b}  D. Dumont and J. Zeng,
   Further results on the Euler and Genocchi numbers,
   Aequationes Math. {\bf 47}, 31--42 (1994).

% \bibitem{Dumont_98}  D. Dumont and J. Zeng,
%    Polyn\^omes d'Euler et fractions continues de Stieltjes-Rogers,
%    Ramanujan J. {\bf 2}, 387--410 (1998).

% \bibitem{Dyachenko_14}  A. Dyachenko, Total nonnegativity of infinite
%    Hurwitz matrices of entire and meromorphic functions,
%    Complex Anal. Oper. Theory {\bf 8}, 1097--1127 (2014).
% 
% \bibitem{Dzhumadildaev_14}  A. Dzhumadil'daev and D. Yeliussizov,
%    Stirling permutations on multisets,
%    European J. Combin. {\bf 36}, 377--392 (2014).
% 
% \bibitem{Edelman_80}  P.H. Edelman, Chain enumeration and noncrossing
%    partitions, Discrete Math. {\bf 31}, 171--180 (1980).

\bibitem{Ehrenborg_00b}   R. Ehrenborg and E. Steingr\'{\i}msson,
   Yet another triangle for the Genocchi numbers,
   European J. Combin. {\bf 21}, 593--600 (2000).

\bibitem{Elvey-Price-Sokal_wardpoly}  A. Elvey Price and A.D. Sokal,
   Phylogenetic trees, augmented perfect matchings,
   and a Thron-type continued fraction (T-fraction) for the Ward polynomials,
   %% arXiv:2001.01468 [math.CO]
   Electron. J. Combin. {\bf 27}(4), article P4.6 (2020).

\bibitem{Enestrom_13}  G. Enestr\"om,
   {\em Die Schriften Eulers chronologisch nach den Jahren geordnet,
    in denen sie verfa{\ss}t worden sind}\/,
   Jahresbericht der Deutschen Mathematiker-Vereinigung
   (Teubner, Leipzig, 1913).

\bibitem{Eu_21}  S.-P. Eu, T.-S. Fu, H.-H. Lai and Y.-H. Lo,
   Gamma-positivity for a refinement of median Genocchi numbers,
   arXiv:2103.09130 [math.CO].

\bibitem{Euler_1755}  L. Euler, Methodis summandi superior ulterius promota,
   Chapter~7 of
   {\em Institutiones Calculi Differentialis cum eius Usu in Analysi
   Finitorum ac Doctrina Serierum}\/
   [Foundations of Differential Calculus, with Applications to
    Finite Analysis and Series]
   (Academiae Imperialis Scientiarum Petropolitanae, Saint Petersburg, 1755),
   pp.~479--514.
   %% https://books.google.co.uk/books?id=sYE_AAAAcAAJ&printsec=frontcover&source=gbs_ge_summary_r&cad=0#v=onepage&q&f=false
   Reprinted in {\em Opera Omnia}\/, ser.~1, vol.~10, pp.~368--395.
   Latin original available at
    \url{http://eulerarchive.maa.org/pages/E212.html};
    English translation available at
    \url{https://www.agtz.mathematik.uni-mainz.de/algebraische-geometrie/van-straten/euler-kreis-mainz/}

\bibitem{Euler_1760}  L. Euler, De seriebus divergentibus,
   Novi Commentarii Academiae Scientiarum Petropolitanae {\bf 5}, 205--237
   (1760);
   reprinted in {\em Opera Omnia}\/, ser.~1, vol.~14, pp.~585--617.
   [Latin original and English and German translations available at
    \url{http://eulerarchive.maa.org/pages/E247.html}]

\bibitem{Euler_1788}  L. Euler, De transformatione seriei divergentis
   $1 - mx + m(m+n)x^2 - m(m+n)(m+2n)x^3 + \hbox{etc.}$
   in fractionem continuam,
   Nova Acta Academiae Scientarum Imperialis Petropolitanae
     {\bf 2}, 36--45 (1788);
   reprinted in {\em Opera Omnia}\/, ser.~1, vol.~16, pp.~34--46.
   %% Reference copied from Dumont and Kreweras \cite{Dumont_88}
   %%   and then corrected according to arXiv:1201.6687
   [Latin original and English and German translations available at
    \url{http://eulerarchive.maa.org/pages/E616.html}]
   %% [German translation available at http://arxiv.org/abs/1201.6687]

% \bibitem{Fallat_11}  S.M. Fallat and C.R. Johnson,
%    {\em Totally Nonnegative Matrices}\/
%    (Princeton University Press, Princeton NJ, 2011).

\bibitem{Feigin_12}  E. Feigin, The median Genocchi numbers,
   $q$-analogues and continued fractions,
   European J. Combin. {\bf 33}, 1913--1918 (2012).

\bibitem{Flajolet_80}  P. Flajolet, Combinatorial aspects of continued
   fractions,  Discrete Math. {\bf 32}, 125--161 (1980).

% \bibitem{Flajolet_82}  P. Flajolet, On congruences and continued fractions
%    for some classical combinatorial quantities,
%    Discrete Math. {\bf 41}, 145--153 (1982).

% \bibitem{Flajolet_09}  P. Flajolet and R. Sedgewick,
%   {\em Analytic Combinatorics}\/ (Cambridge University Press, Cambridge, 2009).

\bibitem{Foata_08}  D. Foata and G.-N. Han,
   {\em Principes de combinatoire classique: Cours et exercices corrig\'es}\/
   (Universit\'e Louis Pasteur, Strasbourg, D\'epartement de math\'ematique,
    2008).
   Available on-line at
   \url{https://irma.math.unistra.fr/~foata/AlgComb.pdf}

\bibitem{Foata_71}  D. Foata and M.-P. Sch\"utzenberger,
   Nombres d'Euler et permutations alternantes,
   Technical report, University of Florida, 1971, 71 pp.
   Available on-line at
   \url{http://www.emis.de/journals/SLC/books/foaschuetz1.html}

\bibitem{Foata_90}  D. Foata and D. Zeilberger,
   Denert's permutation statistic is indeed Euler-Mahonian,
   Stud. Appl. Math. {\bf 83}, 31--59 (1990).

% \bibitem{Fomin_01}  S. Fomin, Loop-erased walks and total positivity,
%    Trans. Amer. Math. Soc. {\bf 353}, 3563--3583 (2001).
% 
% \bibitem{Fomin_10}  S. Fomin, Total positivity and cluster algebras,
%    in {\em Proceedings of the International Congress of Mathematicians}\/,
%    vol.~II,
%    edited by R.~Bhatia, A.~Pal, G.~Rangarajan, V.~Srinivas and M.~Vanninathan
%    (Hindustan Book Agency, New Delhi, 2010), pp.~125--145.
% 
% \bibitem{Fomin_forthcoming}  S. Fomin, L. Williams and A. Zelevinsky,
%    {\em Introduction to Cluster Algebras}\/,
%    forthcoming book;
%    preliminary draft of Chapters~1--5 posted at
%    arXiv:1608.05735 [math.CO] and arXiv:1707.07190 [math.CO]
%    at arXiv.org.
% 
% \bibitem{Fomin_99}  S. Fomin and A. Zelevinsky,
%    Double Bruhat cells and total positivity,
%    J. Amer. Math. Soc. {\bf 12}, 335--380 (1999).
% 
% \bibitem{Fomin_00}  S. Fomin and A. Zelevinsky,
%    Total positivity: tests and parametrizations,
%    Math. Intelligencer {\bf 22}, no.~1, 23--33 (2000).
% 
% \bibitem{Francon_78}  J. Fran\c{c}on, Histoires de fichiers,
%    RAIRO Informat. Th\'eor. {\bf 12}, 49--62 (1978).
% 
% \bibitem{Frank_56}  E. Frank,
%    A new class of continued fraction expansions for the ratios of
%    hypergeometric functions,
%    Trans. Amer. Math. Soc. {\bf 81}, 453--476 (1956).

\bibitem{Fusy_15}  E. Fusy and E. Guitter,
   Comparing two statistical ensembles of quadrangulations:
   A continued fraction approach,
   %% preprint, arXiv:1507.04538 [math.CO] at arXiv.org,
   Ann. Inst. Henri Poincar\'e D {\bf 4}, 125--176 (2017).

% \bibitem{Gandhi_70}  J.M. Gandhi, 
%    A conjectured representation of Genocchi numbers,
%    Amer. Math. Monthly {\bf 77}, 505--506 (1970).

% \bibitem{Gantmacher_02}  F.R. Gantmacher and M.G. Krein,
%    {\em Oscillation Matrices and Kernels and Small Vibrations of
%        Mechanical Systems}\/
%    %% edited and with a preface by Alex Eremenko
%    (AMS Chelsea Publishing, Providence RI, 2002).
%    Based on the second Russian edition, 1950.
% 
% \bibitem{Gantmakher_37}  F. Gantmakher and M. Krein,
%    Sur les matrices compl\`etement non n\'egatives et oscillatoires,
%    Compositio Math. {\bf 4}, 445--476 (1937).
%    %% Available at archive.numdam.org/article/CM_1937__4__445_0.pdf
% 
% \bibitem{Gasca_96}  M. Gasca and C.A. Micchelli, eds.,
%    {\em Total Positivity and its Applications}\/
%    (Kluwer, Dordrecht, 1996).
% 
% \bibitem{Gasper_04}  G. Gasper and M. Rahman,
%    {\em Basic Hypergeometric Series}\/, 2nd ed.
%    (Cambridge University Press, Cambridge--New York, 2004).
%    %% Available on-line at
%    %% http://ebooks.cambridge.org/chapter.jsf?bid=CBO9780511526251
% 
% \bibitem{Gauss_1813}  C.F. Gauss, Disquisitiones generales circa seriem
%    infinitam $1 + {\alpha \beta \over 1 . \gamma} x
%       + {\alpha(\alpha+1)\beta(\beta+1) \over 1 . 2 . \gamma (\gamma+1)} xx
%       + {\alpha(\alpha+1)(\alpha+2)\beta(\beta+1)(\beta+2) \over 1 . 2 . 3 . \gamma (\gamma+1) (\gamma+2)} x^3 + \hbox{etc.}$,
%    Commentationes Societatis Regiae Scientiarum Gottingensis Recentiores,
%    Classis Mathematicae {\bf 2} (1813).
%    [Reprinted in C.F. Gauss, {\em Werke}\/, vol.~3
%    (Cambridge University Press, Cambridge, 2011),
%    pp.~123--162.]
%    % https://www.cambridge.org/core/books/werke/0075304640500FEF9DA26221FCF83EC4

\bibitem{Genocchi_1852}  A. Genocchi,
   Intorno all'espressione generale de' numeri bernulliani,
   Annali di Scienze Mathematiche e Fisiche (Roma)
   {\bf 3}, 395--405 (1852).
   %% https://books.google.co.uk/books?id=Y_A3AAAAMAAJ&printsec=frontcover&source=gbs_ge_summary_r&cad=0#v=onepage&q&f=false

% \bibitem{Gessel_78a} I. Gessel, A note on Stirling permutations,
%    unpublished manuscript, August~1978, cited in \cite{Park_94a}.
% 
% \bibitem{Gessel_private}  I. Gessel, private communication, 27~May 2015.
%    %% letter_from_gessel.May_27_15
% 
% \bibitem{Gessel_16}  I.M. Gessel, Lagrange inversion,
%    J. Combin. Theory A {\bf 144}, 212--249 (2016).
% 
% \bibitem{Gessel-Seo_06}  I.M. Gessel and S. Seo,
%    A refinement of Cayley's formula for trees,
%    Electron. J. Combin. {\bf 11}, no.~2, \#R27 (2006).
% 
% \bibitem{Gessel_78}  I. Gessel and R.P. Stanley, Stirling polynomials,
%    J. Combin. Theory A {\bf 24}, 24--33 (1978).
% 
% \bibitem{Gessel-Viennot_89}  I.M. Gessel and X.G. Viennot,
%    Determinants, paths, and plane partitions,
%    Brandeis University preprint (1989).
%    Available on-line at
%    \url{http://people.brandeis.edu/~gessel/homepage/papers/} or
%    \url{http://xavierviennot.org/xavier/articles.html}
% 
% %% \bibitem{Goldschmidt_16}  C. Goldschmidt,
% %%   A short introduction to random trees,
% %%   Mongolian Math. J. {\bf 20}, 51--72 (2016).
% 
% \bibitem{Goulden_83}  I.P. Goulden and D.M. Jackson,
%    {\em Combinatorial Enumeration}\/ (Wiley, New York, 1983).
%    Reprinted by Dover, Mineola NY, 2004.
% 
% \bibitem{Gouyou-Beauchamps_98}  D. Gouyou-Beauchamps,
%    Construction of $q$-equations for convex polyominoes,
%    paper presented at 10th International Conference on
%    Formal Power Series and Algebraic Combinatorics (FPSAC '98),
%    available on-line at
%    \url{http://www-igm.univ-mlv.fr/~fpsac/FPSAC98/articles.html}
% 
% \bibitem{Graham_94}  R.L. Graham, D.E. Knuth and O. Patashnik,
%    {\em Concrete Mathematics: A Foundation for Computer Science}\/,
%    2nd ed.~(Addison-Wesley, Reading, Mass., 1994).
% 
% \bibitem{Hackl_18}  B. Hackl, C. Heuberger and H. Prodinger,
%    Ascents in non-negative lattice paths,
%    preprint, arXiv:1801.02996 [math.CO] at arXiv.org.
% 
% \bibitem{Haglund_12}  J. Haglund and M. Visontai,
%    Stable multivariate Eulerian polynomials and generalized Stirling
%    permutations, European J. Combin. {\bf 33}, 477--487 (2012).
% 
% \bibitem{Haiman_94}  M.D. Haiman, Conjectures on the quotient ring by
%    diagonal invariants, J. Algebraic Combin. {\bf 3}, 17--76 (1994).

\bibitem{Han_18}  G.-N. Han and J.-Y. Liu,
   Combinatorial proofs of some properties of tangent and Genocchi numbers,
   European J. Combin. {\bf 71}, 99--110 (2018).

\bibitem{Han_99a}  G.-N. Han and J. Zeng,
   $q$-polyn\^omes de Gandhi et statistique de Denert,
   Discrete Math. {\bf 205}, 119--143 (1999).

\bibitem{Han_99b}  G.-N. Han and J. Zeng,
   On a $q$-sequence that generalizes the median Genocchi numbers,
   Ann. Sci. Math. Qu\'ebec {\bf 23}, 63--72 (1999).

% \bibitem{He_15}  T.-X. He, Matrix characterizations of Riordan arrays,
%    Lin. Alg. Appl. {\bf 465}, 15--42 (2015).
% 
% \bibitem{Heine_1847}  E. Heine, Untersuchungen \"uber die Reihe
%  $\displaystyle
%     1 + {(1-q^\alpha)(1-q^\beta) \over (1-q)(1-q^\gamma)} \cdot x
%     + {(1-q^\alpha)(1-q^{\alpha+1})(1-q^\beta)(1-q^{\beta+1})
%        \over (1-q)(1-q^2)(1-q^\gamma)(1-q^{\gamma+1})} \cdot x^2 + \ldots$,
%  J. reine angew. Math. {\bf 34}, 285--328 (1847).
%  Available on-line at
%  \url{http://www.digizeitschriften.de/main/dms/img/?PPN=GDZPPN002145758}
% 
% \bibitem{Holtz_03}  O. Holtz, Hermite--Biehler, Routh--Hurwitz,
%    and total positivity, Lin. Alg. Appl. {\bf 372}, 105--110 (2003).
% 
% \bibitem{Janson_11}  S. Janson, M. Kuba and A. Panholzer,
%    Generalized Stirling permutations, families of increasing trees
%    and urn models,
%    J. Combin. Theory A {\bf 118}, 94--114 (2011).
% 
% \bibitem{Jones_80}  W.B. Jones and W.J. Thron,
%    {\em Continued Fractions: Analytic Theory and Applications}\/
%    (Addison-Wesley, Reading MA, 1980).

% \bibitem{Josuat-Verges_10}  M. Josuat-Verg\`es,
%    Generalized Dumont--Foata polynomials and alternative tableaux,
%    S\'eminaire Lotharingien de Combinatoire {\bf 64}, article B64b (2010),
%    17~pp. 

\bibitem{Josuat-Verges_18}  M. Josuat-Verg\`es,
   A $q$-analog of Schl\"afli and Gould identities on Stirling numbers,
   %% preprint, arXiv:1610.02965 [math.CO] at arXiv.org.
   Ramanujan J. {\bf 46}, 483--507 (2018).

% \bibitem{Karlin_68}  S. Karlin, {\em Total Positivity}\/
%    (Stanford University Press, Stanford CA, 1968).
% 
% \bibitem{Karlin_59}  S. Karlin and J. McGregor, Coincidence probabilities,
%     Pacific J. Math. {\bf 9}, 1141--1164 (1959).
% 
% \bibitem{Kemperman_82}  J.H.B. Kemperman,
%    A Hurwitz matrix is totally positive,
%    SIAM J. Math. Anal. {\bf 13}, 331--341 (1982).
% 
% \bibitem{Khovanskii_63}  A.N. Khovanskii, {\em The Application of Continued
%    Fractions and their Generalizations to Problems in Approximation Theory}\/,
%    translated from the Russian by P.~Wynn
%    (Noordhoff, Groningen, 1963).

% \bibitem{Kim_00}  D. Kim and J. Zeng,
%    On a continued fraction formula of Wall,
%    Ramanujan J. {\bf 4}, 421--427 (2000).

\bibitem{Knuth_98}  D.E. Knuth, {\em The Art of Computer Programming}\/,
   vol.~3, 2nd ed.\ (Addison-Wesley, Reading MA, 1998).

% \bibitem{Krattenthaler_contiguous}  C. Krattenthaler, A systematic list of
%    two- and three-term contiguous relations for basic hypergeometric series,
%    posted 25~August 2002,
%    \url{http://www.mat.univie.ac.at/~kratt/artikel/contrel.ps.gz}
%    %% Index of directory at http://www.mat.univie.ac.at/~kratt/artikel/
% 
% \bibitem{Krattenthaler_HYP}  C. Krattenthaler, HYP
%    [Manual for the {\sc Mathematica} package HYP],
%    posted 20~September 2003,
%    \url{http://www.mat.univie.ac.at/~kratt/hyp_hypq/hypm.pdf}
%    %% Index of directory at http://www.mat.univie.ac.at/~kratt/hyp_hypq/
% 
% \bibitem{Kreweras_72}  G. Kreweras,
%    Sur les partitions non crois\'ees d'un cycle,
%    Discrete Math. {\bf 1}, 333--350 (1972).
% 
% \bibitem{Kuba_09}  M. Kuba and A.L. Varvak,
%    On path diagrams and Stirling permutations,
%    preprint, arXiv:0906.1672v2 [math.CO] at arXiv.org.

\bibitem{Kuznetsov_94}  A.G. Kuznetsov, I.M. Pak and A.E. Postnikov,
   Increasing trees and alternating permutations,
   Uspekhi Mat. Nauk {\bf 49}, no.~6, 79--110 (1994)
   [= Russian Math. Surveys {\bf 49}, no.~6, 79--114 (1994)].

% \bibitem{Lam_84}  T.Y. Lam, An introduction to real algebra,
%    Rocky Mountain J. Math. {\bf 14}, 767--814 (1984).
% 
% \bibitem{Lambert_1768}  J.H. Lambert, M\'emoire sur quelques propri\'et\'es
%    remarquables des quantit\'es transcendentes circulaires et logarithmiques,
%    M\'emoires de l'Acad\'emie Royale des Sciences de Berlin
%    {\bf 17}, 265--322 (1768).
%    Available on-line at \url{http://www.kuttaka.org/~JHL/L1768b.html}

\bibitem{Lazar_20}  A.L. Lazar,
   The homogenized Linial arrangement and its consequences in enumerative
   combinatorics,
   Ph.D.~thesis, University of Miami (August~2020),
   \url{https://scholarship.miami.edu/discovery/delivery/01UOML_INST:ResearchRepository/12367619000002976?l#13367618990002976}
   %% \url{https://na-st01.ext.exlibrisgroup.com/01UOML_INST/upload/1641907374181_axl416S20.pdf?Expires=1641907494&Signature=dKYl1yGWclT7OzIT8OXsssvT6lQw2FE999ZPPfac08~YVASo~QWpCNlZwq-ZUUBW701NNldzSr6d9R35F9eOoMY~35OCsIV5m-iO9SRznjih58fW6-K~-VZ-iIeqkiwUi~84E7eg7WfZtxsGWYe~hbMBw50myFsMcCby9PkGOsBM8kD9LmDieL7O8Nuk5OZFcyolBRReFnkpEP2HPuKrnKh6jYfwDmPal94PQ5Ai-L3HhhemLi-XE6zKy-nIAwV8cU9IQ2gHFwz~FkMGwoW9k2zkM3e~Pjmjfprxt0l68hct19iovK5ZbabVl23RxQsDdS77tbgl5IscmFBELJhsSQ__&Key-Pair-Id=APKAJ72OZCZ36VGVASIA}

\bibitem{Lazar_22}  A. Lazar and M.L. Wachs,
   The homogenized Linial arrangement and Genocchi numbers,
   %% arXiv:1910.07651 [math.CO],
   Combin. Theory {\bf 2}, issue~1, paper no.~2 (2022), 34 pp.

% %% \bibitem{LeGall_05}  J.-F. Le Gall,  Random trees and applications,
% %%    Probab. Surv. {\bf 2}, 245--311 (2005).
% 
% \bibitem{Lenczewski_13}  R. Lenczewski and R. Sa\l{}apata,
%    Multivariate Fuss-Narayana polynomials and their application to
%    random matrices,
%    Electron. J. Combin. {\bf 20}, no.~2, article \#P41 (2013).

% \bibitem{Leroux_90}  P. Leroux, Reduced matrices and $q$-log-concavity
%    properties of $q$-Stirling numbers,
%    J. Combin. Theory A {\bf 54}, 64--84 (1990).

% \bibitem{Liang_16}  H. Liang, L. Mu and Y. Wang,
%    Catalan-like numbers and Stieltjes moment sequences,
%    Discrete Math. {\bf 339}, 484--488 (2016).
% 
% \bibitem{Lindstrom_73}  B. Lindstr\"om, On the vector representations of
%    induced matroids, Bull. London Math. Soc. {\bf 5}, 85--90 (1973).
% 
% \bibitem{Lis_17}  M. Lis, The planar Ising model and total positivity,
%    J. Stat. Phys. {\bf 166}, 72--89 (2017).
% 
% \bibitem{Lorentzen_92}  L. Lorentzen and H. Waadeland,
%    {\em Continued Fractions with Applications}\/
%    (North-Holland, Amsterdam, 1992).
% 
% \bibitem{Lothaire_97}  M. Lothaire, {\em Combinatorics on Words}\/,
%    2nd ed.\ (Cambridge University Press, Cambridge, 1997).

\bibitem{Lucas_1877}  \'E. Lucas, Sur les th\'eor\`emes de Binet
   et de Staudt concernant les nombres de Bernoulli,
   Nouvelles Annales de Math\'ematiques (2${}^e$ s\'erie)
   {\bf 16}, 157--160 (1877).

% % \bibitem{Luke_69}  Y.L. Luke,
% %    {\em The Special Functions and Their Approximations}\/, vol.~I
% %    (Academic Press, New York--London, 1969).
% 
% \bibitem{Lusztig_94}  G. Lusztig, Total positivity in reductive groups,
%    in {\em Lie Theory and Geometry}\/,
%    edited by J.-L.~Brylinski, R.~Brylinski, V.~Guillemin and V.~Kac
%    (Birkh\"auser Boston, Boston MA, 1994), pp.~531--568.
% 
% \bibitem{Lusztig_98}  G. Lusztig, Introduction to total positivity,
%    in {\em Positivity in Lie Theory: Open Problems}\/,
%    edited by J.~Hilgert, J.D.~Lawson, K.-H.~Neeb and E.B.~Vinberg
%    (de Gruyter, Berlin, 1998), pp.~133--145.
%    %% Available at http://dedekind.mit.edu/~gyuri/papers/pod.ps
% 
% \bibitem{Lusztig_08}  G. Lusztig, A survey of total positivity,
%    Milan J. Math. {\bf 76}, 125--134 (2008).
% 
% \bibitem{Macdonald_95}  I.G. Macdonald, {\em Symmetric Functions and
%    Hall Polynomials}\/, 2nd ed. (Clarendon Press, Oxford, 1995).
% 
% \bibitem{Marshall_08}  M. Marshall, {\em Positive Polynomials and
%     Sums of Squares}\/, Mathematical Surveys and Monographs \#146
%     (American Mathematical Society, Providence RI, 2008).
% 
% \bibitem{Mlotkowski_12}  W. M\l{}otkowski,
%    Probability measures corresponding to Aval numbers,
%    Colloq. Math. {\bf 129}, 189--202 (2012).
% 
% \bibitem{Mu_17a}  L. Mu, J. Mao and Y. Wang,
%    Row polynomial matrices of Riordan arrays,
%    Lin. Alg. Appl. {\bf 522}, 1--14 (2017).
% 
% \bibitem{Nica_06}  A. Nica and R. Speicher,
%    {\em Lectures on the Combinatorics of Free Probability}\/,
%    London Mathematical Society Lecture Note Series \#335
%    (Cambridge University Press, Cambridge, 2006).
%    %% Available on-line at
%    %% http://dx.doi.org/10.1017/CBO9780511735127

\bibitem{NIST}  F.W. Olver, D.W. Lozier, R.F. Boisvert and C.W. Clark, eds.,
   {\em NIST Handbook of Mathematical Functions}\/
   (U.S.~Department of Commerce, Washington DC, 2010).
   Available on-line at \url{http://dlmf.nist.gov}

\bibitem{OEIS}  The On-Line Encyclopedia of Integer Sequences,
   published electronically at \url{http://oeis.org}

\bibitem{Oste_15}  R. Oste and J. Van der Jeugt,
   Motzkin paths, Motzkin polynomials and recurrence relations,
   Electron. J. Combin. {\bf 22}, no.~2, \#P2.8 (2015).

\bibitem{Pan_21}  Q. Pan and J. Zeng,
   Cycles of even-odd drop permutations and continued fractions of
   Genocchi numbers,
   arXiv:2108.03200 [math.CO].

% \bibitem{Pan_15}  R. Pan, Algorithmic solution to Problem~1
%    (and linear extensions of general one-level grid-like posets),
%    \url{http://www.math.ucsd.edu/~projectp/problems/p1.html}
%    (24~February 2015)
%    and
%    \url{http://www.math.ucsd.edu/~projectp/problems/solutions/OneLevelGridPoset.pdf}
%    (28~June 2016).

% \bibitem{Park_94a}  S.-K. Park, The $r$-multipermutations,
%    J. Combin. Theory A {\bf 67}, 44--71 (1994).
% 
% \bibitem{Park_94b}  S.-K. Park, Inverse descents of $r$-multipermutations,
%    Discrete Math. {\bf 132}, 215--229 (1994).
% 
% \bibitem{Peart_00}  P. Peart and W.-J. Woan,
%    Generating functions via Hankel and Stieltjes matrices,
%    J. Integer Seq. {\bf 3}, article 00.2.1 (2000).
% 
% \bibitem{Perron}  O. Perron, {\em Die Lehre von den Kettenbr\"uchen}\/
%    (Teubner, Leipzig, 1913).
%    %% Available on-line at https://archive.org/details/dielehrevondenk00perrgoog
%    Second edition: Teubner, Leipzig, 1929; reprinted by Chelsea, New York, 1950.
%    Third edition, 2 vols.: Teubner, Stuttgart, 1954, 1957.
%    %% Third edition available on-line at
%    %% http://link.springer.com/book/10.1007/978-3-663-12289-0
%    %% http://link.springer.com/book/10.1007%2F978-3-663-01496-6
%    %% but one has to pay
% 
% \bibitem{Petersen_15}  T.K. Petersen, {\em Eulerian Numbers}\/
%    (Birkh\"auser--Springer, New York--Heidelberg, 2015).
% 
% \bibitem{Petreolle-Sokal_contiguous}  M. P\'etr\'eolle and A.D. Sokal,
%    Three-term contiguous relations for all generalized hypergeometric
%    functions ${}_r F_s$ and basic hypergeometric functions ${}_r \phi_s$,
%    in preparation.

\bibitem{latpath_SRTR}  M. P\'etr\'eolle, A.D. Sokal and B.-X. Zhu,
   Lattice paths and branched continued fractions:
       An infinite sequence of generalizations
       of the Stieltjes--Rogers and Thron--Rogers polynomials,
       with coefficientwise Hankel-total positivity,
   preprint (2018), arXiv:1807.03271 [math.CO] at arXiv.org,
   to appear in the Memoirs of the American Mathematical Society.

% \bibitem{Pinkus_10}  A. Pinkus, {\em Totally Positive Matrices}\/
%    (Cambridge University Press, Cambridge, 2010).
%    %% Available on-line at
%    %% http://ebooks.cambridge.org/ebook.jsf?bid=CBO9780511691713
% 
% \bibitem{Pitman_06}  J. Pitman, {\em Combinatorial Stochastic Processes}\/
%    (Ecole d'Et\'e de Probabilit\'es de Saint-Flour XXIII -- 2002),
%    Lecture Notes in Mathematics \#1875 (Springer-Verlag, Berlin, 2006).
% 
% \bibitem{Prestel_01}  A. Prestel and C.N. Delzell,
%    {\em Positive Polynomials: From Hilbert's 17th Problem to Real Algebra}\/
%    (Springer-Verlag, Berlin, 2001).
% 
% \bibitem{Prodinger_16}  H. Prodinger, Returns, hills, and $t$-ary trees,
%    J. Integer Seq. {\bf 19}, article 16.7.2 (2016), 8~pp.
% 
% \bibitem{Rainville_60}  E.D. Rainville, {\em Special Functions}\/
%    (Macmillan, New York, 1960).
% 
% \bibitem{Ramanathan_87}  K.G. Ramanathan,
%     Hypergeometric series and continued fractions,
%     Proc. Indian Acad. Sci. (Math. Sci.) {\bf 97}, 277--296 (1987).

% \bibitem{Randrianarivony_94}  A. Randrianarivony,
%    Polyn\^omes de Dumont-Foata g\'en\'eralis\'es,
%    S\'eminaire Lotharingien de Combinatoire {\bf 32}, article B32d (1994),
%    12~pp.

\bibitem{Randrianarivony_97}  A. Randrianarivony,
   Fractions continues, $q$-nombres de Catalan et $q$-polyn\^omes
   de Genocchi,
   Europ. J. Combin. {\bf 18}, 75--92 (1997).

% \bibitem{Randrianarivony_94b}  A. Randrianarivony and J. Zeng,
%    Sur une extension des nombres d'Euler et les records des permutations
%    alternantes,
%    J. Combin. Theory A {\bf 68}, 86--99 (1994).

\bibitem{Randrianarivony_96}  A. Randrianarivony and J. Zeng,
   Une famille de polyn\^omes qui interpole plusieurs suites classiques
   de nombres,
   Adv. Appl. Math. {\bf 17}, 1--26 (1996).

\bibitem{Randrianarivony_96b}  A. Randrianarivony and J. Zeng,
   Some equidistributed statistics on Genocchi permutations,
   Electron. J. Combin. {\bf 3}, no. 2, Research Paper \#22 (1996), 11 pp.

% \bibitem{Rattan_14}  A. Rattan, Parking functions and related
%    combinatorial structures, Ph.D.~thesis, University of Waterloo (2014).

% \bibitem{Riordan_73}  J. Riordan and P.R. Stein,
%    Proof of a conjecture on Genocchi numbers,
%    Discrete Math. {\bf 5}, 381--388 (1973).

% \bibitem{Roblet_94}  E. Roblet,
%    Une interpr\'etation combinatoire des approximants de Pad\'e,
%    Th\`ese de doctorat,
%    Universit\'e Bordeaux~I (1994).
%    Reprinted as Publications du Laboratoire de Combinatoire
%    et d'Informatique Math\'ematique (LACIM) \#17,
%    Universit\'e du Qu\'ebec \`a Montr\'eal (1994).
%    Available on-line at \url{http://lacim.uqam.ca/en/les-parutions/}
% 
% \bibitem{Roblet_96}  E. Roblet and X.G. Viennot,
%    Th\'eorie combinatoire des T-fractions et approximants de Pad\'e
%    en deux points,
%    Discrete Math. {\bf 153}, 271--288 (1996).

\bibitem{Rogers_07}  L.J. Rogers, On the representation of certain asymptotic
   series as convergent continued fractions,
   Proc. London Math. Soc. (series 2) {\bf 4}, 72--89 (1907).

% \bibitem{Roman_15}  S. Roman, {\em An Introduction to Catalan Numbers}\/
%    (Birkh\"auser--Springer, Cham--Heidelberg--New York, 2015).
% 
% \bibitem{Salvy_94}  B. Salvy and P. Zimmermann,
%    Gfun: a Maple package for the manipulation of generating and holonomic
%    functions in one variable,
%    ACM Trans. Math. Software {\bf 20}, 163--177 (1994).
% 
% \bibitem{Schmudgen_17}  K. Schm\"udgen, {\em The Moment Problem}\/
%    (Springer, Cham, 2017).
% 
% \bibitem{Schoenberg_53}  I.J. Schoenberg and A. Whitney,
%    On P\'olya frequency functions. III. The positivity of translation
%    determinants with an application to the interpolation problem
%    by spline curves,
%    Trans. Amer. Math. Soc. {\bf 74}, 246--259 (1953).
% 
% \bibitem{Shapiro_91}  L.W. Shapiro, S. Getu, W.J. Woan and L.C. Woodson,
%    The Riordan group,
%    Discrete Appl. Math. {\bf 34}, 229--239 (1991).
% 
% \bibitem{Shareshian_10}  J. Shareshian and M.L. Wachs,
%    Eulerian quasisymmetric functions,
%    Adv. Math. {\bf 225}, 2921--2966 (2010).
% 
% \bibitem{Shohat_43}  J.A. Shohat and J.D. Tamarkin,
%    {\em The Problem of Moments}\/
%    (American Mathematical Society, New York, 1943).
% 
% \bibitem{Simion_00a}  R. Simion, Noncrossing partitions,
%    Discrete Math. {\bf 217}, 367--409 (2000).
% 
% \bibitem{Simon_98}  B. Simon, The classical moment problem as a
%    self-adjoint finite difference operator,
%    Adv. Math. {\bf 137}, 82--203 (1998).
% 
% \bibitem{Skandera_03}  M. Skandera, Introductory notes on total positivity
%    (June 2003), available at
%    \url{http://www.math.lsa.umich.edu/~fomin/565/intp.ps}
% 
% \bibitem{Slater_66}  L.J. Slater, {\em Generalized Hypergeometric Functions}\/
%    (Cambridge University Press, Cambridge, 1966).

% \bibitem{Sokal_flajolet}  A.D. Sokal, Coefficientwise total positivity
%    (via continued fractions) for some Hankel matrices of combinatorial
%    polynomials, talk at the S\'eminaire de Combinatoire Philippe Flajolet,
%    Institut Henri Poincar\'e, Paris, 5 June 2014;
%    transparencies available at
%    \url{http://semflajolet.math.cnrs.fr/index.php/Main/2013-2014}

\bibitem{Sokal_eulernumbers}  A.D. Sokal,
  The Euler and Springer numbers as moment sequences,
  Expositiones Mathematicae {\bf 38}, 1--26 (2020).
  %% arXiv:1804.04498 [math.CO]
  %% https://doi.org/10.1016/j.exmath.2018.08.001

\bibitem{Sokal_totalpos}  A.D. Sokal, Coefficientwise total positivity
   (via continued fractions) for some Hankel matrices of combinatorial
   polynomials, in preparation.

% \bibitem{Sokal_unpub_2015}  A.D. Sokal, unpublished, May~2015.
% 
% \bibitem{Sokal_alg_contfrac}  A.D. Sokal, A simple algorithm for expanding
%    a power series as a continued fraction, in preparation.

\bibitem{Sokal-Zeng_masterpoly}  A.D. Sokal and J. Zeng,
    Some multivariate master polynomials
    for permutations, set partitions, and perfect matchings,
    and their continued fractions,
    %% preprint (March~2020), arXiv:2003.08192 [math.CO] at arXiv.org.
    Adv. Appl. Math. {\bf 138}, 102341 (2022).

% \bibitem{Song_05}  C. Song, The generalized Schr\"oder theory,
%    Electron. J. Combin. {\bf 12}, paper \#R53 (2005), 10~pp. 
% 
% \bibitem{Speicher_94}  R. Speicher, 
%    Multiplicative functions on the lattice of noncrossing partitions
%    and free convolution, Math. Ann. {\bf 298}, 611--628 (1994).
% 
% \bibitem{Sprugnoli_94}  R. Sprugnoli, Riordan arrays and combinatorial sums,
%    Discrete Math. {\bf 132}, 267--290 (1994).

% \bibitem{Stanley_86}  R.P. Stanley, {\em Enumerative Combinatorics}\/,
%       vol. 1 (Wadsworth \& Brooks/Cole, Monterey, California, 1986).
%       Reprinted by Cambridge University Press, 1999.

% \bibitem{Stanley_97}  R.P. Stanley, Parking functions and noncrossing
%    partitions, Electron. J. Combin. {\bf 4}, no.~2, paper \#R20 (1997), 14 pp.

% \bibitem{Stanley_99}  R.P. Stanley, {\em Enumerative Combinatorics}\/,
%       vol.~2 (Cambridge University Press, Cambridge--New York, 1999).

\bibitem{Stanley_10}  R.P. Stanley, A survey of alternating permutations,
   in: {\em Combinatorics and Graphs}\/,
   edited by R.A.~Brualdi, S.~Hedayat, H.~Kharaghani, G.B.~Khosrovshahi
   and S.~Shahriari,
   Contemporary Mathematics \#531
   (American Mathematical Society, Providence, RI, 2010), pp.~165--196.
   %% http://www-math.mit.edu/~rstan/papers/altperm.pdf


% \bibitem{Stanley_15}  R.P. Stanley, {\em Catalan Numbers}\/
%    (Cambridge University Press, New York, 2015).
% 
% \bibitem{Stanley_18}  R.P. Stanley and Y. Wang, Some aspects of
%    $(r,k)$-parking functions, J. Combin. Theory A {\bf 159}, 54--78 (2018).
% 
% \bibitem{Stembridge_91}  J.R. Stembridge,
%    Immanants of totally positive matrices are nonnegative,
%    Bull. London Math. Soc. {\bf 23}, 422--428 (1991).

\bibitem{Stieltjes_1889}  T.J. Stieltjes, Sur la r\'eduction en fraction
   continue d'une s\'erie proc\'edant selon les puissances descendantes
   d'une variable,
   Ann. Fac. Sci. Toulouse {\bf 3}, H1--H17 (1889).
   %% https://eudml.org/doc/72571
   %% http://www.numdam.org/item?id=AFST_1889_1_3__H1_0

% \bibitem{Stieltjes_1894}  T.J. Stieltjes, Recherches sur les fractions
%     continues, Ann. Fac. Sci. Toulouse {\bf 8}, J1--J122 (1894)
%     and {\bf 9}, A1--A47 (1895).
%     %% https://eudml.org/doc/72663  AND  https://eudml.org/doc/72665
%     [Reprinted, together with an English translation,
%      in T.J. Stieltjes, {\em \OE{}uvres Compl\`etes/Collected Papers}\/
%      (Springer-Verlag, Berlin, 1993), vol.~II, pp.~401--566 and 609--745.]

% \bibitem{Thoma_64}  E. Thoma, Die unzerlegbaren, positiv-definiten
%    Klassenfunktionen der abz\"ahlbar unendlichen, symmetrischen Gruppe, 
%    Math. Z. {\bf 85}, 40--61 (1964).
% 
% \bibitem{Varvak_04}  A.L. Varvak, Encoding properties of lattice paths,
%    Ph.D.~thesis, Brandeis University, May 2004.
%    Available on-line at
%    \url{http://people.brandeis.edu/~gessel/homepage/students/varvakthesis.pdf}
% 
\bibitem{Viennot_81}  G. Viennot, Interpr\'etations combinatoires des nombres
   d'Euler et de Genocchi, S\'eminaire de Th\'eorie des Nombres (Bordeaux),
   Ann\'ee 1980--81, expos\'e no.~11 (1981).
   %% https://www.jstor.org/stable/44165433

\bibitem{Viennot_83}  G. Viennot, Une th\'eorie combinatoire des polyn\^omes
   orthogonaux g\'en\'eraux, Notes de conf\'erences donn\'ees
   \`a l'Universit\'e du Qu\'ebec \`a Montr\'eal,
   septembre-octobre 1983.
   Available on-line at
   \url{http://www.xavierviennot.org/xavier/polynomes_orthogonaux.html}

% % \bibitem{Viennot_85}  G. Viennot, A combinatorial theory for general
% %    orthogonal polynomials with extensions and applications,
% %    in {\em Polyn\^omes Orthogonaux et Applications}\/
% %    (Lecture Notes in Mathematics \#1171),
% %    edited by C.~Brezinski {\em et al.}\/
% %    (Springer-Verlag, Berlin, 1985), pp.~139--157.
% 
% \bibitem{Visontai_13}  M. Visontai, Some remarks on the joint distribution
%    of descents and inverse descents,
%    Electron. J. Combin. {\bf 20}, no.~1, paper \#P52 (2013).

% \bibitem{Wall_40}  H.S. Wall,
%    Continued fractions and totally monotone sequences,
%    Trans. Amer. Math. Soc. {\bf 48}, 165--184 (1940).

% \bibitem{Wall_48}  H.S. Wall, {\em Analytic Theory of Continued Fractions}\/
%    (Van Nostrand, New York, 1948).
%    Reprinted by the American Mathematical Society, Providence RI, 2000.
% 
% \bibitem{wikipedia}  Wikipedia, Gauss's continued fraction,
%    \url{http://en.wikipedia.org/wiki/Gauss%27s_continued_fraction}
% 
% \bibitem{Woan_01}  W.-J. Woan, Hankel matrices and lattice paths,
%    J. Integer Seq. {\bf 4}, article 01.1.2 (2001).
% 
% \bibitem{Yan_15}  C.H. Yan, Parking functions, in
%    {\em Handbook of Enumerative Combinatorics}\/,
%    edited by M.~Bona (CRC Press, Boca Raton FL, 2015), pp.~835--893.

\bibitem{Zeng_96}  J. Zeng,
   Sur quelques propri\'et\'es de sym\'etrie des nombres de Genocchi,
   Discrete Math. {\bf 153}, 319--333 (1996).

% \bibitem{Zhu_13}  B.-X. Zhu, Log-convexity and strong $q$-log-convexity
%    for some triangular arrays,
%    Adv. Appl. Math. {\bf 50}, 595--606 (2013).
% 
% \bibitem{Zhu_14}  B.-X. Zhu, Some positivities in certain triangular arrays,
%    Proc. Amer. Math. Soc. {\bf 142}, 2943--2952 (2014).
% 
% \bibitem{Zhu_17a}  B.-X. Zhu,
%    Log-concavity and strong $q$-log-convexity for Riordan arrays
%    and recursive matrices,
%    Proc. Roy. Soc. Edinburgh A {\bf 147}, 1297--1310 (2017).
%    %%  https://doi.org/10.1017/S0308210516000500
% 
% \bibitem{Zhu_18}  B.-X. Zhu,
%    Total positivity, continued fractions and Stieltjes moment sequences,
%    in preparation.



\end{thebibliography}
\end{document}